\documentclass[11pt]{article}%{amsart}
\usepackage{mathrsfs}
\usepackage{amsmath,amssymb}
\usepackage{amsthm,bm}
\usepackage[colorlinks,linkcolor=blue,anchorcolor=blue,citecolor=green]{hyperref}
\usepackage{graphics,graphicx}
\usepackage{color}
\usepackage{enumerate}
\usepackage[title]{appendix}
\usepackage{subfigure}

\usepackage{caption}

\allowdisplaybreaks

\usepackage{tikz}
\usepackage{indentfirst}

\newtheorem{Thm}{Theorem}[section]
\newtheorem{theorem}[Thm]{Theorem}

\newtheorem{definition}[Thm]{Definition}

\newtheorem{lemma}[Thm]{Lemma}

\newtheorem{corollary}[Thm]{Corollary}

\newtheorem{proposition}[Thm]{Proposition}

\newtheorem{remark}[Thm]{Remark}

\newtheorem{example}[Thm]{Example}
\def\R{{\mathbb{R}}}
\def\E{{\mathbb{E}}}
\def\P{{\mathbb{P}}}

\def\Z{{\mathbb{Z}}}
\def\b{{\mathfrak{b}}}
\def\dd{{\mathfrak{d}}}
\def\pp{{\mathbf{p}}}
\def\qq{{\mathbf{q}}}
\def\vv{{\mathcal{V}}}
\def\zz{{\mathbf{z}}}
\def\xx{{\mathbf{x}}}
\def\yy{{\mathbf{y}}}
\def\nuu{\mbox{\boldmath$\nu$\unboldmath}}
\def\tauu{\mbox{\boldmath$\tau$\unboldmath}}
\title{Large deviation principle for empirical measures of once-reinforced random walks on finite graphs
	\footnote{Y. Liu is supported by CNNSF (No. 11731009, No.11926327, No. 12231002) and Center for Statistical Science, PKU.
		K. Xiang is supported by CNNSF (No.~12171410) and by Hu Xiang Gao Ceng Ci Ren Cai Ju Jiao Gong Cheng-Chuang Xin Ren Cai (No. 2019RS1057).}}
\author{
	Xiangyu Huang\footnote{CMS, Chongqing University, Chongqing 401331, China. \emph{xiangyuhuang077@gmail.com}}
	\and
	Yong Liu\footnote{LMAM, School of Mathematical Sciences, Peking University, Beijing 100871, China. \emph{liuyong@math.pku.edu.cn}}
	\and
	Kainan Xiang\footnote{School of Mathematics and Computational Science, Xiangtan University, Xiangtan City 411105, Hunan Province, China. \emph{kainan.xiang@xtu.edu.cn}}}
\date{}

\usepackage{caption}
\usepackage{amsthm,amsmath,amsfonts,amssymb,mathrsfs}
\usepackage{graphics,graphicx}
\usepackage{bm}

\usepackage[square,numbers]{natbib}

\usepackage[title]{appendix}
\usepackage{subfigure}

\usepackage{tikz}
\usepackage{indentfirst}

\begin{document}
	\maketitle
	\begin{abstract}
		A $\delta$ once-reinforced random walk ($\delta$-ORRW) on connected graph is a self-interacting random walk which moves to its neighbors at each step  according to the weights of the edges at that time, where the weights are $1$ on edges that have not been traversed and $\delta$ otherwise. In this paper, we prove a large deviation principle for empirical measures of $\delta$-ORRWs on finite connected graphs using a  modified weak convergence approach. The rate function of the large deviation principle exhibits a phase transition at the $\delta=1$.\\

	\end{abstract}

	\tableofcontents

\section{Introduction}\label{sec 1}
\setcounter{equation}{0}{\color{black}
	
	\noindent 
	Suppose that $G=(V,E)$\label{G=(V,E)} is a connected locally finite undirected graph with vertex set $V$ and edge set $E$.  If two vertices $u,v\in V$ are adjacent, we write $u\sim v$\label{u sim v} and denote the corresponding edge by $uv=\{u,v\}$. For any subgraph $G'=(V',E')\subseteq G$, we define $\partial E'$ to be the set of all edges $e'\in E'$ for which  there exists an edge  $e\in E\setminus E'$ adjacent to $e'$. We denote the vertex set and the edge set of $G'$ by $V_{G'}$ and $E_{G'}$ \label{VG-EG}.
	
	The edge reinforced random walk (ERRW) $X=(X_n)_{n\geq 0}$\label{X=(X_n)_{n>0}} on $G$ is a self-interacting non-Markovian random walk  introduced by Coppersmith and Diaconis \cite{CD1987} in 1987. This model captures the idea that a random walker may prefer to traverse edges that have been visited previously. Mathematically, let
	$(\mathscr{F}_n)_{n\ge 0}$ be the natural filtration  of $X$. Then, the ERRW $X$ is defined by 
	\[\P\left(X_{n+1}=u\left|X_n=v,\mathscr{F}_n\right.\right)=\left\{
	\begin{aligned}
		&\frac{w_n(uv)}{\mathop{\sum}_{u'\sim v}w_n(u'v)}, & u\sim v,\\
		&0,&\text{otherwise},
	\end{aligned}\right.
	\]
	where $w_n(e)$ is the  weight of edge $e$ at time $n$. In the original model, 
	$$
	w_0(e)=1,\ w_{n+1}(e)=w_n(e)+\delta {\bm 1}_{\{X_nX_{n+1}=e\}},\ e\in E,\,n\geq 0,
	$$ 
	for some $\delta>0$.

	The study of the ERRW is pretty challenging because its future depends on its entire past trajectory. 
	The ERRW does satisfy some partial exchangeability which has been used to obtain several results, such as the recurrence/transience on lattice graphs $\mathbb{Z}^d$ for all integers $d\ge1$ in \cite{ACK2014, P2019, SC-TP2015, SC-ZX2019}, and invariant principle on $b$-regular trees for $b\ge 4$ in \cite{CTU2020}. The ERRW is closely related to the continuous reinforced process studied in \cite{MP-SP2003} (see also \cite{PR2007}), the random Schr\"{o}dinger operator in \cite{SC-ZX2019}, and the $\mathbb{H}^2$ and $\mathbb{H}^{2|2}$ spin systems in \cite{SC-TP2015} (see also \cite{BR-HT2021}). The ERRW can be applied in other scientific fields, such as  biology (\cite{CE-PM2008, SA-OH1997}) and social networks (\cite{KD-TP2016, PR-SB2004, SB-PR2009, VS2001}).

	\vskip2mm
	
	In 1990, Davis introduced a variant of the ERRW: the once reinforced random walk (ORRW) in \cite{Dav1990}. 
	For the ORRW, all initial edge weights are set to $1$, once an edge is traversed, its weight changes permanently to a constant $\delta>0$. That is, for all  $n\ge 0$ and $e\in E$, 
	\begin{equation}\label{weight}
		w_n(e)=1+(\delta-1)\cdot {\bm 1}_{\{N(e,n)>0 \}}=\left\{\begin{array}{ll}
			1 & {\rm if}\ N(e,n)=0,\\
			\delta & {\rm if}\ N(e,n)>0,
		\end{array}
		\right.
	\end{equation}
	where $N(e,n):=|\{i<n:X_iX_{i+1}=e \}|$ is the number of times edge $e$ has been traversed  up to time $n$, and $|A|$ denotes the cardinality of any set $A$. If $\delta<1$, $X$ is negatively reinforced,  if $\delta>1$, it is positively reinforced. We refer to the  ORRW as the ORRW with parameter $\delta$ or $\delta$-ORRW. 
	
	When $\delta < 1$, the reinforcement is negative,  so that  the walk is self-repelling; that is, it prefers to explore vertices that have not been visited before rather than returning to previously visited ones. In contrast, when $\delta > 1$, the reinforcement is positive, making the walk self-attracting;  it prefers to return to  vertices that have been visited before instead of  exploring  new ones.
	
	Compared to ERRWs, the ORRW does not satisfy the partial exchangeability. This makes the study of ORRWs is more challenging than that of ERRWs.

	\subsection{Overview of results}
	\noindent 
	
	Let $G=(V,E)$ be a finite connected graph with $\vert E\vert\geq 1$ and let $X$ be a $\delta$-ORRW 
	on $G$. We define the empirical measure process $(L^n)_{n\geq 1}$ of the $\delta$-ORRW $X$ by
	\begin{align}
		L^n(A)=\frac{1}{n}\sum_{i=0}^{n-1}\dd_{X_i}(A)\text{ for }A\subseteq V,\label{eq-emprical-ORRW}
	\end{align}
	where $\dd_x$\label{hat{delta}_x} is the Dirac measure centered at $x\in V$. Since $G$ is finite and connected,  the $\delta$-ORRW will almost surely visit every edge of $G$ in a finite time (as will be confirmed in Proposition \ref{exponentially integrable}); that is, the edge cover time 
	$$ 
	C_E=\inf\{t:\ \forall e\in E,\ \exists\,s\le t,\ X_{s-1}X_{s}=e \}
	$$\label{cover time}
	is almost surely  finite. Note that after all edges have been traversed,  $\delta$-ORRW reduces to a simple random walk (SRW) on $G$. Thus, by the ergodic theorem for the SRW,  the empirical measure process $(L^n)_{n\geq 1}$ converges weakly to the unique invariant probability of the SRW on $G$.
	The purpose of this paper is to study the large deviation principle (LDP) for the empirical measure process $(L^n)_{n\geq 1}$.
	
	\vskip2mm
	
	Now we recall the definitions of the large deviation principle (LDP) and a good rate function from \cite[Section 1]{BD2019}, \cite[Section 1.1]{DE1997} and \cite[Section 1.2]{DZ1998}:
	
	Let $\mathscr{X}$ be a Polish space, $I:\mathscr{X}\mapsto [0,\infty]$ a function on $\mathscr{X}$, and $(\xi_n)_{n\geq 1}$ a sequence of $\mathscr{X}$-valued random variables. $I$ is said to be  a good rate function if 
	\begin{itemize}
		\item[{(i)}] $I$ is lower semicontinuous, i.e., $\liminf_{y\to x}I(y)\ge I(x)$ for all $x\in \mathscr{X}$;
		\item[{(ii)}]  each level set $\{x:\ I(x)\le M \}$ is compact for all $M>0$; and
		\item[{(iii)}] there exists some $x$ with  $I(x)<\infty$. 
	\end{itemize}
	We say the sequence of  $(\xi_n)_{n\geq 1}$ satisfies a  LDP with a good rate function $I$ if  the following hold:
	\begin{itemize}
		\item[{(a)}] For every closed subset $C\subseteq \mathscr{X}$, $\limsup_{n\to\infty}\frac{1}{n}\log\P(\xi_n\in C)\le -\inf_{x\in C}I(x)$; and 
		\item[{(b)}] For every open subset $O\subseteq \mathscr{X}$, $\liminf_{n\to\infty}\frac{1}{n}\log\P(\xi_n\in O)\ge -\inf_{x\in O}I(x)$.
	\end{itemize}
	
	\vskip2mm

	When $\vert V\vert=2,$ every $\delta$-ORRW $X=(X_n)_{n\geq 0}$ on $G$ is simply  a SRW and the LDP for $(L^n)_{n\ge 1}$ is  trivial. Therefore, throughout this paper we always assume $\vert V\vert\geq 3$.

	For any Polish space $\mathscr{X}$, let $\mathscr{P}(\mathscr{X})$\label{mathscr{X}} be the set of all probability measures on $\mathscr{X}$ equipped with the weak convergence topology. Our main results are summarized as follows:
	\begin{theorem}\label{thm 1.1}
		If $\vert V\vert\geq 3$ (equivalently $\vert E\vert\geq 2$), 
		the empirical measure process $(L^n)_{n\geq 1}$ of the $\delta$-ORRW satisfies an LDP with a good rate function $I_{\delta}: \mathscr{P}(V) \mapsto [0,\infty]$. 
		For any $\nu\in \mathscr{P}(V)$, the function $\delta\to I_{\delta}(\nu)$ is continuous and  decreasing in $[1, \infty)$, constant on $ (0,1]$; and for some $\nu\in \mathscr{P}(V)$, $I_{\delta}(\nu)$ is not differentiable at $\delta=1$. For  all $1\leq\delta_1<\delta_2$, $I_{\delta_1}\not=I_{\delta_2}$.
	\end{theorem}
	
	For the details on the rate function, 
	see Theorems \ref{I_delta} and \ref{phase for rate function}.  	
	When $\delta>1$, $I_\delta$ ($\delta>1$) the rate function $I_\delta$ depends on the starting vertex of the 	$\delta$-ORRW (see Example \ref{ldp on urn},  Example \ref{ldp on path} and Example \ref{ldp on path-2} in Section \ref{sec 3-vertex}).  By contrast,  $I_1$ (the rate function of the SRW) does not depend on the starting point.

	\vskip2mm In this paper, we use the weak convergence approach (see Dupuis and Ellis's book \cite{DE1997} and Budhiraja and Dupuis's book \cite{BD2019}) to study the LDP. Two main ingredients in our arguments are as follows. 
	
	First, we establish  a more general variational formula (see Theorem \ref{estimate for rate})--the generalized Laplace principle. This formula yields  the LDP for empirical measures of the ORRWs on finite graphs, see Theorem \ref{I_delta}. It also enables us to obtain precise probability estimates for various statistics, such as the critical exponent for the exponential integrability of the cover time (which will be discussed in a separate paper, see \cite{HLX2025}). This general variational formula remains valid for general discrete Markov processes on finite graphs. It provides a tool for studying the tail distribution of cover times for discrete Markov processes, as well as estimates for some other stopping times. For a more detailed explanation, please refer to the end of Section \ref{sec 2.16}.
	
	Second, to overcome essential difficulties caused  by 
	the lack of Markovian  property 
	of the ORRWs $(X_n)_{n\geq 0}$, we lift $(X_n)_{n\geq 0}$ to the line digraph and adapt  the weak convergence approach originally developed for the LDP of empirical measures of homogeneous Markov chains. For further details, please refer to Sections \ref{sec 2.2} and \ref{sec 2.3}.

	\subsection{Literature review} 
	\subsubsection{Background on ERRWs and ORRWs}
	\noindent   We first review some results on ERRWs. An important topic in the study of ERRWs is the recurrence/transience phase transition with respect to the reinforcement parameter $\delta$. In 1988 \cite{Pe1988}, Pemantle initiated the study of ERRWs and discovered a recurrence/transience phase transition on a binary tree with a critical point $\delta_c\approx 4.29$. His analysis relies on the observation that the distribution of a certain simple random walk coincides with that of a P\'{o}lya urn. The P\'{o}lya urn  can also be interpreted as the ERRW on the star graph with vertices $\{v_0,v_1,\dots,v_d\}$ and  edge  $\{v_0v_i:i=1,\dots,d\}$ (see Figure \ref{polya urn graph}). Moreover,  Pemantle proved  a criterion for the recurrence/transience of ERRWs on infinite trees.

	In 2000, Takeshima \cite{T2000} proved a recurrence/transience phase transition (with respect to the initial weights rather than the reinforcement factor) of the ERRW on the half-line. In 2005, Merkel and Rolles \cite{MR2005} 
	proved that when $a>3/4$ and the weight function of the ERRW at each edge $e$ is
	$$
	w_n(e)=a+N(e,n),
	$$
	the ERRW on the ladder $\mathbb{Z}\times\{0,1\}$ is recurrent. 
	
	Angel, Crawford and Kozma (2014) \cite{ACK2014} showed that, when the initial weights are sufficiently small, ERRWs on graphs with bounded degrees are recurrent, and that, when the initial weights are sufficiently large, ERRWs on non-amenable graphs are transient. They also identified a constant $a<a_0(d)$ such that the ERRW is recurrent for $a<a_0(d)$ on $\Z^d$ when $d\ge 3$. Subsequently, Disertori, Sabot and Tarr\`{e}s (2015) \cite{DST2015} found a constant $a_1(d)$ such that the ERRW is transient for $a>a_1(d)$ on $\Z^d$ for $d\ge3$. 
	
	Sabot and  Tarr\`{e}s \cite{SC-TP2015} coupled the ERRW with the vertex-reinforced jump process (VRJP) with a random environment following a Gamma distribution. 
	They showed that the VRJP in this random environment has the same distribution as the ERRW 
	and proved
	a “magic formula” to describe this distribution. 
	This magic formula offers a general technique for studying ERRWs. The framework is linked to a supersymmetric hyperbolic sigma model known as the $\mathbb{H}^{2|2}$ model, which is a special case of an important model related to the Arboreal Gas and the Anderson model. For further details, we refer to the survey by Bauerschmidt and Helmuth (2021) \cite{BR-HT2021} on recent developments in this area.

	Building on these works, 
	Poudevigne (2024) \cite{P2019}
	proved the monotonicity of the recurrence on $\Z^d$ for $d\ge2$ with respect to the parameter  
	$a$, which implies that the recurrence/transience phase transition for the ERRW on  $\Z^d$, $d\ge 3$, is unique.
	By
	extending the results of Sabot,  Tarr\`{e}s and Zeng (2017) \cite{STZ2017} from finite to infinite graphs, Sabot and Zeng \cite{SC-ZX2019} proved in 2019 that the ERRW is recurrent on $\Z^2$ for all $a>0$.  For additional results, we refer the reader to    \cite{DSZ2010, MR2006, Se2010}.

	\vskip2mm
	Now we review literature on ORRWs. As mentioned before, ORRWs are more difficult to study due to the loss of Markov property and partial exchangeability. For example, determining their recurrence or transience on infinite graphs can be quite challenging. The first result on trees was obtained by Durrett, Kesten and Limic \cite{DKL2002} in 2002. They proved that, when for $\delta>1$, the transience of the $\delta$-ORRW on $d$-regular trees with $d\ge 3$. They also derived a law of large numbers and an invariance principle for the $\delta$-ORRW. 
	
	Dai (2005) \cite{Da2005} showed that, when $\delta>1$, the $\delta$-ORRW on Galton-Watson trees conditioned on non-extinction is transient. Collevecchio, Holmes and Kious (2018) \cite{CHK2017} extended this result to a more general multiplicative ORRW and proved its recurrence under certain conditions on the reinforcement.
	
	The first complete result on the recurrence/transience phase transition for the ORRW was obtained by Kious and Sidoravicius (2018) \cite{KS2018}. They established a sharp phase transition on $\Z^d$-like trees (a class of trees whose growth is polynomial and resembles that of  $\Z^d$). Specifically, they proved that the ORRW is recurrent for $\delta>\log_2d$ and transient for $\delta<\log_2d$. Subsequently, Collevecchio, Kious and Sidoravicius (2020) \cite{CKS2017} extended this result to general trees. They established a sharp recurrence/transience phase transition with a critical point characterized by the branching-ruin number.

	Recurrence for all $\delta>0$ is known to hold on lines with parallel edges, including   $\mathbb{Z}$ (Vervoort (2002) \cite{Ve2002}), and on $\mathbb{Z}\times\{0,1\}$ (Vervoort (2002) \cite{Ve2002}, Sellke (2006) \cite{Se2006}, Huang, Liu, Sidoravicius and Xiang (2021) \cite{HLSX2021}). Recall the $\delta$-ORRW on ladders $\mathbb{Z}\times\{0,1,\ldots,d\}$ with $d\geq 2$ is recurrent for $\delta\in \left(1-1/(d+1),1+1/(d-1)\right)$ as shown by Sellke \cite{Se2006} in 1994 (published in 2006) and Vervoort \cite{Ve2002} in 2002. More generally, Kious, Schapira and Singh  (2021) \cite{KSS2018} proved that the $\delta$-ORRW on ladder graphs $\mathbb{Z}\times\Gamma$, where  $\Gamma$ is  a finite connected graph,  is recurrent for  sufficiently large $\delta$. 
	
	To date, there is no recurrence/transience result for the $\delta$-ORRW on $\mathbb{Z}^d$ for all $d\ge 2$ and $\delta>0$. A conjecture by Sidoravicius (see \cite[p.\,2122]{KS2018}) states that every $\delta$-ORRW on $\mathbb{Z}^2$ is recurrent, while the on $\mathbb{Z}^d$ with $d\geq 3$ it undergoes a phase transition: recurrent for large $\delta$ and transient for small $\delta$.

	\subsubsection{Related works on the LDP}
	\noindent 
	Now we recall some previous large deviation results for Markov processes. 
	The first LDP for empirical measures was given by Sanov \cite{S1957} in 1957. 
	He proved that the empirical measures of an i.i.d. sequence with a common distribution $\rho$ on a Polish space $\mathscr{X}$ satisfy an LDP with the rate function $R(\cdot\|\rho)$.
	Here, for any probability $\gamma$ on $\mathscr{X}$, $R(\gamma)\|\rho)$ is the relative entropy defined by
	\begin{equation}\label{reltv entrpy}
		R(\gamma\|\rho)=\left\{\begin{array}{cl}
			\int_{\mathscr{X}}\frac{{\rm d}\gamma}{{\rm d}\rho}\log{\frac{{\rm d}\gamma}{{\rm d}\rho}}\ {\rm d}\rho & {\rm if}\ \gamma\ll\rho\\
			\infty  & {\rm otherwise}
		\end{array}
		\right..
	\end{equation}

	In 1975, the LDPs for empirical measures of discrete and continuous time parameter Markov processes were investigated  by Donsker and Varadhan in the seminal papers \cite{DV1975_1,DV1975_2,DV1977}.
	They proved that the empirical measure of 
	a Markov process on a compact metric space satisfies the LDP in \cite{DV1975_2}.  
	In the case of a discrete time Markov process with transition probability $P$, the  rate function is given by 
	$$
	I(\mu)=-\inf_{u\in\mathscr{U}_1}\int_{\mathscr{X}}\log{\frac{Pu}{u}(x)}\ \mu({\rm d}x)
	$$
	for any probability measure $\mu$ on the state space, 
	where $\mathscr{U}_1$ denotes the collection of all positive functions $u$ on the state space $\mathscr{X}$, see \cite[Theorem 1]{DV1975_2}.

	The classical studies in \cite{DV1975_1,DV1975_2,S1957}  rely heavily on the Markov property.
	When $\delta\not=1$, the $\delta$-ORRW is not Markovian, so one can not apply the results of \cite{DV1975_1,DV1975_2,S1957} to get the LDP for empirical measures of the $\delta$-ORRW.
	There are some LDP results for non-Markovian processes;  for instance, an LDP of renewal processes was obtained  in \cite{LMZ2011,MZ2014,MZ2016}. However, these results are based on the Markov property of an  embedding chain, which essentially falls within the framework of Donsker and Varadhan. Therefore, they are also not directly applicable to ORRWs.

	To the best of our knowledge, the only existing LDP-type result for  ORRW was given by Zhang \cite{Z2014} in 2014. 
	We now describe the main result of  \cite{Z2014}. Let $b\geq 2$ be an integer and let $\mathbf{T}$ be an infinite rooted tree   in which the root $\varrho$ has $b$ neighbors and every other vertex has $b+1$ neighbors. 
	For any vertex $x\in\mathbf{T}$, denote by $h(x)$ the distance between $x$ and $\varrho$ (i.e., the height of $x$ relative to $\varrho$). 
	Let $\delta>1$ and let $(X_n)_{n\geq 0}$ be a $\delta$-ORRW on $\mathbf{T}$.  There is a positive constant $s=s(\delta)$ such that 
	$$
	\lim\limits_{n\rightarrow\infty}\frac{h(X_n)}{n}=s\ a.s..
	$$   
	Zhang proved that for any $\varepsilon>0$ there exists  a constant   $\beta=\beta(\delta,b,\varepsilon)\in (0,\infty)$ satisfying 
	$$
	\lim_{n\to\infty}-\frac{1}{n}\log{\P\left(h(X_n)\ge(s+\varepsilon)n\right)}=\beta;
	$$
	and for   $0<\varepsilon<s$,
	$$
	0<\liminf_{n\to\infty}-\frac{1}{n}\log{\P\left(h(X_n)\le(s-\varepsilon)n\right)}\leq \limsup_{n\to\infty}-\frac{1}{n}\log{\P\left(h(X_n)\le(s-\varepsilon)n\right)}<\infty.
	$$
	
	Compared with \cite{Z2014}, we deal with general finite graphs rather than trees and we provide a variational representation of the rate function for the LDP of the empirical distribution of ORRWs, rather than focusing on the distance of the walk from the root. 
	
	Recently, Budhiraja and Waterbury \cite{BW2022} and Budhiraja, Waterbury and Zoubouloglou \cite{BW2023} established   LDPs for the empirical measures of a class of self-interacting Markov chains (such as reinforced random walks) on finite spaces $\{1,2,\ldots,d\}$ using the weak convergence approach. The classical LDP result of Donsker and Varadhan (1975) \cite{DV1975_2} on empirical measures of Markov processes is a specific case of these LDPs.  However,
	unlike the Donsker--Varadhan rate function for  the empirical measures of Markov processes, the rate functions obtained in \cite{BW2022, BW2023} are typically non-convex and expressed through dynamical variational formulas with an infinite-horizon discounted objective function. 
	The reinforced random walks considered in \cite{BW2022, BW2023} do not include the ORRWs. Although we also employ  the weak convergence approach, our techniques differ from those in  \cite{BW2022, BW2023}.  Their work relies heavily on the linearity of the weight function of ERRWs, which allows ERRWs to be viewed as a mixture of two Markov processes or as a random walk in a random environment (see \cite{SC-TP2015}). 
	The techniques of \cite{BW2022, BW2023} do not work for ORRWs, we need to develop a different method to study LDPs for ORRWs.

	\vskip2mm
	
	We use the weak convergence approach to  prove the LDP for the empirical measures of ORRWs. The weak convergence approach relies on the following two key points:
	
	\begin{itemize}
		\item[(i)] the Laplace principle implies the large deviation principle (see Definition \ref{LP} and Theorems \ref{thm-LP-LDP}), and
		\item[(ii)] the logarithmic Laplace functional admits a variational representation in terms of the relative entropy (see Proposition \ref{variational representation}).
	\end{itemize} 
	
	\begin{definition}\label{LP}
		Let $\mathscr{X}$ be a Polish space and $C_b(\mathscr{X})$ the set of all bounded continuous functions on 
		$\mathscr{X}$, and $(\xi_n)_{n\geq 1}$ a sequence of $\mathscr{X}$-valued random variables. We say that $(\xi_n)_{n\geq 1}$ satisfies the Laplace principle (LP)  
		if for some good rate function $I$ on $\mathscr{X}$, the following holds for any $h\in C_b(\mathscr{X})$:
		\begin{equation}\label{(2.1)}
			\lim\limits_{n\to\infty}\frac{1}{n}\log\E\left[\exp{\{-nh(\xi_n) \}}\right]=-\inf_{x\in\mathscr{X}}\{I(x)+h(x)\}.
		\end{equation}	
	\end{definition}
	
	\begin{theorem}[\mbox{\cite[Theorem 1.2.3]{DE1997}}]\label{thm-LP-LDP}
		The LP implies the LDP with the same rate function. Namely, if $I$ is a good rate function on Polish space $\mathscr{X}$ and
		$(\xi_n)_{n\geq 1}$ is a sequence of $\mathscr{X}$-valued random variables such that \eqref{(2.1)} holds for any $h\in C_b(\mathscr{X})$, then $(\xi_n)_{n\geq 1}$ satisfies the LDP with the rate function $I$.
	\end{theorem}
	
	\begin{remark} 
		\rm In fact, the LDP also implies the LP (\cite[Theorem 1.2.1]{DE1997}). 
	\end{remark}
	\begin{proposition}[\mbox{\cite[Proposition 1.4.2]{DE1997}}]\label{variational representation}
		Given Polish space $\mathscr{X}$ and $\theta\in\mathscr{P}(\mathscr{X}).$ Then
		\begin{equation}
			-\log\int_\mathscr{X} e^{-h}\ {\rm d}\theta=\inf\limits_{\gamma\in\mathscr{P}(\mathscr{X})}\left\{
			R(\gamma\|\theta)+\int_\mathscr{X} h\ {\rm d}\gamma \right\},\ h\in C_b(\mathscr{X}),\label{(2.2)}
		\end{equation}
		where the infimum is uniquely attained at $\gamma_0$ with
		$
		\frac{{\rm d}\gamma_0}{{\rm d}\theta}=e^{-h}\cdot\frac{1}{\int_{\mathscr{X}}e^{-h}\ {\rm d}\theta}.
		$
	\end{proposition}
	
	One can use Theorem \ref{thm-LP-LDP} and Proposition \ref{variational representation} to get the large deviation principle for the empirical measure of the Markov chain. For further details, see Section \ref{sec 2.2}.

	\subsection{Outline of the paper}
	\noindent This paper is organized as follows. In Section \ref{sec 2}, we present the main results of the paper first. Then we give  the  core ideas of the paper by recalling the weak convergence approach and our adaptation of it,  and explaining how the adapted weak convergence approach is applied in our study.  The last part of Section \ref{sec 2} summarizes notation used in this paper.
	Section \ref{sec 3} contains the proof of the LDP for the empirical measures of ORRWs and discusses several interesting properties of the rate function. In Section \ref{sec 5}, we simplify the rate functions on trees and provide  explicit expressions for some specific graphs. Finally, in Section \ref{sec final}, 
	%RS we present supplementary results 
	we give some auxilliary results 
	and prove some propositions, lemmas, and theorems that appear in the main text.

	\vskip2mm

	After the submission of this 
	paper and during 
	the revision process
	we noticed several arXiv preprints on ORRWs. Collevecchio and Tarr\`{e}s in \cite{CT2025} showed that the ORRWs are transient on non-amenable graphs for small $\delta$. A key ingredient in their proof is an
	LDP estimate for the empirical measures of ORRWs on finite graphs.  In fact, our LDP result in this paper implies    theirs.  
	
	For ORRWs on $\mathbb{Z}^d$ for $d\ge6$ and $\delta$ sufficiently small, Elboim and Kozma in \cite{EK2026} established transience and an invariance principle by verifying the existence of a relaxed time and coupling the ORRW to a concatenation of many independent shorter processes. For the ORRWs on $\mathbb{Z}^d$ for $2\le d\le 5$, the recurrence/transience results are still unknown for all $\delta\neq 1$.
	
	Very recently, Hu, Ma, Song and Wang in \cite{HMSW2026} studied  the asymptotic behavior of the range of ORRWs on the half-line.

	\section{Main results and preliminaries}\label{sec 2}
	\setcounter{equation}{0}
	\noindent 	In Section 2.1, we present the full statements of the main theorems (see Theorem \ref{I_delta} and Theorem  \ref{phase for rate function}). In Section 2.2, we lift $\delta$-ORRW $X$ to the line digraph of $G$, so that the pair consisting of the lifted process and its empirical measure process becomes a Markov chain, which makes it possible to apply the weak convergence approach. Section 2.3 introduces the generalized Laplace principle and explains the motivation for introducing it. To help our readers  better understand our method, we briefly outline in Section 2.4 the key steps in the proof of the large deviation principle for empirical measures of homogeneous Markov chains on finite state spaces using the weak convergence approach. In Section 2.5, we give an intuitive explanation of the main modifications we make to the weak convergence method in order to prove the LDP for the empirical measures of ORRWs. Finally, in Section 2.6, we provide an index of frequently used notations in the paper for the reader’s convenience.

	\subsection{Main results}\label{sec 2.1}
	\noindent  
	Consider a $\delta$-ORRW $X=(X_n)_{n\ge 0}$ on a finite connected graph $G=(V,E)$ with  initial position $X_0=x_0\in V$. Let $\mathfrak{b}:=|E|\ge 2$. Denote by  $\mathbb{P}_{x_0}$  the law of the $X$ starting at $x_0$, and by $\mathbb{E}_{x_0}$  the corresponding expectation under $\mathbb{P}_{x_0}$. Let $\omega$ be the sample path of the process, which is also an element of the probability space.

	To help clarifying the somewhat intricate mathematical formulation of the theorems and also the key ideas behind their proofs, we first describe the structure of $(X_n)_{n\ge 0}$. 
	\begin{itemize}
		\item[(i)] $X_0=x_0$.   If there are $k_0$ edges incident to $x_0$, then at time $1$, $X$ selects one of these edges, say $x_0v_1$, with probability $\frac{1}{k_0}$, and move to $v_1$. 
		Thus, $X_1=v_1$.    Note that $(V_1,E_1):=\big(\{x_0, v_1\}, \{x_0v_1\}\big)$ is a subgraph of $G$.

		\item[(ii)]   If $v_1$ is incident to $k_1$ edges other than $x_0v_1$, then at time $2$, $X$ has two possible type of options:  
		\begin{itemize}
			\item it retunes to $x_0$ via the edge $x_0v_1$ with probability $\frac{\delta}{\delta+k_1}$;\vskip1mm
			\item otherwise, it moves to a vertex $v_2$ by  traversing a different edge incident to $v_1$, $v_1v_2$ with probability $\frac{1}{\delta+k_1}$.  
		\end{itemize} 
		In both cases we denote the position of $X$ at time 2 as $X_2=v_2$.  Note that $\big(\{x_0, v_1, v_2\}, \{x_0v_1, v_1v_2\}\big)$ is a subgraph of $G$, and $(V_1, E_1)$ is a subgraph of $\big(\{x_0, v_1, v_2\}$, $\{x_0v_1, v_1v_2\}\big)$.
		\item[(iii)]  Suppose vertex $v_2$ is incident to $k_2'$ edges. Among these, let $k_2''$ denote the number of edges that have been traversed, and let $k_2:=k_2'-k_2''$ denote the number of edges that have not been traversed. At time 3, $X$ again has two possible type of options:
		\begin{itemize}
			\item  With probability $\frac{1}{k_2''\delta+k_2}$, it selects 
			each edge 
			from the $k_2''$ previously traversed edges incident to $v_2$ and moves back along the chosen edge to its other endpoint;    
			\vskip1mm
			\item With probability $\frac{1}{k_2''\delta+k_2}$, it selects each edge from the $k_2$ untraversed edges incident to $v_2$,  and moves to the vertex at the other end of the chosen edge.
		\end{itemize}  
		In either case we denote the position of $X$ at time $3$ as $X_3=v_3$. The 
		vertices and edges visited up to time $3$  form the subgraph  $\big(\{x_0, v_1, v_2, v_3\}$, $\{x_0v_1, v_1v_2, v_2v_3\}\big)$ of $G$, which naturally contains the subgraph up to time $2$,  $\big(\{x_0, v_1, v_2\}, \{x_0v_1, v_1v_2\}\big)$.
		
		\item[(iv)] The process continues inductively. For each subsequent time step $n\ge 4$, $X$ at its current vertex makes an analogous choice between backtracking along a used edge and  moving forward along a new edge, with probabilities following the same rule as described for times $2$ and $3$.
	\end{itemize}
	
	Thus, $X$ can be decomposed according to the following inductive scheme.
	\begin{itemize}
		\item[(a)] Set $\tau_1:=1$; at this time $X$  traverses its first new edge.
		\item[(b)] Define $\tau_2>\tau_1$ as the stopping time at which $X$ first leaves the subgraph 
		$$
		(E_1, V_1):=\big(\{x_0,v_1\}, \{x_0v_1\}\big).
		$$
		During $[1, \tau_2-1]$, $X$ remains within this subgraph. At the transition $\tau_2-1\to \tau_2$, it moves along a new edge to a vertex we denote by $v_2'$.
		Observe that throughout the whole interval $[1, \tau_2-1]$, the edge $x_0v_1$ has  weight $\delta$.    Hence, the behavior of $X$ on  $[\tau_1, \tau_2-1]$ coincides with the transition probabilities of a simple random walk on this one-edge subgraph.
		\item[(c)]	Define $\tau_3>\tau_2$ as the stopping time when $X$ first leaves the subgraph 
		$$
		(E_2, V_2):=\big(\{x_0,v_1,v_2'\}, \{x_0v_1, v_1v_2'\}\big).
		$$ 
		During $[\tau_2, \tau_3-1]$), $X$ stays inside  this subgraph. At the transition $\tau_3-1\to \tau_3$, it moves along a new edge to a vertex we denote by $v_3'$.
		In the time  interval $[1, \tau_3-1]$, the edges $x_0v_1$ and $v_1v_2'$ both have weight $\delta$.    Therefore,  the behavior of $X$ on  $[\tau_2, \tau_3-1]$ coincides with the transition probabilities of a simple random walk on this two-edge subgraph.

		\item[(d)]    Proceeding inductively we obtain an increasing sequence of stopping times $\tau_1<\tau_2<\dots<\tau_{\mathfrak{b}-1}$ with $\tau_{\mathfrak{b}}$ being the edge cover time $C_E$. By construction, at time $C_E$ every  edge have weight $\delta$, After the cover time,  $X$ simply follows the transition probability of a simple random walk on the whole graph $G$.
	\end{itemize}

	The decomposition outlined above is central to the proofs of our main theorem and to the representation of the large deviation rate function. To state the theorems  precisely, we now formally introduce the following concepts and notations. 
	\begin{definition}[Renewal time of $\delta$-ORRW $X$] \label{def of renewal time} \  Set $\tau_1:=1$,
		$$
		\tau_k:=\inf\big\{j>\tau_{k-1}: X_{j-1}X_j\notin \{X_{i-1}X_i:i<j\}\big\}.
		$$
	\end{definition}
	
	$\tau_k$ is a stopping time with respect to $(\mathscr{F}_n)_{n\ge 0}$,  and represents the 
	$k$-th time $X$ traverses a new edge. In particular, $\tau_{\mathfrak{b}}$ coincides with the edge cover time $C_E$.
	
	\vskip1mm
	%Based on the decomposition of $X$ via renewal times, we now introduce the collection of all possible %sequences of edge sets that satisfy the following structural properties.
	\begin{definition} [Renewal subsets of the edge set] \ Let $\mathscr{E}$ \label{mathscr{E}} denote the collection of all sequences of edge subsets $\{E_k\}_{1\le k\le \mathfrak{b}}$   satisfying the following conditions:
		\begin{itemize}
			\item[\bf{A1}]$E_k\subset E_{k+1}\subseteq E,\ 1\leq k<\mathfrak{b}$.
			\item[\bf{A2}]$|E_k|=k,\ 1\leq k\leq \mathfrak{b}$.
			\item[\bf{A3}]There exists $v\in V$ such that $x_0v\in E_1$.
			\item[\bf{A4}]The unique edge in $E_{k+1} \setminus E_k$ is adjacent to some edge in $E_k$.
		\end{itemize}
		%{\color{pink} Each $E_k$ thus induces a connected graph on its vertices; we denote this graph by $G_k:=(V_k,E_k)$. ??}
	\end{definition}
	Intuitively, $\mathscr{E}$ collects all sequences of edge sets that can be generated by the edge-traversal process 
	$(X_nX_{n+1})_{n\geq 0}$. More precisely, for any fixed sample path $\omega$, if we define 
	\[E_k(\omega)=\{X_{n-1}(\omega)X_{n}(\omega):\ n\le{\tau}_{k}(\omega)\},\]
	then $\{E_k(\omega)\}_{k=1}^\mathfrak{b}$ belongs to $\mathscr{E}$.
	
	\vskip1mm
	
	We now define two types of  transition probabilities  associated with  an edge  subset 
	$E'\subset E$. Let $G':=(V',E')$ denote the graph induced by $E'$ on the vertices incident to its edges.  
	
	The first one describes the transition probability that, starting from a vertex $x\in V'$, $X$ moves along to an edge $xy$
	(which may or may not belong to $E'$) and reaches vertex $y$. 
	\vskip1mm
	
	\begin{definition} \label{Def for p_u on G}For a subset $E'\subset E$, define
		\[
		p_{E'}(x,y):=\left\{
		\begin{aligned}
			&\frac{g(\delta,E',xy)}{\sum_{z\sim x}g(\delta,E',xz)},\ &y\sim x\\ 
			&0,\ &\text{\rm otherwise},
		\end{aligned}	
		\right.,
		\]
		where 
		$$
		g(\delta,E',xy):=\delta\cdot {\bf{1}}_{\{xy\in E' \}} + {\bf{1}}_{\{xy\notin E' \}}.
		$$
		In particular, $p_E$  corresponds to the transition probability of the simple random walk (SRW) on $G$, which we denote simply by $p$.		
		
	\end{definition}
	
	{\color{black}
		$p_{E'}$ gives  the transition probability of the ORRW when its set of traversed edges is exactly $E'$. 
	}
	
	The second one represents  the transition probability of a Markov chain on $G'=(V', E')$ that, at each step, moves from its current position only to adjacent vertices in $G'$.  \vskip1mm
	
	\begin{definition}  For a subset $E'\subset E$, denote by $\mathscr{T}(E')$ the collection of transition probabilities $q: V'\times V$ $\mapsto [0,1]$ satisfying 
		$$
		\sum_{z\in V}q(x,z)=1,\  \ q(x,y)>0,\ \text{only if }\, xy\in E'.
		$$
		Thus, every $q\in \mathscr{T}(E')$ corresponds to the transition probabilities of a Markov chain on $V'$ that moves only along edges in $E'$ 
		at each step. 
		In particular,  $\mathscr{T}_G:=\mathscr{T}(E)$  \label{mathscr{T}_G} is the set of transition probabilities of Markov chains on the vertex set $V$ of $G$ that  moves from the current vertex to a neighboring vertex at each step.
	\end{definition}
	
	For any Polish space $\mathscr{X}$, let $\mathscr{P}(\mathscr{X})$\label{mathscr{X}} be the set of all probability measures on $\mathscr{X}$ equipped with the weak convergence topology.
	\begin{definition} For any $\nu\in \mathscr{P}(V)$, define $\mathscr{A}(\nu)$ \label{A(nv)} as the set of all quadruples 
		$$
		(\nu_k,r_k,E_k, q_{k})_{1\leq k\leq \mathfrak{b}}
		$$ 
		satisfying the following conditions: 
		\begin{itemize}
			\item[(1)] {\textbf{Edge-set sequences:}}  $\{E_k\}_{1\leq k\leq \mathfrak{b}}\in\mathscr{E}$;
			\item[(2)] {\textbf{Proportions:}} $r_k\ge 0$ for all $k$, $\sum_{k=1}^\mathfrak{b}r_k=1$;
			\item[(3)] {\textbf{Transition probabilities and invariant probability measures:}} For each $k$, 
			$$
			q_k\in\mathscr{T}(E_k),\  \nu_k\in \mathscr{P}(V_k)\ \text{\ \rm such that}\  \nu_k q_k=\nu_k,
			$$ where
			$V_k$ is the set of all vertices of $E_k$, (i.e., $(V_k,E_k)$ denotes the graph induced by $E_k$ on the vertices incident to its edges);
			\item[(4)] {\textbf{Decomposition of $\nu$:}}\  $\nu:=\sum_{k=1}^{\mathfrak{b}}r_k \nu_k$. 
		\end{itemize}
		$\mathscr{A}(\nu)$ collects all representations  of $\nu$
		as a convex combination of invariant distribution $\nu_k$ for Markov chains whose transition probability belongs to $\mathscr{T}(E_k)$. Equivalently,
		\begin{align}
			&{\mathscr{A}}(\nu)=\Big\{ (\nu_k,r_k,E_k,{q}_{k})_{1\leq k\leq \b}:\ \{E_k\}_{1\leq k\leq \b}\in\mathscr{E},\ r_k\ge 0,\ \sum_{k=1}^{\b}r_k=1,\nonumber\\
			&\hskip 1.8cm \nu_k\in \mathscr{P}(V_k),\ 
			{q}_k\in\mathscr{T}(E_k),\ \text{such that }\nu_k{q}_k=\nu_k \text{ and } \sum_{k=1}^{\b}r_k \nu_k=\nu  \Big\}.\label{hat_A}
		\end{align}
		
	\end{definition}

	We are now in a position to state the main results. Recall the definition of the empirical measure process $(L^n)_{n\ge 1}$ given in (1.2). Our main theorems on the LDP for $(L^n)_{n\geq 1}$ are 
	Theorems \ref{I_delta} and  Theorem \ref{phase for rate function}.

	\begin{theorem}[LDP for $(L^n)_{n\geq 1}$] \label{I_delta}
		The empirical measure process $(L^n)_{n\geq 1}$ of the $\delta$-ORRW $X=(X_n)_{n\ge 0}$ on a finite connected graph $G=(V,E)$ satisfies an LDP  with the good rate function $I_\delta$ given by 
		\begin{equation}
			\label{simpler expression of I}
			I_{\delta}(\nu)=\inf_{(\nu_k,r_k,E_k,{q}_k)_{k}\in{\mathscr{A}}(\nu)}
			\sum_{k=1}^\b r_k\int_{V}R\big({q}_{k}(x,\cdot)\|{p}_{E_k}(x,\cdot)\big)\ \nu_{k}({\rm d}x),\ \nu\in\mathscr{P}(V),
		\end{equation}
		where $R(\cdot\|\cdot)$ denotes the relative entropy. 
		
		Moreover,	 for any $\delta>0$, $I_\delta(\nu)<\infty$ if and only if $\nu$ is an invariant distribution for some $q\in \mathscr{T}_G$.
	\end{theorem}
	
	The rate function $I_\delta$ admits an intuitive interpretation: it can be viewed as a convex combination of the rate functions of the simple random walks on the subgraphs $(V_k,E_k)$ for $k=1,2,\dots,\mathfrak{b}$, where the coefficients $\{r_k\}_{1\le k\le \mathfrak{b}}$ represent the proportion of time spent by $(X_n)_{n\ge 0}$ on $(V_k, E_k)_{1\le k\le \mathfrak{b}}$ 
	before the edge cover time $C_E$.
	
	The proportions $\{r_k\}_{1\le k\le \mathfrak{b}}$ and the relative entropy  $\{R(q_k\|p_{E_k})\}_{1\le k\le \mathfrak{b}}$ depend intricately on the structure of the sequence $(V_k,E_k)_{1\le k\le \mathfrak{b}}$ as well as the parameter $\delta$. Hence characterizing the minimizer of the variational representation of $I_\delta$ in (\ref{simpler expression of I})  is highly nontrivial.
	
	Although an explicit expression for $I_\delta$ remains elusive, the variational representation of $I_\delta$ in \eqref{simpler expression of I} still enables us to derive several interesting analytic properties, which are 
	stated in Theorem \ref{phase for rate function} below. In Section 4, we will
	obtain simpler and more explicit formulas for $I_\delta$ on certain specific graphs.

	\begin{theorem}[Analytic property of $I_\delta$]\label{phase for rate function}
		The rate function $I_\delta$ given by \eqref{simpler expression of I} satisfies the following properties:
		\begin{itemize}
			\item[(a)]  For any $\delta\in (0,1)$, $I_\delta(\mu)\equiv I_1(\mu)$ on $\mathscr{P}(V)$;
			\item[(b)]  $I_\delta$ is monotonically decreasing in $\delta\ge 1$, i.e., for all  $\delta_2>\delta_1\ge 1$,
			$$
			I_{\delta_1}(\mu)\ge I_{\delta_2}(\mu),\ \forall\,\mu\in \mathscr{P}(V);
			$$
			and  there exists some $\mu_0\in \mathscr{P}(V)$ (depending on $\delta_1$ and $\delta_2$) such that the inequality is strict: $I_{\delta_1}(\mu_0)> I_{\delta_2}(\mu_0)$.
			\item[(c)]   $I_\delta$ is uniformly continuous in $\delta$ in the sense that
			$$
			\lim\limits_{|\delta_1-\delta_2|\to 0}\sup\limits_{\mu\in
				\mathscr{P}(V):\, I_1(\mu)<\infty}|I_{\delta_1}(\mu)-I_{\delta_2}(\mu)|=0.
			$$
			\item[(d)] For some $\mu\in\mathscr{P}(V)$, the function $\delta\mapsto I_\delta(\mu)$ is not differentiable at $\delta=1$. % This is part \emph{(d)} of Theorem \ref{phase for rate function}.
		\end{itemize}
		
	\end{theorem} 
	
	\begin{remark}
		We conjecture that $I_\delta$ is differentiable and convex in $\delta>1$.
	\end{remark} 
	
	%%%%%%%%%%%%%%%%%%%%%%%%%%
	%%%%%%%%%%%%%%%%%%%%%%%%%%
	%%%%%%%%%%%%%%%%%%%%%%%%%%
	%%%%%%%%%%%%%%%%%%%%%%%%%% 
	
}

\subsection{Lifted process on line digraph}\label{sec 2.15}
\noindent{\color{black}  If $(X_n)_{n\ge 0}$ is a Markov chain, then $(X_n, L^n)_{n\ge 0}$ is also a Markov process with respect to the filtration generated by $X$ (see Proposition \ref{X-L-Mrk}). This property is crucial in the weak convergence approach to LDPs for empirical measures of Markov chains (see detailed discussion in Section  \ref{sec 2.2}). {\color{black} The Markov property of $\left(X_j,\frac{j}{n}L^j\right)_{0\le j\le n}$  is also noted in \cite[Section 4.2]{DE1997}, though no proof is provided there. }

	In the case of the ORRW $(X_n)_{n\ge 0}$, for two different historical paths,   the pair  $(X_n, L^n)$  may occupy the same vertex and share the same empirical measure at time $n$,   yet the transition probabilities from $(X_n, L^n)$ to $(X_{n+1}, L^{n+1})$ may differ.  Hence,   $(X_n, L^n)_{n\ge 0}$ is not a Markov process with respect to the filtration generated by $X$ in general,  
	see Example \ref{ORRW-non_Mrv}  in Section \ref{5-pf sec 2.2}.
	
	To overcome this difficulty, we introduce a new technique:  we lift $X$ from $G$ to  its line digraph $\vec{l}(G)$ (for the original definition of line digraph given by Harary and Norman, refer to \cite{HN1960}), obtaining a vertex-reinforced random walk $\mathcal{Z}:=(\mathcal{Z}_n)_{n\ge 1}$ such that the two-component process $(\mathcal{Z}_n, \mathcal{L}^n)_{n\ge 1}$ ---\,where $\mathcal{L}:=(\mathcal{L}^n)_{n\ge 1}$ is the empirical process of $\mathcal{Z}$ ---\,becomes a Markov process with respect to the filtration generated by $\mathcal{Z}$, see Proposition \ref{Z-L-Mrk} below. This construction makes the weak convergence approach  applicable. 
	\vskip2mm

	We now introduce the line digraph. To conveniently describe the relationship between the empirical measure on directed edges and the empirical measure on vertices of the graph $G$, we adopt the following notation.
	
	\vskip1mm
	
	Let $\overrightarrow{uv}$ denote the directed edge from $u$ to $v$.  For any subset $E'\subseteq E$, define
	\begin{equation}
		\vec{E'}:=\{ \overrightarrow{uv}: u,v\in V, uv\in E' \}. \label{def-directed-edge-set}
	\end{equation}

	%Set a line digraph $\vec{l}(G)$\label{S} associated with $G$ as follows:
	
	\begin{definition}[Line digraph associated with $G$] \label{line digraph}
		
		For a finite connected graph $G=(V,E)$, let $\vec{G}=(V,\vec{E})$ be the directed graph obtained by assigning two opposite directions to each edge of $G$. 
		The line digraph  of  $G$ is $\vec{l}(G)=(V_{\vec{l}(G)},E_{\vec{l}(G)})$\label{S} with
		\begin{itemize}
			\item[(a)] $V_{\vec{l}(G)}:=\vec{E}$.
			\item[(b)] For $\overrightarrow{u_1v_1},\overrightarrow{u_2v_2}\in V_{\vec{l}(G)}$, we write ${\overrightarrow{u_1v_1}\to\overrightarrow{u_2v_2}}$\label{z_1 to z_2}, if $v_1=u_2$. In this case a directed edge is formed from $\overrightarrow{u_1v_1}$ to $\overrightarrow{u_2v_2}$, denoted by $\overrightarrow{\overrightarrow{u_1v_1}\overrightarrow{u_2v_2}}$.  $E_{\vec{l}(G)}$ is the set of these directed edges.
		\end{itemize}
	\end{definition}
	
	We will also use the following concepts related to $\vec{l}(G)$. If $\mathbf{z}:=\overrightarrow{uv}\in V_{\vec{l}(G)}$, we define the head $\mathbf{z}^+$ and the tail $\mathbf{z}^-$ by 
	\begin{equation}
		\mathbf{z}^+=v,\ \mathbf{z}^-=u.\label{head_tail}
	\end{equation}
	For $A\subseteq \vec{E}= V_{\vec{l}(G)}$,  we define
	$\partial A$ \label{partial A} to be the set of all vertices $\mathbf{z}'\in A$ such that there exists some $\mathbf{z}\in V_{\vec{l}(G)}\setminus A$ with $\mathbf{z}'\to \mathbf{z}$.  The set $\partial A$ is the boundary vertex set of $A$.
	
	\vskip1mm
	\begin{definition}[Lifted process of $X$] The process $\mathcal{Z}:=(\mathcal{Z}_n)_{n\ge 0}$  defined by
		\begin{equation}\label{Z-n}
			\mathcal{Z}_n:=\overrightarrow{X_nX_{n+1}},  \ n\ge 0,
		\end{equation}
		is called the lifted process of $(X_n)_{n\ge 0}$ on the line digraph $\vec{l}(G)$.
	\end{definition}
	It follows from 
	\eqref{head_tail},  for every $n\geq 0$,
	\begin{equation}
		\mathcal{Z}_n^-=X_{n}\ \text{ and }\ \mathcal{Z}_n^+=X_{n+1}. 
	\end{equation}
	Note that
	\eqref{head_tail} induces 
	the following mapping from $\mathscr{P}(V_{\vec{l}(G)})$ to $\mathscr{P}(V)$.
	\begin{definition} Define $T: \mathscr{P}(V_{\vec{l}(G)}) \mapsto \mathscr{P}(V)$ by 
		\begin{equation}
			T(\hat{\mu})(V')=\hat{\mu}(\mathbf{z}:\mathbf{z}^-\in V'),  \ \ \hat{\mu}\in \mathscr{P}(V_{\vec{l}(G)}), \ V'\subset V.
		\end{equation}
	\end{definition}
	\begin{proposition}\label{continuity of map T}
		The mapping  $T: \mathscr{P}(V_{\vec{l}(G)}) \mapsto \mathscr{P}(V)$ is continuous with respect to  the weak convergence topology.
	\end{proposition}
	The proof of this Proposition 
	%RS is provided 
	will be given
	in Section \ref{5-pf sec 2.2}. \vskip2mm

	\eqref{head_tail} also induces a natural  mapping  from $V_{\vec{l}(G)}$ to the undirected edge set $E$,  as well as a mapping  from $\mathscr{P}(V_{\vec{l}(G)})$ to $\mathscr{P}(E)$, see the definition below.
	\begin{definition}\label{def-z|E-mu|E}\ \  For $\mathbf{z}\in V_{\vec{l}(G)}$, define its edge projection by 
		\begin{equation}
			\mathbf{z}|_E:=\mathbf{z}^-\mathbf{z}^+. \label{z|_E}
		\end{equation}
		For a probability measure $\hat{\mu}\in \mathscr{P}(V_{\vec{l}(G)})$, define its edge-projected measure $\hat{\mu}|_E\in \mathscr{P}(E)$ by 
		\begin{equation} \label{edge mu}
			\hat{\mu}|_E(E'):=\hat{\mu}(\{ \mathbf{z}:\mathbf{z}|_E\in E' \}), \ \ E'\subset E.
		\end{equation} 
	\end{definition}

	\vskip1mm
	
	We can now give the transition probabilities of $\mathcal{Z}$ using that of $X$.
	\begin{definition}[Transition probability on $\vec{l}(G)$]\label{Def for p_u}  \label{pp_E}
		For every edge set $E'\subset E$, define
		\begin{align}
			\pp_{E'}(\mathbf{z}_1,\mathbf{z}_2;\delta):=
			\left\{
			\begin{aligned}
				&\frac{{\bm 1}_{\{\mathbf{z}_2|_E\notin E' \}}+\delta {\bm 1}_{\{\mathbf{z}_2|_E\in E' \}}}{\sum_{\mathbf{z}\leftarrow \mathbf{z}_1}{\bm 1}_{\{\mathbf{z}|_E\notin E' \}}+\delta {\bm 1}_{\{\mathbf{z}|_E\in E' \}}}, \ &\mathbf{z}_1\to \mathbf{z}_2,\\
				&0, \ &{\rm otherwise,}
			\end{aligned}\right.\label{p_E}
		\end{align}
		where the sum runs over all $\mathbf{z}\in \vec{l}(G)$ with $\mathbf{z}\leftarrow \mathbf{z}_1$.
		
		For every measure $\hat{\mu}$ on $V_{\vec{l}(G)}$ satistying $\text{supp}(\hat{\mu}|_E)=E'$, set 
		\begin{equation}\label{pp_mu}
			\pp_{\hat{\mu}}(\mathbf{z}_1,\mathbf{z}_2;\delta):=\pp_{E'}(\mathbf{z}_1,\mathbf{z}_2;\delta),
		\end{equation} or equivalently 
		\begin{align}
			\pp_{\hat{\mu}}(\mathbf{z}_1,\mathbf{z}_2;\delta):=
			\left\{
			\begin{aligned}
				&\frac{{\bm 1}_{\{\hat{\mu}|_E(\mathbf{z}_2|_E)=0 \}}+\delta {\bm 1}_{\{\hat{\mu}|_E(\mathbf{z}_2|_E)>0 \}}}{\sum_{\mathbf{z}\leftarrow \mathbf{z}_1}{\bm 1}_{\{\hat{\mu}|_E(\mathbf{z}|_E)=0 \}}+\delta {\bm 1}_{\{\hat{\mu}|_E(\mathbf{z}|_E)>0 \}}}, \ &\mathbf{z}_1\to \mathbf{z}_2,\\
				&0, \ &{\rm otherwise.}
			\end{aligned}\right.\label{p_u}
		\end{align}
		
		When the parameter $\delta$ is fixed, we often write
		$\pp_{E'}(\mathbf{z}_1,\mathbf{z}_2;\delta)$ and $\pp_{\hat{\mu}}(\mathbf{z}_1,\mathbf{z}_2;\delta)$ simply as $\pp_{E'}(\mathbf{z}_1,\mathbf{z}_2)$ and $\pp_{\hat{\mu}}(\mathbf{z}_1,\mathbf{z}_2)$, respectively.\end{definition}
	
	{\color{black}
		Similar to the transition probability $p_{E'}$ for any edge subset $E'$,  $\pp_{E'}$ gives the transition probability of the lifted process $\mathcal{Z}$ when the edge projection (defined in \eqref{z|_E}) of the traversed vertex set of $\mathcal{Z}$ equals $E'$.
	}
	
	\begin{definition}\label{mathscr{T}_S}
		$\mathscr{T}_{\vec{l}(G)}$ is defined to be the collection of all transition probabilities $\qq:\vec{l}(G)\times\vec{l}(G)\mapsto[0,1]$ such that
		\[
		\sum_{\zz'\in\vec{l}(G)}\qq(\zz,\zz')=1 \text{ for all }\zz\in 
		%RS \vec{l}(G).
		\vec{l}(G),
		\]
		%RS Moreover, 
		and that, 
		for all  $\zz,\zz'\in V_{\vec{l}(G)}$, $\qq(\zz,\zz')\neq 0$ only if $\zz \to \zz'$. That is, 
		%RS $\qq$ is 
		$\mathscr{T}_{\vec{l}(G)}$ is the collection of all
		transition probabilities of Markov chains on the vertex set $V_{\vec{l}(G)}$  that  move from the current vertex to a neighboring vertex along the directed edge at each step.
	\end{definition}	  
	
	\vskip1mm
	\begin{definition} [Transition probability of $\mathcal{Z}$]   
		For any fixed sample path $\omega$, let	
		$$
		{E}_n^\omega=\left\{e\in E:\ \exists\, m\le n,\ \mathcal{Z}_m|_E(\omega)=e\right\}
		$$
		be the set of edges traversed by $X$ up to time $n\geq 1$. Then, for every $n\geq 1$,  the  transition probability  from $\mathcal{Z}_n$ to $\mathcal{Z}_{n+1}$ is given by  $\pp_{{E}_n^\omega}$;  that is,
		\begin{equation}
			\P(\mathcal{Z}_{n+1}=\zz_2\,|\,\mathcal{Z}_n=\zz_1,\mathscr{F}_n^{\mathcal{Z}})(\omega)=\left\{\begin{aligned}
				&\pp_{{E}_n^\omega}(\zz_1,\zz_2), &\zz_1\to \zz_2,\\
				&0, &\text{otherwise},
			\end{aligned}\right.\label{transition probability}
		\end{equation}
		where $(\mathscr{F}_n^{\mathcal{Z}})_{n\ge 1}$ is the natural filtration 
		%RS generated by the history of $\mathcal{Z}$ up to time $n$.  
		of $\mathcal{Z}$. 
		$\mathscr{F}_n^\mathcal{Z}$ is 
		equal to $\mathscr{F}_{n+1}$ generated by the history of $X$ up to time $n+1$. 
	\end{definition}
	\vskip2mm
	$(\mathcal{Z}_n)_{n\geq 0}$ can be regarded as a vertex once-reinforced random walk on the line digraph $\vec{l}(G)$ with the following property: for every pair of vertices $u,v\in V$ with $u\sim v$, the weights of the vertices  $\overrightarrow{uv}$ and $\overrightarrow{vu}$ are reinforced simultaneously whenever either of them is  traversed (see Figure \ref{state space 1} for a simple illustration).

	\begin{figure}
		\centering\includegraphics[scale=1]{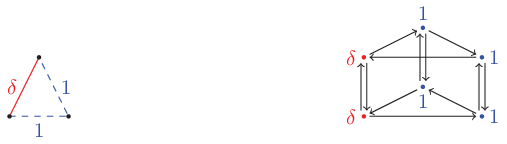}
		\caption{{\small If the red edge is traversed, then the weight on this edge becomes $\delta$. This means that one of the red point in the figure on the right hand side has been traversed, and that the weights on both red points becomes $\delta$.  }} \label{state space 1}
	\end{figure}

	\begin{definition}[Empirical measure of $\mathcal{Z}$]  The empirical measure of $\mathcal{Z}$ is defined by
		\begin{eqnarray}\label{eq-empirical-lifted}	
			\mathcal{L}^n(A)=\frac{1}{n}\sum_{i=0}^{n-1}\dd_{\mathcal{Z}_i}(A),\ A\subseteq V_{\vec{l}(G)},\ n\geq 1,
		\end{eqnarray}
		where $\dd_\zz$ is the Dirac measure concentrated on $V_{\vec{l}(G)}$ at $\zz$. By convention, $\mathcal{L}^0=0$.
	\end{definition}
	
	\vskip1mm
	\begin{proposition}\label{Z-L-Mrk}
		The two-component process
		$(\mathcal{Z}_n,\mathcal{L}^n)_{n\geq 0}$ is a non-homogeneous Markov process with respect to the  filtration  $\big(\mathscr{F}^{\mathcal{Z}}_n\big)_{n\ge 0}$.  Its the transition probability is given by
		\begin{equation}\label{Z_n,L^n}
			\P\left(\left. \mathcal{Z}_{n+1}=\zz',\mathcal{L}^{n+1}=\hat{\nu}\,\right|\,\mathcal{Z}_n=\zz,\mathcal{L}^n=\hat{\mu}\right)=\pp_{\hat{\mu}+\dd_\zz}\left(\zz,\zz'\right)\cdot \dd_{\frac{n\hat{\mu}+\dd_\zz}{n+1}}\left(\hat{\nu}\right),\ n\geq 0,
		\end{equation}
		where $\dd_{\hat{\nu}}$ denotes the Dirac measure on the space $\mathscr{P}(V_{\vec{l}(G)})$\label{mathscr{P}(T)} concentrated at $\hat{\nu}$. 
	\end{proposition}

	\vskip1mm
	
	The proof of Proposition \ref{Z-L-Mrk} is presented in Section \ref{5-pf sec 2.2}.
	\vskip1mm
	\begin{definition}[Renewal time of $\mathcal{Z}$]
		The times at which  $\mathcal{Z}$ traverses a previously unvisited edge, defined by ${\tauu}_1:=0$ and
		\begin{equation}
			{\tauu}_k=\inf\left\{j>{\tauu}_{k-1}:\ Z_j|_E\notin \{Z_i|_E:i<j \} \right\}.\label{def of renewal time Z}
		\end{equation}
		are called the renewal times of $\mathcal{Z}$.  
	\end{definition}
	
	For any $k\ge 1$, $\tauu_k$ is a stopping time with respect to $(\mathscr{F}_n^\mathcal{Z})_{n\ge 0}$.
	The renewal times $\tauu_k$ are related to $\tau_k$ defined in Definition \ref{def of renewal time} via
	$$
	\tauu_k=\tau_k-1 \text{ for }1\le k\le \b.
	$$

	The decomposition of $X$ through the renewal times $\tau_k\   (k=1,\dots,\b)$ induces a 
	%RS corresponding 
	decomposition of  $\mathcal{Z}$ via $\tauu_k\  (k=1,\dots, \b)$. Specifically, for a fixed 
	$m\in[1,\b]$ and conditioned on $\mathscr{F}^\mathcal{Z}_{\tauu_m}$, the segment  $(\mathcal{Z}_n)_{n\in[\tauu_m,\tauu_{m+1}-1]}$ forms a time-homogeneous Markov chain  with the transition probability $\pp_{E_m}$, 
	$$
	\pp_{E_m}(\zz_1,\zz_2;\delta)=\left\{\begin{aligned}
		&\frac{\delta}{\sum_{\mathbf{z}\leftarrow \mathbf{z}_1}{\bm 1}_{\{\zz|_E\notin E_m \}}+\delta {\bm 1}_{\{\zz|_E\in E_m \}}}, \ &\mathbf{z}_1\to \mathbf{z}_2,\, \zz_2|_E\in E_m\\
		&\frac{1}{\sum_{\mathbf{z}\leftarrow \mathbf{z}_1}{\bm 1}_{\{\zz|_E\notin E_m \}}+\delta {\bm 1}_{\{\zz|_E\in E_m \}}}, \ &\mathbf{z}_1\to \mathbf{z}_2,\, \zz_2|_E\notin E_m,\\
		&0, \ &{\rm otherwise.}
	\end{aligned}
	\right.
	$$

	\vskip1mm
	When we prove  our Laplace principle, 
	we need to modify the method presented in  \cite[Chapter 8]{DE1997} to adapt to the situation of dividing the whole trajectory into several parts.

	%%%%%%%%%%%%%%%%%%%%%%%%%%
	%%%%%%%%%%%%%%%%%%%%%%%%%%
	%%%%%%%%%%%%%%%%%%%%%%%%%%
	%%%%%%%%%%%%%%%%%%%%%%%%%%
}

\subsection{Generalized Laplace principle}\label{sec 2.16}
\noindent 
In this subsection, we first present a generalized Laplace principle for the lifted process $(\mathcal{Z}_n, \mathcal{L}^n)_{n\ge 1}$. Then, as applications of the generalized Laplace principle, we get  
a large deviation principle for the empirical measure and provide exponential decay estimates for tail probabilities of certain stopping times, such as the cover time.
Even for the case $\delta=1$, (i.e., the Markov chain setting), to the best of our knowledge, these tail probability estimates are new. 
Using our exponential decay estimates for the tail probabilities, we can precisely characterize the way the exponential integrability of these quantities depends on the parameter  $\delta$ and on the structure of the underlying graph  (a detailed discussion can be found  in arXiv version of our  paper \cite{HLX2025}).
\vskip1mm

%The generalized Laplace principle considered here concerns the limiting form of a restricted logarithmic Laplace functional analogous to \eqref{(2.1)}. Informally, the restricted logarithmic Laplace functional is defined by taking the logarithmic Laplace functional of the empirical measure of $\mathcal{Z}$ restricted to a  given subgraph.

We first give a definition. 
Recall that $G':=(V',E')$ denotes the subgraph induced by the edge set $E'$ on the vertices incident to its edges.

\begin{definition}   
	%RS Let
	We define
	\begin{eqnarray}\label{eq-edge-subset}
		\mathscr{S}:=\Big\{E':\ \exists\,\text{connected}\ G'=(V',E')\subseteq G=(V,E)\ \text{and}\ \exists\,v\in V\ \text{with}\
		x_0v\in E'\Big\},
	\end{eqnarray}
	where $x_0$ is the starting vertex of $\delta$-ORRW $X$. 
	%RS Equivalently, 
	$\mathscr{S}$ consists of all edge subsets $E'$ of $E$ for which the induced subgraph $G'$ is connected and contains the starting vertex $x_0$. 
\end{definition}

\begin{definition}[Decreasing subset of $\mathscr{S}$]  A non-empty subset $\mathscr{S}_0\subset \mathscr{S}$ \label{mathscr{S}_0} is called {\bf decreasing} if for any $\hat{E},\hat{\hat{E}}\in\mathscr{S}$ with $\hat{E}\subseteq \hat{\hat{E}}$, the condition $\hat{\hat{E}}\in\mathscr{S}_0$ implies $\hat{E}\in\mathscr{S}_0$.
\end{definition}
$\mathscr{S}_0$ 
represents the family of edge sets that can be traversed by $X$  along some sample path starting from $x_0$ up to a given time. 

\vskip1mm
\begin{definition}
	Define
	\begin{eqnarray}\label{eq-closed-subset}
		\mathcal{C}l(\mathscr{S}_0):=\left\{\hat{\mu}\in\mathscr{P}(V_{\vec{l}(G)}):\ \text{supp}(\hat{\mu}|_E)\subseteq E'\ \text{for some}\ E'\in \mathscr{S}_0 \right\}.
	\end{eqnarray}
\end{definition}
$\mathcal{C}l(\mathscr{S}_0)$ is a closed subset of $\mathscr{P}(V_{\vec{l}(G)})$ under the weak convergence topology (see Lemma \ref{closed set under topology}).

\vskip2mm

We are now ready to introduce the restricted logarithmic Laplace functional.

\begin{definition} [Restricted logarithmic Laplace functional]
	For every bounded  continuous function $h$ {\color{blue} in $\mathcal{C}l(\mathscr{S}_0)$}, define
	\begin{equation}
		W^n_{\mathscr{S}_0}(\zz):=\left\{\begin{array}{ll}
			-\frac{1}{n}\log{{\E}_\zz\left\{\exp\left[-nh(\mathcal{L}^n)\right]\mathbf{1}_{\{\mathcal{L}^n\in\mathcal{C}l(\mathscr{S}_0)\}}\right\}},&    {\P}_\zz(\mathcal{L}^n\in\mathcal{C}l(\mathscr{S}_0))>0\\
			\infty,&  {\P}_\zz(\mathcal{L}^n\in\mathcal{C}l(\mathscr{S}_0))=0    \end{array}\right., \label{eq-W_S0^n}
	\end{equation}
	where ${\E}_\zz$ and ${\P}_\zz$ denote expectation and probability conditioned on $\mathcal{Z}=\zz$.
	The functional $W^n_{\mathscr{S}_0}(\zz)$ is called the restricted logarithmic Laplace functional.
\end{definition}

To give an explicit expression for the rate function of the generalized Laplace principle concerning
$W^n_{\mathscr{S}_0}(\zz)$, 
we need to introduce a few more concepts.

\begin{definition}
	Let $\mathscr{S}_0$ be an arbitrary non‑empty decreasing subset of $\mathscr{S}$. For any $\hat{\mu}\in\mathscr{P}(V_{\vec{l}(G)})$, we define $\hat{\mathscr{A}}(\hat{\mu},\mathscr{S}_0)$ as the set of all quadruples $$
	(\hat{\mu}_k,r_k,E_k,\pp_k)_{1\leq k\leq \b}
	$$
	satisfying the following conditions:
	\begin{itemize}
		\item[(1)] {\textbf{Edge-set sequences:}}    $\{E_k\}_{1\leq k\leq \mathfrak{b}}\in\mathscr{E}$;
		\item[(2)] {\textbf{Proportions:}} $r_k\ge 0$, $r_l=0$ for $E_l\notin \mathscr{S}_0$, $\sum_{k=1}^\b r_k=1$;
		\item[(3)] {\textbf{Transition probabilities and invariant  measures:}} For each $k$, 
		$$
		\qq_k\in \mathscr{T}_{\vec{l}(G)},\ \hat{\mu}_k\in  \mathscr{P}(V_{\vec{l}(G)}),\  {\rm{supp}}(\hat{\mu}_k)\subseteq \vec{E}_k, 
		$$
		and $\hat{\mu}_k\qq_k=\hat{\mu}_k$ (i.e., $\hat{\mu}_k$ is invariant with respect to $\qq_k$).
		\item[(4)] {\textbf{Decomposition of $\hat{\mu}$:}} $\hat{\mu}:=\sum_{k=1}^\b r_k\hat{\mu}_k$.
	\end{itemize}
	Equivalently, 
	\begin{align}
		\hat{\mathscr{A}}(\hat{\mu},\mathscr{S}_0):=\Bigg\{& (\hat{\mu}_k,r_k,E_k,\qq_k)_{1\leq k\leq \b}:\ \{E_j\}_{1\leq j\leq \b}\in\mathscr{E},\nonumber\\ &\hskip 5mm r_k\ge 0,\ r_l=0\ \text{for}\ E_l\notin\mathscr{S}_0,\ \sum_{k=1}^{\b}r_k=1,\nonumber\\
		&\hskip 5mm  \qq_k\in \mathscr{T}_{\vec{l}(G)},\ \hat{\mu}_k\in  \mathscr{P}(V_{\vec{l}(G)}),\  {\rm{supp}}(\hat{\mu}_k)\subseteq \vec{E}_k \text{ \ such that \ } \hat{\mu}_k\qq_k=\hat{\mu}_k\Bigg\}. \label{A}
	\end{align}   
	
	We define the rate function 
	\begin{align}
		{\Lambda}_{\delta,\mathscr{S}_0}(\hat{\mu}):=\inf_{(\hat{\mu}_k,r_k,E_k,\qq_k)_{k}\in\hat{\mathscr{A}}(\hat{\mu},\mathscr{S}_0)}
		\sum_{k=1}^{\b}r_k\int_{V_{\vec{l}(G)}}R(\qq_{k}\|\pp_{E_k})\ {\rm d}\hat{\mu}_{k}.\label{rate function}
	\end{align}
	In the special case $\mathscr{S}_0=\mathscr{S}$, we denote $\hat{\mathscr{A}}(\hat{\mu}, \mathscr{S})$ simply  by  $\hat{\mathscr{A}}(\hat{\mu})$ and $\Lambda_\delta(\cdot)$ by $\Lambda_{\delta,\mathscr{S}}(\cdot)$.
\end{definition}
\vskip2mm
For  $\zz\in V_{\vec{l}(G)}$, we define  \begin{equation}\label{eq-renew-subset}
	\mathscr{E}_\zz:=\{\{E_k\}_{1\le k\le \b}\in\mathscr{E}:\ E_1=\{\zz|_E\} \}.
\end{equation} 
$\mathscr{E}_\zz$  consists of all sequences of edge subsets (belonging to $\mathscr{E}$) whose first element is exactly the single edge $\zz|_E$. 

\begin{definition} For $\zz\in V_{\vec{l}(G)}$,  we define   
	\begin{equation}
		\hat{\mathscr{A}}_\zz(\hat{\mu},\mathscr{S}_0):=\left\{(\hat{\mu}_{k},r_k,E_k,\qq_k)_{1\leq k\leq \b}\in\hat{\mathscr{A}}(\hat{\mu},\mathscr{S}_0):\ \{E_k\}_{1\le k\le \b}\in \mathscr{E}_\zz \right\},
		\label{A_z}
	\end{equation}
	i.e., the subset of $\hat{\mathscr{A}}(\hat{\mu},\mathscr{S}_0)$, in which the sequence of edge sets starts with the single edge $\zz|_E$. We also define 
	\begin{align}
		{\Lambda}^\zz_{\delta,\mathscr{S}_0}(\hat{\mu}):=\inf_{(\hat{\mu}_k,r_k,E_k,\qq_k)_{k}\in\hat{\mathscr{A}}_\zz(\hat{\mu},\mathscr{S}_0)}
		\sum_{k=1}^{\b}r_k\int_{V_{\vec{l}(G)}}R(\qq_{k}\|\pp_{E_k})\ {\rm d}\hat{\mu}_{k}.\label{Lambda^z}
	\end{align}
\end{definition}

We now present state a generalized Laplace principle for $W^n_{\mathscr{S}_0}(\zz)$, which a central result of this article.
\begin{theorem}[Generalized Laplace principle] \label{estimate for rate}
	For every  bounded continuous function $h$ on $\mathscr{P}(V_{\vec{l}(G)})$,
	\begin{equation}
		\lim_{n\to\infty}W^n_{\mathscr{S}_0}(\zz) = \inf_{\hat{\mu}\in \mathcal{C}l(\mathscr{S}_0)}\left\{\Lambda^\zz_{\delta,\mathscr{S}_0} (\hat{\mu})+h(\hat{\mu})\right\},\label{our variational representation}
	\end{equation}
\end{theorem}	

The proof of Theorem \ref {estimate for rate} will be given in Section \ref{sec 3.4}. We now give some properties of the rate function ${\Lambda}_{\delta,\mathscr{S}_0}^\zz$.

\begin{proposition}\label{prop domain of rate function}
	For all $\delta>0$,
	${\Lambda}_{\delta,\mathscr{S}_0}^\zz(\hat{\mu})<\infty$ if and only if $\hat{\mu}\in\mathcal{C}l(\mathscr{S}_0)$, and is invariant for some $\qq\in \mathscr{T}_{\vec{l}(G)}$.
	Further, ${\Lambda}_{\delta,\mathscr{S}_0}(\hat{\mu})<\infty$ if and only if $\hat{\mu}\in\mathcal{C}l(\mathscr{S}_0)$, and is invariant for some $\qq\in \mathscr{T}_{\vec{l}(G)}$.
\end{proposition}

\vskip1mm
\begin{proposition}\label{lower semicontinuous}
	${\Lambda}_{\delta,\mathscr{S}_0}$ (resp. ${\Lambda}^\zz_{\delta,\mathscr{S}_0}$) is lower semicontinuous,  and $\left\{\hat{\mu}:\ {\Lambda}_{\delta,\mathscr{S}_0}(\hat{\mu})\le M  \right\}$ (resp. $\big\{\hat{\mu}:\ {\Lambda}^\zz_{\delta,\mathscr{S}_0}(\hat{\mu})\le M \big\}$) is compact for all $0<M<\infty$.
\end{proposition}

\vskip2mm
The proofs of the above two propositions will be given in Section \ref{pf of sec 2.16}. 
Now, we give three applications of the generalized Laplace principle. 	\\

\noindent{\bf(i){\emph{\ The first application of the generalized Laplace principle}}}	

Using Theorem \ref{estimate for rate}, we can get the following LDP for the empirical measures of $(\mathcal{Z}_n)_{n\geq 0}$ under $\P_{x_0}$ on $\vec{l}(G)$.

\begin{theorem}\label{ldp for Z}
	The empirical measure process $\left(\mathcal{L}^n\right)_{n\geq 0}$ of $(\mathcal{Z}_n)_{n\geq 0}$ with reinforcement factor $\delta>0$ satisfies Laplace principlr with the rate function $\Lambda_\delta$; that is
	\[
	\lim_{n\to\infty}\frac{1}{n}\log \E_{x_0}\exp\{-nh(\mathcal{L}^n)\}=-\inf_{\hat{\mu}\in \mathcal{C}l(\mathscr{S})}\{\Lambda_\delta(\hat{\mu})+h(\hat{\mu})\},\ h\in C_b\big(\mathscr{P}(V_{\vec{l}(G)})\big).
	\]
	$\Lambda_\delta$  is also the rate function of the LDP of $\left(\mathcal{L}^n\right)_{n\geq 1}$.
\end{theorem}

The proof of Theorem \ref{ldp for Z} will be provided in Section \ref{sec 3.4}.\vskip2mm

It follows Proposition~\ref{continuity of map T} that the mapping	 $T: \mathscr{P}(V_{\vec{l}(G)}) \mapsto \mathscr{P}(V)$ is continuous. Now applying 
the contraction principle, we immediately obtain a large deviation principle for the empirical measures of $X$.

\begin{theorem}\label{rate function 0}
	The empirical measure process $(L^n)_{n\geq 1}$ of $\delta$-ORRW $(X_n)_{n\ge 0}$ on a finite connected graph $G=(V,E)$ satisfies an LDP  
	with good  rate function $I_\delta$ given by 
	\begin{equation}
		\label{rate func.I}
		I_{\delta}(\mu)=\inf\limits_{\hat{\mu}\in\mathscr{P}(V_{\vec{l}(G)}):\, T(\hat{\mu})=\mu}\Lambda_\delta(\hat{\mu}),\ \mu\in\mathscr{P}(V).
	\end{equation}
\end{theorem}

In Section \ref{sec 3.4}, we will use the results above to prove 
Theorem  \ref{I_delta} and  obtain a more concise expression for the rate function $I_\delta$ which depends only on the graph  $G$ and the process $X$.

Applying Theorem \ref{ldp for Z} and the contraction principle, and following an argument analogous to that of Theorem \ref{rate function 0}, we can also obtain the LDP  for the empirical measures defined on edges and for the pair empirical occupation measures of $X$. 
See Corollary  \ref{coro-ldp-edge}.

%This application constitutes the main content of our paper.
\vskip2mm

\noindent{\bf(ii){\emph{\ The second  application of the generalized Laplace principle}}}	

Applying Theorem \ref{estimate for rate}, we 
can get the following precise exponential tail estimates of the cover time
$C_E$,
\begin{equation}\label{exp-decay 1}
	\lim_{n\to\infty}\frac{1}{n}\log \mathbb{P}_{x_0}\left(C_E>n\right)=-\inf\limits_{\hat{\mu}\in Cl(\mathscr{S})}{\Lambda}_{\delta}(\hat{\mu})=:\alpha_c^1(\delta).
\end{equation}
Combining \eqref{exp-decay 1}  with detailed analysis,  
we can show that $\alpha_c^1(\delta)$ is the critical exponent for the exponential integrability of $C_E$, i.e.,
$$
\begin{array}{ll}
	\mathbb{E}_{x_0}\big(e^{\alpha C_E}\big)
	<\infty,&  \alpha< \alpha_c^1(\delta),\\
	\mathbb{E}_{x_0}\big(e^{\alpha C_E}\big)
	=\infty,& \alpha\ge \alpha_c^1(\delta).
\end{array}
$$

This implies that the asymptotic behaviors of $\alpha_c^1(\delta)$ as $\delta\to 0+$ and $\delta\to \infty$ depend on the structure of the graph $G$.

%Compare with  the exponential integrability result for  $C_E$ given in Proposition \ref{exponentially integrable} via the stochastic domination approach in Proposition \ref{domination}, the generalized Laplace principle reveals  finer information about ORRWs..  

Let
$$
\tau_{\mathscr{S}_0}:= \inf\big\{n\ge 1,  \{e\in E: \exists m<n,  X_mX_{m+1}=e\}\notin \mathscr{S}_0\big\},
$$
be  the stopping time at which $X$ leaves the  decreasing subset $\mathscr{S}_0$. We can also determine the 
the critical exponent for the exponential integrability of $\tau_{\mathscr{S}_0}$ through the generalized Laplace principle.

These results will be shown in our forthcoming paper  \cite{HLX2025}.
%%%%%%%%%%%%%%%%%%%%%%%%%%

\vskip2mm
\noindent{\bf(iii){\emph{\ The third application of the generalized Laplace principle}}}	

Taking $\delta=1$, we obtain a generalized LP for the  simple random walk  on $G$. 
\begin{theorem}\label{thm-rate-function-delta=1} 
	If $\delta=1$, the rate function $\Lambda_{1,\mathscr{S}_0}$ of the generalize Laplace principle can be expressed as
	\begin{align}
		{\Lambda}_{1,\mathscr{S}_0}(\hat{\mu})=\inf_{(\hat{\mu}, \{ E_k \}_{1\le k\le \b}, \qq)\in\hat{\mathscr{A}}(1,\hat{\mu},\mathscr{S}_0) }\int_{\vec{E}_l}R(\qq \| \pp){\rm d} \hat{\mu},\label{Lambda_1}
	\end{align}
	where 
	\begin{align*}
		\hat{\mathscr{A}}(1,\hat{\mu},\mathscr{S}_0):=&\ \Big\{ (\hat{\mu}, \{ E_k \}_{1\le k\le \b}, \qq): \  \{E_k\}_{1\le k\le \b}\in \mathscr{E}, \exists\, l,\  E_l\in \mathscr{S}_0,\  E_{l+1}\notin \mathscr{S}_0,\\ 
		&\ \ \ \ \ \ \ \ \ \ \ \ \ \ \ \ \ \ \ \  {\mathrm{supp}(\hat{\mu}|_E)\subset E_l,\  \hat{\mu}\qq=\hat{\mu}}\Big\}.
	\end{align*} 
\end{theorem}

This generalized Laplace principle  not only implies  the LDP for the empirical measures of Markov chain, but also  yields the critical exponent for the exponential integrability of $C_E$ and $\tau_{\mathscr{S}_0}$.
To the best of our knowledge, this critical exponent is also new for the simple random walk  on a finite connected graph.  The proof of Theorem \ref{thm-rate-function-delta=1} is provided  in Section \ref{pf of sec 2.16}.

%%%%%%%%%%%%%%%%%%%%%%%%%%%%%%%%%%%%%%%%%
%%%%%%%%%%%%%%%%%%%%%%%%%%%%%%%%%%%%%%%%%
%%%%%%%%%%%%%%%%%%%%%%%%%%%%%%%%%%%%%%%%%
%%%%%%%%%%%%%%%%%%%%%%%%%%%%%%%%%%%%%%%%%

\subsection{Brief introduction to weak convergence approach to homogeneous Markov chains}\label{sec 2.2}
\noindent  In this subsection, we 
%RS briefly outline 
sketch
the basic idea of the weak convergence approach developed by Dupuis and Ellis \cite{DE1997} for establishing LDP of the empirical measures of Markov chains. For simplicity, we restrict our discussion to Markov chains on a finite state space 
$\mathcal{S}$.

\begin{theorem}[\mbox{\cite[a corollary of Propositions 8.3.3 and 8.6.1]{DE1997}}] \label{thm-ldp-markov}
	For an irreducible Markov chain $X:=(X_n)_{n\ge 0}$ 
	on a finite state space $\mathcal{S}$ with transition probibility $p(x,{\rm d}y):=p(x,y)\dd_y$, the corresponding empirical measure $(L^n)_{n\ge 1}$ satisfies the LDP with a good rate function $I(\cdot)$ given by
	$$
	I(\mu)=\inf_{q\in\mathscr{T}_{\mathcal{S}}: \mu q=\mu}\int_{S}R\big(q(x,\cdot)\|p(x,\cdot)\big){\ }
	\mu({\rm d}x).
	$$
\end{theorem}	
\vskip2mm
\noindent{\bf(i){\emph{\   Solving the dynamic programming equation with respect to the logarithmic Laplace functional}}}	

%RS Let  $\mathbb{E}_x$ denote the expectation conditioned on $X_0=x$,  in view of  \eqref{(2.1)}, define 
Define
$$
W^n(x):=-\frac{1}{n}\log{\mathbb{E}_x\{\exp[-nh(L^n)] \}},\ x\in \mathcal{S}.
$$ 
We first give a computable representation for $W^n(x)$. Recall that the LP is defined in Definition \ref{LP}.

\begin{itemize}
	\item[\bf (i-1)]  \emph{Construction of  dynamic programming equation.} 
	
	For $n\geq 1$ and $0\leq j\leq n,$ let $\mathscr{M}_{j/n}(G)$ be all measures $\nu$ on $\mathcal{S}$ with 
	%RS total mass 
	$\nu(\mathcal{S})=j/n$,  equipped  with the weak convergence topology. Define
	$$
	L_0^n:=0,\ L_{j+1}^n:=L_j^n+\frac{1}{n}\dd_{X_j}, \ j=0,1,\dots,{n-1}, 
	$$
	so that $L^n=L_n^n$. The two-component process $(X_n, L^n)_{n\ge 0}$ forms a Markov Chain.  
	Define 
	$$
	W^n(i,x,\mu):=-\frac{1}{n}\log{\E_{i,x,\mu}\{\exp[-nh(L^n)] \}},
	$$ where $\E_{i,x,\mu}$ is the expectation conditioned on $(X_i,L^n_i)=(x,\mu)$ with $\mu\in \mathscr{M}_{i/n}$. 
	It follows from 
	\cite[Section 4.2]{DE1997} 
	%RS showed 
	that $W^n(i,x,\mu)$ satisfies
	$$
	W^n(i,x,\mu)=-\frac{1}{n}\log \int_\mathcal{S}\exp\left[-nW^n\left(i+1,y,\mu+\frac{1}{n}\dd_x\right)\right]\ p(x,{\rm d}y).
	$$
	
	Now applying Proposition \ref{variational representation} then yields the following dynamic programming equation:
	\begin{equation}\label{2.2dpequ}
		\left\{
		\begin{aligned}
			W^n(i,x,\mu)=&\,\inf_{\nu\in\mathscr{P}(\mathcal{S})}\left\{\frac{1}{n}R(\nu(\cdot)\|p(x,\cdot))+\int_\mathcal{S}         W^n\left(i+1,y,\mu+\frac{1}{n}\dd_x\right)\ \nu({\rm d}y) \right\}\\
			W^n(0,x,0)=&\, W^n(x)\\
			W^n(n,x,\mu)=&\, h(\mu)
		\end{aligned}\right..		
	\end{equation}
	
	\item[\bf (i-2)] \emph{Solving the  dynamic programming equation.}
	
	For each $0\le j\le n-1$, let $\nuu_j^n:=\nuu_j^n({\rm d}y|x,\mu)$\label{nu_j^n_0}  denote a transition probability from $\mathcal{S}\times\mathscr{M}_{j/n}(\mathcal{S})$ to $\mathcal{S}$.  The sequence  $\{\nuu_j^n,j=0,\dots,n-1  \}$ is called an admissible control sequence. Given such an admissible control sequence, we define $\big(\overline{X}_j^n,\ \overline{L}_j^n\big)_{0\leq j\leq n}$ on a probability space $\left(\overline{\Omega},\overline{\mathscr{F}},\overline{\P}_x\right)$ inductively by 
	$$
	\overline{X}_0^n=x,\ \overline{L}_0^n=0,
	$$
	for  $j\in\{0,\dots,n-1 \}$, 
	$$
	\overline{L}_{j+1}^n=\overline{L}_j^n+\frac{1}{n}\dd_{\overline{X}_j^n}\in \mathscr{M}_{(j+1)/n}(\mathcal{S}),
	$$
	where $\overline{X}_{j+1}^n$ 
	follows the distribution $\nuu_j^n\big(\cdot\big|\overline{X}_j^n,\overline{L}_j^n\big)$.
	By \cite[Theorem 1.5.2]{DE1997}, the unique solution to the  dynamic programming equation  (\ref{2.2dpequ}) is 
	\begin{equation*}
		\mathcal{V}^n(i,x,\mu):=\, \inf_{\stackrel{\{\nuu_j^n\}:\nuu_j^n:\mathcal{S}\times\mathscr{M}_{j/n}(\mathcal{S})\mapsto\mathcal{S}}{\text{ for }0\le j\le n-1}}\overline{\E}_{i,x,\mu}\Big\{\frac{1}{n}\sum_{j=i}^{n-1}R\big(\nuu_j^n(\cdot\big|\overline{X}^n_{j},
		\overline{L}^n_{j})\big\|p(\overline{X}^n_{j},\cdot)\big)+h\big(\overline{L}^n_n\big) \Big\},
	\end{equation*}
	where $\overline{\E}_{i,x,{\mu}}$ stands for expectation with respect to $\big(\overline{X}_j^n,\ \overline{L}_j^n\big)_{i\leq j\leq n}$
	given $\big(\overline{X}_i^n,\overline{L}^n_i\big)=(x,\mu).$ Every $\vv^n(i,x, \mu)$ is called a minimal cost function. 
	Applying \cite[Theorem 4.4.2]{DE1997}, and setting $i=0$ and $\mu=0$, we get the following computable representation:
	\begin{equation}
		W^n(x)=\mathcal{V}^n(x)
		:=\inf_{\stackrel{\{\nuu_j^n\}:\nuu_j^n:\mathcal{S}\times\mathscr{M}_{j/n}(\mathcal{S})\mapsto\mathcal{S}}{\text{ for }0\le j\le n-1}}\overline{\E}_x\Big\{\frac{1}{n}\sum_{j=0}^{n-1}
		R\big(\nu_j^n(\cdot\big|\overline{X}^n_{j},\overline{L}^n_{j})\big\|p(\overline{X}^n_{j},\cdot)\big)+h\big(\overline{L}^n_n\big) \Big\}. \label{(2.5)}
	\end{equation}

\end{itemize}

\noindent{\bf(ii){\emph{\ The upper bound of  Laplace principle}}}

We want to proof
\begin{equation}\label{uppbdd-LP}
	\liminf\limits_{n\to\infty}W^n(x)=\liminf\limits_{n\to\infty}\vv^n(x)\ge\inf\limits_{\mu\in \mathscr{P}(V)}\{I(\mu)+h(\mu) \},
\end{equation}
where $I(\cdot)$ is the good rate function of the LDP for the empirical measures defined in Theorem \ref{thm-ldp-markov}. 

\vskip1mm
Note that for probability measures $\alpha$, $\beta$ and $\gamma$ on $\mathcal{S}$, 
\begin{equation}\label{entropy prdct}
	R(\beta\|\gamma)=R(\alpha\times \beta\|\alpha\times\gamma).
\end{equation}
For any probability measure $\alpha$
and transition probability $q$ on $\mathcal{S}$, we define
\begin{equation}\label{trans prob prdct}
	\alpha\otimes q(A\times B):=\int_{A\times B}\alpha({\rm d}x)q(x,y){\rm d}y=\int_A q(x, B)\alpha({\rm d}x),
\end{equation}
where $\alpha\times \beta$ and $\alpha\times\gamma$ are the product measures 
%RS on $\mathcal{S}\times\mathcal{S}$, and $q(x,y)$ is a transition probability on $\mathcal{S}$.
on $\mathcal{S}\times\mathcal{S}$.

Let $\varepsilon$ be a given arbitrary  positive number. 
For each $n\ge  1$,  choose an admissible control sequence $\{\nuu_j^n\}_{0\le j\le n-1}$ satisfying \begin{align}
	\vv^n(x)+\varepsilon\ge &\, \overline{\E}_x\Bigg\{\frac{1}{n}\sum_{j=0}^{n-1}
	R\Big(\nuu_j^n\big(\cdot\big|\overline{X}^n_{j},\overline{L}^n_{j}\big)\Big\|p\big(\overline{X}^n_{j},\cdot\big)\Big)
	+h\big(\overline{L}^n_n\big) \Bigg\}\nonumber\\
	&\ \ \ \ \text{by}\ \eqref{entropy prdct} \ \ \text{and }\ \eqref{trans prob prdct} \nonumber\\
	= &\,      \overline{\E}_x\Bigg\{ \frac{1}{n} \sum_{j=0}^{n-1}R\Big(\dd_{\overline{X}^n_j}({\rm d}x)\times \nu^n_j\big({\rm d}y|\overline{X}^n_{j},\overline{L}^n_{j}\big)\Big\|\dd_{\overline{X}^n_j}({\rm d}x)\otimes p(x,y){\rm d}y\Big)+ h(L_n^n)\Bigg\} \nonumber  \\
	&\ \ \ \ \big(\text{applying Jensen's inequality to the convex function \  } R(\cdot\|\cdot) \ \ \nonumber \\
	&\ \ \ \ \text{ \ \ and using  the time homogeneity of\ } (\overline{X}^n_j)_{j=0}^{n-1} \big) \label{homo} \\   \ge  &\,  \overline{\E}_x\Bigg\{R\Big( \frac{1}{n} \sum_{j=0}^{n-1}\dd_{\overline{X}^n_j}({\rm d}x)\times \nu^n_j\big({\rm d}y|\overline{X}^n_{j},\overline{L}^n_{j}\big)\Big\| \frac{1}{n} \sum_{j=0}^{n-1}\dd_{\overline{X}^n_j}({\rm d}x)\otimes p(x,y){\rm d}y\Big)+ h(L_n^n)\Bigg\}   \nonumber    \\  
	=&\,  \overline{\E}_x\big\{R\big(\nuu^n\big\|\overline{L}^n_n\otimes p\big)+h(\overline{L}^n_n)\big\}, \nonumber   
\end{align}
where 
$$
\nuu^n:=\frac{1}{n}\sum_{j=0}^{n-1}\dd_{\overline{X}^n_{j}}\times {\nuu}_j^n\big(\cdot\big|\overline{X}^n_{j},\overline{L}^n_{j}\big).
$$ 
The asymptotic behaviour of  the pair $\big(\nuu^n, \overline{L}^n\big)_{n\ge 1}$ as $n\to \infty$ plays a key role in establishing the upper bound.

By \cite[Theorem 8.2.8]{DE1997}, there exist a subsequence of $\big(\nuu^{n'}, \overline{L}^{n'}\big)$, a pair $\left(\nuu,\overline{L}\right)$, and a transition probability $\overline{q}:$ $\mathcal{S}\mapsto \mathcal{S}$ such that
for $\overline{\P}_x$-almost surely $\overline{\omega}\in \overline{\Omega}$, 
$$
\big(\nuu^{n'}(\overline{\omega}),\overline{L}^{n'}(\overline{\omega})\big)\Rightarrow\left(\nuu(\overline{\omega}),\overline{L}(\overline{\omega})\right),\ \nuu(\overline{\omega})=\overline{L}(\overline{\omega})\otimes \overline{q}(\overline{\omega}),
$$ 
where "$\Rightarrow$"\label{Rightarrow}
denotes  weak convergence of probability measure (for each sample path $\overline{\omega}$). In particular, $\overline{L}$ is invariant under $\overline{q}$, i.e., $\overline{L}\overline{q}=\overline{L}$, where 
$$
\overline{L}\overline{q}(\cdot)=\int_\mathcal{S} \overline{L}({\rm d}x)\overline{q}(x,\cdot).
$$ 
Combining this weak convergence, Jensen’s inequality and Fatou’s lemma, we get
\begin{align*}
	%\label{}
	\liminf_{n\to \infty} W^{(n)}(x) +\varepsilon&=  \liminf_{n\to \infty} \vv^{(n)}(x) +\varepsilon   \\
	&  \ge  \overline{\E}_x \Big\{ \int_\mathcal{S} R\big( \overline{q}(x, \cdot)\big\| p(x, \cdot) \big)\overline{L}({\rm d}x)+ h(\overline{L})\Big\}\\
	& \ge \inf\limits_{\mu\in \mathscr{P}(\mathcal{S})}\{I(\mu)+h(\mu) \}.
\end{align*}
Since $\varepsilon>0$ is arbitrary, the upper bound follows. For complete details we refer to  \cite[Section 8.3]{DE1997}.  

\vskip2mm 

\noindent{\bf(iii){\emph{\ The lower bound of Laplace principle}}}

We want to prove 
$$
\limsup\limits_{n\to\infty}W^n(x)=\limsup\limits_{n\to\infty}\vv^n(x)\le\inf\limits_{\mu\in \mathscr{P}(V)}\{I(\mu)+h(\mu)\}.
$$
Given any $h\in C_b(\mathscr{P}(\mathcal{S}))$ and $\varepsilon>0$,  choose a probability measure $\gamma\in \mathscr{P}(\mathcal{S})$ such that
$$
I(\gamma)+h(\gamma)\le\inf_{\mu\in \mathscr{P}(\mathcal{S})}\{I(\mu)+h(\mu)\}+\varepsilon<\infty.
$$ 
We need to find a probability measure $\gamma^*\in \mathscr{P}(\mathcal{S})$ satisfying the following conditions:
\begin{itemize}
	\item[\bf{B1}]  There exists a transition probability $q^*(x, y)$ on $\mathcal{S}$ such that $\gamma^*$ is the unique invariant measure of $q^*$; 
	\item[\bf{B2}]  $\|\gamma-\gamma^*\|<\varepsilon$, where $\|\cdot\|$ is the total variation norm of a finite signed measure;  
	\item[\bf{B3}] {\color{black} \begin{equation}\label{S2.2-1}
			h(\gamma^*)\le h(\gamma)+\varepsilon,\ \  I(\gamma^*)\le \int_{\mathcal{S}}R(q^*\|p)\ {\rm d}\gamma^*\le I(\gamma)<\infty.
	\end{equation} }
\end{itemize} 

It follows from \cite[Lemma 8.6.3]{DE1997} that we may take
\begin{equation}\label{lambda*}
	\gamma^*:=\left(1-\frac{\varepsilon}{2}\right)\gamma+\frac{\varepsilon}{2} \mu^*,
\end{equation} 
where $\mu^*$  is the unique ergodic  measure $\mu^*$ of $p(x,y)$ and is supported on the whole state 
space.

Select the admissible control sequence $\big(\overline{X}^n_{j}\big)_{j=0}^{n-1}$ associated with $q^*$. Combining  \eqref{(2.5)}-\eqref{S2.2-1},  the ergodicity of $\big(\overline{X}^n_{j}\big)_{j=0}^{n-1}$, and the $L^1$ ergodic theorem \cite[Theorem A.4.4]{DE1997},  we obtain  \begin{equation}\label{lowbdd}
	\begin{aligned}
		\limsup_{n\to \infty} W^n(x)   &= \limsup_{n\to \infty} \vv^n(x)  \\
		&\le  \lim_{n\to \infty} \overline{\E}_x\Bigg\{\frac{1}{n}\sum\limits_{j=0}^{n-1}R\Big(q^*\big(\overline{X}^n_{j},\cdot\big)\Big\|
		p\big(\overline{X}^n_{j},\cdot\big)\Big)+h\big(\overline{L}^n_n\big) \Bigg\}\\
		&= \int_V R( q^*(x,\cdot)\|p(x,\cdot))\gamma^*({\rm d}x)+h(\gamma^*)\\
		&\le I(\gamma)+h(\gamma)+\varepsilon\\
		&\le \inf_{\mu\in \mathscr{P}(V)}\{I(\mu)+h(\mu)\}+2\varepsilon.
	\end{aligned}
\end{equation} 
Since $\varepsilon>0$ is arbitrary, the lower bound follows. For complete details, see \cite[Section 8.6]{DE1997}.
\begin{remark}
	In the lower bound estimate, the time homogeneity of the transition probabilities of $(\overline{X}^n_{j})_{j=0}^{n-1}$ is the key point that allows us to apply the $L^1$ ergodic theorem.
\end{remark}

%%%%%%%%%%%%%%%%%%%%%%%%%%%%%
%%%%%%%%%%%%%%%%%%%%%%%%%%%%%
%%%%%%%%%%%%%%%%%%%%%%%%%%%%%
%%%%%%%%%%%%%%%%%%%%%%%%%%%%%

\subsection{Modified weak convergence approach}\label{sec 2.3}
\noindent

\noindent {\color{black} 
	\noindent{\bf(i)\ {\emph{Dynamic programming for a restricted logarithmic Laplace functional}}}		
	
	\vskip1mm
	To prove Theorem \ref{estimate for rate}, we 
	first establish the dynamic programming for $(\mathcal{Z}_n)_{n\ge 0}$ on the line  digraph $\vec{l}(G)$.

	For all $n\geq 1$ and $0\leq j\leq n,$ we let $\mathscr{M}_{j/n}(\vec{l}(G))$ be the set of all measures $\hat{\mu}$ on $V_{\vec{l}(G)}$ with $\hat{\mu}(V_{\vec{l}(G)})=j/n$,
	and endow it with the weak convergence topology. 
	%RS Note an admissible control sequence $\{\nuu_j^n,j=0,\dots,n-1\}$ means that
	Note that $\{\nuu_j^n,j=0,\dots,n-1\}$ being an admissible control sequence means that
	\begin{equation}
		\mbox{each}\ \nuu_j^n:=\nuu_j^n\left(\left.{\rm d}\zz'\right|\zz,\hat{\mu}\right)\ \mbox{is a transition probability from}\ V_{\vec{l}(G)}\times\mathscr{M}_{j/n}(\vec{l}(G))\ \mbox{to}\ V_{\vec{l}(G)}.\label{def of nu_j^n}
	\end{equation}

	Given such an admissible control sequence and 
	%RS $\zz\in V_{\vec{l}(G)}$ arbitrarily. Define 
	$\zz\in V_{\vec{l}(G)}$, we define
	$\big(\overline{\mathcal{Z}}_j^n,\overline{\mathcal{L}}_j^n\big)_{0\leq j\leq n}$\label{overline{Z},overline{L}} 
	%RS by induction on some probability space $(\overline{\Omega},\overline{\mathscr{F}},\overline{\P}_\zz)$\label{overline{P}_z}: 
	on some probability space $(\overline{\Omega},\overline{\mathscr{F}},\overline{\P}_\zz)$\label{overline{P}_z} inductively as follows:
	%RS Set $\overline{\mathcal{Z}}_0^n=\zz,\ \overline{\mathcal{L}}_0^n=0$. For 
	$\overline{\mathcal{Z}}_0^n:=\zz,\ \overline{\mathcal{L}}_0^n:=0$, and for 
	all $j=0,1,\dots,n-1$, define
	\[
	\overline{\mathcal{L}}_{j+1}^n:=\overline{\mathcal{L}}_{j}^n+\frac{1}{n}\dd_{\overline{\mathcal{Z}}_j^n}\in\mathscr{M}_{(j+1)/n}(\vec{l}(G)),
	\]
	and let $\overline{\mathcal{Z}}_{j+1}^n$ 
	%RS be driven by 
	be distributed according to
	$\nuu_j^n\big(\cdot\big|\overline{\mathcal{Z}}_j^n,\overline{\mathcal{L}}_j^n\big)$ given $\big(\overline{\mathcal{Z}}_i^n,\ \overline{\mathcal{L}}_i^n\big)_{0\leq i\leq j}$. 
	%RS In the rest part of the paper, 
	In the remainder of this paper, we 
	write $\overline{\mathcal{L}}_n^n$ as $\overline{\mathcal{L}}^n$.

	Recall the transition probability $\pp_{\hat{\mu}}(\zz, \cdot)$ in Definition \ref{Def for p_u} and $\mathcal{C}l(\mathscr{S}_0)$ in \eqref{eq-closed-subset}. For 
	any measurable function $h:\mathscr{P}(V_{\vec{l}(G)})\mapsto \mathbb{R}$ which is bounded and continuous in $\mathcal{C}l(\mathscr{S}_0)$, we define the  minimal cost function by
	\begin{equation} 
		\vv^n_{\mathscr{S}_0}(\zz):=\inf_{\{\nuu_j^n\}_j:\,\overline{\mathcal{L}}^n\in\mathcal{C}l(\mathscr{S}_0)}
		\overline{\E}_\zz\bigg\{\frac{1}{n}\sum_{j=0}^{n-1}R\Big(\nuu_j^n\big(\cdot\big|\overline{\mathcal{Z}}^n_{j},   \overline{\mathcal{L}}^n_{j}\big)
		\big\|\pp_{\overline{\mathcal{L}}^n_{j+1}}\big(\overline{\mathcal{Z}}^n_{j},\cdot\big)\Big)+h\big(\overline{\mathcal{L}}^n\big) \bigg\},\label{(2.6)}
	\end{equation}
	where $\overline{\E}_\zz$\label{overline{E}_z} is the expectation under $\overline{\P}_\zz$.
	%RS Specifically, write $\vv^n_{\mathscr{S}}(\zz)$ as $\vv^n(\zz)$. 
	We will write $\vv^n_{\mathscr{S}}(\zz)$ as $\vv^n(\zz)$.  Since $\pp_{\hat{\mu}}\left(\zz,\cdot\right)$ is supported on those vertices in $V_{\vec{l}(G)}$  connected to $\zz$,  if $\nuu_j^n=\nuu_j^n\left(\cdot\left|\zz,\hat{\mu}'\right.\right)$ assigns positive probability to some vertex $\zz'$ which is not  adjacent to $\zz$, then the relative entropy
	\[
	R(\nuu_j^n(\cdot|\zz,\hat{\mu}')\|\pp_{\hat{\mu}'}(\zz,\cdot)) 
	\]
	becomes infinite; in that case the infimum in \eqref{(2.6)} can not be attained.  Therefore, when evaluating  $\mathcal{V}_{\mathscr{S}_0}^n(\zz)$, we only need to consider sequences $\nuu_j^n$ that are supported on vertices adjacent to $\zz$ for very measure $\hat{\mu}'$ on $V_{\vec{l}(G)}$.  Namely, under any such $\{\nuu_j^n\}_{0\le j\le n-1}$, the process $\big(\overline{\mathcal{Z}}_j^n\big)_{0\le j\le n}$ is a stochastic process moving to the neighbors at each step.
	
	\vskip2mm
	The quantity $W_{\mathscr{S}_0}^n(\zz)$ defined in \eqref{eq-W_S0^n} has the following variational representation:
	\begin{lemma}\label{Theorem 4.2.2}
		{\color{black}For any measurable function $h:\mathscr{P}(V_{\vec{l}(G)}) \mapsto \mathbb{R}$ which is  bounded and  continuous in $\mathcal{C}l(\mathscr{S}_0)$, $W^n_{\mathscr{S}_0}(\zz)$  coincides with minimal cost function $\vv^n_{\mathscr{S}_0}(\zz)$.}
	\end{lemma}

	We will write $W^n_{\mathscr{S}}(\zz)$ as $W^n(\zz)$. 
	This variational representation is a critical tool in our  approach. 
	The proof of Lemma \ref{Theorem 4.2.2} will be given in Section \ref{5-pr-2.5-lem1}. The proof is more complicated, since $W^n_{\mathscr{S}_0}(\zz)$ is a restricted  logarithmic Laplace functional.

	% In particular, write $W^n_{\mathscr{S}}(\zz)$ as $W^n(\zz)$. This variational representation is a critical tool in our  approach. Our proof of Lemma \ref{Theorem 4.2.2} is given by modifying the proof of \eqref{(2.5)}, stated in Section \ref{5-pr-2.5-lem1}. {\color{black}  Since we are considering  the restricted setting,  the proof is more involved. We  introduce Lemma \ref{Lemma of attainment condition} to guarantee that, under the attainment condition of the dynamic programming equation \eqref{equation of dynamic programming}, the process $\overline{L}_j^n$ remains in  $\mathcal{C}l(\mathscr{S}_0)$.  This property is essential for the proof because we must handle the constraint $\overline{\mathcal{L}}^n\in \mathcal{C}l(\mathscr{S}_0)$ in the infimum condition of \eqref{(2.6)}.}
	
	\begin{remark} \rm   
		
		Instead of the LP in Definition \ref{LP},  we consider the restricted logarithmic Laplace functional in Lemma  \ref{Theorem 4.2.2}, which formally can be regarded as
		$$
		W^n_{\mathscr{S}_0}(\zz)=-\frac{1}{n}\log{\overline{\E}_\zz\left\{\exp\left[-n\big(h\left(\mathcal{L}^n\right)+\infty \mathbf{1}_{\{\mathcal{L}^n\in\mathcal{C}l(\mathscr{S}_0)\}}\big)\right]\right\}}.
		$$
		After showing that $W^n_{\mathscr{S}_0}(\zz)$ 
		%RS equals a 
		is equal to the
		restricted minimal cost function $\vv^n_{\mathscr{S}_0}(\zz)$, we can obtain a generalized LDP for the empirical measures that belong to $\mathcal{C}l(\mathscr{S}_0)$ (Theorem \ref{estimate for rate}). 
	\end{remark}
	
	The weak convergence approach in \cite{DE1997} can not be applied directly because  $\mathcal{Z}$ is not Markovian. However, the decomposition of $\mathcal{Z}$ via its the renewal times  leads to  the following observation:    $\big\{\overline{\mathcal{Z}}_j^n\big\}_{0\le j\le n}$ behaves as a time homogeneous Markov chain whenever it remains inside a subgraph whose edges have all been traversed. This property allows us to split the whole trajectory into several segments, each confined to a fixed subgraph. We can then analyze each segment separately. To be more specific,
	
	\begin{definition}[Renewal time for the admissible control sequence]
		We denote the renewal times for the admissible control sequence by
		\begin{equation}
			\overline{\tauu}_k^n:=\inf\Big\{j>\overline{\tauu}_{k-1}^n:\ \overline{\mathcal{Z}}_j^n\big|_E\notin
			\big\{\overline{\mathcal{Z}}_i^n\big|_E:i<j\big\} \Big\}\wedge n,\ 2\leq k\leq \b,\label{renewal time}
		\end{equation}
		where $\overline{\tauu}_1^n:=0$. For convenience, let $\overline{\tauu}_{\b+1}^n=n.$
	\end{definition}
	
	For each $1\le k\le \b$   we analyze  the segment $(\overline{\mathcal{Z}}_j^n)_{\overline{\tauu}_k^n\le j<\overline{\tauu}_{k+1}^n}$ separately; see Figure \ref{schedule line 0}.
	\begin{figure}[htbp]
		\centering{\includegraphics{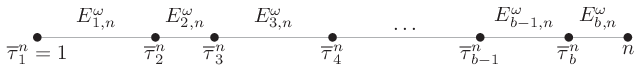}}
		\caption{\small{The process $\overline{\mathcal{Z}}_j^n|_E$ stays in $E_k$ during time $\overline{\tauu}_{k}^n$ and $\overline{\tauu}_{k+1}^n-1$}}
		\label{schedule line 0}
	\end{figure}    
	For a fixed sample path $\omega\in \overline{\Omega}$\   \footnote{ For convenience, provided no confusion arises, we shall continue to denote an element of $\overline{\Omega}$ by the symbol $\omega$.},  let    \begin{equation}
		E_{k,n}^{\omega}:=\big\{e\in E: \exists\, j<\overline{\tauu}_{k+1}^n, \overline{\mathcal{Z}}_j^n(\omega)|_E=e\big\},\ 1\le k\le \b,\label{subgraph}
	\end{equation}
	be the subgraph in which $\overline{\mathcal{Z}}_j^n|_E$ remains for  $j\in[\overline{\tauu}_{k}^n,\overline{\tauu}_{k+1}^n-1]$.  Given the control sequence  $\{\nuu_j^n\}_{0\le j\le n-1}$, the set $E_{k,n}^{\omega}$ depends on both $\omega$ and $n$.

	Since
	every edge in $ E_{k,n}^{\omega}$  has already been traversed by the stopping time $\overline{\tauu}_k^n$,  we have
	$$
	\pp_{\overline{\mathcal{L}}_j^n}=\pp_{E_{k,n}^{\omega}}, \text{ for }  j\in[\overline{\tauu}_{k}^n,\overline{\tauu}_{k+1}^n-1].
	$$ 
	which means that, conditionally on the history of $\overline{\mathcal{Z}}^n$ up to $\overline{\tauu}_k$, the segment  $(\overline{\mathcal{Z}}_j^n)_{\overline{\tauu}_k^n\le j<\overline{\tauu}_{k+1}^n}$ is  a time homogeneous Markov chain. This property allows us to apply the weak convergence technique for Markov chains separately on each interval  $[\overline{\tauu}_{k}^n,\overline{\tauu}_{k+1}^n-1]$ for $1\le k\le \b$.

	Figure \ref{graph of upper} illustrates the trajectory: the requirement that   $\overline{\mathcal{Z}}_\cdot^n$ stays inside $\vec{E}_{k,n}^{\omega}$ is equivalent to requiring that its tail   $\big(\overline{\mathcal{Z}}_\cdot^n\big)^+$ (the tail of $\overline{\mathcal{Z}}_\cdot^n$ defined in (\ref{head_tail}))  remains inside the undirected subgraph  $G_{k,n}^{\omega}=(V_{k,n}^{\omega},E_{k,n}^{\omega})$, where $V_{k,n}^{\omega}$ is the vertex set of $E_{k,n}^{\omega}$.

	\begin{figure}[htbp]
		%\captionsetup[subfigure]{labelformat=empty}
		\subfigure[]{\begin{minipage}[t]{0.4\textwidth}
				\centering{\includegraphics[scale=0.4]{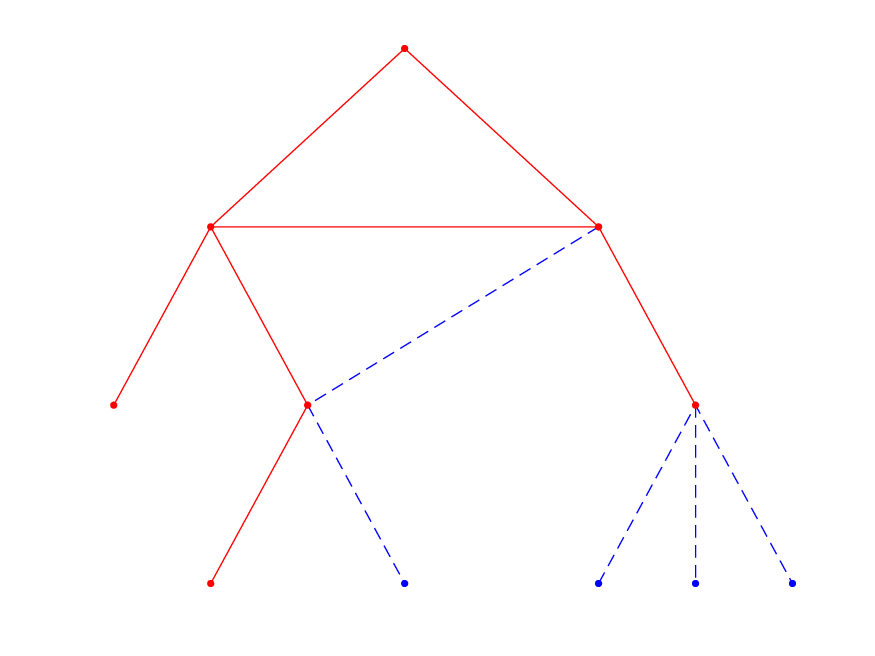}}
		\end{minipage}}
		\subfigure[]{\begin{minipage}[t]{0.4\textwidth}
				\centering{\includegraphics[scale=0.4]{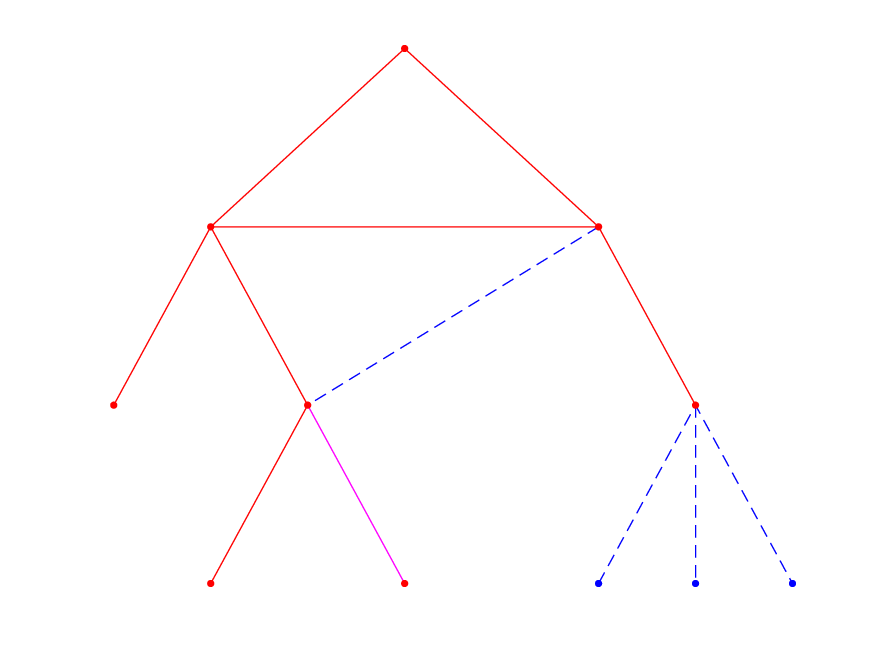}}
		\end{minipage}}
		\caption{\small{  While $\big(\overline{\mathcal{Z}}^n_.\big)^+$ moves along the solid edges, it evolves as a Markov process. (b) At a renewal time it traverses the first dashed edge from the left in (a), then remains on the solid edges shown in (b) for a subsequent interval, during which the process again exhibits Markov process.}}
		\label{graph of upper}
	\end{figure}
	
	\noindent{\bf(ii){\emph{\ Modified upper estimate of generalized Laplace principle}}} 
	
	We want to  prove 
	\begin{equation}\label{uppbdd-GLP}
		\liminf_{n\to\infty}W^n_{\mathscr{S}_0}(\zz)= \liminf_{n\to\infty}\vv^n_{\mathscr{S}_0}(\zz)\ge \inf_{\hat{\mu}\in\mathcal{C}l(\mathscr{S}_0)}\left\{{\Lambda}^\zz_{\delta,\mathscr{S}_0}(\hat{\mu})+h(\hat{\mu})\right\}.    
	\end{equation}
	Following  the same approach as in the part {\bf (ii)} in Section \ref{sec 2.2},  let $\varepsilon>0$ be given.  For each $n\ge  1$,  choose an admissible control sequence $\{\nuu_j^n\}_{0\le j\le n-1}$    such that
	$$
	\vv^n_{\mathscr{S}_0}(\zz)+\varepsilon\ge     \overline{\E}_\zz\Big\{\frac{1}{n}\sum_{j=0}^{n-1}R\Big(\nuu_j^n\big(\cdot\big|\overline{\mathcal{Z}}^n_{j},\overline{\mathcal{L}}^n_{j}\big)
	\Big\|\pp_{\overline{\mathcal{L}}^n_{j+1}}\big(\overline{\mathcal{Z}}^n_{j},\cdot\big)\Big)+h\big(\overline{\mathcal{L}}^n\big) \Big\}.
	$$
	Using \eqref{entropy prdct} and \eqref{trans prob prdct},  we obtain
	\begin{align*}
		&R\Big(\nuu_j^n\big({\rm d}\yy\left|\overline{\mathcal{Z}}_j^n,\overline{\mathcal{L}}_j^n\right.\big)\big\|\pp_{\overline{\mathcal{L}}_{j+1}^n} \big(\overline{\mathcal{Z}}_j^n,{\rm d}\yy\big)\Big)\\
		&=R\Big(\dd_{\overline{\mathcal{Z}}_j^n}({\rm d}\xx)\times\nuu_j^n\big({\rm d}\yy\left|\overline{\mathcal{Z}}_j^n,\overline{\mathcal{L}}_j^n\right.\big)\big\|\dd_{\overline{\mathcal{Z}}_j^n}({\rm d}\xx)\times \pp_{\overline{\mathcal{L}}_{j+1}^n}\big(\overline{\mathcal{Z}}_j^n,{\rm d}\yy\big)\Big)\\
		&=R\left(\dd_{\overline{\mathcal{Z}}_j^n}({\rm d}\xx)\times\nuu_j^n\big({\rm d}\yy\big|\overline{\mathcal{Z}}_j^n,\overline{\mathcal{L}}_j^n\big) \big\|\dd_{\overline{\mathcal{Z}}_j^n}({\rm d}\xx)\otimes \pp_{\overline{\mathcal{L}}_{j+1}^n}(\xx,{\rm d}\yy)\right).
	\end{align*}
	Using the the decomposition structure of $(\mathcal{Z}_n, \mathcal{L}^n)$, we can
	spilt the sum  {  $\sum_{j=0}^{n-1}(\cdot)$} into several segments  $\sum_{j=\overline{\tauu}_k^n}^{\overline{\tauu}_{k+1}^n-1}(\cdot)$ for $k=1,\dots,\b$ so that in the $k$-th segment the process  $\overline{\mathcal{Z}}_j^n$ remains inside  subgraph $\vec{E}_{k,n}^{\omega}$. Set
	$$
	n_k=n_k(n):=\overline{\tauu}^n_{k+1}-\overline{\tauu}^n_{k},\ 1\leq k\leq \b.
	$$   
	Then
	\begin{align*}
		&\overline{\E}_\zz\Big\{ \frac{1}{n}\sum_{j=0}^{n-1}R\Big(\nuu_j^n\big(\cdot\big|\overline{\mathcal{Z}}_j^n,\overline{\mathcal{L}}_j^n\big)\big\|\pp_{\overline{\mathcal{L}}_{j+1}^n}
		\big(\overline{\mathcal{Z}}_j^n,\cdot\big)\Big)+h\big(\overline{\mathcal{L}}^n\big) \Big\}\\
		&=\overline{\E}_\zz\Big\{ \frac{1}{n}\sum_{j=0}^{n-1}R\Big(\dd_{\overline{\mathcal{Z}}_j^n}({\rm d}\xx)\times\nuu_j^n\big({\rm d}\yy\big|\overline{\mathcal{Z}}_j^n,\overline{\mathcal{L}}_j^n\big)\big\|\dd_{\overline{\mathcal{Z}}_j^n}({\rm d}\xx)\otimes \pp_{\overline{\mathcal{L}}_{j+1}^n}(\xx,{\rm d}\yy)\Big)+h\big(\overline{\mathcal{L}}^n\big) \Big\}\\
		&= \overline{\E}_\zz\Bigg\{ \sum_{k=1}^{\b}\frac{n_k}{n}\frac{1}{n_k}\sum_{j=\overline{\tauu}^n_k}^{\overline{\tauu}^n_{k+1}-1}
		R\Big(\dd_{\overline{\mathcal{Z}}_j^n}({\rm d}\xx)\times\nuu_j^n\big({\rm d}\yy\big|\overline{\mathcal{Z}}_j^n,\overline{\mathcal{L}}_j^n\big)\big\|\dd_{\overline{\mathcal{Z}}_j^n}({\rm d}\xx)\otimes {\pp}_{E_{k,n}^\omega}(\xx,{\rm d}\yy)\Big)+h\big(\overline{\mathcal{L}}^n\big)\Bigg\}\\
		&\ \ \ \ \big(\text{applying Jensen's inequality to the convex function \  } R(\cdot\|\cdot) \ \ \nonumber \\
		&\ \ \ \ \text{ \ \ and using  the time homogeneity of\ } \overline{Z}^n_j  \text{ on } \vec{E}_{k,n}^{\omega} \text{ during } \ [\overline{\tauu}_k^n, \overline{\tauu}_{k+1}^n-1] \,\big) \\
		&\ge \overline{\E}_\zz\Bigg\{ \sum_{k=1}^{\b}\frac{n_k}{n}R\Big(\frac{1}{n_k}\sum_{j=\overline{\tauu}^n_k}^{\overline{\tauu}^n_{k+1}-1}\dd_{\overline{\mathcal{Z}}_j^n}({\rm d}\xx)\times\nuu_j^n\big({\rm d}\yy\big|\overline{\mathcal{Z}}_j^n,\overline{\mathcal{L}}_j^n\big)\Big\|\frac{1}{n_k}\sum_{j=\overline{\tauu}^n_k}^{\overline{\tauu}^n_{k+1}-1}
		\dd_{\overline{\mathcal{Z}}_j^n}({\rm d}\xx)\otimes {\pp}_{E_{k,n}^\omega}(\xx,{\rm d}\yy)\Big)\\
		&\hskip 11.5cm+h\big(\overline{\mathcal{L}}^n\big)\Bigg\},
	\end{align*}
	where we set $n_l=0$ whenever $E_{l,n}^\omega\notin\mathscr{S}_0$.
	
	Therefore, for  ORRW, we need  to study the  the tightness of the the following sequence of quadruples	\begin{equation}\label{4-terms upbdd}
		\Bigg(\frac{1}{n_k}\sum_{j=\overline{\tauu}^n_k}^{\overline{\tauu}^n_{k+1}-1}\dd_{\overline{\mathcal{Z}}^n_{j}}\times {\nuu}_j^n\big(\cdot\big|\overline{\mathcal{Z}}^n_{j},\overline{\mathcal{L}}^n_{j}\big),\frac{1}{n_k}\sum_{j=\overline{\tauu}^n_k}^{\overline{\tauu}^n_{k+1}-1}\dd_{\overline{\mathcal{Z}}^n_{j}},\frac{n_k}{n}, \mathbf{1}_{E_{k,n}^\omega}\Bigg)_{1\le k\le \b,\, n\ge1},
	\end{equation}
	which satisfies the following conditions: 
	\begin{itemize}
		\item[\bf{C1}]  $\sum_{k=1}^\b n_k/n=1$;
		\item[\bf{C2}]   the support of $\frac{1}{n_k}\sum_{j=\overline{\tauu}^n_k}^{\overline{\tauu}^n_{k+1}-1}\dd_{\overline{\mathcal{Z}}^n_{j}}$ is a subset of $E_{k,n}^\omega$;
		\item[\bf{C3}]	$\{E_{k,n}^\omega\}_{1\le k\le \b}\in\mathscr{E}_\zz$ (recall $\mathscr{E}_\zz$ from \eqref{eq-renew-subset}).
	\end{itemize}	
	Here the indicator function
	\begin{equation}
		{\bf 1}_{E_{k,n}^\omega}(e)=\left\{\begin{aligned}
			&1, &e\in E_{k,n}^\omega,\\
			&0, &e\notin E_{k,n}^\omega,
		\end{aligned}\right.\label{indicator of graph}
	\end{equation}
	completely characterizes  the support of the empirical measure $\frac{1}{n_k}\sum_{j=\overline{\tauu}^n_k}^{\overline{\tauu}^n_{k+1}-1}\dd_{\overline{\mathcal{Z}}^n_{j}}$. 
	
	To handle the convergence of the edge set $\{E^\omega_{k,n}\}_{1\le k\le \b}$, we equip $\{0,1\}^E$ with 	 the product topology, which induces a natural notion of convergence for such sequences.
	
	\vskip1mm Let $\big(\nuu_k,\overline{\mathcal{L}}_k,\mathscr{R}_k, \mathbf{1}_{E_{k,\infty}^\omega}\big)_{0\le k\le \b-1}$ be a limiting quadruple. Then,  by an   argument similar to that  in  part {\bf (ii)} in Section \ref{sec 2.2}, we obtain the inequality
	\begin{align*}
		\liminf_{n\to\infty}\vv^n_{\mathscr{S}_0}(\zz)+\varepsilon \ge &\ \overline{\E}_\zz\Big\{\sum_{k=1}^{\b}\mathscr{R}_k\int_{V_{\vec{l}(G)}} R\big( \overline{\qq}_k\big\| \pp_{E_{k,\infty}^\omega}\big)\ {\rm d}\overline{\mathcal{L}}_k+h\big(\overline{\mathcal{L}}\big)\Big\},%\
	\end{align*} 	
	which, combined with further analysis, yields he upper bound  \eqref{uppbdd-GLP}.

	The existence and properties of $\big(\nuu_k,\overline{\mathcal{L}}_k,\mathscr{R}_k, \mathbf{1}_{E_{k,\infty}^\omega}\big)_{0\le k\le \b-1}$	are established in Lemma \ref{weak convergence}, Theorem \ref{weak convergence app} and Proposition \ref{convergence of Lp}.  
	The detailed proof of the upper bound estimate for the generalized Laplace principle  \eqref{uppbdd-GLP} will be given	in Section 3.1.

	\vskip2mm
	
	\noindent{\bf(iii){\emph{\ Modified lower bound of generalized Laplace principle}}}
	
	We want to  prove
	\begin{equation}\label{Or-lowbdd}
		\limsup_{n\to\infty}\vv^n_{\mathscr{S}_0}(\zz)\le \inf_{\hat{\mu}\in\mathcal{C}l(\mathscr{S}_0)}\{{\Lambda}_{\delta,\mathscr{S}_0}^\zz(\hat{\mu})+h(\hat{\mu}) \}.
	\end{equation}
	Given any bounded continuous function $h$ on $\mathscr{P}(V_{\vec{l}(G)})$ and $\varepsilon>0$, there exists  $\hat{\gamma}\in\mathcal{C}l(\mathscr{S}_0)$ such that
	\begin{equation}\label{Or-lowbdd1}
		{\Lambda}_{\delta,\mathscr{S}_0}^\zz(\hat{\gamma})+h(\hat{\gamma})\le \inf_{\hat{\mu}\in\mathcal{C}l(\mathscr{S}_0)}\{{\Lambda}_{\delta, \mathscr{S}_0}^\zz(\hat{\mu})+h(\hat{\mu})\}+\varepsilon<\infty.
	\end{equation}
	Our goal is then to  construct  an admissible control sequence $\{\nuu_j^n\}_{0\le j\le n-1}$ satisfying
	\begin{equation}\label{low-bd-constr0}
		\lim_{n\to\infty}\overline{\E}_\zz\Big\{ \frac{1}{n}\sum_{j=0}^{n-1}R\Big(\nuu_j^n\big(\cdot\big|\overline{\mathcal{Z}}_j^n,\overline{\mathcal{L}}_j^n\big)\Big\|
		\pp_{\overline{\mathcal{L}}_{j+1}^n}\big(\overline{\mathcal{Z}}_j^n,\cdot\big)\Big)+h\big(\overline{\mathcal{L}}^n\big) \Big\}\le {\Lambda}_{\delta,\mathscr{S}_0}^\zz(\hat{\gamma})+h(\hat{\gamma})+\varepsilon.
	\end{equation}
	
	The construction of such an admissible control sequence is technically involved. 
	In this part, we outline the intuitive idea and explain the main modifications required. The complete proof is provided in Section \ref{sec 3.3}.
	
	Our strategy of constructing $\{\nuu_j^n\}_{0\le j\le n-1}$ is designed to meet the following five components:

	\begin{itemize}
		\item[\textbf{D1}]
		By  the lower semicontinuity of the relative entropy and the compactness of the state space $\vec{l}(G)$, we can choose a sequence of  $\{E_k\}_{1\le k\le \b}\in\mathscr{E}_\zz$ together with transition probabilities $\qq_k$, probability measures $\hat{\gamma}_k$ and proportions $r_k$, $k=1,\dots,\b$ such that 
		$$
		\sum_{k=1}^\b r_k\hat{\gamma}_k=\hat{\gamma}, \ \hat{\gamma}_k \qq_k=\hat{\gamma}_k, \ \text{\rm{supp}}(\hat{\gamma}_k|_E)\subseteq E_k,
		$$ 
		with the convention $r_l=0$ whenever  $E_l\notin\mathscr{S}_0$ and
		\begin{equation}\label{Or-lowbdd2}
			\sum_{k=1}^\b r_k\int_{V_{\vec{l}(G)}}R( \qq_k\| \pp_{E_k})\ {\rm d}\hat{\gamma}_{k} \le {\Lambda}_{\delta,\mathscr{S}_0}^\zz(\hat{\gamma}) +\varepsilon.
		\end{equation}
		Here we emphasize that the quadruples $(E_k, \qq_k,\hat{\gamma}_k,r_k)_{k=1,\dots, \b}$ are deterministic.
		
		\item[\textbf{D2}]  To simplify the construction, we fix deterministic integer sequences   $\{s_k(n)\}_{1\le k\le \b+1}$ as the renewal times, with 
		$$
		s_1(n)=0, \  s_{\b+1}(n)=n, \  \lim\limits_{n\to\infty}\frac{s_{k+1}(n)-s_k(n)}{n}=r_k.
		$$  
		That is, we set $\overline{\tauu}_{k}^n=s_k(n), 1\le k\le \b$ and require 
		$$
		\lim\limits_{n\to\infty}\frac{\overline{\tauu}^n_{k+1}-\overline{\tauu}^n_{k}}{n}=r_k,\   \ \sum_{k=1}^\b r_k=1.
		$$

		\item[\textbf{D3}]  $\frac{1}{n}\sum_{j=1}^n R\Big(\nuu_j^n\big(\cdot\big|\overline{\mathcal{Z}}_j^n,\overline{\mathcal{L}}_j^n\big)\Big\|\pp_{\overline{\mathcal{L}}_{j+1}^n}
		\Big)\to \sum_{k=1}^\b r_k\int_{V_{\vec{l}(G)}} R\left(\qq_k\|\pp_{E_k}\right)\ {\rm d}\hat{\gamma}_k$ as $n\to\infty$.
		\item[\textbf{D4}]   $h\big(\overline{\mathcal{L}}^n\big)\to h(\hat{\gamma})$ as $n\to\infty$, with probability 1.
		\item[\textbf{D5}]  $\left.\overline{\mathcal{Z}}^n_{\overline{\tauu}_k^n}\right|_E\in E_k\setminus E_{k-1}$.
	\end{itemize}
	Combining (\ref{Or-lowbdd1}) with (\ref{low-bd-constr0}), we obtain
	\begin{align*}
		\limsup_{n\to\infty}\vv^n_{\mathscr{S}_0}(\zz) &\le \lim_{n\to\infty} \overline{\E}_\zz\Bigg\{\frac{1}{n}\sum_{k=1}^\b \sum_{j=s_k}^{s_{k+1}-1}R\Big({\nuu}_j^n\big(\cdot\big|\overline{\mathcal{Z}}^n_{j},\overline{\mathcal{L}}_j^n\big)\Big\|\
		\pp_{E_k}\big(\overline{\mathcal{Z}}^n_{j},\cdot\big)\Big)+h\big(\overline{\mathcal{L}}^n\big)\Bigg\}\\
		&= \ \sum_{k=1}^\b r_k\int_{V_{\vec{l}(G)}}R(\qq_k\|\pp_{E_
			k})\ {\rm d}\hat{\gamma}_k+h(\hat{\gamma})\\
		&\le\   {\Lambda}_{\delta,\mathscr{S}_0}^\zz(\hat{\gamma})+h(\hat{\gamma})+\varepsilon     \\
		&\le   \inf_{\hat{\mu}\in\mathcal{C}l(\mathscr{S}_0)}\{{\Lambda}_{\delta, \mathscr{S}_0}^\zz(\hat{\mu})+h(\hat{\mu})\}+2\varepsilon.   
	\end{align*}	
	This completes the proof of the lower bound estimate \eqref{Or-lowbdd}.
	\vskip1mm
	{\color{black}Nevertheless, \textbf{D2}, \textbf{D3},  \textbf{D4} and \textbf{D5} impose strong mutual constraints. Consequently, $\{\nuu_j^n\}_{0\le j\le n-1}$ must be constructed carefully to satisfy all of them simultaneously.
		
		\vskip2mm
		It should be noted that our construction method for the admissible control sequence here differs from that used in part {\bf (iii)} of Section 	\ref{sec 2.2}.  Specifically, there are two difficulties when applying the method of \cite[Section 8.6]{DE1997}.
		\begin{itemize}
			\item  An intuitive idea is that, once $\qq_k$ and $\hat{\gamma}_k$ are given in {\bf D1}, we first choose an admissible control sequence ${\nuu}_j^n=\qq_k$ for $j\in[\overline{\tauu}_k^n, \overline{\tauu}_{k+1}^n-1]$ to satisfy {\bf D3}, and then  adjust further $\nuu_j^n$ to meet the requirement of {\bf D2}, {\bf D4} and {\bf D5}.
			
			However, in this setting, the probability measure $\hat{\gamma}_k$ may not be ergodic with respect to the transition probability $\qq_k$.  Similar to the lower bound estimate of Laplace principle (see \eqref{lowbdd}), in order to apply  $L^1$ ergodic theorem, we would need an ergodic measure.   If we follow the procedure outlined in  \cite[Section 8.6]{DE1997} 	to obtain an ergodic measure $\hat{\gamma}^*_k$  on $\vec{E}_k$ (see also  \eqref{lambda*}),    then $\hat{\gamma}^*_k$  would have  to satisfy conditions \textbf{B1},\textbf{\ B2} and \textbf{B3} from  Part {\bf (iii)} of Section \ref{sec 2.2}  as well as the support constraint  ${\rm supp}(\hat{\gamma}_k^*|_E)\subseteq E_k$.	Unfortunately, constructing such a $\hat{\gamma}^*_k$ 
			%RS that simultaneously fulfills all four requirements in our context is very difficult.	
			satisfying all four requirements in our context is highly nontrivial.	
			
			\item  
			To devise a practical construction of 	$\hat{\gamma}^*_k$ and $\qq^*_k$ such that	
			\begin{equation}\label{lambda*etp}
				\int_{V_{\vec{l}(G)}}R(\qq^*_k\| \pp_{E_k})\ {\rm d}\hat{\gamma}^*_{k}\le \int_{V_{\vec{l}(G)}}R( \qq_k\| \pp_{E_k})\ {\rm d}\hat{\gamma}_{k},
			\end{equation}
			i.e., to satisfy the second term in  (\ref{S2.2-1}),	it would be convenient if $\hat{\gamma}_k$, $\hat{\gamma}^*_k$, $\qq^*_k$ and $\pp_{E_k}$ were all supported on the same subgraph $\vec{E}_k$. However, this is impossible because the topological support of $\pp_{E_k}$ is the entire state space by the definition of  $\pp_{E_k}$ (see Definition \ref{Def for p_u} and \eqref{p_E}).   
		\end{itemize}
		
		In summary, constructing an ergodic measure $\hat{\gamma}_k^*$ 	that is supported  on the edge subset $E_k$ and satisfying \eqref{lambda*etp} is not a feasible technical approach. Hence we must seek a different route. 
		
		\vskip2mm
		Our approach is to consider each irreducible positive recurrent classes of the transition probability $\qq_k$ specified  in \textbf{D1}. Restricted to such an irreducible positive recurrent class, $\hat{\gamma}_k$ is ergodic for $\qq_k$. This allows us to apply the method in \cite[Section 8.6]{DE1997} separately within each irreducible positive recurrent class (see Lemma \ref{ergodic theorem}), without have to construct a new measure via convex  combination as  \eqref{lambda*}. The remain question is to ensure the connection path among each irreducible positive recurrent class (see Figure \ref{absorbing set} and \eqref{Or-lowbdd-m2}),  and to guarantee the connection path between the fixed $E_k$ and $E_{k+1}$ for $1\le k\le \b-1$ (see Lemma \ref{lower technique}).
		
		Specifically, the corresponding admissible control sequence is constructed according to the following aspects.   
		\begin{itemize}
			\item Let {\color{black} $S_{k,i}, i=1,\cdots,j_k$} denote the irreducible positive-recurrent  classes of the  transition probability $\qq_k$ on $V_{\vec{l}(G)}\times V_{\vec{l}(G)}$.
			If  $\hat{\gamma}_k$ is restricted to $S_{k,i}$, then   it is ergodic with respect with to $\qq_k$  on $S_{k,i}$ on that class.   
			
			\item Inside each irreducible  positive recurrent  class of $\qq_k$, we set $\nuu_j^n=\qq_k$.
			\item The process remains inside each  irreducible  positive recurrent class long enough for the ergodic property to take effect, then moves to the next  irreducible  positive recurrent class.		
			\item At time $\overline{\tauu}_{k+1}$,  the processes $\overline{\mathcal{Z}}^n$ leaves $\vec{E}_k$ and enters $\overrightarrow{\vec{E}_{k+1}\setminus \vec{E}_k}$ exactly.
		\end{itemize}
		The complete construction of the admissible control sequence $\{\nuu_j^n\}_{0\le j\le n-1}$ is  in given detail in Section \ref{sec 3.3}. 
	}
	
	\subsection{Notations}\label{notations}
	\noindent  Since our proof is involved and the notation is rather intricate, we collect in this subsection some frequently used symbols for the reader's convenience.
	\begin{table}[htbp]
		\centering
		\begin{tabular}{l l l}
			\hline
			\ &\textbf{Graph Structure}  & \ \\ \hline
			\textbf{Symbol} & \textbf{Meaning} & \textbf{Page} \\ \hline
			
			$G=(V,E)$ & finite connected graph with vertex set $V$ and edge set $E$ & \pageref{G=(V,E)}\\ \hline
			$u \sim v$ & $u,v\in V$ are adjacent & \pageref{u sim v}\\ \hline
			$V_{G'},E_{G'}$ & the vertex set and the edge set of graph $G'$ &\pageref{VG-EG} \\\hline
			$\mathscr{E}$ & \begin{tabular}[c]{@{}l@{}}the collection of subset sequences $\{E_k\}_{1\le k\le \b}$ with \\ $E_k\subseteq E$ and $\b=\vert E\vert$\end{tabular} & \pageref{mathscr{E}}\\ \hline
			$\vec{E}$ & \begin{tabular}[c]{@{}l@{}}directed edge set generated by orienting two directions \\ of each edge in $E$\end{tabular}  &\pageref{def-directed-edge-set}   \\\hline
			$\vec{l}(G)$ & line digraph with vertex set $V_{\vec{l}(G)}=\vec{E}$ constructed from $G$ & \pageref{S}\\ \hline
			$\zz\to \zz'$ & there is a directed edge from $\zz$ to $\zz'$ for all $\zz,\zz'\in \vec{E}$ & \pageref{z_1 to z_2}\\ \hline
			
			$\zz^+,\zz^-$ & the head and the tail of $\zz$ &\pageref{head_tail}\\\hline
			
			$\partial A$ & \begin{tabular}[c]{@{}l@{}}the boundary of directed subgraph of $\vec{l}(G)$ with the vertex \\ set $A\subseteq V_{\vec{l}(G)}$\end{tabular}  &\pageref{partial A}  \\\hline
			
			$\zz|_E$ & the edge projection of $\zz$ &\pageref{def-z|E-mu|E} \\\hline
			
			$\mathscr{S}$ & \begin{tabular}[c]{@{}l@{}}the collection of edge sets of all connected subgraphs of $G$ \\ that contain the starting point of the $\delta$-ORRW $X$\end{tabular} & \pageref{eq-edge-subset}\\ \hline
			$\mathscr{S}_0$ & a non-empty decreasing subset of $\mathscr{S}$ %((here decreasing means for all $E_1,E_2\in\mathscr{S}$ with $E_1\subseteq E_2$, if $E_2\in\mathscr{S}_0$, then $E_1\in\mathscr{S}_0$)),
			& \pageref{mathscr{S}_0}\\ \hline
			$\mathscr{E}_\zz$ & the subset of $\mathscr{E}$ with $E_1=\{\zz|_E\}$ and $\zz\in V_{\vec{l}(G)}$ & \pageref{eq-renew-subset}\\\hline
			
		\end{tabular}
	\end{table}		
	\begin{table}[htbp]
		\centering
		\begin{tabular}{l l l}	
			\hline			
			
			\ &\textbf{Stochastic Process} &\ \\\hline
			\textbf{Symbol} & \textbf{Meaning} & \textbf{Page} \\ \hline
			
			$X=(X_n)_{n\geq 0}$ & $\delta$-ORRW on $G$ & \pageref{X=(X_n)_{n>0}}\\ \hline
			$(L^n)_{n\geq 0}$ & empirical measure process of $X$ & \pageref{eq-emprical-ORRW} \\\hline
			$\mathscr{T}_G,\mathscr{T}_{\vec{l}(G)}$ & the collection of transition probabilities on $V$ and $V_{\vec{l}(G)}$ & \pageref{mathscr{T}_G}, \pageref{mathscr{T}_S}\\\hline
			$\tau_k$ & Renewal time of $X$ &\pageref{def of renewal time} \\\hline
			$\mathcal{Z}=(\mathcal{Z}_n)_{n\geq 0}$ & $\delta$-ORRW on $\vec{l}(G)$ & \pageref{Z-n}\\ \hline
			
			$(\mathcal{L}^n)_{n\geq 0}$ & empirical measure process of $\mathcal{Z}$ & \pageref{eq-empirical-lifted}\\ \hline
			${\tauu}_k$ & Renewal time of $\mathcal{Z}$ &\pageref{def of renewal time Z} \\\hline
			$\left(\overline{\mathcal{Z}}_j^n,\overline{\mathcal{L}}_j^n\right)_{0\leq j\leq n}$ & \begin{tabular}[c]{@{}l@{}}the process generated by
				$\{\nuu_j^n\}_{0\le j\le n-1}$ and its empirical \\ measure on $\vec{l}(G)$\end{tabular} & \pageref{overline{Z},overline{L}}\\\hline
			
		\end{tabular}
	\end{table}
	\begin{table}[htpb]
		\centering
		\begin{tabular}{l l l}
			\hline
			\ &\textbf{Measure and transition probability}&\ \\ \hline
			\textbf{Symbol} & \textbf{Meaning} & \textbf{Page} \\ \hline
			
			$\dd_x$ & the Dirac measure centered at $x$ on a topology space & \pageref{hat{delta}_x}\\ \hline
			$\mathscr{P}(\mathscr{X})$ & \begin{tabular}[c]{@{}l@{}}the set of all probability measures on a Polish space $\mathscr{X}$, \\ equipped with
				the weak convergence topology\end{tabular} & \pageref{mathscr{X}} \\ \hline
			$p_{E'}, \pp_{E'}$ & \begin{tabular}[c]{@{}l@{}}the transition probabilities of ORRWs $X$ and $\mathcal{Z}$ for the \\ traversed edge set $E'$\end{tabular} &\pageref{Def for p_u on G},\pageref{p_E} \\\hline
			$\hat{\mu}|_E$	& \begin{tabular}[c]{@{}l@{}} edge-projected measure from $\mathscr{P}(V_{\vec{l}(G)})$ to $\mathscr{P}(E)$\end{tabular} &\pageref{edge mu}\\ \hline
			$Cl(\mathscr{S}_0)$ & \begin{tabular}[c]{@{}l@{}}the collection of probability measures, whose first \\ marginal measures are supported on some edge set in $\mathscr{S}_0$\end{tabular} & \pageref{eq-closed-subset}\\ \hline
			
			$\{\nuu_j^n\}_{0\le j\le n-1}$ & the admissible control sequence & \pageref{nu_j^n_0},\pageref{def of nu_j^n}\\\hline
			$\Rightarrow$ & \begin{tabular}[c]{@{}l@{}}the weak convergence of random probability measures \\ (for each sample path $\omega$)\end{tabular} & \pageref{Rightarrow}\\ \hline
			$\overline{\P}_\zz,\overline{\E}_\zz$ & \begin{tabular}[c]{@{}l@{}}the probability generated by $\{\nuu_j^n\}_{0\le j\le n-1}$ conditioned on \\ the starting point $\zz\in V_{\vec{l}(G)}$, and the expectation of $\overline{\P}_\zz$\end{tabular} & \pageref{overline{P}_z}\\\hline
			$\xrightarrow[]{\mathcal{D}}$ & \begin{tabular}[c]{@{}l@{}}the weak convergence of random elements, where \\ the elements are measures equipped with the  weak \\ convergence topology\end{tabular} &\pageref{xrightarrow_D} \\ \hline
			
		\end{tabular}
	\end{table}
	\begin{table}[htpb]
		\centering
		\begin{tabular}{l l l}
			\hline
			\ &\textbf{Rate Function} &\ \\\hline
			\textbf{Symbol} & \textbf{Meaning} & \textbf{Page} \\ \hline
			
			$\mathscr{A}(\nu)$ & the collection of sequence $(\nu_k,r_k,E_k,q_k)$ & \pageref{A(nv)} \\\hline
			
			$I_\delta(\nu)$ & \begin{tabular}[c]{@{}l@{}}the rate function of the LDP for empirical measure process \\ $(L^n)_{n\geq 0}$  on $\mathscr{P}(V)$\end{tabular} & \pageref{I_delta}, \pageref{rate function 0}\\ \hline

			$\hat{\mathscr{A}}(\hat{\mu},\mathscr{S}_0)$ & the collection of sequences $(\hat{\mu}_k,r_k,E_k,\qq_k)$ & \pageref{A}\\\hline
			
			$\Lambda_{\delta,\mathscr{S}_0}$ & the general rate function on $\mathscr{P}(V_{\vec{l}(G)})$ & \pageref{rate function}\\\hline
			
			$\hat{\mathscr{A}}_\zz(\hat{\mu},\mathscr{S}_0)$ & the subset of $\hat{\mathscr{A}}(\hat{\mu},\mathscr{S}_0)$ satisfying $\{E_k\}_{1\le k\le \b}\in\mathscr{E}_\zz$ and $\zz\in V_{\vec{l}(G)}$ & \pageref{A_z}\\\hline
			$\Lambda_{\delta,\mathscr{S}_0}^\zz$ & \begin{tabular}[c]{@{}l@{}}the general rate function on $\mathscr{P}(V_{\vec{l}(G)})$ conditioned on \\ the starting point $\zz\in V_{\vec{l}(G)}$\end{tabular} & \pageref{Lambda^z}\\\hline
			$W^n_{\mathscr{S}_0}(\zz)$ & restricted logarithmic Laplace functional & \pageref{Theorem 4.2.2}\\\hline
			$\vv^n_{\mathscr{S}_0}(\zz)$ & minimal cost function & \pageref{(2.6)}\\\hline

		\end{tabular}	
	\end{table}

	\section{LDP and properties of the rate function}\label{sec 3} {\color{black}
		\noindent   This section is divided into three parts. The first part is devoted to proving the generalized Laplace principle, i.e.,Theorem \ref{estimate for rate}, which is the fundamental theorem of this paper. Section \ref{sec 3.2} and Section \ref{sec 3.3} establish, respectively, the upper bound and the lower bound of the generalized Laplace principle, thereby completing the proof of Theorem \ref{estimate for rate}.
		
		The second part, presented in Section \ref{sec 3.4}, proves the main result of this paper, Thereom \ref{I_delta}. 
		To facilitate reading, we provide a schematic diagram of the proof in Figure \ref{outline_of_the_proof_of_LDP-1} below.
		
		\begin{figure}[htbp]
			\centering{\includegraphics[scale=0.75]{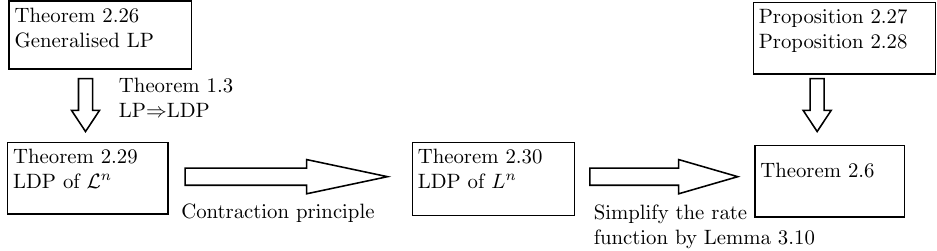}}
			\caption{\small{This is schematic diagram of the proof of Theorem \ref{I_delta}.}}
			\label{outline_of_the_proof_of_LDP-1}
		\end{figure}

		The third part, in Section \ref{sec 3.5}, establishes the analytic properties of the rate function 
		$I_\delta$, i.e., Theorem \ref{phase for rate function}. We first study the analytic properties of the rate function $\Lambda_\delta$ for 
		$(\mathcal{Z}_n)_{n\ge 1}$and its empirical measure process $(\mathcal{L}^n)_{n\ge 1}$. Then,   using the relationship between $G$ and $\vec{l}(G)$ as well as the connection between $((\mathcal{Z}_n, \mathcal{L}^n))_{n\ge 1}$ and $(X_n,L^n)_{n\ge 0}$, we obtain Theorem \ref{phase for rate function}.

		\subsection{The upper bound of generalized Laplace principle}\label{sec 3.2}
		\noindent      The analysis in part {\bf (ii)} of Subsection \ref{sec 2.3} shows that the asymptotic behavior of the quadruple \eqref{4-terms upbdd} with respect to 
		$n$ plays a key role in proving the upper bound. We therefore first present a theorem and related properties concerning the limit of the quadruple. 
		
		\vskip2mm
		
		The following lemma depicts the convergence behaviour of the first two terms appearing  \eqref{4-terms upbdd}.
		
		\begin{lemma}\label{weak convergence}
			Let $\overline{\mathcal{Z}}_0^n={\zz}$ for $n\geq 1,$ and $\{m(n)\}_{n\ge 1}, \{k(n)\}_{n\ge 1}$ be two sequences of stopping times on probability space $\left(\overline{\Omega},\overline{\mathscr{F}},\overline{\mathbb{P}}_{\zz}\right)$ such that $m(n)\le k(n)\le n,\ n\ge 1$.
			Then every subsequence of
			$$
			\Bigg(\frac{1}{k(n)-m(n)}\sum_{j=m(n)}^{k(n)-1}\dd_{\overline{\mathcal{Z}}^n_{j}}\times {\nuu}_j^n\big(\cdot\big|\overline{\mathcal{Z}}^n_{j},\overline{\mathcal{L}}^n_{j}\big),\frac{1}{k(n)-m(n)}\sum_{j=m(n)}^{k(n)-1}\dd_{\overline{\mathcal{Z}}^n_{j}}\Bigg)_{n\ge 1}
			$$
			has a sub-subsequence converging in distribution to some $\left(\nuu,\overline{\mathcal{L}}\right)$ with the following properties. %for $\omega\in\{\omega\in\overline{\Omega}:\lim_{n\to\infty}\frac{k-m}{n}=0 \}^c$.
			\begin{enumerate}
				\item[(a)] There exists a probability space $\left(\overline{\Omega},\overline{\mathscr{F}},\overline{\P}_\zz\right)$\footnote{We obtain the convergence by Skorohod's representation theorem. Here the probability space is actually a new space different from $\left(\overline{\Omega},\overline{\mathscr{F}},\overline{\P}_\zz\right)$ defined in the first paragraph of Subsection \ref{sec 2.3}. In order to facilitate reading, we still write this new space as $\left(\overline{\Omega},\overline{\mathscr{F}},\overline{\P}_\zz\right)$ in the rest of this paper. For more details of the usage of Skorohod's representation theorem, see \cite[Lemma 2.5.1]{DE1997}.} on which the limiting objects
				$$\nuu=\nuu({\rm d}\zz_1\times {\rm d}\zz_2\,\vert\,\omega)\ \mbox{and}\ \overline{\mathcal{L}}=\overline{\mathcal{L}}({\rm d}\zz_1\,\vert\,\omega)$$
				can be realized respectively as a {\color{black} random probability measure} on $V_{\vec{l}(G)}\times V_{\vec{l}(G)}$ and a {\color{black} random probability measure} on $V_{\vec{l}(G)}$. Moreover, for $\overline{\P}_\zz$-almost surely  $\omega\in\overline{\Omega}$, $\overline{\mathcal{L}}({\rm d}\zz_1\,|\,\omega)$ is the first marginal of $\nuu({\rm d}\zz_1\times {\rm d}\zz_2\,|\,\omega)$, i.e.,
				{\color{black}\begin{equation}
						\overline{\mathcal{L}}(A\,|\,\omega)=\nuu\big(A\times V_{\vec{l}(G)}\,|\,\omega\big),\ \forall\, A\subseteq V_{\vec{l}(G)}.\label{prop (a)}
				\end{equation}}
				\item[(b)] There is a transition probability  $\bar{\nu}({\rm d}\zz_2\,|\,\zz_1)=\bar{\nu}({\rm d}\zz_2\,|\,\zz_1,\omega)$ on $V_{\vec{l}(G)}$ such that $\overline{\P}_\zz$-a.s.\,for $\omega\in\overline{\Omega}$,
				\begin{equation}
					\nuu(A\times B\,|\,\omega) = \int_A \bar{\nu}(B\,|\,\zz_1,\omega)\ \overline{\mathcal{L}}({\rm d}\zz_1\,|\,\omega),\ \forall\, A,B\subseteq V_{\vec{l}(G)}.\label{prop (b)}
				\end{equation}
				\item[(c)] For each $\omega\in\overline{\Omega}$, define $\overline{\qq}(\zz_1,{\rm d}\zz_2\,|\,\omega):=\bar{\nu}({\rm d}\zz_2\,|\,\zz_1,\omega)$ is a transition probability  on $V_{\vec{l}(G)}$. Then  
				$\overline{\P}_\zz$-almost surely, 
				$$
				\lim_{n\to\infty}\frac{k(n)-m(n)}{n}=0,
				$$ or
				%(i.e., there is a null set $A'$ w.r.t $\overline{\mathbb{P}}_z$ s.t. for all $\omega\in\{\omega\in\overline{\Omega}:\liminf_{n\to\infty}\frac{k(n)-m(n)}{n}=0 \}\setminus A'$),	
				\begin{equation}\label{prop (c)}
					\begin{aligned}
						\overline{\mathcal{L}}(B\,|\,\omega)&=\int_{V_{\vec{l}(G)}} \overline{\qq}(\zz_1,B\,|\,\omega)\ \overline{\mathcal{L}}({\rm d}\zz_1\,|\,\omega)\\
						&=\int_{V_{\vec{l}(G)}} \bar{\nu}(B\,|\,\zz_1,\omega) \overline{\mathcal{L}}({\rm d}\zz_1\,|\,\omega)=\nuu(V_{\vec{l}(G)}\times B\,|\,\omega),\ B\subseteq V_{\vec{l}(G)}.
					\end{aligned}		
				\end{equation}
				In other words, for $\overline{\P}_\zz$-a.s. $\omega\in\overline{\Omega}$, at least one of the following event  occurs:
				\begin{itemize}
					\item $\lim_{n\to\infty}\frac{k(n)-m(n)}{n}\to 0$; \vskip1mm
					\item $\overline{\mathcal{L}}({\rm d}\zz_1\,|\,\omega)$ is an invariant measure of $\overline{\qq}(\zz_1,{\rm d}\zz_2\,|\,\omega)$ and coincides with the second marginal of $\overline{\mathcal{L}}({\rm d}\zz_1\,|\,\omega)\ \mbox{is the second marginal of}\ \nuu({\rm d}\zz_2\times {\rm d}\zz_1\,|\,\omega)$.
				\end{itemize}	 
				In the sequel we suppress the explicit dependence on $\omega$ in the notation for  $\nuu({\rm d}\zz_1\times {\rm d}\zz_2)$, $\overline{\mathcal{L}}({\rm d}\zz_1)$, $\bar{\nu}({\rm d}\zz_2\,|\,\zz_1)$ and $\overline{\qq}(\zz_1,{\rm d}\zz_2)$.  The relation displayed in part (b) can then be written concisely as
				$$
				\nuu({\rm d}\zz_1\times {\rm d}\zz_2)=\overline{\mathcal{L}}({\rm d}\zz_1)\otimes \bar{\nu}({\rm d}\zz_2\,|\,\zz_1)=\overline{\mathcal{L}}({\rm d}\zz_1)\otimes \overline{\qq}(\zz_1,{\rm d}\zz_2).
				$$
			\end{enumerate}
		\end{lemma}

		Lemma \ref{weak convergence} generalizes \cite[Theorem 8.2.8]{DE1997} (one sees this by setting $k(n)=n$, $m(n)=0$). A key difference in Lemma \ref{weak convergence} is that (\ref{prop (c)}) does not hold $\overline{\mathbb{P}}_\zz$-a.s.\,for $\omega\in\overline{\Omega}$.  The proof of Lemma \ref{weak convergence} which follows arguments similar to those of  \cite[Theorem 8.2.8]{DE1997}, is given in Section \ref{appendix 2}.
		
		\vskip2mm

		\begin{theorem}\label{weak convergence app}
			Let $\overline{\mathcal{Z}}_0^n={\zz},\ n\geq 1$ be initial conditions on $V_{\vec{l}(G)}$. Then every subsequence of
			$$
			\Bigg(\frac{1}{n_k}\sum_{j=\overline{\tauu}^n_k}^{\overline{\tauu}^n_{k+1}-1}\dd_{\overline{\mathcal{Z}}^n_{j}}\times {\nuu}_j^n\big(\cdot\big|\overline{\mathcal{Z}}^n_{j},\overline{\mathcal{L}}^n_{j}\big),\frac{1}{n_k}\sum_{j=\overline{\tauu}^n_k}^{\overline{\tauu}^n_{k+1}-1}\dd_{\overline{\mathcal{Z}}^n_{j}},\frac{n_k}{n}, \mathbf{1}_{E_{k,n}^\omega}\Bigg)_{1\le k\le \b,\, n\ge1}
			$$
			has a sub-subsequence such that for all $1\leq k\leq \b$,
			$$
			\Bigg(\frac{1}{n_k}\sum_{j=\overline{\tauu}^n_{k}}^{\overline{\tauu}^n_{k+1}-1}\dd_{\overline{\mathcal{Z}}^n_{j}}\times {\nuu}_j^n\big(\cdot\big|\overline{\mathcal{Z}}^n_{j},\overline{\mathcal{L}}^n_{j}\big),\frac{1}{n_k}\sum_{j=\overline{\tauu}^n_{k}}^{\overline{\tauu}^n_{k+1}-1}\dd_{\overline{\mathcal{Z}}^n_{j}},\frac{n_k}{n}, \mathbf{1}_{E_{k,n}^\omega}\Bigg)_{n\ge1}
			$$
			converges in distribution to some limit $\big(\nuu_k,\overline{\mathcal{L}}_k,\mathscr{R}_k, \mathbf{1}_{E_{k,\infty}^\omega}\big)$ with the following properties.
			\begin{itemize}
				\item[(a)] There is a probability space $\left(\overline{\Omega},\overline{\mathscr{F}},\overline{\P}_\zz\right)$ on which, for all $k\in[1,\b]$, $\nuu_k\text{ and }\overline{\mathcal{L}}_k$ can be realized as a {\color{black} random probability measure}  on $V_{\vec{l}(G)}\times V_{\vec{l}(G)}$ and a {\color{black} random probability measure} on $V_{\vec{l}(G)}$ respectively, and $\left(\nuu_k,\overline{\mathcal{L}}_k\right)$ satisfies \eqref{prop (a)}-\eqref{prop (b)}; moreover $\overline{\P}_{\zz}$-almost surely,
				$$
				\lim_{n\to\infty}\frac{\overline{\tauu}_{k+1}^n-\overline{\tauu}_{k}^n}{n}=0, \ \text{or } \big(\nuu_k,\overline{\mathcal{L}}_k\big) \ \text{satisfies\ }  \eqref{prop (c)}.
				$$
				\item[(b)] $\mathscr{R}_k$ is a $[0,1]$-valued random variable, and
				$$
				\sum_{k=1}^\b \mathscr{R}_k\overline{\mathcal{L}}_k=\overline{\mathcal{L}},\ \sum_{k=1}^\b \mathscr{R}_k=1,
				\ \overline{\P}_\zz \mbox{-a.s. for}\ \omega\in\overline{\Omega},
				$$
				where $\overline{\mathcal{L}}$ is the limit of $\overline{\mathcal{L}}^n$ satisfying \eqref{prop (a)}, \eqref{prop (b)} and \eqref{prop (c)}. That is, $\overline{\mathcal{L}}(\cdot|\omega)$ is invariant with respect to transition probability $\overline{\qq}(\cdot,\cdot|\omega)\in\mathscr{T}_{\vec{l}(G)}$.
				\item[(c)] $\overline{\mathbb{P}}_\zz$-a.s. for $\omega\in\overline{\Omega}$, $\{E_{k,\infty}^{\omega}\}_{1\le k\le \b}\in\mathscr{E}_\zz$, and $\text{\rm supp}(\overline{\mathcal{L}}_k|_E)\subseteq E^\omega_{k,\infty}$.
			\end{itemize}
		\end{theorem}
		
		\vskip2mm
		\begin{proposition}\label{convergence of Lp}
			Let $\{\hat{\mu}_n \}_{n\ge1}$ be a sequence of probability measures on $V_{\vec{l}(G)}$, and $\{\mathbf{1}_{E_n}\}$ a sequence of indicator functions defined in \eqref{indicator of graph}. Suppose $\hat{\mu}_n\Rightarrow \hat{\mu}$ and $\mathbf{1}_{E_n}\Rightarrow \mathbf{1}_{E_\infty}$. Then
			$$\hat{\mu}_n\otimes \pp_{E_n}\Rightarrow\hat{\mu}\otimes \pp_{E_\infty}.$$
		\end{proposition}
		
		\begin{remark}
			\rm Theorem \ref{weak convergence app} demonstrates the limit properties of \eqref{4-terms upbdd}, which are similar to {\bf C1}, {\bf C2} and {\bf C3} in Part {\bf{(ii)}} in Section \ref{sec 2.3}. These properties are crucial for obtaining a proper rate function of the LDP for empirical measures of ORRWs, as indicated \eqref{Lambda^z}.
			
			{\color{black} The proof of Theorem \ref{weak convergence app} is given in Section \ref{appendix 2}. Part (a) of the theorem shows that, for each segment separated by the stopping times, the limits $\nuu_k$ and $\overline{\mathcal{L}}_k$ not only satisfy the relations in  \eqref{prop (a)}-\eqref{prop (b)}, but also that $\overline{L}_k$ is an invariant measure for some random transition  probability $\hat{\qq}$ whenever the relative length of that segment does not vanish in the limit. Part (b) describes how the different segments are linked and verifies that the convex combination of the $\overline{\mathcal{L}}_k$ still meets the conditions of Lemma \ref{weak convergence}. We also emphasize that part (c) of the theorem is essentially a property of the support of the empirical measures.  It captures the phenomenon that the process remains confined to the subgraph spanned by the already visited vertices until it jumps to a previously unvisited vertex. Indeed, for every fixed $n$, this property holds trivially for the discrete time process  $\big\{\overline{\mathcal{Z}}_j^n\big\}_{0\le j\le n}$. By Theorem \ref{weak convergence app}, any limit points of  the sequence in  \eqref{4-terms upbdd} inherits the same property.}
		\end{remark}

		\begin{proof}[Proof of the upper bound of generalized Laplace principle.] \ Recall the part {\bf (ii)} in Section \ref{sec 2.3}. 
			
			\noindent{\small\bf  \emph{Step 1}}.  Let $\varepsilon>0$ be given. For each $n\ge 1$,  we choose an admissible control sequence $\{\nuu_j^n\}_{0\le j\le n-1}$    such that
			\begin{align*}
				&\,  \vv^n_{\mathscr{S}_0}(\zz)+\varepsilon\\
				\ge &\, \overline{\E}_\zz\Bigg\{ \sum_{k=1}^{\b}\frac{n_k}{n}R\Bigg(\frac{1}{n_k}\sum_{j=\overline{\tauu}^n_k}^{\overline{\tauu}^n_{k+1}-1}\dd_{\overline{\mathcal{Z}}_j^n}({\rm d}\xx)\times\nuu_j^n\big({\rm d}\yy\big|\overline{\mathcal{Z}}_j^n,\overline{\mathcal{L}}_j^n\big)\Big\|\frac{1}{n_k}\sum_{j=\overline{\tauu}^n_k}^{\overline{\tauu}^n_{k+1}-1}
				\dd_{\overline{\mathcal{Z}}_j^n}({\rm d}\xx)\otimes {\pp}_{E_{k,n}^\omega}(\xx,{\rm d}\yy)\Bigg)\\
				&\hskip 11.5cm+h\big(\overline{\mathcal{L}}^n\big)\Bigg\},
			\end{align*}
			where we set $n_l=0$ whenever $E_{l,n}^\omega\notin\mathscr{S}_0$.   
			
			By Theorem \ref{weak convergence app}, we can extract a further subsequence along which
			\begin{align*}
				&\Bigg(\frac{1}{n_k}\sum_{j=\overline{\tauu}^n_k}^{\overline{\tauu}^n_{k+1}-1}\dd_{\overline{\mathcal{Z}}^n_{j}}\times {\nuu}_j^n\big(\cdot\big|\overline{\mathcal{Z}}^n_{j},\overline{\mathcal{L}}^n_{j}\big),\frac{1}{n_k}\sum_{j=\overline{\tauu}^n_k}^{\overline{\tauu}^n_{k+1}-1}\dd_{\overline{\mathcal{Z}}^n_{j}},\frac{n_k}{n}, \mathbf{1}_{E_{k,n}^\omega}\Bigg)_{1\le k\le \b}\\
				&\xrightarrow{\mathcal{D}} \big(\nuu_k,\overline{\mathcal{L}}_k,\mathscr{R}_k, \mathbf{1}_{E_{k,\infty}^\omega}\big)_{1\le k\le \b}.
			\end{align*}
			Here, "$\xrightarrow{\mathcal{D}}$"\label{xrightarrow_D} denotes weak convergence in distribution of random elements. By Skorohod’s representation theorem we may, without loss of generality, assume that this convergence holds almost surely.
			
			\noindent{\small\bf  \emph{Step 2}}.  Using Propositiom \ref{convergence of Lp}, we have
			\[
			\frac{1}{n_k}\sum_{j=\overline{\tauu}^n_k}^{\overline{\tauu}^n_{k+1}-1}
			\dd_{\overline{\mathcal{Z}}_j^n}({\rm d}\xx)\otimes {\pp}_{E_{k,n}^\omega}(\xx,{\rm d}\yy)\Rightarrow\overline{\mathcal{L}}_k\otimes \pp_{E_{k,\infty}^\omega}\text{ with probability 1},
			\]
			Hence, by Fatou's Lemma,  and the lower semicontinuity of $R(\cdot\|\cdot)$, we obtain that
			\begin{align}
				\liminf_{n\to\infty}\vv^n_{\mathscr{S}_0}(\zz)+\varepsilon&\ge\liminf_{n\to\infty}\overline{\E}_\zz\Big\{ \frac{1}{n}\sum_{j=0}^{n-1}R\big(\nuu_j^n(\cdot|\overline{\mathcal{Z}}_j^n,\overline{\mathcal{L}}_{j+1}^n)
				\big\|\pp_{\overline{\mathcal{L}}_j^n}(\overline{\mathcal{Z}}_j^n,\cdot)\big)+h\big(\overline{\mathcal{L}}^n\big) \Big\}\nonumber\\
				&\ge \overline{\E}_\zz\Big\{\sum_{k=1}^{\b}\mathscr{R}_kR\big(\overline{\mathcal{L}}_k\otimes \overline{\qq}_k\big\|\overline{\mathcal{L}}_k\otimes {\pp}_{E_{k,\infty}^\omega}\big)+h\left(\overline{\mathcal{L}}\right)\Big\}\nonumber\\
				&=\overline{\E}_\zz\Big\{\sum_{k=1}^{\b}\mathscr{R}_k\int_{V_{\vec{l}(G)}} R\big( \overline{\qq}_k\big\| \pp_{E_{k,\infty}^\omega}\big)\ {\rm d}\overline{\mathcal{L}}_k+h\left(\overline{\mathcal{L}}\right)\Big\}.\label{(3.7)}
			\end{align}
			
			\noindent{\small\bf  \emph{Step 3}}.   Notice that
			\begin{align*}
				\Big\{(r_k,f_k)_{1\le k\le \b}:&\ f_k\in\{0,1\}^E,\{\text{supp} (f_k)\}_{1\le k\le \b}\in\mathscr{E}_\zz, \\
				&\ r_k\ge 0, r_l=0\text{ for}\, l\ \text{if}\, \text{supp}(f_l)\notin\mathscr{S}_0, \sum_{k=1}^\b r_k=1\Big\}
			\end{align*}
			is a closed set. {\color{black}
				Here is a short proof of the closedness:
				
				Suppose  $(r_{k,n},f_{k,n})_{1\le k\le \b}\to (r_k,f_k)_{1\le k\le \b}$ as $n\to\infty$ for some $(r_{k,n},f_{k,n})_{1\le k\le \b}$ belongs to  this set. Because the domain is finite and the values of $f_{k,n}$ are discrete, there exists a large integer $N$ such that $f_{k,n}=f_k$ for all $n\ge N$. Consequently, for every index $l$ with  $\text{supp}(f_l)\notin \mathscr{S}_0$, we have $r_{l,n}=0$ whenever $n\ge N$, which forces $r_l=0$. The property $\sum_{k=1}^\b r_k=1$   follows from the fact that  $1=\sum_{k=1}^\b r_{k,n}\to\sum_{k=1}^\b r_k$.\\
			}
			
			Since $n_l=0$ whenever  $E_{l,n}^\omega\notin\mathscr{S}_0$ and $(\frac{n_k}{n},\mathbf{1}_{E_{k,n}^\omega})_{1\le k\le \b}$ lies in  this closed set above,
			the same reasoning as in Theorem \ref{weak convergence app} (c) gives  $\mathscr{R}_l(\omega)=0$ for every $\omega$ with $E_{l,\infty}^\omega\notin\mathscr{S}_0$.  Because $\sum_{k=1}^\b \overline{\mathcal{L}}_k\mathscr{R}_k=\overline{\mathcal{L}}$, 
			this implies $\overline{\mathcal{L}}\in\mathcal{C}l(\mathscr{S}_0)$.

			Recall $\hat{\mathscr{A}}(\hat{\mu},\mathscr{S}_0)$ defined in (\ref{A}) and $\hat{\mathscr{A}}_\zz(\hat{\mu},\mathscr{S}_0)$ defined in  (\ref{A_z}). 
			
			Observe that
			$$
			\left(\overline{\mathcal{L}}_k(\omega),\mathscr{R}_k(\omega),E_{k,\infty}^\omega, \qq_k(\omega)\right)_{k}\in\hat{\mathscr{A}}_\zz
			\left(\overline{\mathcal{L}}(\omega),\mathscr{S}_0\right)\ \mbox{for all}\ \omega,
			$$
			and $\overline{\mathcal{L}}_k\overline{\qq}_k=\overline{\mathcal{L}}_k$ on $\{\omega: \mathscr{R}_k>0 \}$ (this follows from the last statement in Theorem \ref{weak convergence app} (a)).

			From the definition of 
			$\Lambda^\zz_{\delta,\mathscr{S}_0}$ given in
			\eqref{Lambda^z},  and the fact that, for each fixed $\omega$, $\overline{\mathcal{L}}(\cdot|\omega)$ is invariant for some $\overline{\qq}(\cdot,\cdot|\omega)\in\mathscr{T}_{\vec{l}(G)}$ (established in (b) of Theorem \ref{weak convergence app}), we obtain
			\begin{align*}
				\text{(\ref{(3.7)})}
				\ge\overline{\E}_\zz\left\{ \Lambda^\zz_{\delta,\mathscr{S}_0}\left(\overline{\mathcal{L}}\right)+h\left(\overline{\mathcal{L}}\right) \right\}%&\ge\inf_{\omega}\left\{\sum_{k=1}^{d}\mathscr{R}_k\int_{V_{\vec{l}(G)}} R\left( \overline{q}_k\left\| p_{E_{k,\infty}^\omega}\right.\right)\ {\rm d}\overline{\mathcal{L}}_k+h\left(\overline{\mathcal{L}}\right)\right\}\\
				\ge \inf_{\hat{\mu}\in \mathcal{C}l(\mathscr{S}_0)}\{{\Lambda}_{\delta,\mathscr{S}_0}^\zz(\hat{\mu})+h(\hat{\mu})\}.
			\end{align*}
			Since $\varepsilon$ is arbitrary, this completes the proof of the upper bound. 
		\end{proof}

		\subsection{The lower bound of generalized Laplace principle}\label{sec 3.3}	
		\noindent      {\color{black}  For every  $h\in C_b\big(\mathscr{P}(V_{\vec{l}(G)})\big)$ and $\varepsilon>0$,   there exists a $\hat{\gamma}\in \mathcal{C}l(\mathscr{S}_0)$ satisfying \eqref{Or-lowbdd1}. To prove  the lower bound of generalized Laplace principle,   we must construct   an admissible control sequence $\{\nuu_j^n\}_{0\le j\le n-1}$ for which \eqref{low-bd-constr0} holds.  We now implement the heuristic construction outlined in  part {\bf(iii)}  in Section \ref{sec 2.3}.
			
			\vskip2mm
			\noindent {\emph{Constructing of admissible control sequence} $\{\nuu_j^n\}_{0\le j\le n-1}$.}   \vskip1mm
			
			\noindent{\small\bf  \emph{Step 1}}.   Choose the deterministic quadruple $(E_k, \qq_k,\hat{\gamma}_k,r_k)_{k=1,\dots, \b}$ described in part {\bf(iii)} {\bf D1}   in Section \ref{sec 2.3} such that \eqref{Or-lowbdd2} holds.
			
			\vskip1mm
			
			\noindent{\small\bf  \emph{Step 2}}. We select  a deterministic integer sequence $\{s_k=s_k(n)\}_{1\le k\le \b+1}$  with $s_1=0$, $s_{\b+1}=n$ such that
			\begin{itemize}
				\item $s_i(n)< s_j(n)$ for $i<j$;
				\item $\displaystyle\lim\limits_{n\to\infty}\frac{s_{k+1}(n)-s_{k}(n)}{n}=r_k$;
				\item for every $l\in[1,\b]$ with $E_l\notin\mathscr{S}_0$, set $s_l(n)=n$;
				\item there exists a path on $\vec{l}(G)$ staring from $\zz$ such that   stays inside $\vec{E}_k$ during the interval $[s_k(n), s_{k+1}(n)-1)$ and reaches $\overrightarrow{E_{k+1}\setminus E_k}$ exactly at time $s_{k+1}(n)$.
			\end{itemize}     
			
			\vskip1mm
			\noindent{\small\bf  \emph{Step 3}}.   Observe that there exists an integer $N$ large enough such that for all $n>N$, 
			$$
			s_{k+1}(n)-s_k(n)\ge 2\b-1,\ \ \ \text{for every}\ \ k \ \ \text{with}  r_k=\lim\limits_{n\to\infty}(s_{k+1}-s_{k})/n>0.  $$
			For $n\ge N$,  we construct  an admissible control sequence $\{\nuu_j^n\}_{0\le j\le n-1}$ for $n>N$ satisfying (writing  $s_k(n)$ simply as $s_k$) satisfying:
			\begin{itemize}
				\item[\textbf{E1}] For $s_k\le j<s_{k+1}$, $\overline{\mathcal{Z}}_j^n\in \vec{E}_k$.
				\item[\textbf{E2}] For $k\in[1,\b-1]$ and $j\in[{s}_{k+1}-2\b-1,s_{k+1}-1]$ with $s_{k+1}-s_k\ge 2\b+1$, under $\nuu_j^n$, we have $\overline{\mathcal{Z}}_{s_{k+1}}^n\in \overrightarrow{E_{k+1}\setminus E_k}$.
				\item[\textbf{E3}] For $s_k\le j< {s}_{k+1}$ ($1\le k\le \b$),
				$$
				\lim_{n\to \infty}\frac{1}{s_{k+1}-s_k}\sum_{j=s_k}^{s_{k+1}-1}R\Big(\nuu_j^n\big(\cdot\big|\overline{\mathcal{Z}}_j^n,\overline{\mathcal{L}}_j^n\big)\Big\|\pp_{E_k}
				\big(\overline{\mathcal{Z}}_j^n,\cdot\big)\Big)=\int_{V_{\vec{l}(G)}}R(\qq_k\|\pp_{E_k}){\rm d}\hat{\gamma}_k.
				$$
			\end{itemize}
			
			We emphasize that, by \eqref{pp_mu} in Definition \ref{pp_E}, \textbf{E1} alone guarantees 
			\begin{equation}\label{condition of lower bound}
				R\Big(\nuu_j^n\big(\cdot\big|\overline{\mathcal{Z}}_j^n,\overline{\mathcal{L}}_j^n\big)\Big\|\pp_{{\overline{\mathcal{L}}_{j+1}^n}}
				\big(\overline{\mathcal{Z}}_j^n,\cdot\big)\Big)=R\Big(\nuu_j^n\big(\cdot\big|\overline{\mathcal{Z}}_j^n,\overline{\mathcal{L}}_j^n\big)\Big\|\pp_{E_k}
				\big(\overline{\mathcal{Z}}_j^n,\cdot\big)\Big),
			\end{equation}
			while \textbf{E1} together with \textbf{E2} ensures that $\overline{\mathcal{Z}}^n_{j}$ depart from $\vec{E}_k$ and arrive at $\overrightarrow{E_{k+1}\setminus E_k}$  exactly  at time $s_{k+1}$.
			This is one of the key technical points of our construction, which differs from the Markov chain case.

			To achieve \textbf{E1}, \textbf{E2} and \textbf{E3}, we build the admissible control  $\{\nuu_j^n\}_{0\le j\le n-1}$
			on each interval $[s_k, s_{k+1}]$ in two parts: 
			\begin{itemize}
				\item the ergodicity part ($s_k\le j\le s_{k+1}-2\b-2$), which ensures \textbf{E1} and \textbf{E3};
				\item the connection part ($s_k\le s_{k+1}-2\b-1\le j<s_{k+1}$), which ensures \textbf{E1} and \textbf{E2}.\end{itemize}
			See Figure \ref{schedule line 2} for an illustration.
			
			\begin{figure}[htbp]
				\centering{\includegraphics[scale=0.7]{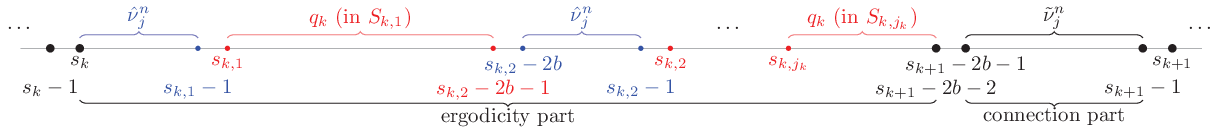}}
				\caption{\small{Ergodicity part} and Connection part}
				\label{schedule line 2}
			\end{figure}
			
			\noindent\textbf{\emph{ {\small Step 3.1}: Ergodicity  part}}.
			\vskip1mm
			We construct an admissible control sequence for the indices $j$ belonging to the interval 
			$$
			[s_k, s_{k+1}-2\b-2],\ \ 1\le k\le \b,
			$$
			so that \textbf{E1} and \textbf{E3} are satisfied.
			
			Since the state space  $\vec{l}(G)$ is finite, the transition probability  $\qq_k$ possesses only finitely many irreducible positive recurrent classes. If $\qq_k$ has a unique  irreducible positive recurrent class--- equivalently, if  $\hat{\gamma}_k$   is the unique invariant measure of $\qq_k$ on $\vec{E}_k$---then the admissible control can be taken simply as $\qq_k$ itself. In this case, the  $L^1$ ergodic theorem \cite[Theorem A.4.4]{DE1997} guarantees that \textbf{E1} and \textbf{E3} hold.

			Currently, we assume that $\qq_k$ has more than one irreducible positive recurrent classes, denote these classes by $S_{k,i}$ for $i=1,\dots,j_k$.    Let $\hat{\gamma}_k|_{S_{k,i}}$ the normalized restriction of $\hat{\gamma}_k$ to $S_{k,i}$, and let $C_{k,i}$ be the corresponding normalizing constants of $\hat{\gamma}_k$ on $S_{k,i}$.  
			
			We choose $\nuu_j^n$ such that the process remains inside $S_{k,i}$ during the interval  $[s_{k,i}, s_{k,i+1}-2\b]$, and  exits $S_{k,i}$ to visit $S_{k,i+1}$ during  $[s_{k,i+1}-2\b, s_{k,i+1}]$, where
			$$
			\lim_{n\to \infty}\frac{s_{k,i+1}-s_{k,i}}{n}=r_k\cdot C_{k,i}, \  \ s_{k,1}=s_k,\ s_{k,j_k+1}=s_{k+1}-1,
			$$
			see Figure \ref{absorbing set}.
			
			Precisely, let   $\varphi_{A}(\cdot,\cdot)$ denote the graph distance on a directed graph $A$. Define
			\begin{equation}\label{Or-lowbdd-m2}
				\nuu_j^n:=\left\{
				\begin{array}{lll}
					\qq_k, &  j\in[s_{k,i},s_{k,i+1}-2\b-1], &i\in[1,j_k],\\
					\hat{\nuu}_j^n, & j \in[s_{k,i+1}-2\b,s_{k,i+1}],& i\in[1,j_k-1],
				\end{array}
				\right.
			\end{equation}
			where {\color{black} for $j\in[s_{k,i+1}-2\b,s_{k,i+1}]$, $i\in[1,j_k-1]$,
				$$
				\hat{\nuu}_j^n(\xx,\yy)=1/n_{\xx,S_{k,i+1}}\ \ \text{if}\  \  \xx\notin S_{k,i+1},\ \  \varphi_{\vec{E}_k}(\xx,S_{k,i+1})-\varphi_{\vec{E}_k}(\yy,S_{k,i+1})=1,
				$$ {\color{black} with $n_{\xx,S_{k,i+1}}$ being the number of vertices $\yy$ satisfying $\varphi_{\vec{E}_k}(\xx,S_{k,i+1})-\varphi_{\vec{E}_k}(\yy,S_{k,i+1})=1$. For  $\xx\in S_{k,i+1}$,  $\hat{\nuu}_j^n(\xx,\cdot)$ is supported uniformly on those $\yy\in S_{k,i+1}$ for which $\qq_k(\xx,\yy)>0$.}
				%where the number of $\yy$ satisfying $\varphi_{\vec{E}_k}(\xx,S_{k,i+1})-\varphi_{\vec{E}_k}(\yy,S_{k,i+1})=1$ is $n_{\xx,S_{k,i+1}}$.
				%And for $\xx\in S_{k,i+1}$,  $\hat{\nuu}_j^n(\xx,\cdot)$ is supported uniformly on those $\yy\in S_{k,i+1}$ for which $\qq_k(\xx,\yy)>0$.
			}
			
			\begin{figure}[htbp]
				\centering{\includegraphics[scale=1]{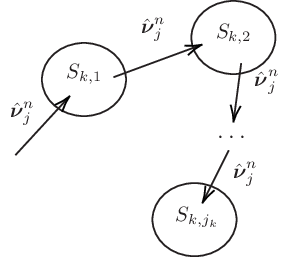}}
				\caption{\small{The process $\overline{\mathcal{Z}}_j^n$ is governed by $\qq_k$ within each irreducible positive recurrent class  $S_{k,i}$ $(i=1,\dots,j_k)$.After spending a sufficiently long time in such a class, the process is steered by $\hat{\nuu}_j^n$ and moves to the next irreducible positive recurrent class.}}
				\label{absorbing set}
			\end{figure}
			
			Furthermore, by  $L^1$ ergodic theorem and the pointwise ergodic theorem \cite[Theorem A.4.4]{DE1997}, we obtain the following lemma:
			
			\begin{lemma}\label{ergodic theorem}
				For $j\in[s_{k,i},s_{k+1,i}-2\b-1]$, the corresponding $\overline{\mathcal{Z}}_j^n$ is driven   by $\qq_k$. Moreover, for all $\zz\in S_{k,i}$, if $\lim\limits_{n\to\infty}\frac{s_{k+1,i}-s_{k,i}}{n}>0$, then as $n\to \infty$,
				\begin{align*}
					\overline{\E}\Bigg\{ \Bigg|\frac{1}{s_{k+1,i}-s_{k,i}-2\b-1}\cdot\sum_{j=s_{k,i}}^{s_{k+1,i}-2\b-1}&R\Big(\qq_k\big(\overline{\mathcal{Z}}_j^n,\cdot\big)\Big\|\pp_{E_k}\big(\overline{\mathcal{Z}}_j^n,\cdot\big)
					\Big) - \\
					&\int_{V_{\vec{l}(G)}}R(\qq_k\|\pp_{E_k})\ {\rm d}(\hat{\gamma}_k|_{S_{k,i}})\Bigg|\ \ \Bigg|\overline{\mathcal{Z}}_{s_{k,i}}^n=\zz \Bigg\}\to 0,
				\end{align*}
				and  conditionally on $\overline{\mathcal{Z}}_{s_{k,i}}^n=\zz$,
				\begin{align}\label{Or-lowbdd-m0}
					\frac{1}{s_{k+1,i}-s_{k,i}-2\b-1}\cdot\sum_{j=s_{k,i}}^{s_{k+1,i}-2\b-1}\dd_{\overline{\mathcal{Z}}_j^n}\Rightarrow\hat{\gamma}_k|_{S_{k,i}}\text{ with probability 1}.
				\end{align}
			\end{lemma}
			
			The proof of Lemma \ref{ergodic theorem} follows a similar line to that of the lower bound argument in \cite[Section 8.6]{DE1997}; details is provided  in  Section \ref{sec pf ergodic theorem}.%{\color{yellow}Appendix \ref{pf of Lem ergdc thm}}. 

			Moreover, because each $S_{k,i}$ is finite, the condition $\overline{\mathcal{Z}}_{s_{k,i}}^n=\zz\in S_{k,i}$ in Lemma \ref{ergodic theorem} can be replaced by an arbitrary initial distribution supported on $S_{k,i}$.  By our construction of  $\hat{\nuu}_j^n$, we guarantee that $\overline{\mathcal{Z}}_{s_{k,i}}^n$ 
			stays inside $S_{k,i}$ whenever $\overline{\mathcal{Z}}_0^n=\zz$ for any $\zz\in V_{\vec{l}(G)}$. Consequently,      
			Lemma \ref{ergodic theorem} remains valid under the initial condition $\overline{\mathcal{Z}}_0^n=\zz$ for every $\zz\in V_{\vec{l}(G)}$.
			
			Applying (\ref{Or-lowbdd-m0}),  summing over all $k$ and $i$, and using the continuity of $h$, we obtain
			\begin{corollary}\label{ergodic theorem-h}
				Under the  condition $\overline{\mathcal{Z}}_0^n=\zz$, we have  $\overline{\mathcal{L}}^n\Rightarrow\hat{\gamma}$ and $h\big(\overline{\mathcal{L}}^n\big)\to h(\hat{\gamma})$ with probability 1.
			\end{corollary}
			
			Because  $s_l=n$ for all $l$ with $E_l\notin\mathscr{S}_0$,  coming the  above two parts, combining the two parts above shows that, under the chosen admissible control sequence,  we obtain that $\overline{\mathcal{L}}^n\in\mathcal{C}l(\mathscr{S}_0)$.

			\begin{remark} 
				The construction of $\nuu_j^n$ described above is quite intuitive and offers an alternative method for verifying the lower bound of the Laplace principle for Markov chains on a finite state space.
			\end{remark}

			\vskip2mm
			
			\noindent\textbf{\emph{ {\small Step 3.2}: Connection  part}}.
			\vskip1mm
			
			We choose the admissible control sequence as $\tilde{\nuu}_j^n$ for indices $j$ belonging  to the connection part
			$$
			[s_{k+1}-2\b-1,s_{k+1}]
			$$
			such that the starting point $\overline{\mathcal{Z}}_{{s}_{k+1}-2\b-1}^n$ is random,  while the terminal point $\overline{\mathcal{Z}}_{s_{k+1}}^n$ is deterministic.
			
			We now state a lemma that guarantees the existence of such an embedding path; the path is deterministic once the starting point is fixed.

			\begin{lemma}\label{lower technique}
				Let $G_0=(V_0,E_0)$ be a finite connected undirected graph with at most  $\b$ edges. 
				Suppose that for a time interval $[s,t]$ with  $s<t-2\b$ and for vertices $u_0,v_0\in V_0$, there exists a path 
				$u_0u_{1}\dots u_{t-s}$ satisfying  $u_{t-s}=v_0$. Then there exists a sequence $(\tilde{p}_n)_{t-2\b\le n\le t}$ of transition probabilities on $G_0$ such that  every Markov chain $(U_n)_{s\leq n<t-2\b}$ on $G_0$ with $U_s=u_0$ can be extended over $[t-2\b, t]$	to a process $(U_n)_{s\leq n\le t}$ with the following properties:
				\begin{itemize}
					\item[(a)] $(U_n)_{s\le n\le t}$ is a Markov Chain;
					\item[(b)] the transition probability of $(U_n)_{t-2\b\le n\le t}$ are given by $(\tilde{p}_n)_{t-2\b\le n\le t}$;
					\item[(c)]   for every  $v\in V_0$ and every $n\in [t-2\b,  t]$,  $\tilde{p}_n(v, v)=0$ and  $\sum_{u\sim v}\tilde{p}_n(v, u)=1$;
					\item[(d)] $U_t=v_0$.
				\end{itemize}
			\end{lemma}
			
			\vskip1mm
			It is worth noting   that $(\tilde{p}_n)_{t-2\b\le n\le t}$ is independent of $(U_n)_{s\leq n<t-2\b}$. The core idea of Lemma \ref{lower technique} is to prescribe a rule that forces $(U_n)_n$ to follow a deterministic trajectory after time  along a deterministic path after time $t-2\b$ and ultimately arrive at $v_0$ at time $t$. The proof of Lemma \ref{lower technique} is provided in Section \ref{pf of Lem 3.8}.
			
			\vskip1mm
			Because the lemma describes a Markov chain on an undirected graph, we must lift the construction to the process $\overline{\mathcal{Z}}_j^n$ on the line digraph $\vec{l}(G)$.
			
			In the ergodicity part, we have constructed $(\overline{\mathcal{Z}}_j^n)_ {s_k\le j\le s_{k+1}-2\b-2}$ with transition probabilities $\nuu_j^n$ given in (\ref{Or-lowbdd-m2}).  
			Recall that $\overline{X}_j^n=(\overline{\mathcal{Z}}_j^n)^-$. Then $\overline{X}_j^n\in G_k$ for  $j\in[s_k,s_{k+1}-2\b-2]$, where $G_k=(V_k,E_k)$ with $V_k=\{v\in V:\ \exists e\in E_k,\ v\in e\}$. Thus $(\overline{X}_j^n)_ {s_k\le j\le s_{k+1}-2\b-2}$ is a Markov chain on $G_k$.   
			
			Applying Lemma \ref{lower technique}, we choose $\tilde{p}_j^n$ for $j\in[s_{k+1}-2\b,s_{k+1}-1]$ and extend $X_j^n$ to the interval $[s_{k+1}-2\b,s_{k+1}-1]$.   These transition probabilities force $\overline{X}_j^n$ to stay inside $G_k$  from time $s_{k+1}-2\b$ up to time $s_{k+1}$, and guarantee that  $\overline{X}_{s_{k+1}}^n=\zz_0^-$ for some fixed $\zz_0\in \overrightarrow{E_{k+1}\setminus E_k}$, 
			regardless of the state of  $\overline{X}_{s_{k+1}-2\b}^n$. 
			
			We then define the admissible control sequence for the connection part by 
			$$
			\tilde{\nuu}_j^n(\zz_1,\zz_2):=\tilde{p}_{j+1}^n(\zz_1^+,\zz_2^+),\ \ j\in[s_{k+1}-2\b-1,s_{k+1}-2],
			$$ 
			and set $\tilde{\nuu}_{s_{k+1}-1}^n(\zz,\zz_0)=1$ for every  $\zz\in V_{\vec{l}(G)}$ with $\zz^+=\zz_0^-$.
			This construction ensures that the connection part of $(\overline{\mathcal{Z}}_j^n)_ {s_k\le j\le s_{k+1}-2\b-2}$  satisfies all the required conditions.

			\begin{proof}[Proof of the lower bound of generalized Laplace principle]
				Combining the ergodicity part and the connection part, and applying Lemma \ref{ergodic theorem}, Corollary \ref{ergodic theorem-h} together with \eqref{condition of lower bound},  we obtain 
				\begin{align*}
					%  & \ \limsup_{n\to\infty}\vv^n_{\mathscr{S}_0}(\zz)\\
					\ \limsup_{n\to\infty}\vv^n_{\mathscr{S}_0}(\zz)\ \le & \lim_{n\to\infty} \overline{\E}_\zz\left\{\frac{1}{n}\sum_{k=1}^\b\Bigg[\sum_{i=1}^{j_k}\sum_{j=s_{k,i}}^{s_{k+1,i}-2\b-1}R\Big(\qq_k\big(\overline{\mathcal{Z}}^n_{j},\cdot\big)
					\Big\|\pp_{E_k}\big(\overline{\mathcal{Z}}^n_{j},\cdot\big)\Big)\right.\\
					&\ \ \ \ \ \ \ \ \ \ \ \ +\sum_{i=1}^{j_k-1}\sum_{j=s_{k,i+1}-2\b}^{s_{k+1,i}}R\Big(\hat{\nuu}_j^n\big(\overline{\mathcal{Z}}^n_{j},\cdot\big)
					\Big\|\pp_{E_k}\big(\overline{\mathcal{Z}}^n_{j},\cdot\big)\Big)\\
					&\left.\ \ \ \ \ \ \ \ \ \ \ \, +\sum_{j=s_{k+1}-2\b-1}^{s_{k+1}-1}
					R\Big(\tilde{\nuu}_j^n\big(\cdot\big|\overline{\mathcal{Z}}^n_{j},\overline{\mathcal{L}}^n_{j}\big)\Big\|\pp_{E_k}
					\big(\overline{\mathcal{Z}}^n_{j},\cdot\big)\Big)\Bigg]+h\big(\overline{\mathcal{L}}^n\big) \right\}\\
					=  &\,  \sum_{k=1}^\b r_k\sum_{i=1}^{j_k}C_{k,i}\int_{V_{\vec{l}(G)}}R(\qq_k\|\pp_{E_k})\ {\rm d}(\hat{\gamma}_k|_{S_{k,i}})+h(\hat{\gamma})\\
					\le&\  {\Lambda}_{\delta,\mathscr{S}_0}^\zz(\hat{\gamma})+h(\hat{\gamma})+\varepsilon\\
					\le&\, \inf_{\hat{\mu}\in\mathcal{C}l(\mathscr{S}_0)}\{{\Lambda}_{\delta,\mathscr{S}_0}^\zz(\hat{\mu})+h(\hat{\mu})\}+2\varepsilon,
				\end{align*}
				where $s_k=\sum_{j=1}^k n_j$ and $\lim\limits_{n\to\infty} \frac{n_k}{n}=r_k$ (for indices $k$ with $r_k=0$, the term $\frac{1}{n}\sum_{j=s_k}^{s_{k+1}-1}(\cdot)$ can be ignored, which allows us to assume $s_{k+1}-s_k>2\b$).  Since $\varepsilon$ is arbitrary, the proof is complete.
			\end{proof}
			
		}

		\subsection{Large deviation principle for empirical measures of   ORRWs}\label{sec 3.4}
		\noindent  {\color{black} In this subsection, we focus on  proving the LDP for empirical measures of   ORRWs, namely Theorem  \ref{ldp for Z}, Theorem \ref{rate function 0} and Theorem  \ref{I_delta}.
			
			\begin{proof}[Proof of Theorem  \ref{ldp for Z}.]  
				By the law of total probability and the finiteness of $\{\zz:\zz^-=x_0\}$, we have 
				\begin{equation}\label{E_v=min}
					\begin{aligned}
						&\ \lim_{n\to\infty}-\frac{1}{n}\log \E_{x_0}\exp\{-nh(\mathcal{L}^n)\}\\
						&=\ \lim_{n\to\infty}-\frac{1}{n}\log \Bigg(\sum_{\zz:\,\zz^-=x_0}\P_{x_0}(\mathcal{Z}_0=\zz)\E_{\zz}\exp\{-  nh(\mathcal{L}^n)\}\Bigg)\\
						&=\min_{\zz:\,\zz^-=x_0}\lim_{n\to\infty}-\frac{1}{n}\log \E_{\zz}\exp\{-nh(\mathcal{L}^n)\}\\
						&= \min_{\zz:\,\zz^-=x_0}\lim_{n\to\infty}W^n(\zz)\\
						&= \min_{\zz:\zz^-=x_0} \inf_{\hat{\mu}\in \mathcal{C}l(\mathscr{S})}\left\{\Lambda^\zz_{\delta,\mathscr{S}} (\hat{\mu})+h(\hat{\mu})\right\}. 
					\end{aligned}
				\end{equation}
				Recall the definition of $\Lambda_{\delta,\mathscr{S}}(\hat{\mu})$ in \eqref{rate function}  and $\Lambda^\zz_{\delta,\mathscr{S}}(\hat{\mu})$ in \eqref{Lambda^z}. We obtain
				$$
				\min_{\zz:\zz^-=x_0} \inf_{\hat{\mu}\in \mathcal{C}l(\mathscr{S})}\big\{\Lambda^\zz_{\delta,\mathscr{S}} (\hat{\mu})+h(\hat{\mu})\big\}= \inf_{\hat{\mu}\in \mathcal{C}l(\mathscr{S})}\left\{\Lambda_{\delta,\mathscr{S}} (\hat{\mu})+h(\hat{\mu})\right\}. 
				$$
				Furthermore, since Laplace principle implies the large deviation principle, see Theorem \ref{thm-LP-LDP}, the proof is complete.
			\end{proof}

			\begin{proof}[Proof of Theorem \ref{rate function 0}.]  For any subset  $V'\subseteq V$, 
				$$
				T(\mathcal{L}^n)(V')=\frac{1}{n}\sum_{j=0}^{n-1}\dd_{\mathcal{Z}_j}(\{\zz:\zz^-\in V'\})=\frac{1}{n}\sum_{j=0}^{n-1}\dd_{\mathcal{Z}_j^-}(V') =L^n(V'),\ n\geq 1.
				$$
				Since the empirical measure process  $\left(\mathcal{L}^n\right)_{n\geq 1}$ of $(\mathcal{Z}_n)_{n\geq 0}$ satisfies the LDP, the continuity of $T$ established in Proposition  \ref{continuity of map T} together with the contraction principle completes the proof.
			\end{proof}

			Applying  Theorem \ref{ldp for Z} and the contraction principle, we obtain the following corollary. 
			
			\begin{corollary}\label{coro-ldp-edge}

				\begin{itemize}
					\item[(a)] Consider the empirical measures  on edges:
					\[
					L_{(1)}^n(e)=\frac{1}{n}\sum_{k=0}^{n-1}\dd_{X_{k}X_{k+1}}(e),\ \text{ for\ }\ n\geq 1,\ e\in E.
					\]
					It satisfies an LDP with a good rate function
					$$
					\inf_{\hat{\mu}\in\mathscr{P}(V_{\vec{l}(G)}):T_1(\hat{\mu})=\nu}\Lambda_\delta(\hat{\mu}),\ \  \nu\in\mathscr{P}(E),
					$$
					where $T_1(\hat{\mu})(e)=\hat{\mu}|_E(e)$ for all $e\in E$.
					\item[(b)] Consider the pair empirical occupation measures introduced in \cite{CVC2022}:
					\[
					L_{(2)}^n(u,v)=\frac{1}{n}\sum_{k=0}^{n-1}\dd_{X_{k}}(u)\dd_{X_{k+1}}(v),\   \text{for }\  n\geq 1,\ u,v\in V.
					\]
					It satisfies an LDP with a good rate function $\Lambda_\delta(T_2(\nu))$ for $\nu\in\mathscr{P}(V\times V)$, where $T_2(\nu)(\zz)=\nu((\zz^-,\zz^+))$ for $\zz\in V_{\vec{l}(G)}$.
				\end{itemize}
			\end{corollary}
			\begin{proof}  Observe that each $T_i(\mathcal{L}^n)=L_{(i)}^n$ is continuous and satisfies $T_i(\mathcal{L}^n)=L_{(i)}^n$ for $i=1,2$, and that $T_2$ is a bijection.  Hence the statement follows directly from the LDP for $\big(\mathcal{L}^n\big)_{n\geq 1}$ established in Theorem \ref{ldp for Z}  
				together with the contraction principle.
			\end{proof}
			
			Now we proceed to prove the main result of this paper, Theorem \ref{I_delta}. The key step is to obtain a concise expression for the rate function $I_\delta$ in \eqref{rate func.I}, which relies on  the following lemma.
			
			\begin{lemma}\label{inf attainment}
				For a fixed $\hat{\mu}\in\mathscr{P}(V_{\vec{l}(G)})$ and a subset $E'\subseteq E$, let $\nu=T(\hat{\mu})$, where $T$ is the mapping defined in  Theorem \ref{rate function 0}. Then 
				\[
				\inf_{\qq\in\mathscr{T}_{\vec{l}(G)}:\,\hat{\mu} \qq=\hat{\mu}}\int_{V_{\vec{l}(G)}}R(\qq\|\pp_{E'})\ {\rm d}\hat{\mu}=\inf_{{q}\in\mathscr{T}_G:\,\forall \zz\in V_{\vec{l}(G)}, \nu(\zz^-){q}(\zz^-,\zz^+)=\hat{\mu}(\zz)}\int_{V}R({q}\|{p}_{E'})\ {\rm d}\nu.
				\]
			\end{lemma}
			
			\vskip2mm
			Observe  that for every  $E'\subseteq E$, and every directed edge $\zz_1\to \zz_2\in  E_{\vec{l}(G)}$,
			\[
			\pp_{E'}(\zz_1,\zz_2)={p}_{E'}(\zz_2^-,\zz_2^+).
			\]
			The main idea for proving Lemma \ref{inf attainment} is to show that the infimum on the left‑hand side is attained by a transition probability  $\qq$ satisfying   $\qq(\zz_1,\zz)=\qq(\zz_2,\zz)$ for all $\zz_1,\zz_2\to \zz$ and for every $1\le k\le \b$.
			This proof is provided  in Section  \ref{pf lem 3.10}.

			\begin{proof}[Proof of Theorem \ref{I_delta}] We first verify the expression \eqref{simpler expression of I} of the rate function $I_\delta(\nu)$ in the following steps.
				
				\noindent{\small\bf  \emph{Step 1.1}}.
				For each $\hat{\mu}\in\mathscr{P}(V_{\vec{l}(G)})$, denote 
				\begin{align*}
					\check{\mathscr{A}}(\hat{\mu}):=& 
					\Bigg\{ (\hat{\mu}_k,r_k,E_k)_{1\leq k\leq \b}:\ \{E_j\}_{1\leq j\leq \b}\in\mathscr{E}; \, r_k\ge 0, \,\sum_{k=1}^{\b}r_k=1; \nonumber\\
					&\hskip 15mm \hat{\mu}_k\in\mathscr{P}(V_{\vec{l}(G)}),\ {\rm{supp}}(\hat{\mu}_k)\subseteq \vec{E}_k,\ \sum_{k=1}^{\b} r_k \hat{\mu}_k=\hat{\mu}\Bigg\}. 
				\end{align*}
				By Theorem \ref{rate function 0} and Lemma \ref{inf attainment}, the rate function $I_\delta$ can be rewritten  as follows: for every $\nu\in\mathscr{P}(V)$,
				\begin{align}
					I_{\delta}(\nu)=&\, \inf\limits_{\hat{\mu}\in\mathscr{P}(V_{\vec{l}(G)}): T(\hat{\mu})=\nu} \Lambda_\delta(\hat{\mu})\nonumber\\
					=&\,  \inf\limits_{\hat{\mu}\in\mathscr{P}(V_{\vec{l}(G)}): T(\hat{\mu})=\nu} \inf_{(\hat{\mu}_k,r_k,E_k,\qq_k)_{k}\in\hat{\mathscr{A}}(\hat{\mu})}\sum_{k=1}^{\b}r_k\int_{V_{\vec{l}(G)}}R(\qq_{k}\|\pp_{E_k})\ {\rm d}\hat{\mu}_{k}\nonumber\\
					=&\, 	 \inf\limits_{\hat{\mu}\in\mathscr{P}(V_{\vec{l}(G)}): T(\hat{\mu})=\nu} \inf_{\stackrel{(\hat{\mu}_k,r_k,E_k)\in \check{\mathscr{A}}(\hat{\mu}):}{q_k\in \mathscr{T}_G,T(\hat{\mu}_k)(\zz^-){q}_{k}(\zz^-,\zz^+)=\hat{\mu}_k(\zz)}}\sum_{k=1}^{\b}r_k\int_{V}R({q}_{k}\|{p}_{E_k})\ {\rm d}T(\hat{\mu}_k)\nonumber\\
					=&\inf_{\stackrel{(\hat{\mu}_k,r_k,E_k)\in \check{\mathscr{A}}(\hat{\mu}):}{q_k\in \mathscr{T}_G,T(\hat{\mu}_k)(\zz^-){q}_{k}(\zz^-,\zz^+)=\hat{\mu}_k(\zz), \sum_{k=1}^\b r_kT(\hat{\mu}_k)=\nu}}\sum_{k=1}^{\b}r_k\int_{V}R({q}_{k}\|{p}_{E_k})\ {\rm d}T(\hat{\mu}_k). \label{step 1 in I_delta proof} 
				\end{align}

				\noindent{\small\bf  \emph{Step 1.2}}. Recall the definition of $\mathscr{A}(\nu)$,
				\begin{align}
					{\mathscr{A}}(\nu)=&\,\Big\{ (\nu_k,r_k,E_k,{q}_{k})_{1\leq k\leq \b}:\ \{E_k\}_{1\leq k\leq \b}\in\mathscr{E},\ r_k\ge 0,\ \sum_{k=1}^{\b}r_k=1,\nonumber\\
					&\hskip 1cm \nu_k\in \mathscr{P}(V_k),\ 
					{q}_k\in\mathscr{T}(E_k),\ \text{such that }\nu_k{q}_k=\nu_k \text{ and } \sum_{k=1}^{\b}r_k \nu_k=\nu  \Big\}.%\label{hat_A}
				\end{align}
				We aim to show that the infimum condition in (\ref{step 1 in I_delta proof}) implies $(\nu_k,r_k,E_k,{q}_k)_{1\le k\le \b}\in {\mathscr{A}}(\nu)$.  Equivalently, define
				\begin{align*}
					\breve{\mathscr{A}}(\nu):=&\, \Big\{  (\hat{\mu}_k,r_k,E_k, q_k)_{1\le k\le \b}: \    (\hat{\mu}_k,r_k,E_k)\in \check{\mathscr{A}}(\hat{\mu}), \, q_k\in \mathscr{T}_G,  \nonumber\\
					&  \hskip1cm T(\hat{\mu}_k)(\zz^-){q}_{k}(\zz^-,\zz^+)=\hat{\mu}_k(\zz),\, \sum_{k=1}^\b r_kT(\hat{\mu}_k)=\nu\Big\},  
				\end{align*}
				we must prove that   for every  $(\hat{\mu}_k,r_k,E_k, q_k)_{1\le k\le \b}\in  \breve{\mathscr{A}}(\nu)$, 
				setting  $\nu_k=T(\hat{\mu}_k)$ yields $(\nu_k,r_k,E_k, q_k)_{1\le k\le \b}\in  {\mathscr{A}}(\nu)$. Consequently, 
				\begin{equation}\label{proof 2.6-1.5}
					{\rm (\ref{step 1 in I_delta proof})}\ge \inf_{(\nu_k,r_k,E_k,{q}_k)_{k}\in{\mathscr{A}}(\nu)}
					\sum_{k=1}^{\b}r_k\int_{V}R({q}_{k}\|{p}_{E_k})\ {\rm d}\nu_{k}.
				\end{equation} 
				% {\color{yellow}
					%\begin{align}\label{proof 2.6-1}
					% \Big\{  (\hat{\mu}_k,r_k,E_k, q_k)_{1\le k\le \b}: &\    (\hat{\mu}_k,r_k,E_k)\in \check{\mathscr{A}}(\hat{\mu}), \, q_k\in \mathscr{T}_G,  \nonumber\\
					%  &\  T(\hat{\mu}_k)(\zz^-){q}_{k}(\zz^-,\zz^+)=\hat{\mu}_k(\zz),\, \sum_{k=1}^\b r_kT(\hat{\mu}_k)=\nu\Big\} \subset {\mathscr{A}}(\nu),
					% \end{align}}
			
			To verify the above claim, we note the following equivalence.
			%\begin{minipage}
			\vskip2mm {\small\emph{
					For $\hat{\mu}_k\in\mathscr{P}(V_{\vec{l}(G)})$ with ${\rm{supp}}(\hat{\mu}_k)\subseteq \vec{E}_k$, set $\nu_k=T(\hat{\mu}_k)$ for $k=1,\dots,\b$. Observe that under the stationary condition: $\nu_k{q}_{k}=\nu_k$,   the statement
					\begin{equation}\label{proof 2.6-2}
						``{q}_{k}\in\mathscr{T}(E_k)\ \ \text{and} \ \ {\rm supp}(\nu_k)\subseteq V_k"
					\end{equation} is equivalent to 
					\begin{equation}\label{proof 2.6-3}
						``\nu_k(x){q}_{k}(x,y)>0\Rightarrow xy\in E_k ".
			\end{equation}}}
			
			\vskip2mm
			
			For every $(\hat{\mu}_k,r_k,E_k. q_k)_{1\le k\le \b}\in \breve{\mathscr{A}}(\nu)$,  the following properties  are satisfied: 
			\begin{itemize}
				\item Summing $\nu_k(\zz^-){q}_k(\zz^-,\zz^+)=\hat{\mu}_k(\zz)$ over all $\zz$ with $\zz^+=y$ for a fixed $y\in V$ yields $\hat{\nu}_k{q}_k=\hat{\nu}_k$;
				\item If  $\nu_k(x){q}_k(x,y)>0$, 	then for the directed edge $\zz$ with $\zz^-=x,\,\zz^+=y$, we have $\hat{\mu}_k(\zz)>0$; hence  $\zz\in{\rm supp}(\hat{\mu}_k)$, which force $xy\in E_k$, i.e., $\nu_k(x){q}_k(x,y)>0\Rightarrow xy\in E_k$.
			\end{itemize}
			Therefore, by  the equivalence between the statement \eqref{proof 2.6-2} and \eqref{proof 2.6-3},  we conclude that $(\nu_k,r_k,E_k, q_k)_{1\le k\le \b}\in  {\mathscr{A}}(\nu)$; consequently \eqref{proof 2.6-1.5} follows. 
			
			\vskip2mm
			\noindent{\small\bf  \emph{Step 1.3}}.  Let $(\nu_k,r_k,E_k,{q}_k)_{1\le k\le \b}\in{\mathscr{A}}(\nu)$ and define $\hat{\mu}_k(\zz)=\nu_k(\zz^-){q}_k(\zz^-,\zz^+)$.  
			Conversely, we aim to verify 
			\begin{equation}\label{proof 2.6-4}
				(\hat{\mu}_k,r_k,E_k,{q}_k)_{1\le k\le \b}\in{\breve{\mathscr{A}}}(\nu).
			\end{equation}
			%{\color{yellow}
				%\begin{align}\label{proof 2.6-4}
				%{\mathscr{A}}(\nu)\subset &\,  \Big\{  (\hat{\mu}_k,r_k,E_k, q_k)_{1\le k\le \b}: \    (\hat{\mu}_k,r_k,E_k)\in \check{\mathscr{A}}(\hat{\mu}), \, q_k\in \mathscr{T}_G,  \nonumber\\
				%&\ \hskip1.8cm T(\hat{\mu}_k)(\zz^-){q}_{k}(\zz^-,\zz^+)=\hat{\mu}_k(\zz), \sum_{k=1}^\b r_kT(\hat{\mu}_k)=\nu\Big\}.
				%\end{align}}
				Equivalently,  every $(\nu_k,r_k,E_k,{q}_k)_{1\le k\le\b}\in{\mathscr{A}}(\nu)$ satisfies the infimum condition in \eqref{step 1 in I_delta proof}, 
				which implies 
				\begin{equation}\label{proof 2.6-5}
					{\rm (\ref{step 1 in I_delta proof})}\le \inf_{(\hat{\nu}_k,r_k,E_k,{q}_k)_{k}\in{\mathscr{A}}(\nu)}
					\sum_{k=1}^{\b}r_k\int_{V}R({q}_{k}\|{p}_{E_k})\ {\rm d}\nu_{k}.
				\end{equation}
				
				For each $(\nu_k,r_k,E_k,{q}_k)_{1\le k\le \b}\in{\mathscr{A}}(\nu)$, define $\hat{\mu}_k(\zz)=\nu_k(\zz^-){q}_k(\zz^-,\zz^+)$. Then the following properties hold: 
				\begin{itemize}
					\item For an arbitrary $x\in V$, summing $\hat{\mu}_k(\zz)=\nu_k(\zz^-){q}_k(\zz^-,\zz^+)$ over all $\zz$ with $\zz^-=x$ and using  $\sum\limits_{y\in V}{q}_k(x,y)=1$  yields $\nu_k(x)=T(\hat{\mu}_k)(x)$.
					\item  For every $\zz\in {\rm supp}(\hat{\mu}_k)$, $\hat{\mu}_k(\zz)=\nu_k(\zz^-){q}_k(\zz^-,\zz^+)>0$ implies 
					$\zz|_{E}\in E_k$; i.e., ${\rm supp}(\hat{\mu}_k|_E)\subset E_k$. Consequently, ${\rm supp}(\hat{\mu}_k)\subseteq \vec{E}_k$.
				\end{itemize}
				Thus, \eqref{proof 2.6-4} holds, and therefore \eqref{proof 2.6-5} is valid. 
				
				Combining the two inequalities \eqref{proof 2.6-1.5} and \eqref{proof 2.6-5} yields the desired expression \eqref{simpler expression of I} for the rate function $I_\delta(\nu)$. 
				
				\vskip2mm
				
				We next prove that  $I_\delta(\nu)<\infty$ if and only if $\nu$ is invariant for some ${q}\in\mathscr{T}_G$.
				\vskip1mm
				\noindent{\small\bf  \emph{Step 2.1}}.  We show that $I_\delta(\nu)<\infty$, then $\nu$ is invariant for some ${q}\in\mathscr{T}_G$. 
				
				\vskip2mm
				
				Since $I_\delta(\nu)<\infty$,  by  Theorem \ref{rate function 0} guarantees the existence of a constant $M<\infty$ such that 
				$$
				\{\hat{\mu}: \Lambda_\delta(\hat{\mu})\le M\}\neq \emptyset.
				$$
				Moreover,  Proposition \ref{lower semicontinuous}  implies that there exists a $\hat{\mu}\in \mathscr{P}(V_{\vec{l}(G)})$ with $T(\hat{\mu})=\nu$  satisfying $I_\delta (\nu)=\Lambda_\delta(\hat{\mu})<\infty$. Hence, Proposition \ref{prop domain of rate function} implies that $\hat{\mu}$ is  invariant for some $\qq\in\mathscr{T}_{\vec{l}(G)}$. That is, for every vertex $x$ of the graph $G$ and  for any $\zz_1,\zz_2$ with $\zz_1^+=\zz_2^-=x$, 
				\[
				\sum_{\zz^-=x}\qq(\zz_1,\zz)=1,\ \ \sum_{\zz^+=x}\qq(\zz,\zz_2)\hat{\mu}(\zz)=\hat{\mu}(\zz_2).
				\]
				Consequently,
				\begin{align*}
					\sum_{\zz_2^-}\hat{\mu}(\zz_2)&=\sum_{\zz_2^-=x}\sum_{\zz_1^+=x}\qq(\zz_1,\zz_2)\hat{\mu}(\zz_1)\\
					&=\sum_{\zz_1^+=x}\sum_{\zz_2^-=x}\qq(\zz_1,\zz_2)\hat{\mu}(\zz_1)\\
					&=\sum_{\zz_1^+=x}\hat{\mu}(\zz_1).
				\end{align*}
				For all adjacent vertices $x,y$ of $G$ and any directed edge $\zz_0$ with $\zz_0^-=x,\zz_0^+=y$, define
				\[
				{q}(x,y)=\frac{\hat{\mu}(\zz_0)}{\sum_{\zz^+=x}\hat{\mu}(\zz)}.
				\]
				Then  ${q}\in\mathscr{T}_G$ and $\nu$ is invariant for ${q}$.\\
				
				\vskip1mm
				\noindent{\small\bf  \emph{Step 2.2}}. We verify that  if $\nu$ is invariant for some ${q}\in\mathscr{T}_G$, then $I_\delta(\nu)<\infty$.
				
				\vskip2mm
				Note that by \cite[Propositions 8.3.3 and 8.6.1]{DE1997},  the rate function of LDP  for  the empirical measure of the simple random walk on $G$, 
				$$
				I_1(\nu)= \inf_{{q}:\nu{q}=\nu}\int_{V}R({q}\| {p}_E)\ {\rm d}\nu,
				$$
				which is finite whenever $\nu$ invariant for some ${q}\in\mathscr{T}_G$.  
				
				Using \eqref{simpler expression of I},  take $(\nu_k, r_k, E_k, q_k)_{1\le k\le \b}\in \mathscr{A}(\nu)$ with  $r_1=\cdots=r_{\b-1}=0$, $r_{\b}=1$ and $\nu_{\b}=\nu$, $q_{\b}=q$.  Then, for every $\delta>0$,
				\[
				I_\delta(\nu)\le I_1(\nu)<\infty.
				\]
			\end{proof}
		}
		
		\subsection{Analytic properties of rate functions $\Lambda_\delta$ and $I_\delta$}\label{sec 3.5}
		\noindent {\color{black}  
			The main purpose of this section is to prove Theorem  \ref{phase for rate function},  which concerns certain analytic properties of the rate function $I_\delta$. We first derive the corresponding properties of the rate function $\Lambda_\delta$  for the LDP of  the empirical measure $(\mathcal{L}^n)_{n\ge 1}$ of the lifted process $(\mathcal{Z}^n)_{n\ge 1}$ on the digraph $\vec{l}(G)$, and then use these to establish the properties of $I_\delta$.

			Recall the rate function $\Lambda_\delta:=\Lambda_{\delta,\mathscr{S}}$ in \eqref{rate function}, which is given as the infimum of a convex combination of the integral of relative entropy $R\big(\qq_k(\xx,\cdot)\|\pp_{E_k}(\xx,\cdot))\big)$ for $k=1,2,\dots,\b$. We begin by  examining the monotonicity of $R(\qq(\xx,\cdot)\|\pp_{E'}(\xx,\cdot))$ with respect to  $\delta$ for a fixed vertex $\xx\in \vec{l}(G)$ and a subset $E'\subseteq E$.

			\begin{lemma}\label{R(|)-R(|)}
				For all positive $\delta_1<\delta_2$, any edge set $E'$, any $\xx\in \vec{E}'$ and any transition probability $\qq(\xx,\cdot)$ supported on $\vec{E}'$, then for $\yy\in \vec{E}'$,
				\begin{equation}\label{3.4-entpy-1}
					\log \frac{\pp_{E'}(\xx,\yy;\delta_1)}{\pp_{E'}(\xx,\yy;\delta_2)}\left\{
					\begin{aligned}
						&<0, &\xx\in\partial(\vec{E}'),\\
						&=0, &\xx\notin\partial(\vec{E}').
					\end{aligned}
					\right.
				\end{equation}
				{\color{black} See the definition of $\partial(\vec{E}')$ in Definition \ref{line digraph} (d).} Moreover, 	
				\begin{equation}
					R\big(\qq(\xx,\cdot)\big\|\pp_{E'}(\xx,\cdot;\delta_2)\big)\le R\big(\qq(\xx,\cdot)\big\|\pp_{E'}(\xx,\cdot;\delta_1)\big). \label{3.4-entpy-2}
				\end{equation}
			\end{lemma}

			\begin{proof} 
				Note that for $\yy\in \vec{E}'$,
				\begin{equation}
					\frac{\pp_{E'}(\xx,\yy;\delta_1)}{\pp_{E'}(\xx,\yy;\delta_2)}=\frac{k(\xx)\delta_2+d(\xx)-k(\xx)}{\delta_2}\cdot\frac{\delta_1}{k(\xx)\delta_1+d(\xx)-k(\xx)},\label{p(delta_1)/p(delta_2)}
				\end{equation}
				where $k(\xx)$ (resp. $d(\xx)$) denotes the number of neighbors of $\xx$ in $\vec{E}'$ (resp. $\vec{l}(G)$). For $\xx\in\partial(\vec{E}')$, we have $k(\xx)<d(\xx)$. Since $\delta_1<\delta_2$,
				\[
				\frac{\pp_{E'}(\xx,\yy;\delta_1)}{\pp_{E'}(\xx,\yy;\delta_2)}<1.
				\]
				For $\xx\notin\partial(\vec{E}')$ it holds that $k(\xx)=d(\xx)$, which implies
				\[
				\frac{\pp_{E'}(\xx,\yy;\delta_1)}{\pp_{E'}(\xx,\yy;\delta_2)}=1.
				\]
				Hence \eqref{3.4-entpy-1} follows.
				
				Furthermore, observe that
				\begin{equation}\label{3.4-entpy-3}
					R\big(\qq(\xx,\cdot)\big\|\pp_{E'}(\xx,\cdot;\delta_2)\big)-R\big(\qq(\xx,\cdot)\big\|\pp_{E'}(\xx,\cdot;\delta_1)\big)=\int_{V_{\vec{l}(G)}}\log{\frac{\pp_{E'}(\xx,\yy;\delta_1)}{\pp_{E'}(\xx,\yy;\delta_2)}}\ \qq(\xx,{\rm d}\yy).
				\end{equation}
				Thus, \eqref{3.4-entpy-1} immediately implies \eqref{3.4-entpy-2}.
			\end{proof}

			Theorem \ref{phase for rate function} is a corollary of the following proposition.
			
			\begin{proposition}\label{phase for rate function 0}
				The rate function $\Lambda_\delta$ satisfies the following properties:
				\begin{itemize}
					\item[(a)]  For all $\delta\in (0,1)$,  $\Lambda_\delta(\hat{\mu})= \Lambda_1(\hat{\mu})$ for every $\hat{\mu}\in  \mathscr{P}(V_{\vec{l}(G)})$.
					\item[(b)]  $\Lambda_\delta$ is monotonically decreasing in $\delta\ge 1$, i.e.,
					$$
					\text{for all} \ \delta_2>\delta_1\ge 1, \   \Lambda_{\delta_1}(\hat{\mu})\ge \Lambda_{\delta_2}(\hat{\mu}),\ \forall\hat{\mu}\in \mathscr{P}(V_{\vec{l}(G)});
					$$ and moreover there exists a $\hat{\mu}_0\in \mathscr{P}(V_{\vec{l}(G)})$ (depending on $\delta_1$ and $\delta_2$) such that $\Lambda_{\delta_1}(\hat{\mu}_0)> \Lambda_{\delta_2}(\hat{\mu}_0)$.
					\item[(c)]  $\Lambda_\delta$ is uniformly continuous in $\delta$ in the sense that
					$$
					\lim_{|\delta_2-\delta_1|\to 0}\sup_{\hat{\mu}\in\mathscr{P}(V_{\vec{l}(G)}):\,\Lambda_1(\hat{\mu})<\infty}|\Lambda_{\delta_2}(\hat{\mu})-\Lambda_{\delta_1}(\hat{\mu})|=0.
					$$
					\item[(d)]  $\Lambda_\delta(\hat{\mu})$ is not differentiable with respect to $\delta$ at $\delta=1$ for some $\hat{\mu}\in\mathscr{P}(V_{\vec{l}(G)})$.
				\end{itemize}
			\end{proposition}
			
			\begin{proof}
				We divide the proof into the following four parts.
				\vskip 2mm
				
				\noindent{\bf(i)\ $\Lambda_\delta=\Lambda_1$} {\bf\emph{for} $\delta\le 1$}.
				
				Denote by $\pp$ the transition probability of the  SRW on $\vec{l}(G)$. For every sequence $\{E_k\}_{1\le k\le \b}$ and $l=1,\dots,\b$, we notice  the fact that $\pp_{E_l}(\xx, \cdot; 1)=\pp(\xx, \cdot)$. By Lemma \ref{R(|)-R(|)}, setting $\delta_1=\delta$ and $\delta_2=1$,  we have
				\[
				R\big(\qq_{l}(\xx,\cdot)\big\|\pp_{E_l}(\xx,\cdot;\delta_2)\big)\ge R\big(\qq_{l}(\xx,\cdot)\big\|\pp_{E_l}(\xx,\cdot;1)\big)=R\big(\qq_{l}(\xx,\cdot)\big\|\pp(\xx,\cdot)\big).
				\]
				Recall the definition of $\Lambda_\delta(\hat{\mu})$ in \eqref{rate function},
				\begin{align*}
					{\Lambda}_{\delta}(\hat{\mu}):=\inf_{(\hat{\mu}_k,r_k,E_k,\qq_k)_{k}\in\hat{\mathscr{A}}(\hat{\mu})}
					\sum_{k=1}^{\b}r_k\int_{V_{\vec{l}(G)}}R(\qq_{k}\|\pp_{E_k})\ {\rm d}\hat{\mu}_{k}.%\label{rate function-3.4}
				\end{align*}
				In particular,	
				\begin{align*}
					{\Lambda}_{1}(\hat{\mu}):=\inf_{(\hat{\mu}_k,r_k,E_k,\qq_k)_{k}\in\hat{\mathscr{A}}(\hat{\mu})}
					\sum_{k=1}^{\b}r_k\int_{V_{\vec{l}(G)}}R(\qq_{k}\|\pp)\ {\rm d}\hat{\mu}_{k}.\label{rate function-3.4}
				\end{align*}		
				For every $(\hat{\mu}_k, r_k, E_k, \qq_k)\in \hat{\mathscr{A}}(\hat{\mu})$, 
				$$
				\sum_{k=1}^\b r_k\int_{V_{\vec{l}(G)}}R(\qq_{k}\|\pp_{E_k})\ {\rm d}\hat{\mu}_{k}\ge \sum_{k=1}^\b r_k\int_{V_{\vec{l}(G)}}R(\qq_{k}\|\pp)\ {\rm d}\hat{\mu}_{k}\ge \Lambda_1(\hat{\mu}),
				$$
				which implies $\Lambda_{\delta}\ge \Lambda_1$.

				To show the reverse inequality $\Lambda_{\delta}\le \Lambda_1$, note that $E_{\b}=E$, $\pp_E=\pp$. Taking  $r_1=\cdots=r_{\b-1}=0$ and $r_{\b}=1$, we obtain 
				\begin{eqnarray*}
					\Lambda_\delta(\hat{\mu})&\le &\inf_{\begin{tiny}\begin{aligned}&(\hat{\mu}_{k},r_k,E_k)_{k}\in \hat{\mathscr{A}}(\hat{\mu},\mathscr{S}),r_{\b}=1\\
								&\qq_{k}\in\mathscr{T}_{\vec{l}(G)}:\hat{\mu}_{E_k}\qq_{k}=\hat{\mu}_{E_k}\end{aligned}\end{tiny}}
					\sum_{k=1}^{\b}r_k\int_{V_{\vec{l}(G)}}R(\qq_{k}\|\pp_{E_k})\ {\rm d}\hat{\mu}_{k}\\
					&=&\inf_{\qq\in\mathscr{T}_{\vec{l}(G)},\hat{\mu} \qq=\hat{\mu}}\int_{V_{\vec{l}(G)}} R(\qq\|\pp)\ \hat{\mu}({\rm d}\xx)=\Lambda_1(\hat{\mu}).
				\end{eqnarray*}
				The last equality follows from \cite[Theorem 8.4.3]{DE1997} or  Theorem \ref{thm-ldp-markov}, because when $\delta=1$,  the process $(\mathcal{Z}_n)_{n\geq 0}$ is exactly the SRW on $\vec{l}(G)$. Consequently, $\Lambda_{\delta}=\Lambda_1$.
				
				\vskip2mm
				\noindent {\bf (ii)\ \emph{Monotonicity for}  $\delta\ge1$}.
				
				For all $1\le \delta_1<\delta_2$, Lemma \ref{R(|)-R(|)} gives
				\begin{equation*}	
					R\big(\qq_{l}(\xx,\cdot)\big\|\pp_{E_l}(\xx,\cdot;\delta_2)\big)\le R\big(\qq_{l}(\xx,\cdot)\big\|\pp_{E_l}(\xx,\cdot;\delta_1)\big).
				\end{equation*}
				Hence $\Lambda_{\delta_2}(\hat{\mu})\le\Lambda_{\delta_1}(\hat{\mu})$ for all $\hat{\mu}$. We now prove the strict monotonicity for some $\hat{\mu}$.  
				
				Note that   for $\xx\in\partial (\vec{E}_l)$, $1\le k(\xx)<d(\xx)\le \b$ for $\xx\in\partial (\vec{E}_l)$,  where $k(\xx)$ (resp. $d(\xx)$) is the number of neighbors of $\xx$ in $\vec{E}_l$ (resp. $\vec{l}(G)$).
				Using the equality \eqref{p(delta_1)/p(delta_2)} from the proof of Lemma \ref{R(|)-R(|)}, for $\xx\in\partial (\vec{E}_l)$ and  $\yy\in \vec{E}_l$,
				\begin{equation}
					\frac{\pp_{E_l}(\xx,\yy;\delta_1)}{\pp_{E_l}(\xx,\yy;\delta_2)}=\frac{k(\xx)\delta_2+d(\xx)-k(\xx)}{\delta_2}\cdot\frac{\delta_1}{k(\xx)\delta_1+d(\xx)-k(\xx)}\le\frac{\delta_1}{\delta_2}\cdot\frac{(\b-1)\delta_2+1}{(\b-1)\delta_1+1}<1.\label{(3.12)}
				\end{equation}
				
				Specifically we choose an edge $e_0=x_0x_1$ adjacent to the starting point $x_0$ of $(X_n)_{n\geq 0}$ and a probability measure $\hat{\mu}_0\in\mathscr{P}(V_{\vec{l}(G)})$ such that 
				$$
				\hat{\mu}_0(\overrightarrow{x_0x_1})=\hat{\mu}_0(\overrightarrow{x_1x_0})=\frac{1}{2}.
				$$
				We now show that the infima  in $\Lambda_{\delta_1},\Lambda_{\delta_2}$ are both attained at the choice $r_1=1$, $E_1=\{e_0\}$. 
				
				Note that
				$
				\pp_{E}(\xx,\yy)\le \pp_{\{e_0\}}(\xx,\yy)\ \mbox{for}\ \xx,\yy\in \{\overrightarrow{x_0x_1},\overrightarrow{x_1x_0}\}.
				$
				For any $(\hat{\mu}_k,r_k,E_k,\qq_k)_{1\le k\le \b}\in\hat{\mathscr{A}}(\hat{\mu}_0)$,   the convexity of the  rate function $\inf_{\qq:\,\hat{\mu} \qq=\hat{\mu}}\int_{V_{\vec{l}(G)}} R(\qq\|\pp)\ {\rm d}\hat{\mu}$ (see \cite[Proposition 8.5.2]{DE1997}) yields
				\begin{align*}
					\sum_{k=1}^\b r_k\int_{V_{\vec{l}(G)}} R(\qq_k\|\pp_{E_k})\ {\rm d}\hat{\mu}_k&\ge \sum_{k=1}^\b r_k\int_{V_{\vec{l}(G)}} R\left(\qq_k\|\pp_{\{e_0\}}\right)\ {\rm d}\hat{\mu}_k\\
					&\ge \sum_{k=1}^\b r_k\inf_{\qq_k:\,\hat{\mu}_k \qq_k=\hat{\mu}_k}\int_{V_{\vec{l}(G)}} R\left(\qq_k\|\pp_{\{e_0\}}\right)\ {\rm d}\hat{\mu}_k\\
					&\ge \inf_{\qq:\,\hat{\mu}_0 \qq=\hat{\mu}_0}\int_{V_{\vec{l}(G)}} R\left(\qq\|\pp_{\{e_0\}}\right)\ {\rm d}\hat{\mu}_0,
				\end{align*}
				which confirms that the infimum is attained at $r_1=1$, $E_1=\{e_0\}$.
				
				\vskip1mm
				
				Since $\hat{\mu}_0(\partial (\{\overrightarrow{x_0x_1},\overrightarrow{x_1x_0}\}))\ge \frac{1}{2}$, and assuming  the infimum  $\Lambda_{\delta_1}$ is attained at $\qq\in\mathscr{T}_{\vec{l}(G)}$, we obtain 
				\begin{align*}
					\Lambda_{\delta_2}(\hat{\mu}_0)-\Lambda_{\delta_1}(\hat{\mu}_0)&\le\int_{V_{\vec{l}(G)}}\Big[R\big(\qq(\xx, \cdot)\big\|\pp_{\{e_0\}}(\xx,\cdot;\delta_2)\big)-
					R\big(\qq(\xx,\cdot)\big\|\pp_{\{e_0\}}(\xx,\cdot;\delta_1)\big)\Big]\ {\rm d}\hat{\mu}_{0}\\
					&=\int_{V_{\vec{l}(G)}}\Bigg[\int_{V_{\vec{l}(G)}} \log{\frac{\pp_{\{e_0\}}(\xx,\yy;\delta_1)}{\pp_{\{e_0\}}(\xx,\yy;\delta_2)}}\ \qq(\xx,{\rm d}\yy)\Bigg]{\rm d}\hat{\mu}_{0}\\
					&\le \log\left\{\frac{\delta_1}{\delta_2}\cdot\frac{(\b-1)\delta_2+1}{(\b-1)\delta_1+1}\right\}\,\hat{\mu}_{0}\big(\partial (\{\overrightarrow{x_0x_1},\overrightarrow{x_1x_0}\})\big)\\
					&\le \frac{1}{2}\log\left\{\frac{\delta_1}{\delta_2}\cdot\frac{(\b-1)\delta_2+1}{(\b-1)\delta_1+1}\right\}\\
					&< 0.
				\end{align*}
				\vskip2mm
				
				\noindent{\bf (iii)\ \emph{Continuity}}.
				
				When studying uniform continuity, we bound $\sup_{\hat{\mu}:\,\Lambda_1(\hat{\mu})<\infty}|\Lambda_{\delta_1}(\hat{\mu})-\Lambda_{\delta_2}(\hat{\mu})|$ by a continuous function of $\delta_1,\delta_2$. Since $\Lambda_\delta=\Lambda_1$ for $\delta\le1$, we only need to prove the uniform continuity in $\delta\ge1$.
				
				Define
				\[
				\mathscr{P}_\Lambda(\vec{l}(G)):=\{\hat{\mu}\in\mathscr{P}(V_{\vec{l}(G)}):\ \Lambda_1(\hat{\mu})<\infty\}.
				\] 
				First, fix $\hat{\mu}\in \mathscr{P}_\Lambda(\vec{l}(G))$, $1\le\delta_1<\delta_2$,   Lemma \ref{R(|)-R(|)} gives 
				\[
				R(\qq_{l}(\xx,\cdot)\|\pp_{E_l}(\xx,\cdot;\delta_2))-R(\qq_{l}(\xx,\cdot)\|\pp_{E_l}(\xx,\cdot;\delta_1))\le 0.
				\]
				Second, using  equality \eqref{p(delta_1)/p(delta_2)} from the proof of Lemma \ref{R(|)-R(|)},
				\[
				\frac{\pp_{E_l}(\xx,\yy;\delta_1)}{\pp_{E_l}(\xx,\yy;\delta_2)}=\frac{k(\xx)\delta_2+d(\xx)-k(\xx)}{\delta_2}\cdot\frac{\delta_1}{k(\xx)\delta_1+d(\xx)-k(\xx)}\ge {\frac{\delta_1}{\delta_2}}.
				\]
				By equality \eqref{3.4-entpy-3} from the proof of Lemma \ref{R(|)-R(|)} , this implies 
				\[
				R(\qq_{l}(\xx,\cdot)\|\pp_{E_l}(\xx,\cdot;\delta_2))-R(\qq_{l}(\xx,\cdot)\|\pp_{E_l}(\xx,\cdot;\delta_1))\ge \log{\frac{\delta_1}{\delta_2}}.
				\]
				Now choose $(\hat{\mu}_{k},r_k,E_k,\qq_{k})_{1\le k\le \b}\in \hat{\mathscr{A}}(\hat{\mu})$ such that
				\[
				\Lambda_{\delta_2}(\hat{\mu})=\sum_{k=1}^{\b}r_k\int_{V_{\vec{l}(G)}}R(\qq_{k}(\xx,\cdot)\|\pp_{E_k}(\xx,\cdot;\delta_2))\ {\rm d}\hat{\mu}_{k}.
				\]
				From the definition of $\Lambda_{\delta_1}(\hat{\mu})$ and the monotonicity of $\Lambda_\delta$,
				\begin{align*}
					0&\ge \Lambda_{\delta_2}(\hat{\mu})-\Lambda_{\delta_1}(\hat{\mu})\\
					&\ge \sum_{k=1}^{\b}r_k\int_{\vec{l}(G)}\left[R(\qq_{k}(\xx,\cdot)\|\pp_{E_k}(\xx,\cdot;\delta_2))-R(\qq_{k}(\xx,\cdot)\|\pp_{E_k}(\xx,\cdot;\delta_1))\right]\ {\rm d}\hat{\mu}_{k}
					\ge\log{\frac{\delta_1}{\delta_2}}.
				\end{align*}
				Since $\hat{\mu}$ is arbitrary, we have shown
				\begin{equation}\label{3.4-entpy-4}
					0\le\sup_{\hat{\mu}\in\mathscr{P}_\Lambda(\vec{l}(G))}|\Lambda_{\delta_2}(\hat{\mu})-\Lambda_{\delta_1}(\hat{\mu})|\le\log{\frac{\delta_2}{\delta_1}}.
				\end{equation}
				By the squeeze theorem, for $\delta_1,\delta_2\ge1$, 
				$$
				\lim\limits_{|\delta_2-\delta_1|\to 0}\sup_{\hat{\mu}\in\mathscr{P}_\Lambda(\vec{l}(G))}|\Lambda_{\delta_2}(\hat{\mu})-\Lambda_{\delta_1}(\hat{\mu})|=0,
				$$ which establishes the uniform continuity.

				\vskip2mm
				
				\noindent  {\bf(iv)\  \emph{Non-differentiability at} $\delta=1$}. 
				
				Since $\Lambda_\delta= \Lambda_1$ for all $\delta\le 1$,  the left derivative at $\delta=1$ is $0$. We now examine the right derivative.

				Let $\hat{\mu}_0$  be the probability measure chosen in the proof of part {\bf (ii)\ (\emph{Monotonicity for}  $\delta\ge1$)}. That is 
				$$
				\hat{\mu}_0(\overrightarrow{x_0x_1})=\hat{\mu}_0(\overrightarrow{x_1x_0})=\frac{1}{2}.
				$$ 
				From the proof of the monotonicity, the infimum of $\Lambda_\delta\ (\delta>1)$ is attained at $r_1=1$ and $E_1=\{e_0\}$. Let $d_i$  denote the degree of $x_i$ for $i=0,1$. Then 
				$$
				\Lambda_{\delta}(\hat{\mu}_0)-\Lambda_{1}(\hat{\mu}_0)=\frac{1}{2}\log{\frac{d_0-1+\delta}{d_0\delta}}+\frac{1}{2}\log{\frac{d_1-1+\delta}{d_1\delta}}.
				$$ 
				Consequently,  the right derivative (with $\hat{\mu}=\hat{\mu}_0$ fixed) is
				\[
				\lim_{\delta\to 1+}\frac{-\Lambda_{1}(\hat{\mu}_0)}{\delta-1}=\frac{1}{2d_0}+\frac{1}{2d_1}-1.
				\]
				Since the graph is connected and contains more than just the two vertices $\{x_0,x_1\}$, the right derivative at 
				$\delta=1$ for $\hat{\mu}_0$ is is strictly negative. Hence, $\Lambda_{\delta}(\hat{\mu}_0)$
				is not differentiable at $\delta=1$.
			\end{proof}

			\begin{proof}[Proof of Theorem \ref{phase for rate function}] Based on Proposition \ref{phase for rate function 0}, we can prove Theorem \ref{phase for rate function} as follows.
				
				\vskip 2mm
				
				\noindent{\bf(i)\ $I_\delta= I_1$} {\bf\emph{for} $\delta\le 1$}.
				
				\vskip 2mm
				Note that Theorem \ref{rate function 0} gives $I_{\delta}(\mu)=\inf_{\hat{\mu}\in\mathscr{P}(V_{\vec{l}(G)}):\, T(\hat{\mu})=\mu}\Lambda_\delta(\hat{\mu}),\ \mu\in\mathscr{P}(V)$,  thus part (a) of Proposition \ref{phase for rate function 0} implies that $I_\delta(\mu)=I_1(\mu)$ for every $\mu\in\mathscr{P}(V)$.
				
				\vskip 2mm
				\noindent {\bf (ii)\ \emph{Monotonicity for}  $\delta\ge1$}.
				\vskip2mm
				
				For $1\le\delta_1<\delta_2$, recall from the proof of part {\bf(ii)} of Proposition \ref{phase for rate function 0} that $\Lambda_{\delta_2}(\hat{\mu}_0)-\Lambda_{\delta_1}(\hat{\mu}_0)<0$ for $\hat{\mu}_0\in\mathscr{P}(V_{\vec{l}(G)})$ satisfying  
				$$
				\hat{\mu}_0(\overrightarrow{x_0x_1})=\hat{\mu}_0(\overrightarrow{x_1x_0})=\frac{1}{2}, 
				$$ 
				where $x_1$ is a vertex in $V$ adjacent to the starting point $x_0$. Let $e_0=x_0x_1$, and define $\nu_0\in\mathscr{P}(V)$ by $\nu_0(x_0)=\nu_0(x_1)=\frac{1}{2}$. Then, $\hat{\mu}_0$ is the unique solution to $T(\hat{\mu})=\nu_0$. By part (b) of Proposition \ref{phase for rate function 0}, we immediately obtain 
				$$
				I_{\delta_2}(\nu_0)=\Lambda_{\delta_2}(\hat{\mu}_0)<\Lambda_{\delta_1}(\hat{\mu}_0)=I_{\delta_1}(\nu_0), \  \ I_{\delta_2}(\nu)\le I_{\delta_2}(\nu)\ \text{for every}\ \  \nu\in\mathscr{P}(V).
				$$  
				
				\vskip 2mm
				
				\noindent{\bf (iii)\ \emph{Continuity}}.
				\vskip 2mm
				
				Define
				$
				\mathscr{P}_I(G):=\{\nu\in\mathscr{P}(V):\ I_1(\nu)<\infty\}.
				$
				For $\nu\in  \mathscr{P}_I(G)$, the preimage  $T^{-1}\{\nu\}$ is a closed set of $\mathscr{P}(\vec{l}(G))$; hence the infimum of $I_{\delta_2}(\nu)$ is attained  at some $\hat{\mu}_\nu\in T^{-1}\{\nu\}$. Consequencly, 
				$$
				0\le I_{\delta_1}(\nu)-I_{\delta_2}(\nu)\le \Lambda_{\delta_1}(\hat{\mu}_\nu)-\Lambda_{\delta_2}(\hat{\mu}_\nu).
				$$
				Using inequality \eqref{3.4-entpy-4} from the proof of the continuity in Proposition \ref{phase for rate function 0}, we obtain
				\[
				0\le\sup_{\nu\in\mathscr{P}_I(G)}|I_{\delta_1}(\nu)-I_{\delta_2}(\nu)|\le \sup_{\nu\in\mathscr{P}_I(G)}|\Lambda_{\delta_1}(\hat{\mu}_\nu)-\Lambda_{\delta_2}(\hat{\mu}_\nu)|\le\log{\frac{\delta_2}{\delta_1}}.
				\] 
				Thus, $I_\delta$ is  uniformly continuous.  
				
				\vskip2mm
				\noindent  {\bf(iv)\  \emph{Non-differentiability at} $\delta=1$}.
				
				\vskip2mm
				
				Fixing $\hat{\mu}=\hat{\mu}_0$, we have $I_\delta(\nu_0)=\Lambda_\delta(\hat{\mu}_0)$, $I_1(\nu_0)=\Lambda_1(\hat{\mu}_0)$. Therefore, using the same reasoning as in the proof of the non-differentiability of $\Lambda_\delta$ at $\delta=1$, we conclude that the left derivative is $0$, while the right derivative is strictly negative. This establishes the non-differentiability of $I_\delta$ at $\delta=1$.
			\end{proof}  
			
		}
		
		\section{Some examples}\label{sec 5}
		
		{\color{black}
			\noindent   In this section, we demonstrate  more concise expressions for the rate functions and highlight several interesting phenomena.
			
			We first show that for the empirical measures of ORRWs on trees, the method developed in Section \ref{sec 2.3} can be applied directly to verify the large deviation principle, without having to pass through the lifted directed graph. This simplification is possible because trees contain no cycles, which allows us to construct a bijection between the edge set and the vertex set (excluding the root).  Consquently, the infimum condition in the rate function can be simplified accordingly.

			Next,  we compute   the rate functions of ORRWs on the three-vertex tree. We observe that the rate functions differ depending on the starting point—a phenomenon that distinguishes ORRWs from the SRW. Furthermore, we derive  that the rate functions on the  star-shaped graph and the  path graph.  Finally, we obtain the rate function for ORRW on the triangle.

			\subsection{Trees}
			
			\noindent If a $\delta$-ORRW $X$ moves on a finite tree, the absence of cycles allows us to construct a bijection between edges and vertices. This bijection makes it possible to determine whether an edge has been traversed directly from the empirical measure    $L^{n+1}$. Consequently, by Definition    \ref{Def for p_u on G},    the transition probability from      $X_n$ to $X_{n+1}$    can be written as  ${p}_{L^{n+1}}$,    which is uniquely determined by   $L^{n+1}$.       Applying the same procedure used to prove Theorem \ref{estimate for rate} and Theorem \ref{rate function 0} of $\mathcal{Z}$ to $X$ on trees,    we obtain the LDP for $L^{n}$ directly.   This simplification, however, does not extend to general graphs containing cycles.
			
			\vskip2mm
			Consider the rate function $I_\delta$ on trees. {\color{black} Recall that the infimum condition ${\mathscr{A}}(\nu)$ in \eqref{hat_A} for $(\nu_k,r_k,E_k,q_k)_{1\le k\le \b}$ is given by: 
				\begin{itemize}
					\item[(a)] $\{E_k \}_{1\le k\le \b}\in\mathscr{E}$;
					\item[(b)] $r_k\ge 0$, $\sum_{k=1}^\b r_k=1$;
					\item[(c)] $\nu_k\in\mathscr{P}(V_k)$, $q_k\in\mathscr{T}(E_k)$, such that $\nu_k q_k=\nu_k$, and $\sum_{k=1}^\b r_k\nu_k=\nu$.
				\end{itemize}
				This condition can be simplified in the following way: }
			
			\begin{itemize}
				
				\item[(a)]   Condition ${q}_{k}\in\mathscr{T}(E_k)$  can be dropped. 
				\vskip1mm
				This follows from the bijection between edges and vertices (except the root) on trees: the relation $\nu_k{q}_{k}=\nu_k$ together with ${\rm supp}(\nu_k)\subseteq V_k$ already implies ${q}_{k}\in\mathscr{T}(E_k)$.

				\item[(b)]  For  an ORRW on trees, the transition probability ${q}_k$ appearing in the rate function  (\ref{simpler expression of I}) is uniquely determined by  $\nu_k$. 
			\end{itemize}
			
			More precisely, we have   \begin{proposition}\label{Induction on tree}
				For all measure $\nu_{V'}$ supported on a subset $V'\subseteq V$,  the equation  $\nu_{V'}{q}_{\nu_{V'}}=\nu_{V'}$ has a unique solution ${q}_{\nu_{V'}}\in\mathscr{T}_G$.
			\end{proposition}
			We will give the proof of Proposition \ref{Induction on tree} in  Section \ref{Proof Prop5.1}.
			
			\vskip2mm
			\begin{remark}\rm Proposition \ref{Induction on tree}  does not hold in general for graphs containing cycles. 
				For instance,  consider the complete graph on vertices  $\{1,2,3\}$ (a triangle) and the uniform measure  $\nu$ with $\nu(i)=1/3,\ i=1,2,3$.  Then any transition probability of the form          
				$$
				{q}(i,i+1)=x \ \  (\text{with} \ \  {q}(3,4)={q}(3,1))
				$$ satisfies $\nu {q}=\nu$ for every $x\in[0,1]$.
			\end{remark}   
			
			\vskip2mm
			
			Based on the properties above, we obtain the following corollary.
			\begin{corollary} Let $X$ be a $\delta$-ORRW on a finite tree $\mathbf{T}=(V,E)$,  and define
				\begin{align*}
					&{\mathscr{A}}_{\bf T}(\nu)=\Big\{ (\nu_k,r_k,E_k)_{1\leq k\leq \b}:\ \{E_k\}_{1\leq k\leq \b}\in\mathscr{E},\ r_k\ge 0,\ \sum_{k=1}^{\b}r_k=1,\nonumber\\
					&\hskip 3.5cm \nu_k\in \mathscr{P}(V),\ {\rm supp}(\nu_k)\subseteq V_k,\ \sum_{k=1}^{\b} \nu_kr_k=\nu\Big\}.
				\end{align*}
				Then the  empirical measure of $X$ satisfies the LDP with rate function
				\begin{equation}
					I_{\delta}(\nu)=\inf_{(\nu_k,r_k,E_k)_{k}\in{\mathscr{A}}_{\bf T}(\nu)}
					\sum_{k=1}^\b r_k\int_{V}R({q}_{\nu_k}\|{p}_{E_k})\ {\rm d}\nu_{k}.\label{rate function tree}
				\end{equation}
			\end{corollary}
			\vskip2mm
			
			Moreover, we can give an alternative representation  of $I_\delta(\nu)$.
			
			\begin{corollary} \label{rate function 1} Let $X$ be a $\delta$-ORRW on a finite tree $\mathbf{T}=(V,E)$,  and for each $k=1,\dots,\b$ let  $G_k=(V_k,E_k)$  be the subgraph induced by $E_k$. Then     the rate function of the LDP for the  empirical measure of $X$ can also be represented as
				\[
				I_\delta(\nu)=\inf_{(\nu_k,r_k,E_k)_{k}\in{\mathscr{A}}_{\bf T}(\nu)}\Big\{-\inf_{u_{k}\in\mathscr{U}_1}\sum_{k=1}^\b r_k\int_{V}\log{\frac{{p}_{E_k}u_{k}}{u_{k}}(x)}\ {\rm d}\nu_{k}\Big\},
				\]
				where $\mathscr{U}_1$ denote the set of all positive functions $u$ on $V$, and ${p}_{E_k}u(x):=\int_V u(y){p}_{E_k}(x,{\rm d}y)$.
			\end{corollary}
			\begin{proof}  When analyzing the rate function, we replace the condition
				$$
				``{q}_k\in\mathscr{T}(E_k),\ \ {\rm supp}(\nu_k)\subseteq V_k {\rm "}\ \ \text{for}\ \ k=1,\dots,\b\ \ \text{in}\ \ (\ref{hat_A})
				$$
				by the single condition 
				$$
				``{\rm supp}(\nu_k)\subseteq V_k{\rm "}.
				$$ 
				Indeed, using the variational representation of relative entropy together with the minimax lemma cite[Appendix 2, Lemma 3.3]{KL1999},  we obtain 
				\begin{align*}
					\inf_{{q}\in\mathscr{T}_G:\nu {q}=\nu}\int_V R\big({q}(x,\cdot)\big\|{p}(x,\cdot)\big)\ \nu({\rm d}x)
					=&\ \inf_{{q}\in\mathscr{T}_G:\nu {q}=\nu}\int_V\sup_{\phi\in C}\{Q(\phi)-\log{P(e^\phi)} \}\ \nu({\rm d}x)\\
					=&\ \sup_{\phi\in C}\inf_{{q}\in\mathscr{T}_G:\nu {q}=\nu}\int_V \big(Q(\phi)-\log{P(e^\phi)}\big)\ \nu({\rm d}x)\\
					=&\ \sup_{\phi\in C}\int_V\big(\phi-\log{P(e^\phi)}\big)\ \nu({\rm d}x)\\
					=&\ \sup_{\phi\in C}\int_V \log{\frac{e^\phi}{P(e^\phi)}}\ \nu({\rm d}x)\\
					=&\ -\inf_{\phi\in C}\int_V \log{\frac{P(e^\phi)}{e^\phi}}\ \nu({\rm d}x),
				\end{align*}
				where $Pf(x)=\int_V f(y){p}(x,{\rm d}y)$, $Qf(x)=\int_V f(y){q}(x,{\rm d}y)$, and $C$ denotes the set of all real-valued functions on $V$.
				
				Setting $u=e^\phi$ yields 
				\[
				\inf_{{q}\in\mathscr{T}_G:\nu {q}=\nu}\int_V R({q}(x,\cdot)\|{p}(x,\cdot))\ \nu({\rm d}x)=-\inf_{u\in\mathscr{U}_1}\int_V\log{\frac{Pu}{u}(x)}\ \nu({\rm d}x).
				\]     
				Applying this identity to each term in the sum appearing in the definition of $I_\delta(\nu)$     gives the desired representation.     
			\end{proof}
			
			\vskip1mm
			
			Although  the expression in Corollary \ref{rate function 1} resembles the classic result of  Donsker and Varadhan in \cite{DV1975_2},  we are unable to employ their method to ORRWs directly,  their method cannot be applied directly to ORRWs. The reason is that Donsker and Varadhan’s approach fundamentally requires the underlying process to be time‑homogeneous and Markovian, whereas the ORRW lacks these properties over the whole time horizon.     
			
			\subsection{Three-vertex tree}\label{sec 3-vertex}
			\noindent In this subsection, we study the ORRW on the three-vertex tree,   one of the  simplest  simple graphs considered in this paper (see this graph in Figure \ref{path-graph-0}).  We not only derive an explicit expression for the rate function of the LDPs for the empirical measures of ORRWs on this graph, but also observe several interesting phenomena.  
			
			\begin{figure}[htbp]  	
				\centering{\includegraphics[scale=0.8]{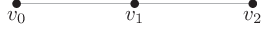}}
				\caption{\label{path-graph-0} \small{  This figure stands for the three-vertex tree $\{0,1,2\}$ with edges $\{0,1\},\{1,2\}$.}}
			\end{figure}
			
			\begin{example} \label{ldp on urn}  Let be a   $\delta$-ORRW
				$X$ start at vertex $v_1$   on  the three-vertex tree $\{v_0,v_1, v_2\}$ with edge set $\{v_0v_1,v_1v_2\}$ (see Figure \ref{path-graph-0}). Then, 	
				the good rate function of the empirical measure process of $X$ is given explicitly by
				\[
				I_{\delta}^{v_1}(\nu)=\left\{
				\begin{array}{ll}
					R(\nu\|{\mu})-R(\nu\|\mathfrak{m}_1(\nu)), &\nu(v_1)=\frac{1}{2}\\
					\infty,&\text{otherwise}
				\end{array}
				\right.,
				\]
				where ${\mu}(v_1)=2{\mu}(v_0)=2{\mu}(v_2)=\frac{1}{2}$, and
				\[
				\mathfrak{m}_1(\nu)(v_i)=\left\{
				\begin{array}{ll}
					\left(\nu(v_i)\vee\frac{\delta-1}{4\delta}\right)\wedge\frac{\delta+1}{4\delta}, &i\neq1\\
					\frac{1}{2}, & i=1
				\end{array}
				\right..
				\]
				The rate function above is obtained by solving the variational problem.	The computation
				is elementary but lengthy; it is provided in Section \ref{sec 5.4}.
				
				\vskip2mm
				This explicit rate function reveals several interesting observations, for instance: \begin{itemize}
					\item let 
					\begin{eqnarray*}
						\nu(v_0)=0,&  \nu(v_1)=\frac{1}{2},&\nu(v_2)=\frac{1}{2}\ \ \  \text{or}\\	
						\nu(v_0)=\frac{1}{2},& \nu(v_1)=\frac{1}{2}, &\nu(v_2)=0.	 
					\end{eqnarray*}	 
					%	\begin{equation}
						%	 \begin{array}{cccc}
							%	  \nu(\{0\})=0,&  \nu(\{1\})=\frac{1}{2},&\nu(\{2\})=\frac{1}{2}&\  \text{or}\\	
							%  	\nu(\{0\})=\frac{1}{2},& \nu(\{1\})=\frac{1}{2}, &\nu(\{2\})=0.&
							%	\end{array}
						%	\end{equation}
					Then, $I_\delta^{v_1}(\nu)$ is not differentiable with respect to $\delta$ at $\delta=1$.  This is a most intuitive example that leads to the discovery of 	 part \emph{(d)} of Theorem \ref{phase for rate function}.
					\item Fix $\delta>1$. Note that $\nu(v_1)=\frac{1}{2}$ if and only if  $I_\delta^{v_1}(\nu)<\infty$. Consequently, the rate function $I_\delta^{v_1}(\nu)$  can be viewed as a function of $\nu(v_0)$. We then have
					$$
					\left\{
					\begin{array}{ll}
						I_\delta^{v_1}(\nu)=I_1(\nu), & \nu(v_0)\in \big[\frac{\delta-1}{4\delta},\frac{\delta+1}{4\delta}\big]\\ 
						I_\delta^{v_1}(\nu)<I_1(\nu), & \nu(v_0)\in  \big[0,\frac{\delta-1}{4\delta}\big)\cup \big(\frac{\delta+1}{4\delta},\frac{1}{2}\big]
					\end{array}\right..
					$$
					Taking $\delta=2$ (see Figure \ref{graph of rate function}), note that  $I_1$  is  exactly the  rate function of SRW on the three-vertex tree. This figure clearly show  the intuitive  difference between the   rate functions of    $\delta$-ORRW  and the SRW.

				\end{itemize}	
				
			\end{example}
			
			\begin{figure}[htbp]
				%\captionsetup[subfigure]{labelformat=empty}
				\centering{\includegraphics[scale=0.6]{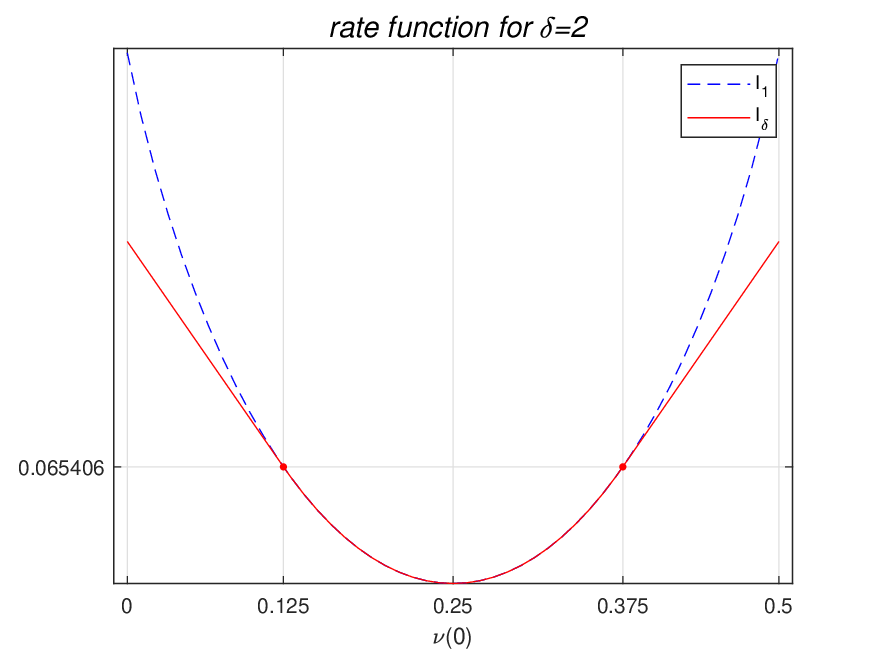}}
				\caption{\label{graph of rate function} \small{ The dashed line is the rate function of the SRW and the solid line is that of the ORRW for $\delta=2$. Here, $\nu(0)$ denotes $\nu(v_0)$. }}
			\end{figure}
			
			\begin{example}\label{ldp on path} Unlike the previous example, here we consider the ORRW $X$ on the three-vertex tree (see Figure \ref{path-graph-0}) starting from vertex $v_0$.  Then, the good rate function of the empirical measure process of $X$ is given  explicitly by
				\[
				I_{\delta}^{v_0}(\nu)=\left\{
				\begin{array}{ll}
					R(\nu\|{\mu})-R(\nu\|\mathfrak{m}_0(\nu)),\ &\nu(v_1)=\frac{1}{2}\\
					\infty,\ &\text{otherwise}
				\end{array}
				\right.,
				\]
				where ${\mu}(v_1)=2{\mu}(v_0)=2{\mu}(v_2)=\frac{1}{2}$, and
				\[
				\mathfrak{m}_0(\nu)(v_i)=\left\{
				\begin{array}{ll}
					\nu(v_0)\wedge\frac{\delta+1}{4\delta},\ &i=0\\
					\frac{1}{2},\ &i=1\\
					\nu(v_2)\vee\frac{\delta-1}{4\delta},\ &i=2
				\end{array}
				\right..
				\]
				
				Similarly, the derivation of  the rate function  is presented  in Section \ref{sec 5.4}, and  resulting $I_\delta^{v_0}$ above exhibits  the following phenomena:
				\begin{itemize}
					\item Let
					\begin{equation*}
						\begin{array}{ccc}
							\nu(v_0)=\frac{1}{2}, &\nu(v_1)=\frac{1}{2}, &\nu(v_2)=0.
						\end{array}
					\end{equation*}
					Then, $I_\delta^{v_0}(\nu)$ is not differentiable with respect to $\delta$ at $\delta=1$.   
					
					\item Let
					\begin{equation*}
						\begin{array}{ccc}
							\nu(v_0)=0, &\nu(v_1)=\frac{1}{2}, &\nu(v_2)=\frac{1}{2}.
						\end{array}
					\end{equation*}
					Then, $I_\delta^{v_0}(\nu)$ in Example \ref{ldp on path} is differentiable with respect to $\delta$ at $\delta=1$. In contrast to Example \ref{ldp on urn},  we emphasize that  the differentiability of this particular  $\nu$ at $\delta=1$ differs from that in Example \ref{ldp on urn}.
					
					\item Fix $\delta>1$. Note that $\nu(v_1)=\frac{1}{2}$ if and only if  $I_\delta^{v_0}(\nu)<\infty$. Consequently, the rate function $I_\delta^{v_0}(\nu)$  can be viewed as a function of $\nu(v_0)$. We then have
					$$
					\left\{
					\begin{array}{ll}
						I_\delta^{v_0}(\nu)=I_1(\nu), & \nu(v_0)\in \big[0,\frac{\delta+1}{4\delta}\big]\\ 
						I_\delta^{v_0}(\nu)<I_1(\nu), & \nu(v_0)\in  \big(\frac{\delta+1}{4\delta},\frac{1}{2}\big]
					\end{array}\right..
					$$
					Choosing  $\delta=2$ (see Figure \ref{graph of rate function 1} (a)), this figure again illustrates the intuitive  difference between the   rate functions of    $\delta$-ORRW  and the SRW.

				\end{itemize}
				
			\end{example}

			\begin{figure}[htbp]
				\subfigure[]{\begin{minipage}[t]{0.4\textwidth}
						\centering{\includegraphics[scale=0.4]{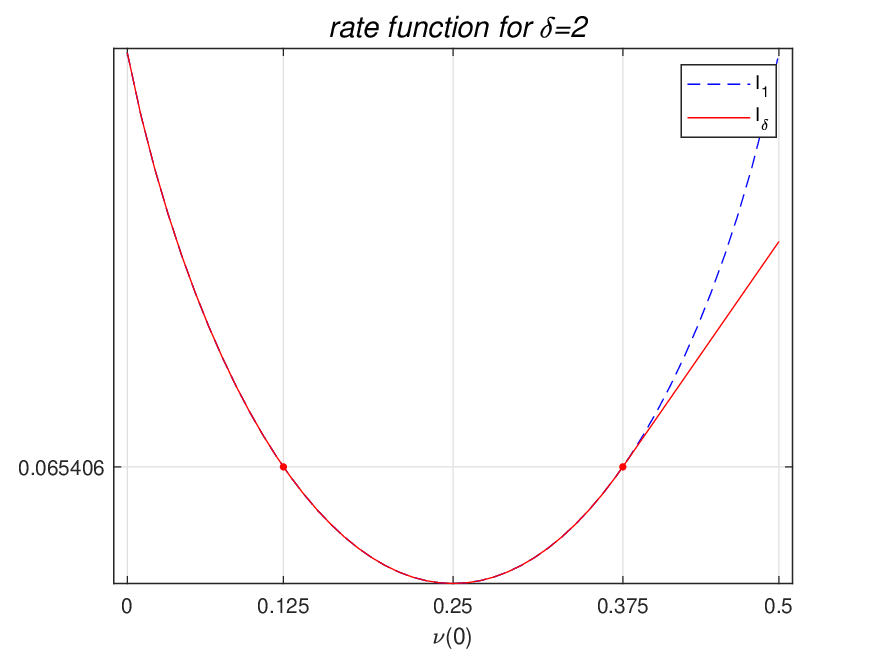}}
				\end{minipage}}
				\subfigure[]{\begin{minipage}[t]{0.7\textwidth}
						\centering{\includegraphics[scale=0.4]{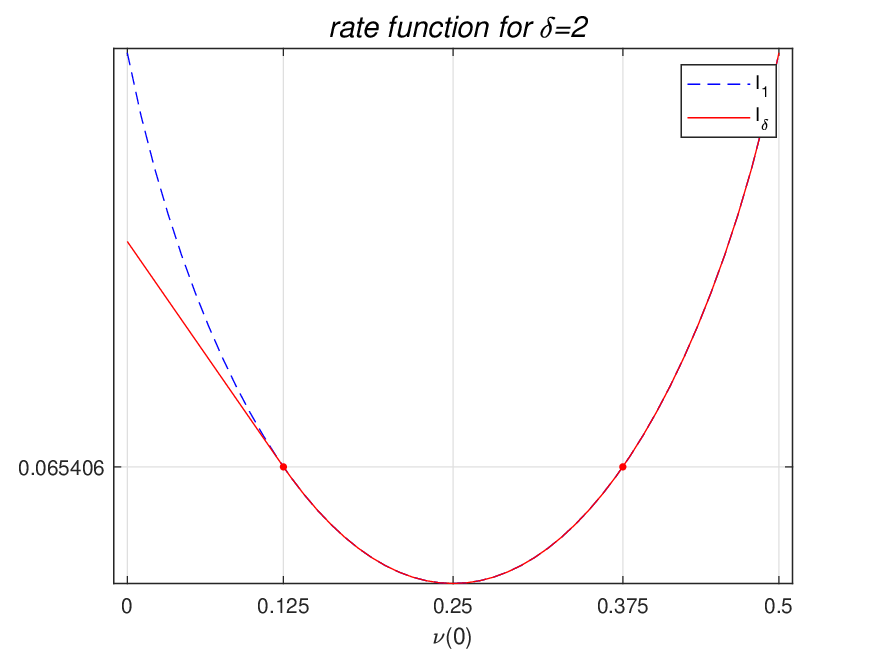}}
				\end{minipage}}
				\caption{\label{graph of rate function 1} \small{(a) is the rate function of the ORRW starting at $v_0$, while (b) is that of starting at $v_2$, where the dashed line is the rate function of the SRW, and the solid line is that of the ORRW for $\delta=2$}. We observe  that the two rate functions are related by  \eqref{3-vertex-1} and coincide  when $\nu(v_0)\in[(\delta-1)/4\delta,(\delta+1)/4\delta]$.}
			\end{figure}
			
			\begin{example}  \label{ldp on path-2} Compared with the previous  Example \ref{ldp on path}, where the ORRW starts from $v_0$, we now consider the ORRW $X$ starting from  $v_2$.
				Then, we obtain the explicit good rate function of the empirical measure process of $X$ as follows:
				\[
				I_{\delta}^{v_2}(\nu)=\left\{
				\begin{array}{ll}
					R(\nu\|{\mu})-R(\nu\|\mathfrak{m}_2(\nu)),\ &\nu(v_1)=\frac{1}{2}\\
					\infty,\ &\text{otherwise}
				\end{array}
				\right.,
				\]
				where ${\mu}(v_1)=2{\mu}(v_0)=2{\mu}(v_2)=\frac{1}{2}$, and
				\[
				\mathfrak{m}_2(\nu)(v_i)=\left\{
				\begin{array}{ll}
					\nu(v_0)\vee\frac{\delta-1}{4\delta},\ &i=0\\
					\frac{1}{2},\ &i=1\\
					\nu(v_2)\wedge\frac{\delta+1}{4\delta},\ &i=2
				\end{array}
				\right..
				\]
				
				Define a mapping  ${\bf m}:$ $\{v_0,v_1, v_2\} \mapsto \{v_0,v_1, v_2\}$ by
				\[
				{\bf m}(v_0)=v_2,\  {\bf m}(v_1)=v_1,\  {\bf m}(v_2)=v_0.
				\]
				For the ORRW $X$ on the three-vertex tree starting at the vertex $v_2$, the transformed process ${\bf m}(X)$ is again  an ORRW on the same tree, now  starting at the vertex $v_0$, see Example \ref{ldp on path}. Consequently,
				\begin{equation}\label{3-vertex-1}
					I_\delta^{v_2}(\nu)= I_\delta^{v_0}(\nu\circ {\bf m}),\ \ \text{for every}\ \ \nu\in\mathscr{P}(\{v_1,v_2,v_3\}),
				\end{equation}
				where $\nu\circ{\bf m}(\{v_i\})=\nu(\{v: {\bf m}^{-1}(v)=v_i\})$,  ${\bf m}^{-1}$ denote the inverse map of ${\bf m}$.
				
			\end{example}
			\vskip2mm	
			
			%\vskip1cm

			\vskip2mm
			{\color{black} It is worth noting that the choice of the starting vertex of the ORRW $X$ affects the rate function of the LDP for its empirical measure when $\delta>1$. This can be seen by comparing Figure \ref{graph of rate function}, Figure \ref{graph of rate function 1} (a) and  Figure \ref{graph of rate function 1} (b).

				In contrast, for the simple random walk ($\delta=1$),  on a finite graph, the rate function 
				$I_1$ does not depend on the starting point. This raises an interesting question: can one run an ORRW on an unknown finite graph or network to extract information about its structure?
				
			}
			
			\vskip2mm
			The  three‑vertex tree is simultaneously a star‑shaped graph and a path graph. We now extend our study to the expressions of the rate function in the LDP for the empirical measure of $X$ on general star‑shaped graphs and path graphs.
			
			\subsection{Star-shaped graphs}\label{sec 5.2}
			\noindent To begin we take a look at the simplest star-shaped graph.  
			\begin{figure}[htbp]
				%\captionsetup[subfigure]{labelformat=empty}
				\centering{\includegraphics[scale=1]{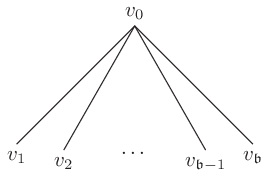}}
				\caption{\label{polya urn graph} \small{This figure indicates a star-shaped graph with $\b$ edges.}}
			\end{figure}
			We consider the star-shaped graph with vertices $\{v_0,v_1,\dots,v_\b\}$ and  edges $\{v_0v_i:i=1,\dots,v_\b \}$, and let $\delta$-ORRW $X$ start at vertex $v_0$ (see Figure \ref{polya urn graph}).
			
			For any subgraph $\{v_0,v_{s_1},\dots,v_{s_k}\}$ with  edge set $E_k=\{v_0v_{s_i}:i=1,\dots,k\}$, the conditions  $\nu_{k}q_{k}=\nu_{k}$, and $\text{\rm{supp}}(\nu_{k})=\{v_0,v_{s_1},\dots,v_{s_k}\}$ imply 
			\begin{align*}
				\nu_{k}(v_0)=&\ \frac{1}{2},\\
				{q}_{k}(v_0,x)=&\  2\nu_{k}(x)\\
				{p}_{E_k}(v_0,y)=&\ \frac{\delta\mathbf{1}_{\{v_{s_1},\dots,v_{s_k}\}}(y)+\mathbf{1}_{\{v_1,\dots,v_\b\} \setminus\{v_{s_1},\dots,v_{s_k}\}}(y)}{k(\delta-1)+\b}.
			\end{align*}
			Consequently, 
			$$
			\int_V R({q}_k\|{p}_{E_k})\ \nu_k({\rm d}x)=\sum_{l=1}^k\nu_k(v_{s_l})\log{2\nu_k(v_{s_l})}+\frac{1}{2}\log{\frac{k(\delta-1)+\b}{\delta}}.
			$$ Thus, by Theorem \ref{I_delta}, the good rate function of the empirical measure process of $X$ is
			\[
			I_{\delta}(\nu)=\left\{
			\begin{array}{ll}\displaystyle
				\inf_{(\nu_k,r_k,s_k)_k\in\mathscr{A}_{\text{star}}(\nu)} \sum_{k=1}^\b r_k\Big[\sum_{l=1}^k\nu_k(v_{s_l})\log{2\nu_k(v_{s_l})}+\frac{1}{2}\log{\frac{k(\delta-1)+\b}{\delta}}\Big], &\nu(v_0)=\frac{1}{2}\\
				\infty,&\text{otherwise}
			\end{array}
			\right.,
			\] where
			\begin{align*}
				\mathscr{A}_{\text{star}}(\nu):=&\ \Big\{(\nu_k,r_k,s_k)_k:\ \sum_{k=1}^\b\nu_k r_k=\nu,\sum_{k=1}^\b r_k=1,r_k\ge 0,
				\nu_k(v_0)=1/2,\\
				&\hskip 3mm \text{\rm{supp}}(\nu_k)=\{v_0,v_{s_1},\dots,v_{s_k}\},\{v_{s_1},\dots,v_{s_\b}\}\ \text{\rm is a permutation of}\ \{v_1,\dots,v_\b\}\Big\}.
			\end{align*}
			
			It is indeed difficult to obtain an explicit expression for $I_\delta$ on a general star‑shaped graph with $\b \ge 3$ as simple as the one in Example \ref{ldp on urn}.

			\subsection{Path graph $\{0,1,\dots,\b \}$}\label{sec 5.3}
			\noindent  In this subsection,  we study the rate function of the ORRW  on the path graph $\{v_0,v_1,\dots,v_\b \}$ with edges $\{v_iv_{i+1}:i=0,\dots,\b-1\}$.
			\begin{figure}[htbp]%\label{path-graph}
				%\captionsetup[subfigure]{labelformat=empty}
				
				\centering{\includegraphics[scale=0.8]{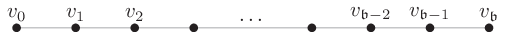}}
				
				\caption{\label{path-graph} \small{This figure indicates a path graph with $\b$ edges.}}
			\end{figure}
			
			Note that in this case $\mathscr{E}=\big\{\{v_iv_{i+1}:0\le i\le k-1 \}_{1\le k\le \b}\big\}$. For each subgraph $G_k=\{v_0,v_1,\dots,v_k\}$ with edge set $E_k=\{v_0v_1,v_1v_2,\dots,v_{k-1}v_k \}$ and a given probability measure   $\nu_k$ on $G_k$, the corresponding transition probability  ${q}_k$ can be determined inductively from the leaf $v_0$ as shown in Proposition \ref{Induction on tree}.
			Specifically,
			\[
			{q}_k(v_i,v_{i-1})=\frac{\sum_{j=1}^i (-1)^{j-1}\nu_k(v_{i-j})}{\nu_k(v_i)},\ {q}_k(v_i,v_{i+1})=1-{q}_k(i,i-1).
			\]
			Since ${q}_k(v_k,v_{k-1})=1$. we obtain the condition $\sum_{j=0}^k (-1)^{j}\nu_k(v_{k-j})=0$. Consequently,
			\begin{align*}
				&\int_{V}R({q}_k\|{p}_{E_k})\ {\rm d}\nu_k\\
				&=\sum_{i=1}^{k-1}\nu_k(v_i)\Bigg[\log{2\frac{\sum_{j=1}^i (-1)^{j-1}\nu_k(v_{i-j})}{\nu_k(v_i)}}+\log{\Big(2-2\frac{\sum_{j=1}^i
						(-1)^{j-1}\nu_k(v_{i-j})}{\nu_k(v_i)}\Big)}\Bigg]\\
				&\ \ \ \ +\nu_k(v_k)\log{\frac{1+\delta}{\delta}},
			\end{align*}
			subject to $\sum_{j=0}^k (-1)^{j}\nu_k(v_{k-j})=0$.

			Thus the good rate function for the empirical measure process of  $X$ on the path graph $\{0,1,\dots,\b\}$ with edges $\{\{i,i+1\},i=0,\dots,\b-1\}$ starting from $v_0$   can be written as  
			\begin{align*}
				I_\delta({\mu})=&\,\inf_{(\nu_k,\lambda_k)_k\in \mathscr{A}_{\text{path}}({\mu})}\Bigg\{\sum_{k=1}^\b \lambda_k\sum_{i=1}^{k-1}\nu_k(v_i)\Big[\log{2\frac{\sum_{j=1}^i (-1)^{j-1}\nu_k(v_{i-j})}{\nu_k(v_i)}}\\
				&\hskip2cm +\log{\Big(2-2\frac{\sum_{j=1}^i (-1)^{j-1}\nu_k(v_{i-j})}{\nu_k(v_i)}\Big)}\Big]
				+\nu_k(v_k)\log{\frac{1+\delta}{\delta}}I_{\{k<d\}}\Bigg\},
			\end{align*}
			where
			\begin{align*}\mathscr{A}_{\text{path}}({\mu}):=\Bigg\{&(\nu_k,\lambda_k)_k:\sum_{k=1}^\b\nu_k\lambda_k={\mu},\sum_{k=1}^\b\lambda_k=1,
				\\
				&\sum_{j=0}^k (-1)^{j}\nu_k(v_{k-j})=0,\text{\rm{supp}}(\nu_k)=\{v_0,\dots,v_k\},\lambda_k\ge0 \Bigg\}.
			\end{align*}
			
			Even for path graphs, we still lack a better computational method to obtain a more concise form of the rate function.
			% \item
			%  \item

			\subsection{Triangles} 
			For ORRWs on the triangle $\{v_1,v_2,v_3\}$ starting at $v_1$,   consider any probability measure $\nu$ on the triangle. Write  $\nu(v_1)=a$, $\nu(v_2)=b$, and $\nu(v_3)=1-a-b$ (the graph is depicted in Figure \ref{rate-function-on-triangle} (a)). For notational convenience, identify site $v_0$ with site $v_3$ and site $v_4$ with site $v_1$. Then the rate function $I_\delta$ is given by 
			\begin{align}
				I_\delta(\nu)%&=I_\delta((a,b,1-a-b))\nonumber\\
				=&\ \inf_{\{E_1,E_2,E_3\}\in\mathscr{E}}\Bigg\{\inf_{(r_1,r_2,r_3,\nu_2,\nu_3)\in\mathscr{A}_{\text{triangle}}(E_1,E_2,E_3)}
				\bigg\{ r_1\log\frac{\delta+1}{\delta}+\nonumber\\
				&\hskip2cm  r_2\Big(\frac{1}{2}\log\frac{\delta+1}{2\delta}+\nu_2(v_1)\log(4\nu_2(v_1))+\nu_2(v_2)\log(4\nu_2(v_2)) \nonumber \\ 
				&\hskip6cm +\nu_2(v_3)\log(4\nu_2(v_3)) \Big)+r_3I_1(\nu_3)  \bigg\}\Bigg\},
			\end{align}
			where $I_1$ is the rate function for $\delta=1$, and
			\begin{align}
				&\mathscr{A}_{\text{triangle}}(E_1,E_2,E_3)\nonumber\\
				:=&\ \ \Big\{ (r_1,r_2,r_3,\nu_2,\nu_3):r_1+r_2+r_3=1,\ r_1,r_2,r_3\ge0, \nu_2,\nu_3\in\mathscr{P}(\{v_1,v_2,v_3\}),\nonumber\\ 
				&\hskip1cm \frac{r_1}{2}+r_2\nu_2(v_1)+r_3\nu_3(v_1)=a,
				\frac{r_1}{2}\mathbf{1}_{\{ v_1v_2\in E_1 \}}+r_2\nu_2(v_2)+r_3\nu_3(v_2)=b, \nonumber\\
				&\hskip5cm  \text{for }i \text{ with }v_iv_{i+1}\notin E_2,\  \nu_2(v_{i-1})=\frac{1}{2}\Big\}.
			\end{align}
			For the rate function $I_1$ and a probability measure $\nu\in\mathscr{P}(\{v_1,v_2,v_3\})$,
			\[
			I_1(\nu)=\inf_{(q_1,q_2,q_3)\in Q(\nu(v_1),\nu(v_2),\nu(v_3))} \sum_{i=1,2,3} \nu(v_i)\big(q_i\log{2q_i}+(1-q_i)\log{2(1-q_i)}\big),
			\]
			where
			\begin{align*}
				Q(a,b,c)=&\ \{(q_1,q_2,q_3): 0\le q_i\le 1\text{ for }i=1,2,3,\\
				&\hskip 1cm b(1-q_2)+cq_3=a, aq_1+c(1-q_3)=b, bq_2+a(1-q_1)=c\}.
			\end{align*}
			\begin{figure}[htbp]
				\subfigure[]{\begin{minipage}[t]{0.4\textwidth}
						\centering{\includegraphics[scale=1.3]{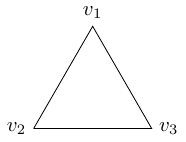}}
				\end{minipage}}
				\subfigure[]{\begin{minipage}[t]{0.4\textwidth}
						\centering{\includegraphics[scale=0.5]{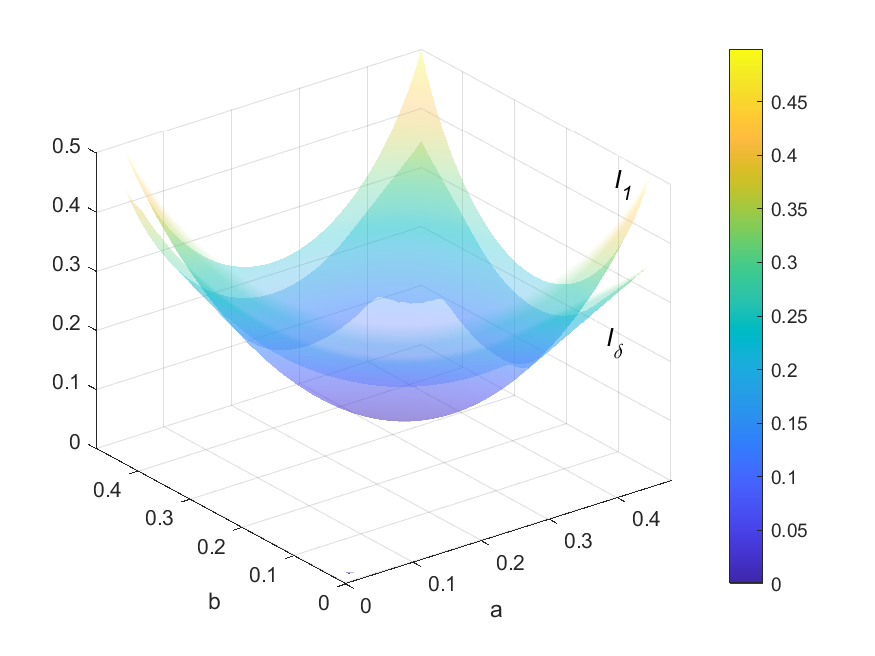}}
				\end{minipage}}
				\caption{\small{(a) illustrates the triangle graph. (b) indicates the rate function on triangle when $\delta=2$. Here $a$ is the measure on the starting point $v_1$, and $b$ is the measure on $v_2$. The rate function is finite if and only if $a+b\ge1/2$ and $0\le a,b\le1/2$.}}
				\label{rate-function-on-triangle}
			\end{figure}
			
			The expression of the rate function $I_\delta$ on the triangle is complicated, and it is difficult to obtain an explicit value for a fixed $\nu$. So far, we have not found an effective method to compute its explicit expression.

			Through numerical computation and simulation, we can plot the rate function as a function of  $\nu(v_1)=a$ and $\nu(v_2)=b$, see Figure \ref{rate-function-on-triangle} (b). This visualization helps to intuitively understand the shape of the rate function and reveals interesting phenomena, such as, 
			the rate function $I_\delta$ is less than $I_1$ if $\delta>1$; in the neighborhood of the zero point $(a,b)=(1/3,1/3)$, the rate function $I_\delta$ equals $I_1$.

		}

		\section{Supplementary Results and Proofs} \label{sec final}
		
		{\color{black} 
			\subsection{Supplement Results in Section \ref{sec 1}}
			\noindent  The following proposition is well‑known and will be used to streamline the proof of Proposition \ref{exponentially integrable}.
			\begin{proposition}\label{domination}
				For two random variables $\xi$ and $\eta$, if $\P(\xi \ge t)\le \P(\eta  \ge t)$ for any $t\in \R$, and $f$ is differentiable and increasing, then $\E f(\xi)\le \E f(\eta)$.
			\end{proposition}
			\begin{proof}   $\E f(\xi)=\int_{-\infty}^\infty f'(t)\P(\xi>t)\ {\rm d}t\le \int_{-\infty}^\infty f'(t)\P(\eta>t)\ {\rm d}t=\E f(\eta)$.
			\end{proof}
			
			\begin{proposition}\label{exponentially integrable}
				For the $\delta$-ORRW $(X_n)_{n\ge 0}$ on finite connected graph $G=(V,E)$ with $|E|=\b\ge 1$, there exists some $\alpha>0$ such that $\E\, e^{\alpha C_E}<\infty$.
			\end{proposition}
			\begin{proof}
				Fix an ordering $\vec{e}_1,\dots,\vec{e}_{2\b}$ of all directed edges of $G$. Let $\mathbb{N}=\{1,2,\cdots\}$, and define 
				$$
				\tilde{\tau}_k=\inf\left\{t>\tilde{\tau}_{k-1}:\ t\in \mathbb{N}_V, \, \exists s\in[\tilde{\tau}_{k-1}+1,t],\ \text{such that}\ \ \overrightarrow{X_{s-1}X_{s}}=\vec{e}_k \right\},\ 1\le k\le 2\b,
				$$
				as a sequence of stopping times, where $\tilde{\tau}_0=0$, and $\tilde{\tau}_{2\b+1}=\infty$. Then $C_E\le \tilde{\tau}_{2\b}$.

				Our strategy is to construct  an i.i.d.\,sequence $\{N_k\}_{1\le k\le 2\b}$ obeying geometric distribution with parameter $p_0\in (0,1]$ (to be determined later) such that for every non-negative integers $t$, 
				\begin{equation}
					\mathbb{P}(\tilde{\tau}_k\ge t) \le \mathbb{P}\left(|V| \sum_{j=1}^k N_j\ge t\right). \label{geometric d control}
				\end{equation}

				Observing that  $f(x)=e^{\alpha x}$ is non-negative, differentiable and increasing,  Proposition \ref{domination} yields,
				$$
				\E e^{\alpha C_E}\le \E e^{\alpha |V|\sum_{k=1}^{2\b}N_k}=\left( \sum_{k=0}^\infty [e^{\alpha |V|}(1-p_0)]^k p_0 \right)^{2\b}.
				$$
				Since  $\left( \sum_{k=0}^\infty [e^{\alpha |V|}(1-p_0)]^k p_0 \right)^{2\b}<\infty$ if and only if $e^{\alpha |V|}(1-p_0)<1$, we obtain 	
				$$
				\E e^{\alpha C_E}<\infty\ \mbox{for every}\ \alpha \in\left(0,\frac{1}{|V|}\log{\frac{1}{1-p_0}}\right).
				$$ 	
				
				To prove \eqref{geometric d control}, we first determine $p_0$. For any vertices $v,u\in V$ and a stopping time $\tau$, by the definition of weight function $w_\tau$ in \eqref{weight} gives 
				$$
				\P(X_{\tau+1}=u|X_\tau=v,\mathscr{F}_\tau)= \frac{w_\tau(uv)}{\sum_{u'\sim v}w_\tau(u'v)} \ge\frac{\delta \wedge 1}{(|V|-1)(\delta\vee 1)},\ \ \mbox{if}\ u\sim v;
				$$
				If $u$ is not adjacent to $v$,   consider the shortest path $vv_1\dots v_{n-1}u$ from $v$ to $u$ (noting $n\le |V|$). Then
				\begin{align*}
					\P(X_{\tau+n}=u|X_\tau=v,\mathscr{F}_\tau)&\,\ge\P(X_{\tau+n}=u, X_{\tau+n-1}=v_{n-1},\dots,X_{\tau+1}=v_1|X_\tau=v,\mathscr{F}_\tau)\\
					&\, \ge\left( \frac{\delta \wedge 1}{(|V|-1)(\delta\vee 1)} \right)^{|V|}.
				\end{align*}	
				Consequently, for any $\tilde{u}\sim u$,
				\begin{equation*}
					\P(X_{\tau+n+1}=\tilde{u},X_{\tau+n}=u |X_\tau=v,\mathscr{F}_\tau)\ge \left( \frac{\delta \wedge 1}{(|V|-1)(\delta\vee 1)} \right)^{|V|+1}. \label{technique of cover time}
				\end{equation*}
				We now set 
				$$
				p_0=\left( \frac{\delta \wedge 1}{(|V|-1)(\delta\vee 1)} \right)^{|V|+1}.
				$$
				By \eqref{technique of cover time}, for any $1\le k\le 2\b$ and vertex  $v\in V$, and any integer $m\in\mathbb{N}$,
				\begin{align}
					&\P\Big(\exists n\in\big[|V| m,|V|(m+1)\big], \overrightarrow{X_{n}X_{n+1}}=\vec{e}_{k}\,\Big|\,X_{|V|m}=v,\mathscr{F}_{|V| m}\Big)\ge p_0,\nonumber\\
					&\P\Big(\forall n\in\big[|V| m,|V|(m+1)\big], \overrightarrow{X_{n}X_{n+1}}\neq\vec{e}_{k}\, \Big|\,X_{|V| m}=v,\mathscr{F}_{|V| m}\Big)<1-p_0.\label{technique of cover time 1}
				\end{align}
				Hence, from \eqref{technique of cover time 1} we obtain for every $m\in\mathbb{N},$
				\[
				\P\left.\left(\frac{\tilde{\tau}_k-\tilde{\tau}_{k-1}}{|V|}\ge m \right| \mathscr{F}_{\tilde{\tau}_{k-1}} \right)\le (1-p_0)^m \text{ for all }k=1,\dots,2\b;
				\]
				and further for any  integers $m_1,\cdots,m_{2\b}\in\mathbb{N}$ and any $1\le k\le 2\b,$
				\[
				\P\left(\frac{\tilde{\tau}_j-\tilde{\tau}_{j-1}}{|V|}\ge m_j,\, j=1,\dots,k\right)\le (1-p_0)^{m_1+\dots+m_k}=\P(N_j\ge m_j,j=1,\dots,k).
				\]
				This implies  \eqref{geometric d control} for any $1\le k\le 2\b$. 
			\end{proof}

			\subsection{Supplement Results and Proofs in Section \ref{sec 2}}
			
			\subsubsection{Supplement Results and Proof in Section \ref{sec 2.15}}\label{5-pf sec 2.2}
			\begin{proposition}\label{X-L-Mrk}
				Let  $(X_n)_{n\ge 0}$ be a discrete-time Markov processes  defined on Probability space $(\Omega,\mathscr{F},\P)$ and taking values in a Polish space $\mathscr{X}$, and let $(L^n)_{n\ge 0}$ be its associated  empirical measure process.  Then the two-component processes $\big((X_n,L^n)\big)_{n\ge 0}$ is  again  Markovian  with respect to the natural filtration $\big\{\mathscr{F}^X_n=\sigma(X_k, 0\le k\le n)\big\}_{n\ge 0}$ generated by $(X_n)_{n\ge 0}$.
			\end{proposition}
			\begin{proof} 
				For each bounded measurable function $f: \mathscr{X}\times\mathscr{P}(\mathscr{X})\mapsto \mathbb{R}$, we aim to show
				\begin{equation}\label{5.1.1-1}
					\E(f(X_{n+1},L^{n+1})|\mathscr{F}^X_n)=\E\Big(f\big(X_{n+1},\frac{1}{n+1} (nL^n+\dd_{X_n})\big)\Big| X_n,L^n\Big).
				\end{equation}
				
				\vskip2mm
				
				\noindent{\small\bf  \emph{Step 1}}.  Define  $\mathscr{F}_n^{(X,L)}:=\sigma\big( (X_k, L^k), 1\le k\le n\big)$. Since $L^{n+1}=\frac{1}{n+1} (nL^{n}+\dd_{X_n})$, it follows that  $L^n$ is measurable with respect to $\sigma(X_k:k\le n)$.  Consequently,  $\big\{\mathscr{F}_n^{(X,L)}\big\}_{n\ge 1}$    coincides with $\big\{\mathscr{F}^X_n\big\}_{n\ge 1}$; i.e., $\big((X_n,L^n)\big)_{n\ge 0}$ is adapted to  $\big\{\mathscr{F}^X_n\big\}_{n\ge 0}$. 
				
				\vskip2mm
				
				\noindent{\small\bf  \emph{Step 2}}.   
				For $\mu\in \mathscr{P}(\mathscr{X})$ and $x\in \mathscr{X}$, define 
				$$
				g(\mu, x,X_n):=\E\Big(f\big(X_{n+1},\frac{1}{n+1} (n\mu+\dd_{x})\big)\Big|X_n\Big).
				$$  By the Markov property of $(X_n)_{n\ge 1}$,  
				\begin{equation}\label{5.1.1-2}
					g(\mu, x,X_n)=\E\Big(f\big(X_{n+1},\frac{1}{n+1} (n\mu+\dd_{x})\big)\Big|X_n\Big)=\E\Big(f\big(X_{n+1},\frac{1}{n+1} (nL^{n}+\dd_{X_n})\big)\Big|\mathscr{F}_n^X\Big).
				\end{equation}  
				Because $(X_n, L^n)$ is $\mathscr{F}^X_n$-measurable, we obtain
				\begin{equation}\label{5.1.1-2}
					\E\Big(f\big(X_{n+1},\frac{1}{n+1} (nL^{n}+\dd_{X_n})\big)\Big|\mathscr{F}_n^X\Big)=g(L^n,X_n, X_n).
				\end{equation}
				\vskip2mm
				
				\noindent{\small\bf  \emph{Step 3}}.   Since  $g(L^n, X_n, X_n)$ is $\sigma(X_n,L^n)$-measurable, we have  \begin{align*}
					g(L^n,X_n,X_n)
					=&\, \E\big(g(L^n, X_n,\omega)\big|X_n, L^n\big)\\
					=&\, \E\Bigg(\E\Big(f\big(X_{n+1},\frac{1}{n+1} (nL^{n}+\dd_{X_n})\big)\Big|\mathscr{F}_n\Big)\Big|X_n,L^n\Bigg)\\
					=&\, \E\Big(f\big(X_{n+1},\frac{1}{n+1} (nL^n+\dd_{X_n})\big)\Big| X_n,L^n\Big).
				\end{align*}
				Hence, \eqref{5.1.1-2} together with the preceding equality yields \eqref{5.1.1-1}. Consequently, $(X_n, L_n)_{n\ge 1}$ is a Markov process.
			\end{proof}

			%\begin{align*}%\label{eq-1-prop-Markov-(X,L)}	
			%       \E(f(X_{n+1},L^{n+1})|\mathscr{F}_n)
			%   =&\,\E\Big(f\big(X_{n+1},\frac{1}{n+1} (nL^{n}+\dd_{X_n})\big)\Big|\mathscr{F}_n\Big)\\
			%   =&\,\E\Big(f\big(X_{n+1},\frac{1}{n+1} (n\mu+\dd_{x})\big)\Big|\mathscr{F}_n\Big)\Big|_{(\mu=L^n, x=X_n)}\\
			%   =&\, \E\Big(f\big(X_{n+1},\frac{1}{n+1} (n\mu+\dd_{x})\big)\Big|X_n\Big)\Big|_{(\mu=L^n, x=X_n)}\\ 
			%   =&\, \E\Bigg(\E\Big(f\big(X_{n+1},\frac{1}{n+1} (n\mu+\dd_{x})\big)\Big|X_n\Big) \Bigg|X_n,L^n\Bigg) \Big|_{(\mu=L^n, x=X_n)} \\
			%  =&\, \E\Bigg(\E\Big(f\big(X_{n+1},\frac{1}{n+1} (n\mu+\dd_{x})\big)\Big|X_n\Big)\Big|_{(\mu=L^n, x=X_n)}  \Bigg|X_n,L^n\Bigg)\\
			%   =&\, \E\Bigg(\E\Big(f\big(X_{n+1},\frac{1}{n+1} (n\mu+\dd_{x})\big)\Big|\mathscr{F}_n\Big)\Big|_{(\mu=L^n, x=X_n)}  \Bigg|X_n,L^n\Bigg)\\
			%  =&\,  \E\Bigg(\E\Big(f\big(X_{n+1},\frac{1}{n+1} (nL^n+\dd_{X_n})\big)\Big|\mathscr{F}_n\Big)  \Bigg|X_n,L^n\Bigg)\\
			%  =&\, \E\Big(f\big(X_{n+1},\frac{1}{n+1} (nL^n+\dd_{X_n})\big)\Big| X_n,L^n\Big),
			%    \end{align*}
		%which completes the proof of Markov property.

		\begin{figure}[htbp]
			\centering{\includegraphics[scale=0.3]{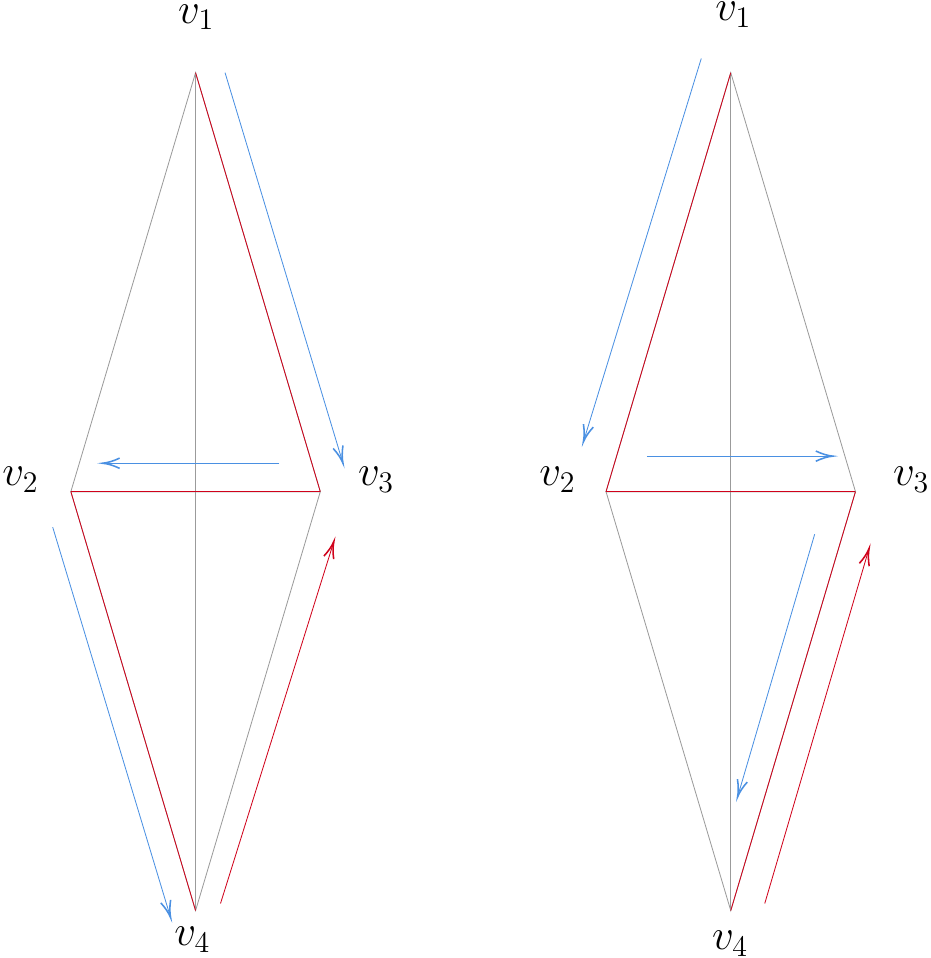}} 
			\caption{\small{ORRW on complete graph}}
			\label{counter-ex}
		\end{figure}
		
		\begin{example}\label{ORRW-non_Mrv} 
			\label{counter-example}    Considering a $\delta$-ORRW $(X_n)_{n\ge0}$ on the complete graph with four vertices $v_1,v_2,v_3,v_4$, as shown in  Figure \ref{counter-ex}. We that this process is not Markovian.
			
			For $\delta\neq1$,
			\begin{align*}
				&\P\big(X_4=v_3\big|X_3=v_4,L^3=\frac{1}{3}\sum_{i=1}^3\dd_{v_i},X_2=v_2,X_1=v_3,X_0=v_1\big)=\frac{1}{2+\delta}\\
				&\neq \P\big(X_4=v_3\big|X_3=v_4,L^3=\frac{1}{3}\sum_{i=1}^3\dd_{v_i},X_2=v_3,X_1=v_2,X_0=v_1\big)=\frac{\delta}{2+\delta}.
			\end{align*}
			This demonstrates that the transition probability depends on the entire history even when the current empirical measure is known; in other words, $(X_n)_{n\ge 0}$ is not a Markov process. 
		\end{example}
		
		\vskip2mm 
		%\begin{proposition}\label{Z-L-Mrk} 
		% Let $(X_n)_{n\ge0}$ be a $\delta$-ORRW  on finite connected graphs $G=(V,E)$, let $(\mathcal{Z}_n)_{n\ge0}$ be  its lifted process on  line digraph $\vec{l}(G)$ and let $(\mathcal{L}^n)_{n\ge0}$ be the  empirical measure processes of $(\mathcal{Z}_n)_{n\ge0}$. Denote by  $\mathscr{F}^{\mathcal{Z}}_n:=\sigma(\mathcal{Z}_k: k\le n)$.  Then, the two-component process $((\mathcal{Z}_n,\mathcal{L}^n))_{n\ge0}$ is a  discrete-time Markov processes with respect to $\{\mathscr{F}^{\mathcal{Z}}_n\}_{n\ge 0}$.
		% \end{proposition}
	
	\begin{proof}[Proof of Proposition \ref{Z-L-Mrk}] For every  bounded measurable function $f:V_{\vec{l}(G)}\times \mathscr{P}(V_{\vec{l}(G)}) \mapsto\mathbb{R}$, we aim to show, 
		\begin{equation}\label{5.1.1-Z-1}
			\E\big(f(\mathcal{Z}_{n+1},\mathcal{L}^{n+1})\big|\mathscr{F}^{\mathcal{Z}}_n)\big) =\E\big(f(\mathcal{Z}_{n+1},\mathcal{L}^{n+1})\big| \mathcal{Z}_{n},\mathcal{L}^{n} \big).
		\end{equation}
		\vskip2mm
		\noindent{\small\bf  \emph{Step 1}}.  Denote by  $\mathscr{F}^{\mathcal{Z}, \mathcal{L}}_n:=\sigma((\mathcal{Z}_k,\mathcal{L}^k):k\le n)$. By the definition of ORRWs $X_n$ and the construction of $\mathcal{Z}_n$, we  know that 
		$$
		\mathscr{F}^{\mathcal{Z}, \mathcal{L}}_n=\mathscr{F}^{\mathcal{Z}}_n=\sigma(X_k:k\le n+1).
		$$ 
		Therefore, 	$\big((\mathcal{Z}_n,\mathcal{L}^n)\big)_{n\ge 0}$ is adapted to  $\big\{\mathscr{F}^\mathcal{Z}_n\big\}_{n\ge 0}$. 
		
		\vskip2mm
		\noindent{\small\bf  \emph{Step 2}}. 	For all bounded measurable function $g:V_{\vec{l}(G)}\mapsto\mathbb{R}$,
		\begin{equation}
			\E(g(\mathcal{Z}_{n+1})|\mathscr{F}^{\mathcal{Z}}_n)=\int_{V_{\vec{l}(G)}} g(\zz) \pp_{n\mathcal{L}^n+\dd_{\mathcal{Z}_n}}(\mathcal{Z}_n, {\rm d}\zz)=\E(g(\mathcal{Z}_{n+1})|\mathcal{Z}_n,\mathcal{L}^n).\label{eq-1-prop-Markov-(Z,L)}
		\end{equation}	
		For every $\mu\in  \mathscr{P}(V_{\vec{l}(G)}$ and $z\in V_{\vec{l}(G)}$,  define  
		$$
		\hat{g}(\mu, z,\omega):= \E\big(f(\mathcal{Z}_{n+1},\frac{1}{n+1}(n\mu+\dd_{z}))\big|\mathscr{F}^{\mathcal{Z}}_n\big),
		$$ 		
		and
		$$
		\check{g}(\mu, z,\omega):= \E\big(f(\mathcal{Z}_{n+1},\frac{1}{n+1}(n\mu+\dd_{z}))\big|\mathcal{Z}_n,\mathcal{L}^n \big).
		$$ 	
		Equality \eqref{eq-1-prop-Markov-(Z,L)} implies that 
		\begin{equation}\label{5.1.1-Z-2}
			\hat{g}(\mu, z,\omega)=\check{g}(\mu, z,\omega).
		\end{equation} 
		
		\vskip2mm		
		\noindent{\small\bf  \emph{Step 3}}. 	Note that $(\mathcal{Z}_n,\mathcal{L}^n)$ is $\mathscr{F}^{\mathcal{Z}}$-measurable and $\mathcal{L}^{n+1}=\frac{1}{n+1}(n\mathcal{L}^n+\dd_{\mathcal{Z}_n})$. Hence
		$$
		\hat{g}(\mathcal{L}^n, \mathcal{Z}_n,\omega):= \E\big(f(\mathcal{Z}_{n+1},\frac{1}{n+1}(n\mathcal{L}^n,+\dd_{\mathcal{Z}_n}))\big|\mathscr{F}^{\mathcal{Z}}_n\big)=\E\big(f(\mathcal{Z}_{n+1},\mathcal{L}^n)\big|\mathscr{F}^{\mathcal{Z}}_n\big).
		$$ 	
		Similarly, 
		$$
		\check{g}(\mathcal{L}^n, \mathcal{Z}_n,\omega):= \E\big(f(\mathcal{Z}_{n+1},\frac{1}{n+1}(n\mathcal{L}^n,+\dd_{\mathcal{Z}_n}))\big|\mathcal{Z}_n, \mathcal{L}^n\big)=\E\big(f(\mathcal{Z}_{n+1},\mathcal{L}^n)\big|\mathcal{Z}_n, \mathcal{L}^n\big).
		$$ 
		Thus,  \eqref{5.1.1-Z-2} yields  \eqref{eq-1-prop-Markov-(Z,L)}. 
	\end{proof}

	\begin{proof}[Proof of Proposition \ref{continuity of map T}]
		Fix $\hat{\mu}_0\in\mathscr{P}(V_{\vec{l}(G)})$ and $\varepsilon>0$. For any $\hat{\mu}\in\mathscr{P}(V_{\vec{l}(G)})$ satisfying  
		$$
		|\hat{\mu}(\mathbf{z})-\hat{\mu}_0(\mathbf{z})|<\varepsilon \text{ \ for all \ } \mathbf{z}\in \mathscr{P}(V_{\vec{l}(G)}), $$
		we have  for every $v\in V$,
		\[
		\left|T(\hat{\mu})(v)-T(\hat{\mu}_0)(v)\right|\le \sum_{\mathbf{z}:\,\mathbf{z}^-=v}\left|\hat{\mu}(\mathbf{z})-\hat{\mu}_0(\mathbf{z})\right|< \b \varepsilon.
		\]
		This implies the continuity of $T$, with respect to the weak‑convergence topology, completing the proof.	\end{proof}

	\subsubsection{Proofs and Supplement Results in Section \ref{sec 2.16}} \label{pf of sec 2.16}
	\begin{lemma}\label{closed set under topology}
		$\mathcal{C}l(\mathscr{S}_0)$ is a closed subset of $\mathscr{P}(V_{\vec{l}(G)})$ with respect to  the weak convergence topology.
	\end{lemma}
	\begin{proof} Assume, for contradiction, that there exists a sequence $\hat{\mu}_n\Rightarrow\hat{\mu}$ in $\mathscr{P}(V_{\vec{l}(G)})$ with $\hat{\mu}_n\in\mathcal{C}l(\mathscr{S}_0)$ but $\hat{\mu}\notin\mathcal{C}l(\mathscr{S}_0)$. Then, for every $E'\in\mathscr{S}_0$, 
		$$
		\text{supp} (\hat{\mu}|_E)\setminus E'\neq\emptyset.
		$$
		Because  $V_{\vec{l}(G)}$ is finite, we may extract a subsequence $\hat{\mu}_{n_k}$ and find sets $E_0, E_1\in\mathscr{S}_0$ with  $E_0\subset E_1$ such that 
		$$
		E_0=\text{\rm supp}(\hat{\mu}_{n_k}|_E), \ \text{for all}\ k.
		$$ 
		Note that $\text{supp}(\hat{\mu}|_E)\setminus E_1\neq \emptyset$. Choose an edge $e\in\text{supp}(\hat{\mu}|_E)\setminus E_1$; then $\hat{\mu}|_E(\{e\})\neq0$ while $\hat{\mu}_{n_k}|_E(\{e\})=0$ for all $k$. This contradicts the weak convergence $\hat{\mu}_{n_k}\Rightarrow\hat{\mu}$. Hence  $\mathcal{C}l(\mathscr{S}_0)$ must be closed.
	\end{proof}
	
	\noindent \begin{proof}[Proof of Proposition \ref{prop domain of rate function}]
		Here we prove that $\Lambda^\zz_{\delta,\mathscr{S}_0}(\hat{\mu})<\infty$ if and only if $\hat{\mu}\in \mathcal{C}l(\mathscr{S}_0)$, and $\hat{\mu}$ is invariant under some $\qq\in\mathscr{T}_{\vec{l}(G)}$. The proof for $\Lambda_{\delta,\mathscr{S}_0}(\hat{\mu})<\infty$ follows a similar argument.\\

		\noindent{\small\bf  \emph{Step 1}}. We prove the implication "$\Rightarrow$": If $\Lambda^\zz_{\delta,\mathscr{S}_0}(\hat{\mu})<\infty$, then $\hat{\mu}\in \mathcal{C}l(\mathscr{S}_0)$ and $\hat{\mu}$ is invariant under  some $\qq\in\mathscr{T}_{\vec{l}(G)}$.
		
		\vskip2mm
		
		\noindent{\small\bf  \emph{Step 1.1}}.  We first show that $\hat{\mu}\in\mathcal{C}l(\mathscr{S}_0)$.\vskip1mm
		
		There exists a sequence of tuples $(\hat{\mu}_k,r_k,E_k,\qq_k)_{1\le k\le \b}\in \hat{\mathscr{A}}_\zz(\hat{\mu},\mathscr{S}_0)$ such that the infimum in \eqref{Lambda^z} is attained at this sequence.  Let $l$ 
		be the smallest integer such that $E_l\notin \mathscr{S}_0$. From the definitions of $\hat{\mathscr{A}}_\zz(\hat{\mu},\mathscr{S}_0)$ in \eqref{A_z} and $\mathcal{C}l(\mathscr{S_0})$ in \eqref{eq-closed-subset}, we conclude that $r_k=0$ for $k\ge l$, and that ${\rm supp}(\hat{\mu}_k)\subseteq \vec{E}_{l-1}$, where $E_{l-1}\in\mathscr{S}_0$. Using  $\sum_{k=1}^\b r_k \hat{\mu}_k=\hat{\mu}$, we obtain   ${\rm supp}(\hat{\mu})\subseteq \vec{E}_{l-1}$; hence $\hat{\mu}\in\mathcal{C}l(\mathscr{S}_0)$. \vskip2mm
		
		\noindent{\small\bf  \emph{Step 1.2}}.  Now we prove by contradiction  that $\Lambda_{\delta,\mathscr{S}_0}^\zz(\hat{\mu})<\infty$ implies the invariance of $\hat{\mu}$ for some $\qq\in\mathscr{T}_{\vec{l}(G)}$.
		
		\vskip1mm
		
		A key observation is that  $\hat{\mu}$ is invariant under some $\qq\in\mathscr{T}_{\vec{l}(G)}$ if and only if 
		$$ 
		\sum_{\zz^-=x}\hat{\mu}(\zz)=\sum_{\zz^+=x}\hat{\mu}(\zz)\  \text{for every vertex}\ \ x\in V,
		$$   
		where $V$ denotes the vertex set of $G$.    We postpone the proof of this observation for a moment.   \vskip1mm
		
		If there exists  a $\hat{\mu}$ that is not invariant under any   $\qq\in\mathscr{T}_{\vec{l}(G)}$ yet satisfies $\Lambda_{\delta,\mathscr{S}_0}^\zz(\hat{\mu})<\infty$, then there is some vertex    $x\in V$ such that 
		$$
		\sum_{\zz^-=x}\hat{\mu}(\zz)\neq\sum_{\zz^+=x}\hat{\mu}(\zz).
		$$ 
		Assume the infimum in \eqref{Lambda^z} is attained at $(\hat{\mu}_k,r_k,E_k,\qq_k)_{1\le k\le \b}$,  where $\sum_{k=1}^{\b}\hat{\mu}_k r_k=\hat{\mu}$ and $\hat{\mu}_k\qq_k=\hat{\mu}_k$ for all $k=1,\dots,\b$. Then  applying the same observation to each invariant measure    $\hat{\mu}_k$ gives
		$$
		\sum_{\zz^-=x}\hat{\mu}_k(\zz)\neq\sum_{\zz^+=x}\hat{\mu}_k(\zz), \ \ \text{for all } k.
		$$ 
		Summing over $k$ with weight $r_k$ yields
		$$
		\sum_{\zz^-=x}\hat{\mu}(\zz)=\sum_{\zz^+=x}\hat{\mu}(\zz),
		$$ 
		which contradicts the initial  assumption. Consequently, for any $\hat{\mu}$ that is not invariant under every    $\qq\in\mathscr{T}_{\vec{l}(G)}$, we have  ${\Lambda}_{\delta,\mathscr{S}_0}^\zz(\hat{\mu})=\infty$.   
		\vskip1mm     
		
		\noindent{\small\bf  \emph{Step 1.3}}. We now verify the observation.

		\begin{itemize}
			\item If $\hat{\mu}$ is invariant under $\qq$, then for every  vertex $x$ of  $G$ and for $\zz_1,\zz_2$ with $\zz_1^+=\zz_2^-=x$, we have
			\[
			\sum_{\zz^-=x}\qq(\zz_1,\zz)=1,\ \sum_{\zz^+=x}\qq(\zz,\zz_2)\hat{\mu}(\zz)=\hat{\mu}(\zz_2).
			\]
			Consequently, 
			\begin{align*}
				\sum_{\zz_2^-=x}\hat{\mu}(\zz_2)&=\sum_{\zz_2^-=x}\sum_{\zz_1^+=x}\qq(\zz_1,\zz_2)\hat{\mu}(\zz_1)\\
				&=\sum_{\zz_1^+=x}\sum_{\zz_2^-=x}\qq(\zz_1,\zz_2)\hat{\mu}(\zz_1)\\
				&=\sum_{\zz_1^+=x}\hat{\mu}(\zz_1).
			\end{align*}
			
			\item Conversely, If 
			$$
			\sum_{\zz^-=x}\hat{\mu}(\zz)=\sum_{\zz^+=x}\hat{\mu}(\zz), \ \text{for all vertex}\ \ 
			x \in V,
			$$  
			define a  transition probability by
			$$
			\qq(\zz_1,\zz_2)=\frac{\hat{\mu}(\zz_2)}{\sum_{\zz^+=x}\hat{\mu}(\zz)}\ \ \text{for all}\ \  	\zz_1,\zz_2\in V_{\vec{l}(G)}\ \ \text{with}	\ \ \zz_1\to \zz_2.
			$$
			Then for every 
			$
			\zz_1\in V_{\vec{l}(G)}$, 
			$$
			\sum_{\zz^-=\zz_1^+}\qq(\zz_1,\zz)=1.
			$$
			Moreover, for any $\zz_2\in V_{\vec{l}(G)}$,
			$$
			\sum_{\zz^+=\zz_2^-}\qq(\zz,\zz_2)\hat{\mu}(\zz)=\sum_{\zz^+=\zz_2^-}\frac{\hat{\mu}(\zz_2)}{\sum_{\zz^+=\zz_2^-}\hat{\mu}(\zz)}\hat{\mu}(\zz)=\hat{\mu}(\zz_2).
			$$
			Hence, $\hat{\mu}$ is invariant under $\qq$.	         			
		\end{itemize}
		
		\noindent{\small\bf  \emph{Step 2}}.    We prove the implication "$\Leftarrow$":    if $\hat{\mu}\in \mathcal{C}l(\mathscr{S}_0)$, and is invariant for some $\qq\in\mathscr{T}_{\vec{l}(G)}$, then $\Lambda_{\delta,\mathscr{S}_0}^\zz(\hat{\mu})<\infty$.      
		
		\vskip1mm  
		Note that
		\[
		\Lambda_{\delta,\mathscr{S}_0}^\zz(\hat{\mu})\le \int_{V_{\vec{l}(G)}}R(\qq\|\pp_{E'})\ {\rm d}\hat{\mu},
		\]
		where $E'\subseteq \mathscr{S}_0$ is the support of $\hat{\mu}|_E$.
		By the definition of $\pp_{E'}$ in Definition \ref{Def for p_u} and because $\qq\in\mathscr{T}_{\vec{l}(G)}$, we have that $\qq(\zz,\cdot)\ll \pp_{E'}(\zz,\cdot)$ for every $\zz\in V_{\vec{l}(G)}$. Since $V_{\vec{l}(G)}$ is finite, it follows  that 
		$$
		\int_{V_{\vec{l}(G)}}R(\qq\|\pp_{E'})\ {\rm d}\hat{\mu}<\infty,
		$$ which completes the proof.      
	\end{proof}
	
	\vskip2mm
	
	\begin{proof}[Proof of Proposition \ref{lower semicontinuous}]\label{appendix 3}
		Recall the definition of ${\Lambda}_{\delta,\mathscr{S}_0}(\hat{\mu})$ in \eqref{rate function},
		\begin{align*}
			{\Lambda}_{\delta,\mathscr{S}_0}(\hat{\mu}):=\inf_{(\hat{\mu}_k,r_k,E_k,\qq_k)_{k}\in\hat{\mathscr{A}}(\hat{\mu},\mathscr{S}_0)}
			\sum_{k=1}^{\b}r_k\int_{V_{\vec{l}(G)}}R(\qq_{k}\|\pp_{E_k})\ {\rm d}\hat{\mu}_{k}.%\label{rate function}
		\end{align*}
		We first present an equivalent representation of ${\Lambda}_{\delta,\mathscr{S}_0}(\hat{\mu})$, which will be convenient for our subsequent analysis.
		\vskip2mm
		For every $\{E_k\}_{1\le k\le \b}\in\mathscr{E}$, define
		\begin{align*}
			{\hat{\mathscr{A}}}|_{\{E_k\}}(\hat{\mu},\mathscr{S}_0):=& \Big\{ (\hat{\mu}_k,r_k)_k:\ \sum_{k=1}^{\b}r_k\hat{\mu}_k=\hat{\mu},\ \text{\rm supp}(\hat{\mu}_k|_E)\subset E_k,\\
			&\hskip 2.4cm r_k\ge 0,\sum_{k=1}^{\b}r_k=1,\ r_l=0\text{ for }E_l\notin \mathscr{S}_0  \Big\}.
		\end{align*}
		Then
		$$
		{\Lambda}_{\delta,\mathscr{S}_0}(\hat{\mu})=\inf_{\{E_k\}\in\mathscr{E}} \Big\{\inf_{(\hat{\mu}_k,r_k)_k\in {\hat{\mathscr{A}}}|_{\{E_k\}}(\hat{\mu},\mathscr{S}_0)}\big\{\sum_{k=1}^\b r_k\inf_{\qq\in\mathscr{T}_{\vec{l}(G)}:\,\hat{\mu}_k \qq=\hat{\mu}_k}\int_{V_{\vec{l}(G)}} R(\qq\|\pp_{E_k})\ {\rm d}\hat{\mu}_k\big\}\Big\}.
		$$
		For brevity, introduce
		\begin{align*}
			\Lambda_{E_k}^0(\hat{\mu}):=&\inf_{\qq\in\mathscr{T}_{\vec{l}(G)}:\,\hat{\mu} \qq=\hat{\mu}}\int_{V_{\vec{l}(G)}} R(\qq\|\pp_{E_k})\ {\rm d}\hat{\mu},\ k=1,\cdots,d,\\
			\Lambda_{\{E_k\}}^1(\hat{\mu}):=& \inf_{(\hat{\mu}_k,r_k)_k\in {\hat{\mathscr{A}}}|_{\{E_k\}}(\hat{\mu},\mathscr{S}_0)}\sum_{k=1}^\b r_k\Lambda_{E_k}^0(\hat{\mu}_k).
		\end{align*}
		That is 
		\begin{equation}\label{5-pr-2.3-1}
			{\Lambda}_{\delta,\mathscr{S}_0}(\hat{\mu})=\inf_{\{E_k\}\in\mathscr{E}} \big\{ \Lambda_{\{E_k\}}^1(\hat{\mu}) \big\}.
		\end{equation}
		
		\noindent{\small\bf  \emph{Step 1}}.  We show for any $\{E_k\}_{1\le k\le \b}\in\mathscr{E}$ and  $\hat{\mu}\in\mathscr{P}(V_{\vec{l}(G)})$, there exist $(\tilde{r}_1,\cdots, \tilde{r}_\b)$ and  $(\tilde{\hat{\mu}}_k)_{1\le k\le\b}$ such that
		\begin{equation}\label{Lem2.4-1}
			\sum_{k=1}^\b \tilde{r}_k \Lambda_{E_k}^0(\tilde{\hat{\mu}}_k)=\Lambda_{\{E_k\}}^1(\hat{\mu}),
		\end{equation}
		where  $\tilde{r}_k\ge 0$,  $\sum_{k=1}^\b\tilde{r}_k=1$, $\sum_{k=1}^\b \tilde{r}_k\tilde{\hat{\mu}}_k=\hat{\mu}$ and $\text{\rm{supp}}(\tilde{\hat{\mu}}_k|_E)\subset E_k$.
		
		\vskip2mm
		
		By the definition of $ \Lambda_{\{E_k\}}^1(\hat{\mu})$, there is a sequence $\big\{(\hat{\mu}_k^n,r_k^{(n)})_k\big\}_{n\ge 1}\subset{\hat{\mathscr{A}}}|_{\{E_k\}}(\hat{\mu},\mathscr{S}_0)$ such that
		$$
		\sum_{k=1}^\b r_k^{(n)} \Lambda_{E_k}^0(\hat{\mu}_k^n)<\Lambda_{\{E_j\}}^1(\hat{\mu})+\frac{1}{n}.
		$$
		Then by the diagonal argument and the lower-semicontinuity of $\Lambda_{E_k}^0(\hat{\mu}_k)$, $(k=1,\dots,\b)$, there is a subsequence such that
		$$
		\big(\hat{\mu}_k^{n'},r_k^{(n')}\big)_{1\le k\le \b}\to \big(\tilde{\hat{\mu}}_k,\tilde{r}_k\big)_{1\le k\le \b}, \ \ \text{as}\ \ n'\to\infty\ \  \mbox{and}\ \sum_{k=1}^\b \tilde{r}_k \Lambda_{E_k}^0(\tilde{\hat{\mu}}_k)\le\Lambda_{\{E_j\}}^1(\hat{\mu}).
		$$
		This implies \eqref{Lem2.4-1}.
		\vskip2mm
		
		Note that by the equivalent representation \eqref{5-pr-2.3-1} and the finiteness of $\mathscr{E}$,  the lower semicontinuity and compactness of ${\Lambda}_{\delta,\mathscr{S}_0}$ are equivalent to the corresponding properties of $\Lambda_{\{E_k\}}^1(\hat{\mu})$. Hence it suffices to prove these properties for $\Lambda_{\{E_k\}}^1(\hat{\mu})$ with a fixed $\{E_k\}_{1\le k\le \b}$.
		
		\vskip2mm
		
		\noindent{\small\bf  \emph{Step 2}}. We verify the lower semicontinuity. \vskip2mm
		
		For an arbitrary sequence $\hat{\mu}^n\Rightarrow \hat{\mu}$, by the result of {\small\bf  \emph{Step 1}},  we can choose $(\hat{\mu}_k^n)_{1\le k\le \b}$ and $(r_k^{(n)})_{1\le k\le \b}$ such that
		$$
		\sum_{k=1}^\b r_k^{(n)} \Lambda_{E_k}^0(\hat{\mu}_k^n)=\Lambda_{\{E_k\}}^1(\hat{\mu}^n).
		$$
		Since $V_{\vec{l}(G)}$ is finite, $\{\hat{\mu}_k^n\}_{n\ge 1}$ is tight.  By the diagonal argument, there exists a subsequence $n'$ such that
		\begin{itemize}
			\item $\hat{\mu}_k^{n'}\Rightarrow \hat{\mu}_k$;
			\item $\lim\limits_{n'\to\infty} \Lambda_{\{E_k\}}^1(\hat{\mu}^{n'})=\liminf\limits_{n\to\infty}\Lambda_{\{E_k\}}^1(\hat{\mu}^n)$;
			\item $r_k^{n'}\to r_k$;
			\item $\hat{\mu}=\sum_{k=1}^\b r_k\hat{\mu}_k$.
		\end{itemize}
		Using the lower semicontinuity of $\Lambda_{E_k}^0$, we obtain
		$$
		\liminf_{n\to\infty}\Lambda_{\{E_k\}}^1(\hat{\mu}^n)\ge\sum_{k=1}^\b r_k\Lambda_{E_k}^0(\hat{\mu}_k)\ge \Lambda_{\{E_k\}}^1(\hat{\mu}).
		$$
		\vskip2mm
		\noindent{\small\bf  \emph{Step 3}}. We show the compactness. 
		\vskip2mm
		
		Note that $\mathscr{P}(V_{\vec{l}(G)})$ is compact in the weak convergence topology because $V_{\vec{l}(G)}$ is finite. Moreover the set  
		$$
		\big\{\hat{\mu}:\Lambda_{\{E_k\}}^1(\hat{\mu})\le M\big\}\subseteq\mathscr{P}(V_{\vec{l}(G)})
		$$ is closed due to  the lower semicontinuity of $\Lambda_{\{E_k\}}^1$.  Consequently, the level set $\big\{\hat{\mu}:\ \Lambda_{\{E_k\}}^1(\hat{\mu})\le M \big\}$ is compact.
	\end{proof}
	
	\vskip2mm
	\begin{remark}  For every $0<\lambda<1$ and $\tilde{\hat{\mu}},\bar{\hat{\mu}}\in\mathscr{P}(V_{\vec{l}(G)})$,  the convexity of $\Lambda_{E_k}^0$ yields 
		\begin{align*}
			\lambda \Lambda_{\{E_k\}}^1(\tilde{\hat{\mu}})+(1-\lambda) \Lambda_{\{E_k\}}^1(\bar{\hat{\mu}})&=\sum_{k=1}^\b \lambda \tilde{r}_k \Lambda_{E_k}^0(\tilde{\hat{\mu}}_k)+(1-\lambda)\bar{r}_k \Lambda_{E_k}^0(\overline{\hat{\mu}}_k)\\
			&\ge \sum_{k=1}^\b (\lambda \tilde{r}_k+(1-\lambda)\bar{r}_k)\Lambda_{E_k}^0\Big(\frac{\lambda \tilde{r}_k\tilde{\hat{\mu}}_k+(1-\lambda)\bar{r}_k\bar{\hat{\mu}}_k}{\lambda \tilde{r}_k+(1-\lambda)\bar{r}_k}\Big)\\
			&\ge \Lambda_{\{E_k\}}^1(\lambda\tilde{\hat{\mu}}+(1-\lambda)\bar{\hat{\mu}}).
		\end{align*}
		Thus  $\Lambda_{\{E_k\}}^1$ is convex for a fixed  
		$\{E_k\}_{1\le k\le \b}$. However, it remains unclear whether the full functional  ${\Lambda}_{\delta,\mathscr{S}_0}(\hat{\mu})$ is convex.
	\end{remark}
	
	\begin{proof}[Proof of Theorem \ref{thm-rate-function-delta=1}]  Note that the relative entropy $R(\cdot\|\cdot)$ is convex. By \eqref{rate function}, \begin{align}
			{\Lambda}_{1,\mathscr{S}_0}(\hat{\mu})&=\inf_{(\hat{\mu}_k,r_k,E_k,\qq_k)_{k}\in\hat{\mathscr{A}}(\hat{\mu},\mathscr{S}_0)}
			\sum_{k=1}^{\b}r_k\int_{V_{\vec{l}(G)}}R(\qq_{k}\|\pp)\ {\rm d}\hat{\mu}_{k}\nonumber\\
			&= \inf_{(\hat{\mu}_k,r_k,E_k,\qq_k)_{k}\in\hat{\mathscr{A}}(\hat{\mu},\mathscr{S}_0)}
			\sum_{k=1}^{\b}r_k R(\hat{\mu}_{k}\otimes\qq_{k} \|\hat{\mu}_{k}\otimes\pp)\nonumber\\
			&\ge \inf_{(\hat{\mu}_k,r_k,E_k,\qq_k)_{k}\in\hat{\mathscr{A}}(\hat{\mu},\mathscr{S}_0)}
			R\left(\sum_{k=1}^{\b}r_k\hat{\mu}_{k}\otimes\qq_{k} \left\|\sum_{k=1}^{\b}r_k\hat{\mu}_{k}\otimes\pp\right)\right..\label{eq-1-proof-Lambda-1}
		\end{align}
		For $l$ such that $E_l\in\mathscr{S}_1$ and $E_{l+1}\notin\mathscr{S}_0$, we have  $r_k=0$ for all $k\ge l+1$. Hence, for a fixed $\{E_k\}_{1\le k\le \b}\in\mathscr{E}$, we obtain  
		$$
		\sum_{k=1}^{\b}r_k\hat{\mu}_{k}=\hat{\mu_l},
		$$
		and  the first marginal of $\sum_{k=1}^{\b}r_k\hat{\mu}_{k}\otimes\qq_{k}$ is $\hat{\mu}_l$. 
		Since for every $A\subset \vec{l}(G)$,
		\[
		\int_{\vec{l}(G)}\sum_{k=1}^{\b}r_k\hat{\mu}_{k}({\rm d}x)\otimes\qq_{k}(x,A)
		=
		\sum_{k=1}^{\b}r_k\hat{\mu}_{k}(A)
		=\hat{\mu}_l(A),
		\]
		there exists a transition probability $\qq\in\mathscr{T}_{\vec{l}(G)}$ satisfying $\hat{\mu}_l\qq=\hat{\mu}_l$ and
		\[
		\hat{\mu}_l\otimes \qq = \sum_{k=1}^{\b}r_k\hat{\mu}_{k}\otimes\qq_{k}.
		\]
		Therefore,
		\begin{align}
			\eqref{eq-1-proof-Lambda-1}
			&= \inf_{(\hat{\mu}_k,r_k,E_k,\qq_k)_{k}\in\hat{\mathscr{A}}(\hat{\mu},\mathscr{S}_0)}
			R(\hat{\mu}_l\otimes\qq\|\hat{\mu}_l\otimes \pp)\nonumber\\
			&\ge \inf_{(\hat{\mu}, \{ E_k \}_{1\le k\le \b}, \qq)\in\hat{\mathscr{A}}(1,\hat{\mu},\mathscr{S}_0) }R(\hat{\mu}\otimes\qq \|\hat{\mu}\otimes \pp)\nonumber\\
			&= \inf_{(\hat{\mu}, \{ E_k \}_{1\le k\le \b}, \qq)\in\hat{\mathscr{A}}(1,\hat{\mu},\mathscr{S}_0) }\int_{\vec{E}_l}R(\qq \| \pp){\rm d} \hat{\mu}.\label{eq-2-proof-Lambda-1}
		\end{align}
		For the reverse inequality, note that the infimum condition in \eqref{eq-2-proof-Lambda-1} is a specific case  in \eqref{Lambda_1} (by setting all $r_k=0$ for $k<l$). This derives that $\Lambda_{1,\mathscr{S}_0}\le \eqref{eq-2-proof-Lambda-1}$, completing the proof of \eqref{Lambda_1}.
		%Specifically for $\mathscr{S}_0=\mathscr{S}$, $l$ is exactly equal to $\b$, and $E_l$ is always equal to $E$. Therefore, by \eqref{eq-2-proof-Lambda-1},
		%$
		%\Lambda(\hat{\mu})=\inf_{\qq\in\mathscr{T}_{\vec{l}(G)}:\hat{\mu}\qq=\hat{\mu}}\int_{\vec{E}}R(\qq \| \pp){\rm d} \hat{\mu}.
		%$
	\end{proof}

	%%%%%%%%%%%%%%%%%%%%%%%%%%%%%%%%%%%%

	\subsubsection{Proof of Lemma \ref{Theorem 4.2.2}} \label{5-pr-2.5-lem1}
	
	\noindent Now, we prove Lemma \ref{Theorem 4.2.2}.  The overall structure of the proof follows the framework of \cite[Theorems 1.5.2 and 8.2.1]{DE1997}.
	\begin{itemize}
		\item First, set $\mathcal{L}_j^n=\frac{j}{n}\mathcal{L}^j$.  Then  $(\mathcal{Z}_j,\mathcal{L}_j^n)_{j=0,1,\dots,n}$ is a homogeneous Markov process with the transition probability
		$$
		\P\big((\mathcal{Z}_{i+1},\mathcal{L}_{i+1}^n)\in {\rm d}\zz_2\times {\rm d}\hat{\nu}\big|(\mathcal{Z}_{i},\mathcal{L}_{i}^n)=(\zz_1,\hat{\mu})\big)=\pp_{\hat{\mu}+\dd_{\zz_1}}(\zz_1,{\rm d}\zz_2)\cdot \dd_{\hat{\mu}+\frac{1}{n}\dd_{\zz_1}}({\rm d}\hat{\nu}).
		$$
		Define
		\begin{equation}\label{6-W}
			W^n_{\mathscr{S}_0}(i,\zz,\hat{\mu})=\left\{
			\begin{array}{ll}
				-\frac{1}{n}\log{\E_{i,\zz,\hat{\mu}}\left\{\exp[-nh(\mathcal{L}^n)]\textbf{1}_{\{\mathcal{L}^n\in\mathcal{C}l(\mathscr{S}_0)\}} \right\}}, & \P_{i,\zz,\hat{\mu}}(\mathcal{L}^n\in\mathcal{C}l(\mathscr{S}_0))>0,  \\
				\infty, &   \P_{i,\zz,\hat{\mu}}(\mathcal{L}^n\in\mathcal{C}l(\mathscr{S}_0))=0,
			\end{array}
			\right.
		\end{equation}
		where $\P_{i,\zz,\hat{\mu}}$ and  $\E_{i,\zz,\hat{\mu}}$ denote the probability and expectation conditioned on $(\mathcal{Z}_i,\mathcal{L}^n_i)=(\zz,\hat{\mu})$, respectively.
		
		We show that $W^n_{\mathscr{S}_0}(i,\zz,\hat{\mu})$   can be expressed via the variational  formula \eqref{equation of dynamic programming} in Lemma \ref{Lemma of dynamic programming equation} below, which is interpreted as a dynamic programming equation.
		
		\item  We then verify that \eqref{equation of dynamic programming} satisfies the \emph{Attainment Condition} stated in Lemma \ref{Lemma of attainment condition} below.
		
		\item We demonstrate that the terminal condition of $W^n_{\mathscr{S}_0}(i,\zz,\hat{\mu})$ is also satisfied by the minimal cost function $\vv^n_{\mathscr{S}_0}(i,\zz,\hat{\mu})$ provided in \eqref{6-dyn-prg}. Finally, using the same argument  as  \cite[Theorem 1.5.2]{DE1997} , we obtain  
		$$
		W^n_{\mathscr{S}_0}(i,\zz,\hat{\mu})=\vv^n_{\mathscr{S}_0}(i,\zz,\hat{\mu}),
		$$ which establishes  Lemma \ref{Theorem 4.2.2}.
		
	\end{itemize}
	
	However,   the restricted logarithmic Laplace functionals  $W^n_{\mathscr{S}_0}(i,\zz,\hat{\mu})$ in our setting differ from those in  \cite{DE1997}. Consequently, we must introduce additional notations related to  
	$\mathcal{C}l(\mathscr{S}_0)$, establish auxiliary lemmas, and adapt the approach of  \cite{DE1997}.

	Define {\color{black}
		\begin{equation}
			\mathbb{S}:=\Big\{(\zz,\hat{\mu}):\ \hat{\mu}+\frac{1}{n}\dd_\zz\in \mathcal{C}l(\mathscr{S}_0)\Big\}, \ S':=\Big\{\zz'\in V_{\vec{l}(G)}:\ (\zz',\hat{\mu}+\frac{1}{n}\dd_\zz)\in\mathbb{S}\Big\}.\label{condition for dynamic programming}
	\end{equation}}
	
	\begin{lemma}\label{Lemma of dynamic programming equation}
		If $(\zz,\hat{\mu})\in V_{\vec{l}(G)}\times \mathscr{P}(V_{\vec{l}(G)})$ belongs to $\mathbb{S}$, 
		then $W^n_{\mathscr{S}_0}(i,\zz,\hat{\mu})$ satisfies the following  equation on $S'$:
		\begin{eqnarray}
			\nonumber W^n_{\mathscr{S}_0}(i,\zz,\hat{\mu})&=&\inf_{\nu\in\mathscr{P}(V_{\vec{l}(G)}):\,\text{\rm supp}(\hat{\nu})\subseteq S'}\Bigg\{\frac{1}{n}R(\hat{\nu}(\cdot)\|\pp_{\hat{\mu}+\dd_\zz}(\zz,\cdot))\\
			&&\hskip 2cm  +\int_{V_{\vec{l}(G)}} W^n_{\mathscr{S}_0}\left(i+1,\zz',\hat{\mu}+\frac{1}{n}\dd_\zz\right)\ \nu({\rm d}\zz') \Bigg\},\label{equation of dynamic programming}
		\end{eqnarray}
		where the infimum is attained at the unique measure $\nuu_i^n$ defined  by
		\begin{equation}
			\nuu_i^n(A|\zz,\hat{\mu}):= \frac{\int_A\exp[-nW^n_{\mathscr{S}_0}(i+1,\zz',\hat{\mu}+\frac{1}{n}\dd_\zz)]\pp_{\hat{\mu}+\dd_\zz}(\zz,{\rm d}\zz')}{\int_{V_{\vec{l}(G)}}\exp[-nW^n_{\mathscr{S}_0}(i+1,\zz',\hat{\mu}+\frac{1}{n}\dd_\zz)]\pp_{\hat{\mu}+\dd_\zz}(\zz,{\rm d}\zz')},\label{unique attainment}
		\end{equation}
		and $\text{\rm supp}( \nuu_i^n(\cdot|\zz,\hat{\mu}))\subseteq S'$.
	\end{lemma}
	
	\vskip2mm
	
	Recall the set $\mathbb{S}$ defined in \eqref{condition for dynamic programming}.  We now prove the attainment condition for  $W^n_{\mathscr{S}_0}(i,\zz,\hat{\mu})$ as follows.
	
	\begin{lemma}\label{Lemma of attainment condition}
		For a fixed index $i$, suppose  $\big(\overline{\mathcal{Z}}_i^n,\overline{\mathcal{L}}_i^n\big)=(\zz,\hat{\mu})\in\mathbb{S}$. Then there exists a unique sequence $\left\{\nuu_j^n\right\}_{j\ge i}$ such that: 
		\begin{itemize}
			\item the measures $\nuu_j^n\left(\cdot\left|\overline{\mathcal{Z}}_j^n,\overline{\mathcal{L}}_j^n\right.\right)$ are well-defined;
			\item   $\overline{\mathcal{L}}_j^n\in\mathcal{C}l(\mathscr{S}_0)$ for every $j\ge i$ and every sample path $\omega$;   
		\end{itemize}
		Moreover, the infimum in the {\color{black}dynamic programming equation} \eqref{equation of dynamic programming} is attained at $\left\{\nuu_j^n\right\}_{j\ge i}$.
	\end{lemma}
	\vskip2mm
	
	We postpone the proofs of the above lemmas  for a moment.   
	
	\vskip2mm
	
	\begin{proof}[Proof of Lemma \ref{Theorem 4.2.2}] Define the minimal cost function {\color{black}
			\begin{align}\label{6-dyn-prg}
				\vv^n_{\mathscr{S}_0}(i,\zz,\hat{\mu}):=&\,\inf_{\stackrel{\{\nuu_j^n\}:\nuu_j^n:V_{\vec{l}(G)}\times\mathscr{M}_{j/n}(V_{\vec{l}(G)})\mapsto V_{\vec{l}(G)}}{\text{ for }i\le j\le n-1}}\overline{\E}_{i,\zz,\hat{\mu}}\Big\{\frac{1}{n}\sum_{j=i}^{n-1}R\big(\nuu_j^n(\cdot|\overline{\mathcal{Z}}^n_{j},\overline{\mathcal{L}}^n_{j})\big\|\pp_{\overline{\mathcal{L}}^n_{j+1}}(\overline{\mathcal{Z}}^n_{j},\cdot)\big)\nonumber\\
				&\, \hskip6cm+h(\overline{\mathcal{L}}^n)+\infty\cdot\textbf{1}_{\{\overline{\mathcal{L}}^n\in\mathcal{C}l(\mathscr{S}_0)^c\}} \Big\},
		\end{align} }
		where $\overline{\E}_{i,\zz,\hat{\mu}}$ denote the expectation under $\nuu_j^n$ conditioned on $(\overline{\mathcal{Z}}_i^n,\overline{\mathcal{L}}^n_i)=(\zz,\hat{\mu})$. In particular, we adopt the convention that $\infty\cdot 0=0$	, which ensures that the definition \eqref{6-dyn-prg} is well‑defined mathematically.	
		
		By definition \eqref{6-dyn-prg} and the convention  $\infty\cdot 0=0$, it follows that
		$$
		\vv^n_{\mathscr{S}_0}(i,\zz,\hat{\mu})=\inf_{\{\nuu_j^n\}:\overline{\mathcal{L}}^n\in\mathcal{C}l(\mathscr{S}_0)}\overline{\E}_{i,\zz,\hat{\mu}}\Big\{\frac{1}{n}\sum_{j=i}^{n-1}R\big(\nuu_j^n(\cdot|\overline{\mathcal{Z}}^n_{j},\overline{\mathcal{L}}^n_{j})\big\|\pp_{\overline{\mathcal{L}}^n_{j+1}}(\overline{\mathcal{Z}}^n_{j},\cdot)\big)+h(\overline{\mathcal{L}}^n) \Big\}.
		$$
		Our goal  is to prove
		$$
		W^n_{\mathscr{S}_0}(i,\zz,\hat{\mu})=\vv^n_{\mathscr{S}_0}(i,\zz,\hat{\mu}) \text{ for } i=0,1,\dots,n.
		$$
		Taking  $i=0, \hat{\mu}(\cdot)=0$ then yields  Lemma \ref{Theorem 4.2.2}.
		
		\vskip2mm

		\noindent{\small\bf  \emph{Step 1}}. \ {\boldmath $W^n_{\mathscr{S}_0}(n,z,\hat{\mu})=\vv^n_{\mathscr{S}_0}(n,z,\hat{\mu})$.}
		
		Setting $i=n$ in (\ref{6-W}) gives
		$$
		W^n_{\mathscr{S}_0}(n,\zz,\hat{\mu})=\left\{
		\begin{array}{ll}
			h(\hat{\mu}), & \hat{\mu} \in \mathcal{C}l(\mathscr{S}_0)\\
			\infty, &  \hat{\mu} \notin \mathcal{C}l(\mathscr{S}_0)
		\end{array}\right.;
		$$
		Similarly, taking $i=n$ in (\ref{6-dyn-prg}) yields
		$$
		\vv^n_{\mathscr{S}_0}(n,\zz,\hat{\mu})=\left\{
		\begin{array}{ll}
			h(\hat{\mu}), & \hat{\mu} \in \mathcal{C}l(\mathscr{S}_0)\\
			\infty, &  \hat{\mu} \notin \mathcal{C}l(\mathscr{S}_0)
		\end{array}\right..
		$$
		Thus, $W^n_{\mathscr{S}_0}$ and  $\vv^n_{\mathscr{S}_0}$ satisfy the same terminal condition.
		\vskip2mm
		
		\noindent{\small\bf  \emph{Step 2}}. Furthermore, for $i<n$,  $\vv^n_{\mathscr{S}_0}(i,\zz,\hat{\mu})=W^n_{\mathscr{S}_0}(i,\zz,\hat{\mu})=\infty$ whenever $(\zz,\hat{\mu})\notin\mathbb{S}$.
		Hence it remains to prove 
		$$
		W^n_{\mathscr{S}_0}(i,\zz,\hat{\mu})=\vv^n_{\mathscr{S}_0}(i,\zz,\hat{\mu})\ \ \text{for}\ \  (\zz,\hat{\mu})\in\mathbb{S}.
		$$ 
		
		\noindent{\small\bf  \emph{Step 2.1}}. {\boldmath $W^n_{\mathscr{S}_0}(i,\zz,\hat{\mu})\le \vv^n_{\mathscr{S}_0}(i,\zz,\hat{\mu})\text{\emph{ for }}(\zz,\hat{\mu})\in\mathbb{S}.$} \vskip1mm
		
		Note that if the process only moves to the adjacent vertices in $\vec{l}(G)$, i.e., 
		$$
		\text{supp}(\hat{\nu}(\cdot))\subseteq\{\zz'\in \vec{l}(G):\zz\to \zz'\},
		$$ 
		then $\frac{1}{n}R(\hat{\nu}(\cdot)\|\pp_{\hat{\mu}+\frac{1}{n}\dd_\zz}(\zz,\cdot))<\infty$. 
		
		Now choose an admissible control sequence $\nuu_j^n$ such that $\overline{\mathcal{L}}^n\in\mathcal{C}l(\mathscr{S}_0)$ for almost surely $\omega$. We verify 
		\begin{equation}\text{supp}\left(\nuu_j^n(\cdot|\overline{\mathcal{Z}}_j^n,\overline{\mathcal{L}}_{j}^n)\right)\subseteq \left\{\zz':\ \big(\zz',\overline{\mathcal{L}}_{j+1}^n\big)\in\mathbb{S}\right\}\label{(7.31)}
		\end{equation} 
		by contradiction. If \eqref{(7.31)} fails, then $\overline{\mathcal{L}}_{j+1}^n\notin \mathcal{C}l(\mathscr{S}_0)$. 
		This contradicts the fact that 
		$$
		\text{supp}\big(\overline{\mathcal{L}}_{j+1}^n\big)\subseteq\text{supp}\big(\overline{\mathcal{L}}^n\big)\ \ \text{and}\ \  \overline{\mathcal{L}}^n\in\mathcal{C}l(\mathscr{S}_0),
		$$
		because $\mathscr{S}_0$ is decreasing.
		
		Consequently, by \eqref{(7.31)} the integral
		\[
		\int_{V_{\vec{l}(G)}} W^n_{\mathscr{S}_0}(j+1,\zz',\overline{\mathcal{L}}_{j+1}^n))\nuu_j^n({\rm d}\zz'|\overline{\mathcal{Z}}_j^n,\overline{\mathcal{L}}_j^n))
		\]
		remains bounded for every $j\ge i$.   Meanwhile, (\ref{(7.31)}) implies that $(\overline{\mathcal{Z}}_j^n,\overline{\mathcal{L}}_{j}^n)\in\mathbb{S}$ for all $j\ge i$, hence $W^n_{\mathscr{S}_0}(j,\overline{\mathcal{Z}}_j^n,\overline{\mathcal{L}}_j^n)$ is bounded. From {\color{black} the dynamic programming equation}  (\ref{equation of dynamic programming}), we obtain 
		\begin{align}
			&\,\frac{1}{n}R\big(\nuu_j^n(\cdot|\overline{\mathcal{Z}}_j^n,\overline{\mathcal{L}}_j^n)\big\|\pp_{\overline{\mathcal{L}}_{j+1}^n}(\overline{\mathcal{Z}}_j^n,\cdot)\big)\nonumber\\ 
			&\ \ \ge \, W^n_{\mathscr{S}_0}(j,\overline{\mathcal{Z}}_j^n,\overline{\mathcal{L}}_{j}^n)-\int_{V_{\vec{l}(G)}} W^n_{\mathscr{S}_0}(j+1,\zz',\overline{\mathcal{L}}_{j+1}^n)\nuu_j^n({\rm d}\zz'|\overline{\mathcal{Z}}_j^n,\overline{\mathcal{L}}_j^n).\label{(7.19)}
		\end{align}
		
		Therefore, using the Markov property and the terminal condition $W_{\mathscr{S}_0}^n(n,\overline{\mathcal{Z}}_n^n,\overline{\mathcal{L}}^n)=h(\overline{\mathcal{L}}^n)$ (which holds because  $\overline{\mathcal{L}}^n\in\mathcal{C}l(\mathscr{S}_0)$  almost surely),
		\begin{align}
			& \overline{\E}_{i,\zz,\hat{\mu}}\Big\{\frac{1}{n}\sum_{j=i}^{n-1}R\big(\nuu_j^n(\cdot|\overline{\mathcal{Z}}^n_{j},\overline{\mathcal{L}}^n_{j})\big\|\pp_{\overline{\mathcal{L}}^n_{j+1}}(\overline{\mathcal{Z}}^n_{j},\cdot)\big)+h(\overline{\mathcal{L}}^n) \Big\}\nonumber\\
			&\ge \overline{\E}_{i,\zz,\hat{\mu}}\Big\{\sum_{j=i}^{n-1}\Big[ W^n_{\mathscr{S}_0}(j,\overline{\mathcal{Z}}_j^n,\overline{\mathcal{L}}_j^n)-\int_{V_{\vec{l}(G)}} W^n_{\mathscr{S}_0}(j+1,\zz',\overline{\mathcal{L}}_{j+1}^n)\nuu_j^n({\rm d}\zz'|\overline{\mathcal{Z}}_j^n,\overline{\mathcal{L}}_j^n) \Big]\nonumber\\
			&\hskip6cm +W^n_{\mathscr{S}_0}(n,\overline{\mathcal{Z}}_n^n,\overline{\mathcal{L}}^n) \Big\}\nonumber\\
			&=\overline{\E}_{i,\zz,\hat{\mu}}\Big\{\sum_{j=i}^{n-1}\Big[ W^n_{\mathscr{S}_0}(j,\overline{\mathcal{Z}}_j^n,\overline{\mathcal{L}}_j^n)-\overline{\E}_{i,\zz,\hat{\mu}}\big\{ W^n_{\mathscr{S}_0}(j+1,\overline{\mathcal{Z}}_{j+1}^n,\overline{\mathcal{L}}_{j+1}^n)|\overline{\mathcal{Z}}_j^n,\overline{\mathcal{L}}_j^n)\big\} \Big]\nonumber\\
			&\hskip6cm +W^n_{\mathscr{S}_0}(n,\overline{\mathcal{Z}}_n^n,\overline{\mathcal{L}}^n) \Big\}\nonumber\\
			&=\overline{\E}_{i,\zz,\hat{\mu}}\Big\{\sum_{j=i}^{n-1}\Big[ W^n_{\mathscr{S}_0}(j,\overline{\mathcal{Z}}_j^n,\overline{\mathcal{L}}_j^n)- W^n_{\mathscr{S}_0}(j+1,\overline{\mathcal{Z}}_{j+1}^n,\overline{\mathcal{L}}_{j+1}^n) \Big]+W^n_{\mathscr{S}_0}(n,\overline{\mathcal{Z}}_n^n,\overline{\mathcal{L}}^n) \Big\}\nonumber\\
			&=W^n_{\mathscr{S}_0}(i,\zz,\hat{\mu}).\label{step of dynamic}
		\end{align}
		Since this inequality holds for every  $\{\nuu_j^n\}$ with $\overline{\mathcal{L}}^n\in\mathcal{C}l(\mathscr{S}_0)$, we obtain
		$$
		\vv^n_{\mathscr{S}_0}(i,\zz,\hat{\mu})\ge W^n_{\mathscr{S}_0}(i,\zz,\hat{\mu}).
		$$
		\vskip2mm
		\noindent{\small\bf  \emph{Step 2.2}}. {\boldmath $W^n_{\mathscr{S}_0}(i,\zz,\hat{\mu})\ge \vv^n_{\mathscr{S}_0}(i,\zz,\hat{\mu})\text{\emph{ for }}(z,\hat{\mu})\in\mathbb{S}.$}	
		\vskip2mm
		Assuming that $\{\overline{\nuu}_j^n\}$is the unique attaining sequence given by Lemma \ref{Lemma of attainment condition}.  Then inequality  \eqref{(7.19)} becomes  an equality, and by the same reasoning as in
		\eqref{step of dynamic} we obtain 
		\[
		\overline{\E}_{i,\zz,\hat{\mu}}\Big\{\frac{1}{n}\sum_{j=i}^{n-1}R\big(\overline{\nuu}_j^n(\cdot|\overline{\mathcal{Z}}^n_{j},\overline{\mathcal{L}}^n_{j})\big\|\pp_{\overline{\mathcal{L}}^n_{j+1}}(\overline{\mathcal{Z}}^n_{j},\cdot)\big)+h(\overline{\mathcal{L}}^n) \Big\}=W^n_{\mathscr{S}_0}(i,\zz,\hat{\mu}).
		\] 
		By Lemma \ref{Lemma of attainment condition}, $\overline{\mathcal{L}}^n\in\mathcal{C}l(\mathscr{S}_0)$;
		consequently, 
		$$
		\vv^n_{\mathscr{S}_0}(i,\zz,\hat{\mu})\le W^n_{\mathscr{S}_0}(i,\zz,\hat{\mu}).
		$$ 
		Together with the opposite inequality already established, this yields  
		$$
		\vv^n_{\mathscr{S}_0}(i,\zz,\hat{\mu})= W^n_{\mathscr{S}_0}(i,\zz,\hat{\mu}).
		$$
	\end{proof}
	
	\begin{proof}[Proof of Lemma \ref{Lemma of dynamic programming equation}] Note that
		\begin{align*}
			\exp\big[-nW^n_{\mathscr{S}_0}(i,\zz,\hat{\mu})\big]&=\E_{i,\zz,\hat{\mu}}\big\{\exp[-nh(\mathcal{L}^n)]\mathbf{1}_{\{\mathcal{L}^n\in\mathcal{C}l(\mathscr{S}_0)\}} \big\}\\ 
			&=\E_{i,\zz,\hat{\mu}}\big\{\E_{i+1,\mathcal{Z}_{i+1},\mathcal{L}^n_{i+1}}\exp[-nh(\mathcal{L}^n)]\mathbf{1}_{\{\mathcal{L}^n\in\mathcal{C}l(\mathscr{S}_0)\}}
			\big\}\\
			&=\E_{i,\zz,\hat{\mu}}\big\{\exp\big[-nW^n_{\mathscr{S}_0}(i+1,\mathcal{Z}_{i+1},\mathcal{L}^n_{i+1})\big]\big\}\\
			&=\int_{V_{\vec{l}(G)}}\exp\big[-nW^n_{\mathscr{S}_0}(i+1,\zz',\hat{\mu}+\frac{1}{n}\dd_\zz)\big]\pp_{\hat{\mu}+\dd_\zz}(\zz,{\rm d}\zz').
		\end{align*}	
		For a fixed $i$, $W^n_{\mathscr{S}_0}(i+1,\cdot,\hat{\mu}+\frac{1}{n}\dd_\zz)$ is not identically infinite when $(\zz,\hat{\mu})\in\mathbb{S}$. In fact,
		$$
		W^n_{\mathscr{S}_0}(i+1,\cdot,\hat{\mu}+\frac{1}{n}\dd_\zz)\neq \infty\ \mbox{if and only if}\ \zz'\in S'.
		$$
		Since $S'$ is finite,  $W^n_{\mathscr{S}_0}(i+1,\cdot,\hat{\mu}+\frac{1}{n}\dd_\zz)$ is bounded on $S'$. Using the variational representation in Proposition \ref{variational representation}, we obtain  the dynamic programming equation \eqref{equation of dynamic programming}, where the infimum is attained at the unique  sequence $\{\nuu_j^n\}$ defined in \eqref{unique attainment}.  This attaining sequence is well-defined because $W^n_{\mathscr{S}_0}(i+1,\cdot,\hat{\mu}+\frac{1}{n}\dd_\zz)$ is not identically infinity; consequently, the denominator in \eqref{unique attainment} is  non‑zero. 
		
		Note that
		\[
		\exp\big[-nW^n_{\mathscr{S}_0}\big(i+1,\zz',\hat{\mu}+\frac{1}{n}\dd_\zz\big)\big]=0\ \mbox{\ if\ }\ \zz'\notin S'.
		\]
		From \eqref{unique attainment}, it follows that $\nuu_i^n(\cdot|\zz,\hat{\mu})$ is supported on $S'$. 
	\end{proof}
	
	% \vskip2mm
	
	\begin{proof}[Proof of Lemma \ref{Lemma of attainment condition}] 
		We show that $\big(\overline{\mathcal{Z}}_{j}^n,\overline{\mathcal{L}}_{j}^n\big)\in\mathbb{S}$ for all $j\ge i$ by induction. Assume this property holds for some $j\ge i$. By Lemma \ref{Lemma of dynamic programming equation}, we have that
		\[
		\text{\rm supp}\big( \nuu_j^n\big(\cdot\big|\overline{\mathcal{Z}}_{j}^n,\overline{\mathcal{L}}_{j}^n\big)\big)\subseteq\big\{\zz':\ \big(\zz', \overline{\mathcal{L}}_{j+1}^n\big)\in\mathbb{S}\big\},
		\]
		which immediately yields  $\big(\overline{\mathcal{Z}}_{j+1}^n,\overline{\mathcal{L}}_{j+1}^n\big)\in\mathbb{S}$.
		Since base case $\big(\overline{\mathcal{Z}}_{j}^n,\overline{\mathcal{L}}_{j}^n\big)\in\mathbb{S}$ is given, the induction is complete. The desired conclusion then follows directly from Lemma\ref{Lemma of dynamic programming equation}.
	\end{proof}
	
	\subsection{Proofs in Section \ref{sec 3}}
	
	\subsubsection{Proof of Lemma \ref{weak convergence}, Theorem \ref{weak convergence app} and Proposition \ref{convergence of Lp}}\label{appendix 2}
	
	\noindent To prove Lemma \ref{weak convergence}, we require the following auxiliary lemma. For convenience, define
	\[\nuu^n_{m,k}:=\frac{1}{k-m}\sum_{j=m}^{k-1}\dd_{\overline{\mathcal{Z}}^n_{j}}\times {\nuu}_j^n(\cdot|\overline{\mathcal{Z}}^n_{j}, \overline{\mathcal{L}}^n_{j}),\ \  \overline{\mathcal{L}}^n_{m,k}:=\frac{1}{k-m}\sum_{j=m}^{k-1}\dd_{\overline{\mathcal{Z}}^n_{j}}.\]
	%where for the non-negative integers $m,k$ with $m=m(n)<k=k(n)\le n$, $\lim_{n\to\infty}k/n$, $\lim_{n\to\infty}m/n$ exist, and $\lim_{n\to\infty}(k-m)/n>0$.
	
	\begin{lemma}\label{estimate for nu}
		Let $\phi$ be a bounded measurable function from $V_{\vec{l}(G)}$ to $\R$. Then for every $\varepsilon>0$, every sequence $\{\zz^n \}$ in $V_{\vec{l}(G)}$ and every $n\ge 4\|g\|_{\infty}/\varepsilon$,
		\begin{align*}
			&\overline{\P}_{\zz^n}\left[\frac{k-m}{n}\left|\int_{V_{\vec{l}(G)}} \phi(\yy)\overline{\mathcal{L}}^n_{m,k}({\rm d}\yy)-\int_{V_{\vec{l}(G)}\times V_{\vec{l}(G)}} \phi(\yy)\nuu^n_{m,k}({\rm d}\xx\times {\rm d}\yy) \right|\ge\varepsilon \right]\\
			&=\overline{\P}_{\zz^n}\left[\frac{k-m}{n}\left|\int_{V_{\vec{l}(G)}} \phi(\yy)\overline{\mathcal{L}}^n_{m,k}({\rm d}\yy)-\int_{V_{\vec{l}(G)}} \phi(\yy)(\nuu^n_{m,k})\big|_2({\rm d}\yy) \right|\ge\varepsilon \right]\le \frac{16\|\phi\|_\infty^2 }{n\varepsilon^2},
		\end{align*}
		where $k,m,\nuu_j^n$ satisfy conditions in Lemma \ref{weak convergence}, and for any Borel sets $A$,
		$$(\nuu^n_{k,m})\big|_2(A)=\nuu^n_{k,m}(G\times A).$$
	\end{lemma} 
	
	We postpone the proof of the above lemma  for a moment.

	\begin{proof}[Proof of Lemma \ref{weak convergence}]  The proofs of part \emph{(a)} and part \emph{(b)} follow the same line of argument as in  \cite[Theorem 8.2.8]{DE1997}.

		\vskip 2mm
		\noindent\emph{(a)}. Indeed, because $V_{\vec{l}(G)}$ is finite, every subsequence of $\{\nuu^n_{m,k}\}_{n\ge 0}$ is tight. By Prohorov’s theorem, we can extract a further subsequence that converges in distribution to some random variable $\nuu$.  
		
		Since the map  
		$$
		\nuu^n_{m,k}\mapsto(\nuu^n_{m,k},(\nuu^n_{m,k})|_E)=(\nuu^n_{m,k},\overline{\mathcal{L}}^n_{m,k})
		$$
		is  continuous,  it follows from  \cite[Theorem A.3.6]{DE1997} that 
		$$
		(\nuu^n_{m,k},\overline{\mathcal{L}}^n_{m,k})\Rightarrow(\nuu,\overline{\mathcal{L}}),
		$$  and $\overline{\mathcal{L}}$ coincides with the first marginal of $\nuu$ with probability 1.
		
		Moreover, because $\{\zz^n\}$ takes value in finite space $V_{\vec{l}(G)}$ (defined in Definition \ref{line digraph}), there exists some $\zz\in V_{\vec{l}(G)}$ such that 
		along a further subsequence
		$$
		(\nuu^n_{m,k},\overline{\mathcal{L}}^n_{m,k},\zz^n)\Rightarrow(\nuu,\overline{\mathcal{L}},\zz).
		$$ 
		Here we note that for different integers $n$ the triples $\nuu^n_{m,k},\overline{\mathcal{L}}^n_{m,k}$ live  on different probability spaces.  
		
		Let $(\overline{\Omega}, \overline{\mathscr{F}},\overline{\P}_\zz)$ be a probability space on which  $\nuu$ and $\overline{\mathcal{L}}$  are defined as, respectively,  a transition probability on $V_{\vec{l}(G)}\times V_{\vec{l}(G)}$ and a transition probability measure on $V_{\vec{l}(G)}$ given $\omega\in \bar{\Omega}$, respectively.

		\vskip 2mm
		\noindent\emph{(b)}. Because $\overline{\mathcal{L}}$ is the first marginal of $\nuu$, \cite[Theorem A.5.6]{DE1997} guarantees the existence of  a  transition  probability $\bar{\nu}({\rm d}\zz'|\zz)=\bar{\nu}({\rm d}\zz'|\zz,\omega)$ on $V_{\vec{l}(G)}$ conditional on $(\zz, \omega)\in V_{\vec{l}(G)}\times \overline{\Omega}$
		such that  (\ref{prop (b)}) holds $\overline{\P}_\zz$-a.s. for $\omega\in\overline{\Omega}$.

		\vskip 2mm
		\noindent\emph{(c)}.  For convenience, we fix a convergent subsequence of $(\nuu^n_{m,k},\overline{\mathcal{L}}^n_{m,k})$, simply denote its index again by $n$.  Since $V_{\vec{l}(G)}$ is finite,  the space of bounded functions on $V_{\vec{l}(G)}$ equipped with the discrete metric  is separable. Let $\Xi$ be a countable dense subset of this bounded function space.

		Because $\overline{\mathcal{L}}^n_{m,k}({\rm d}\yy)$ (resp. $\nu^n_{m,k}({\rm d}\xx\times {\rm d}\yy)$) is supported on $V_{\vec{l}(G)}$ (resp. $V_{\vec{l}(G)}\times V_{\vec{l}(G)}$), Lemma \ref{estimate for nu} implies that for all $\phi\in\Xi$ and $\varepsilon>0$,
		\[
		\overline{\P}_{\zz^n}\left(\frac{k-m}{n}\left(\int_{V_{\vec{l}(G)}} \phi(\yy)\overline{\mathcal{L}}^n_{m,k}({\rm d}\yy)-\int_{V_{\vec{l}(G)}\times V_{\vec{l}(G)}} \phi(\yy)\nuu^n_{m,k}({\rm d}\xx\times {\rm d}\yy)\right)\ge\varepsilon\right)\to 0.
		\]
		Thus by \cite[Theorem A.3.7]{DE1997},
		\[
		\frac{k-m}{n}\left(\int_{V_{\vec{l}(G)}} \phi(\yy)\overline{\mathcal{L}}^n_{m,k}({\rm d}\yy)-\int_{V_{\vec{l}(G)}\times V_{\vec{l}(G)}} \phi(\yy)\nuu^n_{m,k}({\rm d}\xx\times {\rm d}\yy)\right)\Rightarrow 0.
		\]
		Since $(\nuu^n_{m,k},\overline{\mathcal{L}}^n_{m,k})\Rightarrow(\nuu,\overline{\mathcal{L}})$,  Skorohod's representation theorem allows us to  assume that,  with probability 1, the following limit holds for all $\phi\in\Xi$:
		\begin{align}
			\lim_{n\to\infty}\frac{k-m}{n}\left(\int_{V_{\vec{l}(G)}}\phi(\yy)\overline{\mathcal{L}}^n_{m,k}({\rm d}\yy)-\int_{V_{\vec{l}(G)}\times V_{\vec{l}(G)}} \phi(\yy)\nu^n_{m,k}({\rm d}\xx\times {\rm d}\yy)\right)=0.\label{estimate in proof of thm 3.1}
		\end{align}
		By part \emph{(a)},
		\begin{align*}
			\lim_{n\to\infty}\int_{V_{\vec{l}(G)}}\phi(\yy)\overline{\mathcal{L}}^n_{m,k}({\rm d}\yy)=&\,\int_{V_{\vec{l}(G)}}\phi(\yy)\overline{\mathcal{L}}({\rm d}\yy),\\
			\lim_{n\to\infty}\int_{V_{\vec{l}(G)}\times V_{\vec{l}(G)}} \phi(\yy)\nu^n_{m,k}({\rm d}\xx\times {\rm d}\yy)=&\,\int_{V_{\vec{l}(G)}\times V_{\vec{l}(G)}} \phi(\yy)\nuu({\rm d}\xx\times {\rm d}\yy).
		\end{align*}
		Combining these limits and employing the decomposition given in part \emph{(b)} of the theorem, we conclude that $\overline{\P}_{\zz}$-a.s. for every $\phi\in\Xi$ (and  hence for every bounded function), 
		$$
		\lim_{n\to\infty}\frac{k-m}{n}=0
		$$
		or
		\[
		\int_{V_{\vec{l}(G)}}\phi(\yy)\overline{\mathcal{L}}({\rm d}\yy)=\int_{V_{\vec{l}(G)}\times V_{\vec{l}(G)}} \phi(\yy)\nuu({\rm d}\xx\times {\rm d}\yy)=\int_{V_{\vec{l}(G)}\times V_{\vec{l}(G)}} \phi(\yy)\bar{\nu}({\rm d}\yy|\xx)\overline{\mathcal{L}}({\rm d}\xx).
		\]
		In the latter case, by \cite[Theorem A.2.2(b)]{DE1997}, 
		$$
		\overline{\mathcal{L}}(\cdot)=\int_{V_{\vec{l}(G)}}\nuu({\rm d}\xx\times \cdot)=\int_{V_{\vec{l}(G)}}\bar{\nu}(\cdot|\xx)\overline{\mathcal{L}}({\rm d}\xx),
		$$ i.e., for every Borel set $B$
		\[
		\overline{\mathcal{L}}(B)=\int_{V_{\vec{l}(G)}} \bar{\nu}(B|\xx) \overline{\mathcal{L}}({\rm d}\xx)=\int_{V_{\vec{l}(G)}} \overline{\qq}(\xx,B)\overline{\mathcal{L}}({\rm d}\xx)=\nuu|_2(B).
		\]
	\end{proof}
	
	\begin{remark}
		\rm The difference between the part \emph{(c)} of Lemma \ref{weak convergence} and  \cite[Theorem 8.2.8]{DE1997} stems from $(k-m)/n$ in (\ref{estimate in proof of thm 3.1}), which preventd the invariance of $\bar{\nu}({\rm d} \yy|\xx)$ from always holding. Indeed, if	
		\begin{equation}
			\lim_{n\to\infty}\Big(\int_{V_{\vec{l}(G)}}\phi(\yy)\overline{\mathcal{L}}^n_{m,k}({\rm d}\yy)-\int_{V_{\vec{l}(G)}\times V_{\vec{l}(G)}} \phi(\yy)\nuu^n_{m,k}({\rm d}\xx\times {\rm d}\yy)\Big)=0,\label{estimate in proof of thm 8.2.8}
		\end{equation}
		then invariance follows directly. However, when $\lim\limits_{n\to\infty}(k-m)/n=0$, condition (\ref{estimate in proof of thm 3.1}) does not guarantee (\ref{estimate in proof of thm 8.2.8}).
	\end{remark}
	
	\vskip2mm
	
	We now proceed to prove Theorem  Theorem \ref{weak convergence app}.	
	
	\begin{proof}[Proof of Theorem \ref{weak convergence app}] We prove the three parts of the theorem separately.
		\vskip2mm
		
		\emph{(a)}.  Note that $\frac{n_k(\omega)}{n}$ and $\mathbf{1}_{E_{k,n}^\omega}$ (with respect to $\omega$) are random variables taking value in a compact set $[0,1]$ and $\{0,1\}^{E}$, respectively.  For every subsequence of $n$, by Lemma \ref{weak convergence} and a diagonal argument, we can select a common sub-subsequence $n'$ such that for all $k=1,\dots,\b$,
		$$
		\Bigg(\frac{1}{n_k(n')}\sum_{j=\overline{\tauu}^{n'}_k}^{\overline{\tauu}^{n'}_{k+1}-1}\dd_{\overline{\mathcal{Z}}^{n'}_{j}}\times {\nuu}_j^{n'}\big(\cdot\big|\overline{\mathcal{Z}}^{n'}_{j},\overline{\mathcal{L}}^{n'}_{j}\big),\frac{1}{n_k(n')}
		\sum_{j=\overline{\tauu}^{n'}_k}^{\overline{\tauu}^{n'}_{k+1}-1}\dd_{\overline{\mathcal{Z}}^{n'}_{j}},
		\frac{n_k(n')}{n'}, \mathbf{1}_{E_{k,n'}^\omega}\Bigg)_{n'\ge 1}
		$$
		converges in distribution to some $\left(\nuu_k,\overline{\mathcal{L}}_k,\mathscr{R}_k,f^\omega\right)$ possessing  the stated  mentioned in Lemma \ref{weak convergence}. Since $f^\omega$ takes value in $\{0,1\}^{E}$, we may write it as $\mathbf{1}_{E_{k,\infty}^\omega}$, where $E_{k,\infty}^\omega$ denote the support of $f^\omega$.

		%%%%%%%%%%%%%%%%%%%%%%%%%%%
		%%%%%%%%%%%%%%%%%%%%%%%%%%%
		
		\vskip 2mm
		\noindent\emph{(b)}. Setting $k=n,m=0$ in Lemma \ref{weak convergence}, we  obtain the limit of $\overline{\mathcal{L}}^{n'}$ exists, denoted it by $\overline{\mathcal{L}}$. Moreover, $\overline{\mathcal{L}}$ satisfies the properties listed in Lemma \ref{weak convergence}.
		
		Observing that $\frac{n_k}{n}\in [0,1]$, 
		$$
		\sum_{k=1}^\b \frac{n_k}{n}\cdot\frac{1}{n_k}\sum_{j=\overline{\tauu}^n_k}^{\overline{\tauu}^n_{k+1}-1}\dd_{\overline{\mathcal{Z}}^n_{j}}=\overline{\mathcal{L}}^n,\ \ \sum_{k=1}^\b n_k=n,
		$$ 
		and $\{E_{k,n}^{\omega}\}\in\mathscr{E}_\zz$ for fixed $\omega$, we have, as $n\to\infty$, 
		$$
		\mathscr{R}_k\in[0,1] ,\ \ \sum_{k=1}^\b \mathscr{R}_k\overline{\mathcal{L}}_k =\overline{\mathcal{L}},\ \ \sum_{k=1}^\b\mathscr{R}_k=1,\ \ \overline{\P}_\zz \mbox{-a.s. for}\ \omega\in\overline{\Omega}.
		$$
		
		\vskip 2mm
		\noindent\emph{(c)}. For convenience, set 
		$$
		\overline{\mathcal{L}}_k^n=\frac{1}{n_k}\sum_{j=\overline{\tauu}^n_{k}}^{\overline{\tauu}^n_{k+1}-1}\dd_{\overline{\mathcal{Z}}^n_{j}}.
		$$
		Note that $(\overline{\mathcal{L}}_k^n,\mathbf{1}_{E_{k,n}^\omega})_{1\le k\le \b}$ belongs to the  closed set $\mathscr{C}$ defined by
		$$
		\big\{(\hat{\mu}_k,f_k)_{1\le k\le \b}:\hat{\mu}_k\in\mathscr{P}(V_{\vec{l}(G)}),f_k\in\{0,1\}^E,\, \{\text{supp} (f_k)\}_{1\le k\le \b}\in\mathscr{E}_\zz, \,\hat{\mu}_k|_E(\text{\rm supp}(f_k))=1\big\}.
		$$
		For every $(\hat{\mu}_k,f_k)_{1\le k \le \b}\in\mathscr{C}$, the following hold: %the properties in part (c) of Theorem \ref{weak convergence app}, i.e.,
		\begin{itemize}
			\item { $f_k=\{f_k(e),e\in E\}\in\{0,1\}^E$},
			\item $\{\text{\rm supp}(f_k)\}_{1\le k\le \b}\in\mathscr{E}_\zz$,
			\item $\text{\rm supp}(\hat{\mu}_k|_E)\subseteq \text{\rm supp}(f_k)$.
		\end{itemize}
		Thus, to prove part \emph{(c)} Theorem \ref{weak convergence app}  it suffices to show that
		$$ 
		(\overline{\mathcal{L}}_k,\mathbf{1}_{E_{k,\infty}^\omega})_{1\le k\le \b}\in\mathscr{C},\ \ \overline{\mathbb{P}}_\zz\text{-a.s.} \ \text{for}\  \omega\in\overline{\Omega}.
		$$
		
		Let $P_n$ and $P_\infty$ denote the distribution measure of
		$$
		(\overline{\mathcal{L}}_k^n(\cdot|\omega),\mathbf{1}_{E_{k,n}^\omega})_{1\le k\le \b}\ \mbox{and}\ (\overline{\mathcal{L}}_k(\cdot|\omega),\mathbf{1}_{E_{k,\infty}^\omega})_{1\le k\le \b}
		$$
		on $\left(\mathscr{P}(V_{\vec{l}(G)})\times\{0,1\}^E\right)^{\b}$. Since 
		$$
		(\overline{\mathcal{L}}_k^{n'}(\cdot|\omega),\mathbf{1}_{E_{k,n'}^\omega})_{1\le k\le \b}\xrightarrow{\mathcal{D}}(\overline{\mathcal{L}}_k(\cdot|\omega),\mathbf{1}_{E_{k,\infty}^\omega})_{1\le k\le \b},$$ we have
		$$
		1=\limsup_{n'\to\infty}P_{n'}(\mathscr{C})\le P_\infty(\mathscr{C}).
		$$
		Hence for $\overline{\mathbb{P}}_\zz$-a.s. $\omega\in\overline{\Omega}$, 
		$$
		(\overline{\mathcal{L}}_k,\mathbf{1}_{E_{k,n}^\omega})_{1\le k\le \b}\in\mathscr{C}.
		$$
	\end{proof} 	
	\vskip2mm
	\begin{proof}[Proof of Proposition \ref{convergence of Lp}] Since $V_{\vec{l}(G)}$ and $E$ are finite, the weak convergence of $\{\hat{\mu}_n \}_{n\ge1}$ and $\{\mathbf{1}_{E_n}\}$ is equivalent to
		\begin{align*}
			\sup_{\zz\in V_{\vec{l}(G)}}|\hat{\mu}_n(\zz)-\hat{\mu}(\zz)|\to 0, {\rm\ and\ }
			\sup_{e\in E}|\mathbf{1}_{E_n}(e)-\mathbf{1}_{E_\infty}(e)|\to 0.
		\end{align*}
		Because $\mathbf{1}_{E_n}(e)$  take values only in $\{0,1\}$ for every l $e\in E$, there exists an integer $N$ such that for all $n>N$,
		\[
		\mathbf{1}_{E_n}(e)=\mathbf{1}_{E_\infty}(e)\text{, for all } e\in E.
		\]
		Thus, $\mathbf{1}_{E_n}= \mathbf{1}_{E_\infty}$ for $n>N$, which implies $\pp_{E_n}=\pp_{E_\infty}$ for $n>N$. Consequently, for sufficiently large  $n$ and any $\zz_1,\zz_2\in V_{\vec{l}(G)}$,
		\begin{align*}
			(\hat{\mu}_n\otimes \pp_{E_n})(\zz_1,\zz_2)&=(\hat{\mu}_n\otimes \pp_{E_\infty})(\zz_1,\zz_2)\\
			&=\hat{\mu}_n(\zz_1)\cdot \pp_{E_\infty}(\zz_1,\zz_2)\\
			&\to \hat{\mu}(\zz_1)\cdot \pp_{E_\infty}(\zz_1,\zz_2)\\
			&=(\hat{\mu}\otimes \pp_{E_\infty})(\zz_1,\zz_2).
		\end{align*}
		Hence $\hat{\mu}_n\otimes \pp_{E_n}\Rightarrow\hat{\mu}\otimes \pp_{E_\infty}$. 
	\end{proof}
	
	%%%%%%%%%%%%%%%%%%%%%%%%%%%
	%%%%%%%%%%%%%%%%%%%%%%%%%%%
	\vskip2mm
	\begin{proof}[Proof of Lemma \ref{estimate for nu}] Define $\overline{\mathscr{F}}_j^n$ as the $\sigma$-field generated by $\{\overline{\mathcal{Z}}_i^n,\overline{\mathcal{L}}_i^n, i=0,1,\dots,j \}$. Noting that $\overline{\mathcal{Z}}_{j+1}^n$ is generated according to $\nuu_j^n(\cdot|\overline{\mathcal{Z}}_j^n,\overline{\mathcal{L}}_j^n)$,  for each bounded measurable function $\phi$,  we have
		$$
		\overline{\E}_{\zz^n}\Big[\phi(\overline{\mathcal{Z}}_{j+1}^n)-\int_{V_{\vec{l}(G)}}
		\phi(\yy)\nuu_j^n({\rm d}\yy|\overline{\mathcal{Z}}_j^n,\overline{\mathcal{L}}_j^n)\Big|\overline{\mathscr{F}}_j^n \Big]=0,\ \overline{\P}_{\zz^n}\text{-a.s.}
		$$
		In other words, 
		$$
		M_j:=\sum_{i=0}^j\phi(\overline{\mathcal{Z}}_{i+1}^n)-\int_{V_{\vec{l}(G)}}
		\phi(y)\nuu_i^n({\rm d}\yy|\overline{\mathcal{Z}}_i^n,\overline{\mathcal{L}}_i^n)
		$$
		is a martingale. Consequently, its increments are uncorrelated:  for $j<j'\le n$,
		\begin{align*}
			&\overline{\E}_{\zz^n}\Big[ \Big(\phi(\overline{\mathcal{Z}}_{j+1}^n)-\int_{V_{\vec{l}(G)}} \phi(\yy)\nuu_j^n({\rm d}\yy|\overline{\mathcal{Z}}_j^n,\overline{\mathcal{L}}_j^n)\Big)\cdot \Big(\phi(\overline{\mathcal{Z}}_{j'+1}^n)-\int_{V_{\vec{l}(G)}} \phi(\yy)\nuu_{j'}^n({\rm d}\yy|\overline{\mathcal{Z}}_{j'}^n,\overline{\mathcal{L}}_{j'}^n)\Big)\Big]\\
			&=\overline{\E}_{\zz^n}\Big[ \Big(\phi(\overline{\mathcal{Z}}_{j+1}^n)-\int_{V_{\vec{l}(G)}} \phi(\yy)\nuu_j^n({\rm d}\yy|\overline{\mathcal{Z}}_j^n,\overline{\mathcal{L}}_j^n)\Big)\\
			&\ \ \ \ \ \ \cdot \overline{\E}_{\zz^n}\Big(\phi(\overline{\mathcal{Z}}_{j'+1}^n)-\int_{V_{\vec{l}(G)}} \phi(\yy)\nuu_{j'}^n({\rm d}\yy|\overline{\mathcal{Z}}_{j'}^n,\overline{\mathcal{L}}_{j'}^n)\Big|\overline{\mathscr{F}}_{j'}^n\Big)\Big]\\
			&=0.
		\end{align*}
		Therefore, 
		\begin{align*}
			\overline{\E}_{\zz^n}\big( M_j^2  \big)=&\ \overline{\E}_{\zz^n}\Big[\Big(\sum_{i=0}^j \big(\phi(\overline{\mathcal{Z}}_{i+1}^n)-\int_{V_{\vec{l}(G)}} \phi(\yy)\nuu_i^n({\rm d}\yy|\overline{\mathcal{Z}}_i^n,\overline{\mathcal{L}}_i^n)\big)\Big)^2\Big]\\
			=&\ \overline{\E}_{\zz^n}\Big[\sum_{i=0}^j \Big(\phi(\overline{\mathcal{Z}}_{i+1}^n)-\int_{V_{\vec{l}(G)}} \phi(\yy)\nuu_i^n({\rm d}\yy|\overline{\mathcal{Z}}_i^n,\overline{\mathcal{L}}_i^n)\Big)^2\Big]\\
			\le&\ 4n\|\phi\|^2_\infty.
		\end{align*}
		Note that
		\begin{align*}
			&\frac{k-m}{n}\Big(\int_{V_{\vec{l}(G)}} \phi(\yy)\overline{\mathcal{L}}^n_{m,k}({\rm d}\yy)-\int_{V_{\vec{l}(G)}\times V_{\vec{l}(G)}} \phi(\yy)\nuu^n_{m,k}({\rm d}\xx\times {\rm d}\yy)\Big)\\
			&=\frac{1}{n}\sum_{j=m}^{k-1}\Big(\phi(\overline{\mathcal{Z}}_{j+1}^n)-\int_{V_{\vec{l}(G)}} \phi(\yy)\nuu_j^n({\rm d}\yy|\overline{\mathcal{Z}}_j^n,\overline{\mathcal{L}}_j^n)\Big)+\frac{1}{n}(\phi(\overline{\mathcal{Z}}_m^n)-\phi(\overline{\mathcal{Z}}_k^n)),
		\end{align*}
		and $\big|\phi(\overline{\mathcal{Z}}_m^n)-\phi(\overline{\mathcal{Z}}_k^n)\big|\le 2\|\phi\|_{\infty}$. Applying Chebyshev's inequality, we obtain 
		for every $\varepsilon>0$ and every $n\ge 4\|\phi\|_\infty/\varepsilon$,
		\begin{align}
			&\overline{\P}_{\zz^n}\Big[\frac{k-m}{n}\Big|\int_{V_{\vec{l}(G)}} \phi(\yy)\overline{\mathcal{L}}^n_{m,k}({\rm d}\yy)-\int_{V_{\vec{l}(G)}\times V_{\vec{l}(G)}} \phi(\yy)\nuu^n_{m,k}({\rm d}\xx\times {\rm d}\yy) \Big|\ge\varepsilon \Big]\nonumber\\
			&\le \overline{\P}_{\zz^n}\Big[\Big|\frac{1}{n}\sum_{j=m}^{k-1}\big(\phi(\overline{\mathcal{Z}}_{j+1}^n)-\int_{V_{\vec{l}(G)}} \phi(\yy)\nuu_j^n({\rm d}\yy|\overline{\mathcal{Z}}_j^n,\overline{\mathcal{L}}_j^n)\big) \Big|\ge\frac{\varepsilon}{2} \Big]\nonumber\\
			&\le \frac{4}{\varepsilon^2}\overline{\E}_{\zz^n}\Big[\frac{1}{n^2}\Big(\sum_{j=m}^{k-1}\big(\phi(\overline{\mathcal{Z}}_{j+1}^n)-\int_{V_{\vec{l}(G)}} \phi(\yy)\nuu_j^n({\rm d}\yy|\overline{\mathcal{Z}}_j^n,\overline{\mathcal{L}}_j^n)\big) \Big)^2 \Big]\nonumber\\
			&=\frac{4}{\varepsilon^2}\overline{\E}_{\zz^n}\big[\frac{1}{n^2}\left(M_{k}-M_{m} \right)^2 \big].\label{estimaite in lemma 7.5}
		\end{align}
		Since $m,k$ are stopping times with $m\le k\le n$. Hence
		\begin{equation}
			\overline{\E}_{\zz^n}[M_{k}M_{m}]=\overline{\E}_{\zz^n}[M_{m}\overline{\E}_{\zz^n}[M_{k}|\mathscr{F}_{m}]]=\overline{\E}_{\zz^n}[M_{m}^2].\label{(7.11)}
		\end{equation}
		Because  $M_t^2$ is a sub-martingale, using (\ref{(7.11)}), we obtain
		\begin{align*}
			(\ref{estimaite in lemma 7.5})=\frac{4}{\varepsilon^2}\overline{\E}_{\zz^n}\left[\frac{1}{n^2}(M_{k}^2-M_{m}^2) \right]
			\le\frac{4}{\varepsilon^2}\overline{\E}_{\zz^n}\left[\frac{1}{n^2}M_{n}^2 \right]
			\le \frac{16\|\phi\|^2_\infty}{n\varepsilon^2}.
		\end{align*}
		This completes the proof.	
	\end{proof}

	\subsubsection{Proof of Lemma \ref{ergodic theorem}} \label{sec pf ergodic theorem}
	\noindent Recall that $\hat{\gamma}_k|_{S_{k,i}}$ is the the normalized restriction of  $\hat{\gamma}_k$ to  $S_{k,i}$. For  convenience,  we introduce the following notations. Let $\overline{\P}_{\hat{\gamma}_k|_{S_{k,i}}}$  be the probability measure conditioned on  $\overline{\mathcal{Z}}_{s_{k,i}}^n$ having law $\hat{\gamma}_k|_{S_{k,i}}$, i.e.,
	\[
	\overline{\P}_{\hat{\gamma}_k|_{S_{k,i}}}(\cdot)=\int_{V_{\vec{l}(G)}} \overline{\P}(\cdot|\overline{\mathcal{Z}}_{s_{k,i}}^n=\zz) \hat{\gamma}_k|_{S_{k,i}}({\rm d}\zz).
	\]
	Denote by $\overline{\E}_{\hat{\gamma}_k|_{S_{k,i}}}$ the expectation induced by  $\overline{\P}_{\hat{\gamma}_k|_{S_{k,i}}}$.
	We now prove   Lemma \ref{ergodic theorem}.
	\begin{proof}[Proof of Lemma \ref{ergodic theorem}] Recall that $S_{k,i}$ is an irreducible positive recurrent set  for transition probability $\qq_k$. For every  $\zz\in S_{k,i}$, the invariant measure $\hat{\gamma}_k|_{S_{k,i}}$ satisfies that $\hat{\gamma}_k|_{S_{k,i}}(\zz)>0$. 
		
		Define
		\begin{align*}
			D^n=\overline{\E}_{\hat{\gamma}_k|_{S_{k,i}}}\Bigg\{ \Bigg|\frac{1}{s_{k+1,i}-s_{k,i}-2\b-1}\cdot\sum_{j=s_{k,i}}^{s_{k+1,i}-2\b-1}R\big(\qq_k(\overline{\mathcal{Z}}_j^n,\cdot)\big\|\pp_{E_k}(\overline{\mathcal{Z}}_j^n,\cdot)\big)\ \ &\\
			- \int_{V_{\vec{l}(G)}}R\big(\qq_k\|\pp_{E_k}\big)\ {\rm d}(\hat{\gamma}_k|_{S_{k,i}})\Bigg|\Bigg\}.&
		\end{align*}
		Set 
		\begin{align*}
			D^n_\zz=\overline{\E}\Bigg\{ \Bigg|&\frac{1}{s_{k+1,i}-s_{k,i}-2\b-1}\cdot\sum_{j=s_{k,i}}^{s_{k+1,i}-2\b-1}R\big(\qq_k(\overline{\mathcal{Z}}_j^n,\cdot)\|\pp_{E_k}(\overline{\mathcal{Z}}_j^n,\cdot)\big)\\
			&- \int_{V_{\vec{l}(G)}}R(\qq_k\|\pp_{E_k})\ {\rm d}(\hat{\gamma}_k|_{S_{k,i}})\Bigg|\  \Bigg|\overline{\mathcal{Z}}_{s_{k,i}}^n=\zz \Bigg\},
		\end{align*}
		where the expectation is conditioned on $\overline{\mathcal{Z}}_{s_{k,i}}^n=\zz\in S_{k,i}$.
		Because the relative entropy is nonnegative and finite for measures on finite space, the $L^1$ ergodic theorem \cite[Theorem A.4.4]{DE1997} yields
		\[
		\lim_{n\to\infty}D^n=\lim_{n\to\infty}\int_{V_{\vec{l}(G)}}D^n_\zz\hat{\gamma}_k|_{S_{k,i}}({\rm d}\zz)=0.
		\]
		Since $S_{k,i}$ is finite and $\hat{\gamma}_k|_{S_{k,i}}(z)>0$ for every $\zz\in S_{k,i}$, the non-negativity of $D^n_\zz$ gives
		\[
		\lim_{n\to\infty}D^n_\zz=0\text{\ \  for all \  }\zz\in S_{k,i}.
		\]
		We now turn to  the convergence of
		\[
		\frac{1}{s_{k+1,i}-s_{k,i}-2\b-1}\cdot\sum_{j=s_{k,i}}^{s_{k+1,i}-2\b-1}\dd_{\overline{\mathcal{Z}}_j^n}.
		\]
		For every $\zz\in S_{k,i}$, define 
		\[
		A(\zz)=\Big\{ \lim_{n\to\infty}\frac{1}{s_{k+1,i}-s_{k,i}-2\b-1}\cdot\sum_{j=s_{k,i}}^{s_{k+1,i}-2\b-1}\dd_{\{\overline{\mathcal{Z}}_j^n=\zz\}} = \hat{\gamma}_k|_{S_{k,i}}(\zz) \Big\}.
		\]
		By the pointwise ergodic theorem \cite[Theorem A.4.4]{DE1997}, for every $\zz'\in S_{k,i}$
		\[
		\overline{\P}_{\hat{\gamma}_k|_{S_{k,i}}}\big(A(\zz')\big)=\int_{V_{\vec{l}(G)}}\overline{\P}\left(A(\zz')\left|\overline{\mathcal{Z}}_{s_{k,i}}^n=\zz\right.\right)\hat{\gamma}_k|_{S_{k,i}}({\rm d}\zz)=1,
		\]
		which in turn implies  
		$$
		\overline{\P}\Big(A(\zz')\Big|\overline{\mathcal{Z}}_{s_{k,i}}^n=\zz\Big)=1,\  \text{for all}\ \ \zz\in S_{k,i}.
		$$
		Therefore,
		\[
		\overline{\P}\Big(\bigcap_{\zz'\in S_{k,i}} A(\zz')\Big|\overline{\mathcal{Z}}_{s_{k,i}}^n=\zz\Big)=1,\  \text{ for all }\zz\in S_{k,i},
		\]
		from which we conclude that, with probability 1,
		\[
		\lim_{n\to\infty}\Big(\frac{1}{s_{k+1,i}-s_{k,i}-2\b-1}\cdot\sum_{j=s_{k,i}}^{s_{k+1,i}-2\b-1}\dd_{\overline{\mathcal{Z}}_j^n}\Big)(\zz)=\hat{\gamma}_k|_{S_{k,i}}(\zz),\ \text{ for all }\zz\in S_{k,i}.
		\]
		Thus the weak convergence holds almost surely.
	\end{proof}

	\subsubsection{Proof of Lemma \ref{lower technique}} \label{pf of Lem 3.8}
	\noindent Before proving this lemma, we outline the intuitive idea.  The Markov chain $(U_n)_{n\ge s}$ starts  at vertex $u_0$ at time $s$, and moves to an adjacent vertex at each step according to a prescribed transition probability, until it reaches some vertex $u$ at time $t-2\b$. Here $u$ can be any vertex of $G$. We then need specify a rule---independent of the law of $(U_n)_{s\le n\le t-2\b}$---that, starting from $u$ at time at  time $t-2\b$, makes $(U)_n$ follow a deterministic path and arrive at $v_0$ at exactly at time $t$.
	
	A heuristic approach is to let  $(U_n)$ follow  the  shortest path from $u$ to $v_0$  after time $t-2\b$.  Once $(U_n)$ first reaches  $v_0$ at some time $t'$,  it then alternates between $v_0$ and one of its neighbors,  finally stopping  at $v_0$ at time $t$.  Two key points arise in this scheme:
	\begin{itemize}
		\item The time difference $t-t'$ must be  even.  While this forces the length of the shortest path from from $u$ to $v_0$ to be even number as well,  the idea can be implemented directly on trees and bipartite graphs (see {\bf \emph{Case 1}} in the proof). For other graphs, modifications are needed: we must design a path that brings $U_0$ to $v_0$ at some time $t'$ while keeping  $t-t'$ even.(see {\bf \emph{Case 2}} in the proof). Since the shortest path from  $u$ to $v_0$   may not be  unique,  one way to fix a unique choice is to take the shortest path on a fixed spanning tree.   
		\item  Although the path  from time $t-2\b$ to $t$ is chosen deterministically, we must still preserve the Markov property of , but we still need to ensure the Markov property of $(U_n)_{t-2\b\le n\le t}$.  That is, for every   $n\in [t-2\b, t)$ and every state   $U_n=\tilde{u}$, the next position  $U_{n+1}$ is determined solely by $n$ and $\tilde{u}$ with probability $1$.
	\end{itemize}
	
	\begin{proof}[Proof of Lemma \ref{lower technique}] On a tree $\mathbf{T}=(E,V)$  rooted at $\varrho$,  let  $l:V\mapsto \Z$ denote the level of a vertex in $\mathbf{T}$. That is, $l(\varrho)=0$ and $l(v)=k$ if $v$ belongs to the  $k$-th generation from $\varrho$.     Let $f_v$ and $c_v$  be the parent and the children of $v$, respectively.  Denote by $A_v$ the set of the ancestors of $v$.   For two $v$ and $w$,  let $a_{vw}$ be
		their youngest common ancestor of  $v$ and $w$, i.e. $a_{vw}\in A_v\cap A_w$ and 
		$$
		l(a_{vw})=\max\{ l(u'): u'\in A_v\cap A_w\}.
		$$
		
		\vskip2mm
		
		\noindent{\bf \emph{Case 1.}} Suppose $G_0$ is a bipartite graph (a tree is a special bipartite graph). That is,  the vertex set of  $G_0$ can be partitioned into two parts, \textsf{Part I} and \textsf{Part II}, such that no two vertices within the same part are adjacent.  Without loss of generality,  assume $u_0\in \textsf{Part I}$.
		
		Note that there exists a path $u_0u_{1}\dots u_{t-s}$ with $u_{t-s}=v_0$. 
		\begin{itemize}
			\item  If $t-s$ is  even, then $v_0\in \textsf{Part I}$ and consequently $U_{t-2\b}=u\in \textsf{Part I}$; 
			\item If $t-s$ is  odd, then $v_0\in \textsf{Part II}$   and   $U_{t-2\b}=u\in \textsf{Part II}$.
		\end{itemize}
		Therefore, the length of the shortest path from $u$ to $v_0$ is necessarily even.
		
		To guarantee the uniqueness of the path, we fix a spanning tree $\mathbf{T}$ (for instance,  by taking the minimum spanning tree on $G_0$ with  edge weights defined as   $10^{|V|} i+j$ for the edges $v_iv_j\ (1\le i<j)$)  and  then construct  the required path on $\mathbf{ T}$. Concretely, the deterministic path is defined by the following transition probability:
		$$\left\{
		\begin{aligned}
			\tilde{p}(v, f_v)=1, &\   \text{if}\ v \notin A_{v_0} \\
			\tilde{p}(v, c_v)=1, &\  \text{if}\ v\in A_{v_0}
		\end{aligned}
		\right..
		$$
		
		Indeed, the process governed by the above transition probability moves from 
		$u$ at time $t-2\b$ to $v_0$  at time $t$ through the following steps (see Figure \ref{tree path}):
		\begin{enumerate}
			\item[\small\bf  \emph{Step 1}.]   $(U_n)$ walks from $u$ to $a_{uv_0}$ along the  shortest path on $\mathbf{ T}$;
			\item[\small\bf  \emph{Step 2}.]   $(U_n)$ walks from $a_{uv_0}$ to $v_0$ along the shortest path on $\mathbf{ T}$;
			\item[\small\bf  \emph{Step 3}.]   $(U_n)$ alternates  between $v_0$ and its parent $f_{v_0}$. (If $v_0$ is the root $\varrho$ of $\mathbf{ T}$,  $(U_n)$ alternates between $v_0$ and a chosen descendants of $v_0$).
		\end{enumerate}
		
		\begin{figure}
			\centering{\includegraphics[scale=1]{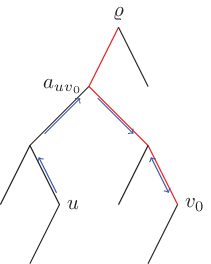}}
			\caption{\small{The red line indicates all possible ancestors of $v_0$.}}
			\label{tree path}
		\end{figure}
		\vskip1mm
		
		\noindent{\bf \emph{Case 2.}}         If $G_0$ is not  bipartite, it contains an odd cycle. Let
		$$
		\mathbf{ C}=w_1w_2\dots w_{2m+1}w_1\ (w_i\neq w_j\text{ if }i\neq j)
		$$
		be such an odd cycle   (see Theorem 11.4 in \cite{Bru2010}).  For each connected component   $G_l\ (1\le l\le L)$ of  $G_0\setminus \mathbf{ C}$   
		denoted by $\mathbf{ C}=w_1w_2\dots w_{2m+1}w_1\ (w_i\neq w_j\text{ if }i\neq j)$ (see Theorem 11.4 in \cite{Bru2010}).  For each connected  component $G_l\ (1\le l\le L)$ of  $G_0\setminus \mathbf{ C}$, choose choose a spanning tree and denote it by     $\mathbf{ T}_l$.    
		Let $A_0=\{1,2,\dots,L\}$.  We now define the family $\{\tilde{\mathbf{ T}}_k:\ 1\le k\le 2m+1\}$ inductively:
		
		\begin{itemize}
			\item Define $A_1=\{ l \in A_0:  \mathbf{ T}_l \text{ is not adjacent to } w_1\text{ by any edge of }E_0\}$.
			\begin{itemize}
				\item  If $A_0\setminus A_1\neq \emptyset$, then for all  $r\in A_0\setminus A_1$ choose a vertex $v_r\in \mathbf{ T}_r$ such that $v_r\sim w_1$.  Let $\tilde{\mathbf{ T}}_1$ be a tree with the vertex set 
				$$
				\{w_1\}\cup(\cup_{r\in A_0\setminus A_1}V(\mathbf{ T}_r))
				$$ and  edge set 
				$$
				\{w_1v_r:\ r\in A_0\setminus A_1\}\cup(\cup_{r\in A_0\setminus A_1}E(\mathbf{ T}_r)).
				$$
				
				\item If $A_0=A_1$, let $\tilde{\mathbf{ T}}_1$ be the  graph (a trivial  tree) with a unique vertex $w_1$ and no edges.
				
			\end{itemize}
			\item For $k\le 2m$, assume that  $A_k$ and  $\tilde{\mathbf{ T}}_k$ have already been defined. Set
			$$
			A_{k+1}=\{r\in A_k:\ \mathbf{ T}_r \text{ is not adjacent to } w_{k+1}\text{ bya any  edge of }E_0\}.
			$$
			\begin{itemize}
				\item If $A_k\setminus A_{k+1}\neq \emptyset$, then for each  $r\in A_k\setminus A_{k+1}$ choose a vertex $v_r\in \mathbf{ T}_r$ such that $v_r\sim w_{k+1}$. Define $\tilde{\mathbf{ T}}_{k+1}$ as the tree with the vertex set 
				$$
				\{w_{k+1}\}\cup(\cup_{r\in A_k\setminus A_{k+1}}V(\mathbf{T}_r))
				$$ and edge set 
				$$
				\{w_{k+1}v_r:\ r\in A_k\setminus A_{k+1}\}\cup(\cup_{r\in A_k\setminus A_{k+1}}E(\mathbf{ T}_r)).
				$$
				
				\item If $A_k=A_{k+1}$, let $\tilde{\mathbf{ T}}_{k+1}$ be the graph (a trivial tree) consisting only of the vertex  $w_{k+1}$ and no edges.
			\end{itemize}
			
			\item  Due to  the connectivity of $G_0$,  every tree $\mathbf{ T}_l$ is adjacent to at least one vertex  $w_k$ with $k\in\{1,\dots,2m+1\}$. Consequently,  $A_{2m+1}=\emptyset$ and the vertex set of $\cup_{k=1}^{2m+1}\tilde{\mathbf{ T}}_k$ coincides with $V_0$.  Take $w_k$ as the root of the tree $\tilde{\mathbf{ T}}_k=(V(\tilde{\mathbf{ T}}_k), E(\tilde{\mathbf{ T}}_k))$.

			Intuitively, the graph $(\cup_{k=1}^{2m+1}\tilde{\mathbf{ T}}_k)\cup \mathbf{ C}$ resembles a key ring if we picture each $\tilde{\mathbf{ T}}_k$ as a key and the odd cycle $\mathbf{ C}$  as the ring (see Figure \ref{cycle tree} (a)). 
		\end{itemize}

		\begin{figure}[htbp]
			\centering
			\subfigure[]{\begin{minipage}[t]{0.3\textwidth}
					\centering{\includegraphics[scale=0.5]{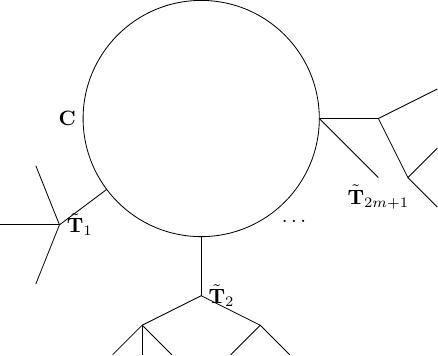}}
			\end{minipage}}
			\subfigure[]{\begin{minipage}[t]{0.3\textwidth}
					\centering{\includegraphics[scale=0.5]{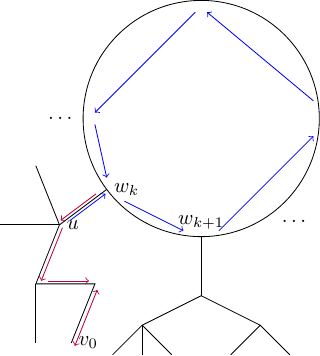}}
			\end{minipage}}
			\subfigure[]{\begin{minipage}[t]{0.3\textwidth}
					\centering{\includegraphics[scale=0.5]{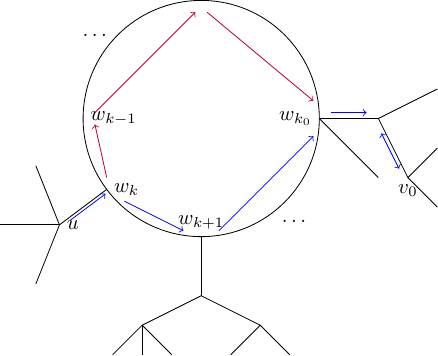}}
			\end{minipage}}
			\caption{}
			\label{cycle tree}
		\end{figure}

		In this case, heuristically, if the length of the shortest path (on the reduced graph) between $u$ and $v_0$ is odd, we can adjust the direction of the path on  $\mathbf{ C}$  to obtain a new path sarisfying our requirements.
		
		If $v$ belongs to $\tilde{\mathbf{ T}}_k$,  denote by  $l_k(v)$, $f_v^k$, $c_v^k$, $A_v^k$ and $a_{vw}^k$ respectively the level, parent, children, ancestors, and youngest common ancestor of 
		$v$ within  $\tilde{\mathbf{ T}}_k$, where we set  $l_k(w_k)=0$.
		\vskip1mm

		Let $U_{t-2\b}=u\in \tilde{\mathbf{ T}}_k$ and $v_0\in \tilde{\mathbf{ T}}_{k_0}$ . We now define the transition probability for all $n\ge t-2\b$:
		
		\begin{itemize}
			\item If $v\in \tilde{\mathbf{ T}}_{k}$ and $n-t+l_{k}(v)+l_{k_0}(v_0)+k_0-k+{\bf 1}_{\{k_0<k\}}\in 2\mathbb{Z}$, then 
			$$
			\left\{
			\begin{aligned}
				\tilde{p}_n(v, f_v^{k})=1, &  \ \ \text{if}\ v\notin A_{v_0}^{k}\cup \mathbf{C}   \\
				\tilde{p}_n(v, a_{vv_0}^{k})=1, &\ \ \text{if}\ v\in  A_{v_0}^{k} \\
				\tilde{p}_n(v, w_{k+1})=1, &\ \ \text{if}\ v=w_{k}
			\end{aligned}
			\right..
			$$
			
			\item If $v\in \tilde{\mathbf{ T}}_{k}$ and $n-t+l_{k}(v)+l_{k_0}(v_0)+k_0-k+{\bf 1}_{\{k_0<k\}}\in 2\mathbb{Z}+1$, then
			$$
			\left\{
			\begin{aligned}
				\tilde{p}_n(v, f_v^{k})=1, &  \ \ \text{if}\ v\in G_0\setminus \mathbf{ C}  \\
				%       \tilde{p}_n(v, a_{vv_0}^{k})=1, &\ \ \text{if}\ v\in  A_{v_0}^{k}  %\\
				\tilde{p}_n(v, w_{k-1})=1, &\ \ \text{if}\ v=w_{k}
			\end{aligned}
			\right..
			$$
		\end{itemize}
		
		\vskip2mm
		
		If $k=2m+1$, we identify $w_{k+1}$ with $w_1$; if $k=1$, we identify $w_{k-1}$ with $w_{w_{2m+1}}$. The term ${\bf 1}_{\{k_0<k\}}$ is used to decide whether, when $U_n=w_{k_0}$, the next step goes to $w_{k_0+1}$ or $w_{k_0-1}$ or to $w_{k_0-1}$.
		Following the transition probability above, we now describe the explicit path of $(U_n)_{t-2\b\le n\le t}$:
		
		\vskip1mm
		
		\noindent{\bf \small \emph{Case 2.1.}}    $k=k_0$.  
		
		Note that $l_{k}(u)+l_{k}(v_0)$ has the same parity as the length of the shortest path of $u$ and $v_0$ on $\tilde{\mathbf{ T}}_{k}$.
		\vskip1mm
		
		\begin{itemize}
			\item If $l_{k}(u)+l_{k}(v_0)$ is  even,  we construct the required path using the method described in {\bf \emph{Case 1}};
			\item If $l_{k}(u)+l_{k}(v_0)$ is  odd,  we proceed as follows::
			\begin{itemize}
				\item[\small\bf  \emph{Step 1}.] $(U_n)$ walks from $u$ to $w_k$ along their shortest path in $\tilde{\mathbf{ T}}_{k}$;\vskip1mm
				
				\item[\small\bf  \emph{Step 2}.] $(U_n)$ traverses the cycle $\mathbf{ C}$ in the order $w_kw_{k+1}\cdots w_{2m+1}w_1\cdots w_k$  (If need,  set $w_{2m+2}=w_1$);\vskip1mm
				
				\item[\small\bf  \emph{Step 3}.] From the vertex $w_{k}$   reached in   {\small \bf  \emph{Step 2}}, $(U_n)$   follows the shortest path in    $\tilde{\mathbf{ T}}_{k}$ to arrive at $v_0$ at time $t'$ 
				(see Figure \ref{cycle tree} (b)).
			\end{itemize}    
		\end{itemize}
		% \end{itemize}
	
	\noindent{\bf \small \emph{Case 2.2.}}  $k\neq k_0$  (see Figure \ref{cycle tree} (c)).
	
	\begin{itemize}
		\item If $l_k(u)+l_{k_0}(v_0)+k-k_0$ is even  and $k<k_0$,   we construct the required path through  the following steps:
		\begin{itemize}
			\item[\small\bf  \emph{Step 1}.] $(U_n)$ walks from $u$ to $w_k$ along their shortest path;\vskip1mm
			\item[\small\bf  \emph{Step 2}.] $(U_n)$ walks along the path $w_kw_{k+1}\cdots w_{k_0}$ on the cycle $\mathbf{ C}$;\vskip1mm
			\item[\small\bf  \emph{Step 3}.] From $w_{k_0}$,  $(U_n)$ proceeds along the shortest path in $\tilde{\mathbf{ T}}_{k_0}$ to reach $v_0$ at time $t' $.
		\end{itemize}
		\item If $l_k(u)+l_{k_0}(v_0)+k-k_0$ is  even and $k>k_0$   we construct the required path through the following steps:
		\begin{itemize}
			\item[\small\bf  \emph{Step 1}.]  $(U_n)$ walks from $u$ to $w_k$ along their shortest path;\vskip1mm
			\item[\small\bf  \emph{Step 2}.]  $(U_n)$ walks along the path $w_kw_{k-1}\cdots w_{k_0}$ in the cycle $\mathbf{C}$;\vskip1mm
			\item[\small\bf  \emph{Step 3}.]  From $w_{k_0}$, $(U_n)$ proceeds along the shortest path in $\tilde{\mathbf{ T}}_{k_0}$ to  reach $v_0$ at time $t' $.
		\end{itemize}
		\item If $l_k(u)+l_{k_0}(v_0)+k-k_0$  odd  and $k<k_0$,   we construct the required path through the following steps:
		\begin{itemize}
			\item[\small\bf  \emph{Step 1}.]  $(U_n)$ walks from $u$ to $w_k$ along their shortest path;\vskip1mm
			\item[\small\bf  \emph{Step 2}.]  $(U_n)$ walks along the path  $w_kw_{k-1}\cdots w_1w_{2m+1}\cdots w_{k_0+1} w_{k_0}$ in the cycle $\mathbf{C}$;\vskip1mm
			\item[\small\bf  \emph{Step 3}.]   From $w_{k_0}$, $(U_n)$ proceeds along the shortest path in $\tilde{\mathbf{ T}}_{k_0}$ to  reach $v_0$ at time $t' $.  
		\end{itemize}
		\item  $l_k(u)+l_{k_0}(v_0)+k-k_0$ is  odd  and $k>k_0$,   we construct the required path through the following steps:
		\begin{itemize}
			\item[\small\bf  \emph{Step 1}.]  $(U_n)$ walks from $u$ to $w_k$ along their shortest path;\vskip1mm
			\item[\small\bf  \emph{Step 2}.]  $(U_n)$ walks along the path $w_kw_{k+1}\cdots w_{2m+1}w_1\cdots w_{k_0-1} w_{k_0}$ in the cycle $\mathbf{C}$;\vskip1mm
			\item[\small\bf  \emph{Step 3}.]  From $w_{k_0}$, $(U_n)$ proceeds along the shortest path in $\tilde{\mathbf{ T}}_{k_0}$ to  reach $v_0$ at time $t' $.  
		\end{itemize}
		Hence, $t-t'$ s necessarily even. In all the cases above, we complete the construction by:
		
		\item[\small\bf  \emph{Step 4}.] $(U_n)$ alternates between between $v_0$ and its parent $f_{v_0}$ ( If $v_0=w_k, k\le 2m$, then  $(U_n)$ alternates  between $w_k$ and $w_{k+1}$; in particular, if $k=2m+1$,  $(U_n)$ it alternates between $w_{2m+1}$ and $w_1$. )
	\end{itemize}
\end{proof}

\subsubsection{Proof of Lemma \ref{inf attainment}}\label{pf lem 3.10}

\noindent Recall that
\begin{align}
	&\sum_{\zz'\to \zz}\qq(\zz',\zz)\hat{\mu}(\zz')=\hat{\mu}(\zz),\label{equ_1 in I}\\
	&\sum_{\zz'\leftarrow \zz}\qq(\zz,\zz')=1.\label{equ_2 in I}
\end{align}
Our aim is to determine the condition under which $\qq$ attains the infimum.
Using the method of Lagrange multipliers, we seek the minimum of
\begin{align}
	F\big(\{\qq(\zz,\zz')\}_{\zz,\zz'\in \vec{l}(G)}, \{\lambda_\zz^1,\lambda_\zz^2 \}_{\zz\in \vec{l}(G)}\big):=&\int_{V_{\vec{l}(G)}}R(\qq\|\pp_{E'})\ {\rm d}\hat{\mu}\nonumber\\
	&\ \ \ +\sum_{\zz\in V_{\vec{l}(G)}}\lambda_\zz^1(\sum_{\zz'\to \zz}\qq(\zz',\zz)\hat{\mu}(\zz')-\hat{\mu}(\zz))\nonumber\\
	&\ \ \ +\sum_{\zz\in V_{\vec{l}(G)}}\lambda_\zz^2(\sum_{\zz'\leftarrow \zz}\qq(\zz,\zz')-1).
\end{align}

\begin{proof}[Proof of Lemma \ref{inf attainment}]   
	
	Note that
	\[
	\int_{V_{\vec{l}(G)}}R(\qq\|\pp_{E'})\ {\rm d}\hat{\mu}_k=\sum_{\zz\in V_{\vec{l}(G)}}\sum_{\zz'\in V_{\vec{l}(G)}}\qq(\zz,\zz')\log{\frac{\qq(\zz,\zz')}{\pp_{E'}(\zz,\zz')}}.
	\]
	\vskip2mm 
	
	\noindent{\small\bf  \emph{Step 1}}.   
	We examine the denominator inside the logarithm.
	Since  $\pp_{E'}(\zz,\zz')={p}_{E'}((\zz')^-,(\zz')^+)$ for every $\zz\to \zz'$, and using (\ref{equ_1 in I}), we obtain 
	\begin{align*}
		\sum_{\zz\in V_{\vec{l}(G)}}\hat{\mu}_k(\zz)\sum_{\zz'\in V_{\vec{l}(G)}}\qq(\zz,\zz')\log{\pp_{E'}(\zz,\zz')}&=\sum_{\zz'\in V_{\vec{l}(G)}}\sum_{\zz\in V_{\vec{l}(G)}}\hat{\mu}(\zz)\qq(\zz,\zz')\log{{p}_{E'}((\zz')^-,(\zz')^+)}\\
		&=\sum_{\zz'\in V_{\vec{l}(G)}}\hat{\mu}(\zz')\log{{p}_{E'}((\zz')^-,(\zz')^+)}.
	\end{align*}
	This quantity does not depend on the variables of $F$. Consequently, it suffices to find the minimum of  
	\begin{align*}
		\tilde{F}&=F+\sum_{\zz\in V_{\vec{l}(G)}}\hat{\mu}(\zz)\sum_{\zz'\in V_{\vec{l}(G)}}\qq(\zz,\zz')\log{\pp_{E'}(\zz,\zz')}\\
		&=\sum_{\zz\in V_{\vec{l}(G)}}\Big[\sum_{\zz'\in V_{\vec{l}(G)}}\qq(\zz,\zz')\log{\qq(\zz,\zz')}+\lambda_\zz^1(\sum_{\zz'\to \zz}\qq(\zz',\zz)\hat{\mu}(\zz')-\hat{\mu}(\zz))
		+\lambda_\zz^2(\sum_{\zz'\leftarrow \zz}\qq(\zz,\zz')-1)\Big].
	\end{align*}
	\vskip2mm
	
	\noindent{\small\bf  \emph{Step 2}} We solve the following  equations:
	\begin{equation}
		\frac{\partial \tilde{F}}{\partial \qq(\zz,\zz')}=0,\ \ \frac{\partial \tilde{F}}{\partial \lambda_\zz^1}=0,\ \ \frac{\partial \tilde{F}}{\partial \lambda_\zz^2}=0.\label{equ for Lagrange}
	\end{equation}
	
	The first equation in (\ref{equ for Lagrange}) gives
	\[
	\hat{\mu}(\zz)(1+\log{\qq(\zz,\zz')})+\lambda_{\zz'}^1\hat{\mu}(\zz)+\lambda_\zz^2=0.
	\]
	Take $\zz'=\zz_1,\zz_2$ and subtract the two equations. We obtain
	\[
	\log{\frac{\qq(\zz,\zz_1)}{\qq(\zz,\zz_2)}}+\lambda_{\zz_1}^1-\lambda_{\zz_2}^1=0.
	\]
	Therefore, for any $\tilde{\zz}$ with $\tilde{\zz}^+=\zz^+$,
	\[
	\log{\frac{\qq(\zz,\zz_1)}{\qq(\zz,\zz_2)}}=\log{\frac{\qq(\tilde{\zz},\zz_1)}{\qq(\tilde{\zz},\zz_2)}},
	\]
	which is equivalent to
	\[
	\frac{\qq(\zz,\zz_1)}{\qq(\tilde{\zz},\zz_1)}=\frac{\qq(\zz,\zz_2)}{\qq(\tilde{\zz},\zz_2)}.
	\]
	Using (\ref{equ_2 in I}), we conclude 
	\begin{equation}
		\qq(\zz,\zz')=\qq(\tilde{\zz},\zz')\ {\rm for\ all\ }\zz,\tilde{\zz}\to \zz'.\label{min point of Lagrange}
	\end{equation}
	This is the condition for the  infimum.
	\vskip2mm
	
	\noindent{\small\bf  \emph{Step 3}}.  Under this condition \eqref{min point of Lagrange}, for $x,y\in V$ we define 
	$${q}(x,y):=\qq(\zz,\zz'), \ \  \text{where}\ \  \zz\to \zz'\ \  \text{and}\ \  (\zz')^-=x,(\zz')^+=y.$$
	Using (\ref{equ_1 in I}), for every $x\in V$
	\[
	\sum_{\zz^+=x}{q}(x,y)\hat{\mu}(\zz)=\hat{\mu}(\zz_0),
	\]
	where $\zz_0^-=x$ and $\zz_0^+=y$.
	
	Summing over all $y$, we obtain that
	\[
	\sum_{\zz^+=x}\hat{\mu}(\zz)=\sum_{\zz^-=x}\hat{\mu}(\zz)=T(\hat{\mu})(x).
	\]
	Set $\nu=T(\hat{\mu})$.
	Then, under the infimum condition (\ref{min point of Lagrange}),
	\begin{align}
		\int_{V_{\vec{l}(G)}}R(\qq\|\pp_{E'})\ {\rm d}\hat{\mu} &= \sum_{\zz\in V_{\vec{l}(G)}}\hat{\mu}(\zz)\sum_{\zz'\in V_{\vec{l}(G)}}\qq(\zz,\zz')\log{\frac{\qq(\zz,\zz')}{\pp_{E'}(\zz,\zz')}}\nonumber \\
		&= \sum_{\zz\in V_{\vec{l}(G)}}\hat{\mu}(\zz)\sum_{y\in V}{q}(\zz^+,y)\log{\frac{{q}(\zz^+,y)}{{p}_{E'}(\zz^+,y)}}\nonumber\\
		&= \sum_{x\in V}\sum_{\zz^+=x}\hat{\mu}(\zz)\sum_{y\in V}{q}(\zz^+,y)\log{\frac{{q}(\zz^+,y)}{{p}_{E'}(\zz^+,y)}}\nonumber \\
		&= \sum_{x\in V}\nu(x)\sum_{y\in V}{q}(x,y)\log{\frac{{q}(x,y)}{{p}_{E'}(x,y)}}\nonumber\\
		&= \int_{V}R({q}\|{p}_{E'})\ {\rm d}\nu. \label{qq=q-min point of Lagrange}
	\end{align}
	
	\vskip2mm
	\noindent{\small\bf  \emph{Step 4}}.  Using the equality \eqref{qq=q-min point of Lagrange} , we obtain
	\begin{align*}
		\inf_{\qq\in\mathscr{T}_{\vec{l}(G)}:\hat{\mu} \qq=\hat{\mu}}\int_{V_{\vec{l}(G)}}R(\qq\|\pp_{E'})\ {\rm d}\hat{\mu}&=\inf_{\stackrel{\qq\in\mathscr{T}_{\vec{l}(G)}:\hat{\mu} \qq=\hat{\mu}}{\qq(\zz,\zz')=\qq(\tilde{\zz},\zz')\ {\rm for}\ \zz,\tilde{\zz}\to \zz'}}\int_{V}R({q}\|{p}_{E'})\ {\rm d}\nu\nonumber\\
		&=\inf_{{q}\in\mathscr{T}_G: \sum_{\zz':(\zz')^+=\zz^-}\hat{\mu}(\zz'){q}(\zz^-,\zz^+)=\hat{\mu}(\zz) }\int_{V}R({q}\|{p}_{E'})\ {\rm d}\nu\\
		&=\inf_{{q}\in\mathscr{T}_G: \nu(\zz^-){q}(\zz^-,\zz^+)=\hat{\mu}(\zz) }\int_{V}R({q}\|{p}_{E'})\ {\rm d}\nu
	\end{align*}
	For any  $x\in V$ with $\nu(x)>0$, 
	$$
	\sum_{y\in V}{q}(x,y)=\sum_{\zz:\zz^-=x}\frac{\hat{\mu}(\zz)}{\nu(\zz^-)}=1,
	$$
	which shows that ${q}\in\mathscr{T}_G$. This completes the proof.     
\end{proof}

\subsection{Proofs and Derivation of the Rate Function in Section \ref{sec 5}}\label{Proof Prop5.1}
\subsubsection{Proof of Proposition \ref{Induction on tree}}\label{Proof Prop5.1}
\begin{proof}[Proof of Proposition \ref{Induction on tree}]  For readability, write ${q}_{\nu_{V'}}={q}_{V'}$. 
	
	Let $f(v)$ and $c(v)$ denote respectively  the parent (if exists) and the  set of children (if exists) of a vertex $v$. For a measure $\nu_{V'}\in\mathscr{P}(V)$ supported on $V'$. The transition probability $q_{V'}$  can be determined inductively starting from the leaves.
	
	\noindent{\small\bf  \emph{Step 1}}.    For a leaf $v$ we have  ${q}_{V'}(v,f(v))=1$ and the  stationarity condition 
	$$
	{q}_{V'}(f(v),v)\nu_{V'}(f(v))=\nu_{V'}(v),
	$$
	which gives 
	$$
	{q}_{V'}(f(v),v)=\frac{\nu_{V'}(v)}{\nu_{V'}(f(v))}.
	$$ 
	\vskip1mm
	\noindent{\small\bf  \emph{Step 2}}.  Now suppose for some vertex 
	$v$ the values ${q}_{V'}(v,u)$ and ${q}_{V'}(u,v)$  are already known for every child   $u\in c(v)$. If $f(v)$ exists,  then the normalization condition forces
	$$
	{q}_{V'}(v,f(v))=1-\sum_{u\in c(v)}{q}_{V'}(v,u).
	$$
	Moreover, when both $f(v)$ and $c(v)$ exist, the  stationarity equation at $v$ reads
	$$
	\nu_{V'}(v)=\nu_{V'}(f(v)){q}_{V'}(f(v),v)+\sum_{u\in c(v)}\nu_{V'}(u){q}_{V'}(u,v),
	$$
	which uniquely determines ${q}_{V'}(f(v),v)$.
	
	\vskip1mm
	\noindent{\small\bf  \emph{Step 3}}. Proceeding inductively from the leaves toward the root, every entry of $q_{V'}$ is uniquely fixed by $\nu_{V'}$.  This completes the proof. 
\end{proof}

\subsubsection{Derivation of the rate function in Example \ref{ldp on urn} and Example \ref{ldp on path}} \label{sec 5.4}
\noindent{\emph{Derivation of the rate function in Example \ref{ldp on urn}. }} 

\vskip2mm

Note that $\mathscr{E}=\{\{E_{v_0},E\},\{E_{v_2},E\} \}$, where $E_{v_i}=\{v_1v_i\}\ (i=0,2)$ and assume for  $\nu\in \mathscr{P}(\{v_0,v_1,v_2\})$, the infimum in the representation of rate function \eqref{rate function tree} is attained at $(\nu_k, r_k, \{E_{v_j},E \})_{k=1,2}$ and for some $j=0$ or $j=2$.  

\vskip2mm

\noindent{\small\bf  \emph{Step 1}}. Let  $\nu_1$ be the probability measures on the three-vertex tree satisfying $\text{\rm{supp}}(\nu_{1})=\{v_1,v_j\}$. 

Let $q_1\in \mathscr{T}_{\{v_1,v_j\}}$ be the unique solution to $\nu_{1}{q}_{1}=\nu_{1}$. Then, we have  ${q}_{1}(v_1,v_j)=1$. This forces  $\nu_{1}(v_j)=\nu_{1}(v_1)=\frac{1}{2}$.

Let $q_2\in \mathscr{T}_{\{v_0,v_1,v_2\}}$ be the unique solution to $\nu_{2}{q}_{2}=\nu_{2}$. Then, we obtain ${q}_2(v_i,v_1)=1$ for each $i=0,2$, it follows that $\nu_2(v_1)=\frac{1}{2}$.
\vskip2mm

Now set $\nu(v_0)=x$, $q_2(v_1,v_0)=y$, $r_1=a$, and $r_2=b$. Then
\[\left\{
\begin{aligned}
	& x=\frac{1}{2} (y\cdot b + a),\\
	& a+b=1,\\
	& \nu_2(v_0)=\frac{1}{2}y,\ \nu_2(v_2)=\frac{1}{2}(1-y).
\end{aligned}
\right.
\]

\noindent{\small\bf  \emph{Step 2}}. If $j=0$, observe  that
\begin{align}
	I_\delta^{v_1}(\nu)&=\inf_{y,b}\Big\{(1-b)\cdot\frac{1}{2}\log{\frac{1+\delta}{\delta}}+b\cdot\frac{1}{2}\big(y\log{2y}+(1-y)\log{2(1-y)}\big)   \Big\}\nonumber\\
	&=\inf_{b}\Big\{(1-b)\cdot\frac{1}{2}\log{\frac{1+\delta}{\delta}}+\frac{1}{2}\big((2x-1+b)\log{\frac{4x-2+2b}{b}}+(1-2x)\log{\frac{2-4x}{b}}\big)\Big\},\label{5-3-3-vertex-1}
\end{align}
where $\frac{2x-1+b}{b}\in[0,1]$. Let $g(b)$ denote the expression inside  the infimum of $I_\delta^{v_1}(\nu)$. Then
\begin{align*}
	2\cdot\frac{\partial g(b)}{\partial b}=\log{\left(2\cdot\frac{2x-1+b}{b}\right)}-\log{\frac{1+\delta}{\delta}}.
\end{align*}
Thus, 
$$
\frac{\partial g(b)}{\partial b}>0 \Leftrightarrow \frac{1-2x}{b}<\frac{\delta-1}{2\delta}.
$$
\vskip1mm

\begin{itemize}
	\item
	For $\delta\le 1$, $\frac{1-2x}{b}<\frac{\delta-1}{2\delta}$ never  holds for any $b>0$; consequently, the infimum in \eqref{5-3-3-vertex-1} is attained at  $b=1$. Hence 
	$$
	I_\delta^{v_1}(\nu)=I_1(\nu)=R(\nu\|{\mu}).
	$$
	
	\item For $\delta>1$,  
	$$
	b>\frac{2\delta}{\delta-1}(1-2x)>1-2x.
	$$ 
	\begin{itemize}
		\item If $\frac{2\delta}{\delta-1}(1-2x)< 1$, i.e., $x>\frac{\delta+1}{4\delta}$, the infimum the infimum in \eqref{5-3-3-vertex-1} is attained at 
		$$
		b=\frac{2\delta}{\delta-1}(1-2x),
		$$
		and
		$$
		I_\delta^{v_1}(\nu)=x\log{\frac{1+\delta}{\delta}}+\left(\frac{1}{2}-x\right)\log{\frac{\delta-1}{\delta}}.
		$$
		\item
		If  $\frac{2\delta}{\delta-1}(1-2x)\ge 1$, i.e., $x\le \frac{\delta+1}{4\delta}$, the infimum in \eqref{5-3-3-vertex-1} is attained at 
		$
		b=1.
		$
	\end{itemize}
\end{itemize}
\vskip1mm

\noindent{\small\bf  \emph{Step 3}}.  If $j=2$, the infimum can be obtained in a similar manner. 
\begin{itemize}
	\item For $\delta\le 1$, the infimum in \eqref{5-3-3-vertex-1} is attained at $b=1$. 
	
	\item For $\delta>1$, 
	\begin{itemize}
		\item if $x<\frac{\delta-1}{4\delta}$, the infimum in \eqref{5-3-3-vertex-1}  is attained at $b=\frac{2\delta}{\delta-1}\cdot 2x$,
		\[
		I_\delta^{v_1}(\nu)=\left(\frac{1}{2}-x\right)\log{\frac{1+\delta}{\delta}}+x\log{\frac{\delta-1}{\delta}}.
		\]
		\item if $x\ge \frac{\delta-1}{4\delta}$, the infimum in \eqref{5-3-3-vertex-1}  is attained at $b=1$.
	\end{itemize}
	
\end{itemize}

\noindent{\small\bf \emph{Step 4}}.  Therefore, 
\begin{itemize} 
	\item for $\delta\le1$, $I_\delta^{v_1}(\nu)=R(\nu\|{\mu})$;
	
	\item for $\delta>1$,  using the inequalities
	\begin{eqnarray*}
		&&x\log{\frac{1+\delta}{\delta}}+\left(\frac{1}{2}-x\right)\log{\frac{\delta-1}{\delta}}\le x\log{4x}+\left(\frac{1}{2}-x\right)\log{2(1-2x)}, \ x>\frac{\delta+1}{4\delta},\\
		&&\left(\frac{1}{2}-x\right)\log{\frac{1+\delta}{\delta}}+x\log{\frac{\delta-1}{\delta}}\le x\log{4x}+\left(\frac{1}{2}-x\right)\log{2(1-2x)}, \ x<\frac{\delta-1}{4\delta},
	\end{eqnarray*}
	the rate function can be written as
	$$
	I_\delta^{v_1}(\nu)=\left\{
	\begin{aligned}
		&x\log{\frac{\delta-1}{\delta}}+\left(\frac{1}{2}-x\right)\log{\frac{\delta+1}{\delta}},\ &0\le x\le\frac{\delta-1}{4\delta},\\
		&x\log{4x}+\left(\frac{1}{2}-x\right)\log{2(1-2x)},\ &\frac{\delta-1}{4\delta}\le x\le \frac{\delta+1}{4\delta},\\
		&\left(\frac{1}{2}-x\right)\log{\frac{\delta-1}{\delta}}+x\log{\frac{\delta+1}{\delta}},\ &\frac{\delta+1}{4\delta}\le x\le\frac{1}{2}.
	\end{aligned}
	\right.
	$$
	For all $\delta>0$, it can also be expressed in the compact form
	\[
	I_{\delta}^{v_1}(\nu)=\left\{
	\begin{aligned}
		&R(\nu\|{\mu})-R(\nu\|\mathfrak{m}_1(\nu)),\ &\nu(v_1)=\frac{1}{2},\\
		&\infty,\ &\text{otherwise},
	\end{aligned}
	\right.
	\]
	where ${\mu}(v_1)=2{\mu}(v_0)=2{\mu}(v_2)=\frac{1}{2}$, and  the mapping $\mathfrak{m}': \mathscr{P}(\{v_0,v_1,v_2\})\mapsto \mathscr{P}(\{v_0,v_1,v_2\})$ is defined by
	
	\[
	\mathfrak{m}_1(\nu)(v_i)=\left\{
	\begin{aligned}
		&\left(\nu(v_i)\vee\frac{\delta-1}{4\delta}\right)\wedge\frac{\delta+1}{4\delta},\ &i\neq1,\\
		&\frac{1}{2},\ &i=1.
	\end{aligned}
	\right.
	\]
\end{itemize}
\hfill $\Box$

\vskip2mm

\noindent{\emph{Derivation of the rate function in Example \ref{ldp on path}. }} 

\vskip2mm

In the example, note that $\mathscr{E}=\{\{E_{v_0},E\}\}$ with $E_{v_0}=\{v_0v_1\}$, which coincides with $j=0$ case  in the derivation of Example \ref{ldp on urn}.   Hence, setting $\nu(v_0)=x$ and following  the same derivation procedure, we obtain
$$
I_\delta^{v_0}(\nu)=\left\{
\begin{aligned}
	& x\log{\frac{1+\delta}{\delta}}+\left(\frac{1}{2}-x\right)\log{\frac{\delta-1}{\delta}}, &  \delta>1,\ x>\frac{\delta+1}{4\delta} \\
	&  I_1(\nu)=R(\nu\|{\mu}), &  \text{otherwise} 
\end{aligned}
\right..
$$
Equivalently,
$$
I_\delta^{v_0}(\nu)=\left\{
\begin{aligned}
	&x\log{4x}+\left(\frac{1}{2}-x\right)\log{2(1-2x)},\ &0\le x\le \frac{\delta+1}{4\delta},\\
	&\left(\frac{1}{2}-x\right)\log{\frac{\delta-1}{\delta}}+x\log{\frac{\delta+1}{\delta}},\ &\frac{\delta+1}{4\delta}\le x\le\frac{1}{2}.
\end{aligned}
\right.
$$
For all $\delta>0$,  it can also be written in the compact form
\[
I_{\delta}^{v_0}(\nu)=\left\{
\begin{aligned}
	&R(\nu\|{\mu})-R(\nu\|\mathfrak{m}_0(\nu)),\ &\nu(v_1)=\frac{1}{2},\\
	&\infty,\ &\text{otherwise},
\end{aligned}
\right.
\]
where ${\mu}(v_1)=2{\mu}(v_0)=2{\mu}(v_2)=\frac{1}{2}$, and the mapping $\mathfrak{m}_0: \mathscr{P}(\{v_0,v_1,v_2\})\mapsto \mathscr{P}(\{v_0,v_1,v_2\})$ is defined by
\[
\mathfrak{m}_0(\nu)(v_i)=\left\{
\begin{aligned}
	&\nu(v_0)\wedge\frac{\delta+1}{4\delta},\ &i=0,\\
	&\frac{1}{2},\ &i=1,\\
	&\nu(v_2)\vee\frac{\delta-1}{4\delta},\ &i=2.
\end{aligned}
\right.
\]
\hfill $\Box$
%%%%%%%%%%%%%%%%%%%%%%%%%%%%%%%%%
%%%%%%%%%%%%%%%%%%%%%%%%%%%%%%%%%

}
\noindent 

\noindent 
\vskip1cm

\noindent    {\bf Acknowldgement:} We are deeply grateful to the referees, the editor and Professor  Renming Song for their valuable suggestions and assistance, which have significantly improved the paper.  
%%%%%%%%%%%%%%%%%%%%%%%%%%
%%%%%%%%%%%%%%%%%%%%%%%%%%  

%%%%%%%%%%%%%%%%%%%%%%%%%%%%%%%%%%%	
%%%%%%%%%%%%%%%%%%%%%%%%%%%%%%%%%%%
%%%%%%%%%%%%%%%%%%%%%%%%%%%%%%%%%%%	
%%%%%%%%%%%%%%%%%%%%%%%%%%%%%%%%%%%	

{
\end{document}